\documentclass[3p,times]{elsarticle}

\usepackage{graphicx}
\usepackage{psfrag}
\usepackage[svgnames]{xcolor}
\usepackage{amssymb}
\usepackage{amsmath}
\usepackage{dsfont}
\usepackage{rotating}
\usepackage{bm}
\usepackage{url}
\usepackage[toc]{appendix}
\usepackage{relsize} 
\usepackage{hyperref}
\usepackage{multirow}
\hypersetup{
	colorlinks=true,                          
	linkcolor=DarkRed,
	citecolor=DarkRed,
	urlcolor=DarkRed  }

\newcommand{\Frac}{\displaystyle\frac}

\newcommand{\bea}{\begin{equation}}
	\newcommand{\eea}{\end{equation}}
\newcommand{\beqa}{\begin{eqnarray}}
	\newcommand{\eeqa}{\end{eqnarray}}

\newcommand{\rmd}{{\rm d}}
\newcommand{\rmdt}{{\rm d}\,t}

\newcommand{\tr}{\text{tr}}
\newcommand{\mddt}[1]{\Frac{\rmd {#1}}{\rmdt}}

\newcommand{\pdd}[2]{\Frac{\partial #1}{\partial #2}}

\newcommand{\DIV}[1]{\nabla \cdot #1}

\renewcommand{\d}{\textnormal{d}}
\newcommand{\Id}{\tensor{I}}

\newcommand{\dt}{\Delta t}

\newcommand{\Q}{\mathbf{Q}}

\newcommand{\x}{\mathbf{x}} 
\newcommand{\n}{\mathbf{n}}
\renewcommand{\v}{\mathbf{v}}

\newcommand{\vv}{\mathbf{v}}

\newcommand{\w}{\mathbf{w}}

\newcommand{\svol}{\omega}
\newcommand{\del}{\nabla}
\newcommand{\transpose}{{\mathsf{\mathsmaller T}}}
\newcommand{\csh}{c_\textrm{sh}}
\newcommand{\ch}{c_\textrm{h}}
\newcommand{\Ft}{\tensor{F}_t}
\newcommand{\Fe}{\tensor{F}_e}
\newcommand{\Fp}{\tensor{F}_p}

\newcommand{\Ge}{{\tensor{G}_e}}

\newcommand{\Gei}{\tensor{G}_{e_i}}
\newcommand{\Geio}{\tensor{G}_{e_i(0)}}
\newcommand{\Geioo}{\tensor{G}_{e_i(1)}}
\newcommand{\devGe}{\mathring{{\tensor{G}}_e}}
\newcommand{\devGei}{\mathring{\tensor{G}}_{e_i}}
\newcommand{\devS}{\mathring{\tensor{\sigma}}}
\newcommand{\J}{\mathbf{J}}
\renewcommand{\H}{\mathbf{H}}
\newcommand{\q}{\mathbf{q}}
\newcommand{\cv}{\mathbf{c}}
\newcommand{\sG}{\Theta}
\newcommand{\sH}{\Psi}
\newcommand{\visc}{\mu}
\newcommand{\Gs}{G}
\newcommand{\gradv}{\mathbb{L}}

\newcommand{\Ae}{\tensor{A}_e}

\newcommand{\Mcp}{\tensor{M}_{ir}}

\newcommand{\Mp}{\tensor{M}_{r}}

\newfont{\numerikEleven}{ecrm1000}
\newfont{\numerikTen}{cmss10}
\newfont{\numerikNine}{cmss9}
\newfont{\numerikEight}{cmss8}
\DeclareMathAlphabet{\mathsfbi}{OT1}{\sfdefault}{bx}{sl}
\DeclareMathVersion{sfletters}
\SetSymbolFont{letters}{sfletters}{OML}{ntxsfmi}{b}{it}

\makeatletter
\newcommand{\mathbfsbilow}[1]{%
	\text{\mathversion{sfletters}$\m@th#1$}%
}
\DeclareRobustCommand{\tensor}[1]{%
	\begingroup
	\ifcat\noexpand #1\relax
	\edef\greek@test{\detokenize{#1}}%
	\edef\greek@test{\expandafter\@cdr\greek@test\@nil}%
	\edef\greek@test{\expandafter\@car\greek@test\@nil}%
	\edef\x{\the\lccode\expandafter`\greek@test}%
	\edef\y{\number\expandafter`\greek@test}%
	\ifnum\x=\y\relax
	\mathbfsbilow{#1}%
	\else
	\mathsfbi{#1}%
	\fi
	\else
	\mathsfbi{#1}%
	\fi
	\endgroup
}
\makeatother

\renewcommand{\vec}[1]{\mathbf{#1}}

\newcommand{\Lstar}{\tensor{L}_\ast}
\newcommand{\Pstar}{\vec{P}_\ast}

\newcommand{\de}{\mathrm{d}}

\newcommand{\dev}{\mathrm{dev}}

\newcommand{\abs}{\mathrm{abs}}

\journal{Journal of Computational Physics}
\begin{document}

\begin{frontmatter}	
	
\title{A cell-centered implicit-explicit Lagrangian scheme for a unified model of nonlinear 
continuum mechanics on unstructured meshes}
	
\author[ferrara]{Walter Boscheri$^*$}
\ead{walter.boscheri@unife.it}
\cortext[cor1]{Corresponding author}

\author[trento]{Simone Chiocchetti}
\ead{simone.chiocchetti@unitn.it}

\author[trento]{Ilya Peshkov}
\ead{ilya.peshkov@unitn.it}

\address[ferrara]{Department of Mathematics and Computer Science, University of
	Ferrara, 44121 Ferrara, Italy}	
\address[trento]{Department of Civil, Environmental and Mechanical Engineering, University of
	Trento, 38123 Trento, Italy}
	
\begin{abstract}
	A cell-centered implicit-explicit updated Lagrangian finite volume scheme on unstructured grids 
	is proposed for a unified first-order hyperbolic formulation of continuum fluid and solid mechanics, 
	namely the Godunov-Pehskov-Romenski (GPR) model. The scheme provably respects the stiff relaxation 
	limits of the continuous model at the fully discrete level, thus it is asymptotic preserving. 
	Furthermore, the GCL is satisfied by a compatible discretization that makes use of a nodal 
	solver to compute vertex-based fluxes that are used both for the motion of the computational 
	mesh as well as for the time evolution of the governing PDEs. Second order of accuracy 
	in space is achieved using a TVD piecewise linear reconstruction, while an implicit-explicit 
	(IMEX) Runge-Kutta time discretization allows the scheme to obtain higher accuracy also in time. 
	Particular care is devoted to the design of a stiff ODE solver, based on approximate analytical 
	solutions of the governing equations, that plays a crucial role 
	when the visco-plastic limit of the model is approached. We demonstrate the accuracy and robustness of 
	the scheme on a wide spectrum of material responses covered by the unified continuum model 
	that includes inviscid hydrodynamics, viscous heat conducting fluids, elastic and elasto-plastic 
	solids in multidimensional settings.
\end{abstract}

\end{frontmatter}

\section{Introduction}	\label{sec.introduction}

The motion of a continuous medium (either fluid or solid) can be described using two different types 
of reference frames. The first type exploits frames that are co-moving and co-deforming together 
with the medium and traditionally called Lagrangian frames. The second type is referred to as Eulerian 
frames of reference that are neither co-moving nor co-deforming with the continuous medium. Each of the two types 
has its pros and cons when used for numerical simulations. For example, approaches based on the 
Lagrangian description allow for accurate tracking of material interfaces and thus are 
dominated in the computational solid mechanics, or even fluid mechanics algorithms aiming to compute 
hydrodynamics effects when the vorticity field is not too severe. However, Lagrangian methods 
are impractical for the simulation of phenomena with vorticity-dominated motion or highly distorted 
patterns of the medium because the computational mesh, which moves with the material, becomes too 
distorted and eventually gets tangled. On the other hand, Eulerian schemes are the natural 
choice for simulating complex flows but also provide a framework for modeling interfaces of 
arbitrary complexity and topology changes if equipped with some interface-capturing techniques, 
e.g. \cite{TavelliELDI,Hank2017,Barton2019,Jackson2019,Kemm2020,Busto2020,Barton2020}. 
Nevertheless, due to the diffusion intrinsically embedded into any Eulerian scheme, the errors in 
the interface approximation accumulate over time and Lagrangian schemes are preferable when the 
accuracy of the interface location is critical, e.g. multi-material computations or plasma flows 
for inertial confinement fusion applications \cite{MaireMM1,MaireMM2}. 

The first Lagrangian schemes were proposed in \cite{Neumann1950} using a formulation of the governing 
equations in primitive variables, which was also employed later in \cite{Benson1992,Caramana1998}. 
However, most of the modern Lagrangian finite volume schemes use the conservation form of the 
equations based on the physically conserved quantities like mass, momentum and total energy in 
order to compute shock waves
properly, see e.g. \cite{munz94,Smith1999,phm109,Despres2009,Morgan2019c}. The original approach 
for Lagrangian algorithms makes use of a staggered grid discretization 
\cite{wil1,StagLag,LoubereSedov3D}, where the velocity is defined at the cell interfaces and the 
other variables at the cell centers. Staggered schemes are compatible with the second law of 
thermodynamics thanks to the addition of an artificial viscosity in the internal energy equation 
allowing the dissipation of kinetic energy into internal energy through shock waves, thus total 
energy is not conserved. Although strategies to overcome this drawback have been proposed in 
\cite{Caramana1998,Despres_book_2017}, the development of cell-centered finite volume (FV) 
Lagrangian schemes dedicated to solve the hydrodynamics equations has started to gain visibility 
from the seminal works presented in \cite{DepresMazeran2003,Despres2005,Maire2007}. As a purely 
Lagrangian formulation, no mass flux is allowed across cell interfaces,
and the cells move and deform by considering the effects of the entire neighborhood. Specifically, 
these Lagrangian methods are based on the definition of a nodal solver which  
takes into account multiple one-dimensional Riemann problems occurring across all surrounding cells, 
and eventually uniquely determines the node velocity. Consequently, compatibility 
between the motion of the cell, thus its volume, and 
the cell deformation computed through PDEs involving the discretization of 
the velocity divergence is ensured. This is the so-called Geometrical Conservation Law (GCL) that is
nothing but the discrete Gauss theorem that is satisfied by construction in classical Godunov-type 
finite volume schemes on fixed grids.

Second order extension in space is typically achieved using a piecewise linear reconstruction of the 
conserved variables, while second order time stepping can be obtained relying on a classical 
Runge-Kutta scheme \cite{phmbn09,3DHydroMaire}, a Generalized Riemann problem 
methodology \cite{phm109,CCL2020} following the seminal ideas outlined 
in \cite{Artzi,Raviart.GRP.1,Raviart.GRP.2}, or the ADER strategy \cite{LAM2018}, originally 
proposed in \cite{Toro:2006a,toro3,titarevtoro} for fixed meshes. Cell-centered finite volume 
schemes have been extended to 2D and 3D Lagrangian solid mechanics in \cite{Maire_elasto,Kluth10,CCL2020}, 
while curved meshes are considered in \cite{Depres2012,chengshu5}. 

Higher order of accuracy in space was first achieved in 
\cite{chengshu1,chengshu2,chengshu3,chengshu4}, where a third order accurate  essentially 
non-oscillatory (ENO) reconstruction operator is introduced into a Godunov-type Lagrangian finite 
volume scheme. The mesh velocity is simply computed as the arithmetic average of the 
corner-extrapolated values in the cells adjacent to a mesh vertex and the numerical fluxes across 
element interfaces are solved at the aid of approximate Riemann solvers, thus these schemes cannot 
be regarded as Lagrangian methods \textit{sensu stricto}, since mass flux can in 
principle take place across cells. In the finite element framework, high order Lagrangian schemes 
have been developed in \cite{scovazzi1,scovazzi2} and also in \cite{Dobrev1,Dobrev2,Dobrev3}, who 
solved the equations for Lagrangian hydrodynamics using high order curvilinear finite element 
methods. Displacement-based finite element schemes for simulating engineering large strain 
transient situations can be found in \cite{GOUDREAU1982,Flanaghan_Belytchko_81}. Discontinuous 
Galerkin (DG) methods on polygonal unstructured meshes have been proposed for hydrodynamics in 
\cite{Vilar2,Vilar3}, while the equations for solid mechanics in hypo- and hyperelastic 
formulations have been solved with Lagrangian DG 
schemes in \cite{Morgan2019a,Morgan2019b}.

Lagrangian methods also may differ concerning the formulation of the governing equations. 
More precisely, there do exist total and updated Lagrangian schemes. \emph{Total Lagrangian} schemes 
rely on 
the discretization of the time rate of change of the deformation gradient, its determinant and its 
co-factor and the computations are performed on a fixed mesh. Contrarily, \emph{updated Lagrangian} 
methods physically move the computational mesh and the 
deformation gradient is the Lagrange-Euler mapping which relates the initial and the 
current mesh location. The solid dynamics equations written 
under total Lagrangian formulation have been successfully solved in
\cite{GodPesh2010,Aguirre2014,Bonet2015,Gil2016,Bonet2021,Gil2D_2014,Haider_2018}, where a 
cell-centered finite volume computational framework was employed. Hypoelastic solid mechanics 
models 
\cite{Truesdell55,Bernstein60} have 
also been numerically solved with Lagrangian methods, see e.g.
\cite{wil1,Gavriluk08,Maire_elasto,Sambasivan2013,cheng_jia_jiang_toro_yu_2017}. However, these 
models are not compatible with the second law of thermodynamics, contrarily to the hyperelasticity 
equations. A parallel discussion about hypo- and hyper-elastic models and their resolution can be 
found for instance in \cite{HyperHypo2019}. 

From the physical viewpoint, all contributions listed so far deal with either fluid \emph{or} solid 
mechanics. In the field of fluid mechanics, most of the existing Lagrangian methods are designed 
for the solution of hydrodynamics equations, i.e. ideal fluids. On the other hand, elastic or 
nearly incompressible solids are mainly addressed by moving mesh schemes. An attempt to derive a 
unified formulation of continuum mechanics in first-order hyperbolic form that includes fluid 
mechanics as well as solid mechanics has been very recently proposed in \cite{HPR2016}, referred to 
as Godunov-Peshkov-Romenski (GPR) model. This model makes use of relaxation-type PDEs which govern 
the deformation of the material. The relaxation parameter in the sources determines the 
\emph{inelasticity} time scale of the phenomena under consideration, thus permitting to recover 
either 
ideal fluids or 
ideal elastic solids in the model limits. Furthermore, because of the continuous transition from 
stiff to non-stiff time scale, the GPR model is also capable of capturing the behavior of viscous 
fluids 
as well as elasto-plastic solids. The numerical solution of this model has been first proposed on 
fixed meshes in \cite{DPRZ2016}, and consequently extended to plasma flows in \cite{DPRZ2017}. In 
\cite{LagrangeHPR}, a high order moving mesh scheme is adopted for the solution of the GPR model, 
which is based on the direct Arbitrary-Lagrangian-Eulerian (ALE) framework introduced in \cite{BosARCME} 
and references therein.

In this work, we aim at developing a cell-centered updated Lagrangian finite volume scheme for the 
solution of the GPR model. We consider unstructured meshes, composed of simplex 
control volumes, both in 2D an 3D. The novel algorithm arises from the previous contribution for hydrodynamics 
\cite{LAM2018} and the recent work on Lagrangian hyperelasticity \cite{Boscheri2021}. Moreover, the 
numerical schemes presented in the aforementioned references are exactly recovered by the new 
finite volume method in the limits of the GPR model. Indeed, the asymptotic preserving property of the 
fully discrete scheme is derived to demonstrate the consistency with already existing and 
well-established Lagrangian schemes, e.g. the EUCCLHYD scheme \cite{Maire2007}. Finally, the 
compatible discretization of the GCL allows the asymptotic limit of the distortion tensor to be 
consistently retrieved without resorting to staggered meshes as done in \cite{SIGPR2021}. To the 
knowledge of the authors, this is the first contribution towards the solution of relaxation-type 
systems of hyperbolic conservation laws in the updated Lagrangian framework. Second order extension 
in space and time is devised as well. In order to treat the \emph{arbitrarily stiff} source
terms with which the unified model is endowed, while maintaining the asymptotic properties of the model at the 
discrete level, we employ a semi-analytical time integration scheme 
for the stiff source terms \cite{chiocchettimueller,tavellicrack}, 
which yields a second order discretization of the fully coupled system of PDEs thanks to 
the implicit-explicit coupling obtained at the aid of the class of IMEX Runge-Kutta time stepping techniques 
\cite{AscRuuSpi,BosRus,PR_IMEX,BOSCHERI2021110206}.

The rest of this article is structured as follows. In Section \ref{sec.GPR}, we present the 
governing equations of the GPR model. Section \ref{sec.numscheme} is devoted to detail the 
Lagrangian finite volume scheme, the ODE solver for stiff relaxation sources and the extension 
of the algorithm to second order of accuracy. The asymptotic 
analysis of the fully discrete numerical method is developed in Section \ref{sec.AP}. In Section 
\ref{sec.test}, we show numerical 
convergence studies and a wide range of test cases on 2D and 3D unstructured moving meshes for different 
types of materials and different equations of state. Finally, in Section \ref{sec.concl} we give some 
concluding remarks and an outlook to future research and developments.

\section{Governing equations: the GPR model for continuum mechanics}	\label{sec.GPR}

The unified first-order hyperbolic model of continuum mechanics discussed in this section can 
describe fluid flows and 
deformations of solids in a single system of first-order hyperbolic equations. The model 
originates from the works by Godunov and Romenski in 1970s 
\cite{GodRom1972,God1978,Romenski1979,GodRom2003} 
on the modeling of large elasto-plastic deformations of solids. Similar to the work by Wilkins 
\cite{wil1}, the important feature of the model is that the description is made in an 
Eulerian frame of reference in contrast to the conventional Lagrangian description of solid 
mechanics. Later, in the work by Peshkov and Romenski \cite{HPR2016}, this allows for further 
generalization of the Godunov-Romenski model (which would 
be impossible under the Lagrangian description) towards incorporating the description of viscous fluids 
as a particular case of \emph{``extremely''} inelastic deformations of solids, see also 
\cite{DPRZ2016,HYP2016,Busto2020,SIGPR2021}. Therefore, in this 
paper, we shall refer to this model as the Godunov-Peshkov-Romenski (GPR) model.

In contrast to the conventional Navier-Stokes-based fluid mechanics, in the GPR model, fluids are 
characterized 
not with a strain-rate measure like in the Newton law of viscosity, but with a deformation 
measure and a special procedure of relaxation of tangential stresses similar to the Maxwell model 
of viscoelasticity. Nevertheless, it can be shown via a formal asymptotic analysis that the 
Navier-Stokes Cauchy stress tensor is recovered in the relaxation limit of the GPR 
model \cite{DPRZ2016,SIGPR2021}.  Thus, according to the GPR model, each material element of an 
isotropic continuum is 
defined by an infinitesimal frame $ \Ae $ (local basis triad) that characterizes
deformation and 
orientation of the material particles. The frame field $ \Ae $ is also called the \emph{effective 
elastic distortion} \cite{GodRom2003}, or simply the distortion field.
The relaxed (stress-free) state is defined for each material element individually as a state when 
the frame field $ \Ae $ is orthonormal. In general, such a state may not be reachable globally for 
all the material elements simultaneously, hence resulting in residual stresses in the material\footnote{In 
conventional fluids with zero yield strength, the residual stresses cannot be generated because the 
relaxation process acts until all stresses vanish.}, see Section  \ref{ssec.Shell} and 
\ref{ssec.TaylorBar2D}.

The GPR model can be also seen as a finite-strain elasto-plasticity model of \emph{hyperelastic} 
type 
\cite{HyperHypo2019}. Thus, 
according to the traditional finite-strain elasto-plasticity framework, the total deformation 
gradient $ \Ft $ can be decomposed as $ \Ft = \Fe \Fp $ into the elastic $ \Fe $ and 
plastic $ \Fp $ part. The natural Lagrangian stress, i.e. the first Piola-Kirchhoff stress, is then 
computed as $ \tensor{\Pi} = \pdd{E}{\Ft} $, where $ E $ is the total energy potential, which, in 
fact, must be a function of only the elastic part $ \Fe $ of the deformation gradient, that is $ E=E(\Fe) 
$. Hence, $ \tensor{\Pi} = \pdd{E}{\Ft} = \pdd{E}{\Fe}\Fp^{-\transpose}$, and therefore, to compute 
the stress tensor in the Lagrangian frame, one needs to know any two of the three quantities $ \Ft $, 
$ \Fe $, or $ \Fp $. On the other hand, in the Eulerian frame, the natural Eulerian stress tensor 
is the Cauchy stress $ \tensor{T} = \rho \Ft \pdd{E}{\Ft}^\transpose $, that can be written as
\begin{equation}
	\tensor{T} = \rho \Ft \pdd{E}{\Ft}^\transpose = \rho \Ft \Fp^{-1}\pdd{E}{\Fe}^\transpose
	=
	\rho \Fe \pdd{E}{\Fe}^\transpose.
\end{equation}
Consequently, to compute the stress tensor in the Eulerian frame, one needs to know only the elastic part $ 
\Fe $ of the total gradient $ \Ft $. In other words, we do not need to compute the entire history 
of inelastic deformations stored in $ \Fp $. This opens the great possibility to describe 
\emph{arbitrarily} 
large 
inelastic deformations typical of fluid-type motions. The distortion field $ \Ae $ can also be seen as 
the inverse of $ \Fe $, thus $ \Ae=\Fe^{-1} $.

Apart of being a deformation measure, the frame field $ \Ae $ also encodes the orientational 
degrees of freedom of the material elements, e.g. see the computational results in 
\cite{DPRZ2016,SIGPR2021,nonNewtonian2021}. The orientational degrees of freedom might be relevant 
for turbulence 
modeling \cite{Torsion2019} or anisotropic plasticity models \cite{Rubin2019}.  However, in this 
paper, we ignore these 
degrees of freedom and instead of the frame field $ \Ae $ we consider the metric tensor $ 
\Ge=\Ae^\transpose\Ae $. However, in the context of structure preserving numerical methods, the 
formulation in terms of $ \Ae $ is more preferable because the deformation compatibility condition 
$ \nabla\times\Ae = 0 $ is linear while, in terms of $ \Ge $, the compatibility condition 
(vanishing 
Riemann curvature tensor) is 
nonlinear, e.g. see \cite{GodRom2003}. In addition, we also take into account heat conduction 
effects via a 
first-order hyperbolic formulation proposed by Romenski in \cite{RomenskiMalyshev1987,Rom1989}, which 
characterizes the heat propagation via a relaxation equation for the thermal impulse $ \J $, 
e.g. see \cite{DPRZ2016}. In later papers \cite{SIGPR2021,PTRSA2020,SHTC-GENERIC-CMAT}, we however, 
use slightly 
different hyperbolic heat conduction model that coincides with the one used in this paper in the 
Fourier limit 
(local thermodynamic equilibrium).

For the updated Lagrangian scheme discussed in Section \ref{sec.numscheme}, we write the Eulerian equations 
\cite{DPRZ2016} in the co-moving frame of reference, that is the time derivative is the Lagrangian 
derivative (or material derivative), while the spatial derivatives are the Eulerian ones. The space 
is defined in $\mathbb{R}^d$ with $d=\{2,3\}$ representing the number of space dimensions. The time 
coordinate is given by $t$, while $\x=(x,y,z)$ denotes the spatial position vector. The system of 
governing equations for the unknowns $\Q:=\{ \svol, \vv, E, \J, \Ge \}$
therefore reads
\begin{subequations}
	\label{eqn.cl}
	\begin{align}
	&\rho \mddt{\svol}-\nabla \cdot \vv=0,\label{eqn.cl1}\\[2mm]
	&\rho \mddt{\vv}-\nabla \cdot \tensor{T}=\tensor{0},\label{eqn.cl2}\\[2mm]
	&\rho \mddt{E}-\nabla \cdot (\tensor{T}\vv) +\DIV{\q}=0, 
	\label{eqn.cl3}\\[2mm]
	&\rho \mddt{\J} + \del T =- \frac{\rho \H}{\sH},
	\label{eqn.cl4}\\[2mm]
	&\phantom{\rho} \mddt{\Ge} + \Ge \del\vv + \del\vv^\transpose \Ge = 
	\frac{2}{\rho \sG} \tensor{\sigma}, \label{eqn.cl5}
	\end{align}
\end{subequations}
where $ \rho $ is the mass density, $ \svol = \rho^{-1} $ is the specific volume, $ \vv=(u,v,w) $ is the 
velocity, $ \tensor{T}$ is the 
Cauchy stress 
tensor,
$ E(\rho,p,\vv,\Ge) $ is the total energy
\begin{subequations}\label{eqn.E}
	\begin{equation}
	E = E_h(\rho,p) + E_e(\Ge) + E_{th}(\J) + E_k(\vv),
	\end{equation}
	\begin{equation}\label{eqn.EE}
	E_h = \varepsilon(\rho,p), 
	\qquad
	E_e = \frac{\csh^2}{4}\Vert \devGe \Vert^2,
	\qquad
	E_{th} = \frac{1}{2}\alpha^2 \Vert\J\Vert^2,
	\qquad
	E_k = \frac{1}{2}\Vert\vv\Vert^2,
	\end{equation}
\end{subequations}
with $ E_h = \varepsilon(\rho,p) $ being the pure hydrodynamics part, $ p $ being the hydrodynamics 
pressure defined only from the hydrodynamics equation of state 
$ \varepsilon(\rho,p) $, $ E_e $ being the elastic energy stored in the 
deformed material elements, with $ \csh $ being the shear sound speed which characterizes the 
rigidity of the material elements, and $ \devGe = \Ge - \frac{1}{3}\tr(\Ge)\tensor{I} $ being the 
deviatoric part of 
the metric tensor $ \Ge $, and  $ \Id $ is the identity tensor. 
The third term $ E_{th} $ represents the energy carried by the thermal impulse $ \J $ with $ 
\alpha $ characterizing the velocity of propagation of thermal perturbations (defined later in 
Section \ref{ssec.heat}). The last term $ E_k $ denotes the standard kinetic energy of 
the continuous medium. The source terms in \eqref{eqn.cl4} and \eqref{eqn.cl5} are the algebraic 
relaxation sources defined in Section \ref{ssec.irrevers} and \ref{ssec.heat}.
The temperature $ T $ in \eqref{eqn.cl4} is not a state variable and must be defined from the 
equations of state $ \varepsilon(\rho,p) $, $ \q $ is the heat flux defined later in 
Section \ref{ssec.heat}. The Cauchy stress $ \tensor{T} $ is split into two parts:
\begin{equation}\label{eqn.T}
\tensor{T} = -p \Id + \tensor{\sigma}, 
\qquad
\tensor{\sigma} = -2\rho\Ge \pdd{E}{\Ge} = -\rho \csh^2 \Ge\devGe,
\end{equation} 
where the pressure $ p $ is defined from the hydrodynamics energy $ E_h(\rho,p) $ and the tangential (or viscous)
stress $ \tensor{\sigma} $ results from the definition of the elastic energy $ E_e(\Ge) $. 
In all computational examples for viscous fluids and elasto-plastic solids in 
Section \ref{sec.test}, we use $ 
E_e $ in the form 
\eqref{eqn.EE}, simply because the asymptotic analysis for the Navier-Stokes limit in \cite{DPRZ2016} was 
performed for this particular choice and we aim at demonstrating the asymptotic property for the fully discrete scheme. Nevertheless, other elastic energies could be used.
We remark that the spherical part of $ \tensor{\sigma} $ for our choice of the elastic energy is not zero 
but scales as $ \sim \Vert\devGe\Vert^2 $, i.e. quadratically in $ \devGe  $. For example, for 
fluid flows it is negligibly small, while for hyperelastic solids it might be not small and the 
total pressure should be defined as $ P = -\frac{1}{3} \tr(\tensor{T}) = p - 
\frac{1}{3}\tr(\tensor{\sigma})$.

\subsection{Equation of state for $ E_h $}	\label{ssec.hydroEOS}

In the numerical examples in Section \ref{sec.test}, we use three equations of state (EOS) for the 
hydrodynamics energy $ E_h $. 
For gases, we use the ideal gas EOS: 
\begin{equation}\label{eqn.IG}
	\varepsilon(\rho,p) = \frac{p}{\rho(\gamma-1)}, \qquad T = \frac{\varepsilon}{c_v}, \qquad 
	c_0^2 = \frac{\gamma p}{\rho},
\end{equation}
where $ \gamma = c_p/c_v $ is the ratio of specific heats, $ c_p $ is the specific heat at 
constant pressure, $ c_v $ is the specific heat at constant volume, and $ c_0 $ is the adiabatic 
sound speed.

For liquids and solids, one can use the Mie-Grüneisen equations of state:
\begin{equation}\label{eqn.MG}
	\varepsilon(\rho,p) = \frac{p - \rho_0 c_0^2 f(J)}{\rho_0 \Gamma_0},
	\qquad
	f(J) = \frac{(J-1)(J-\frac12\Gamma_0(J-1))}{(J-s(J-1))^2}, \qquad J =
	\frac{\rho}{\rho_0},
\end{equation}
where $ c_0 = const $ is the adiabatic sound speed, $ \Gamma_0 = const$, $ s = const $ are material 
constants, and $ \rho_0 $ is the reference mass density.

The third type of hydrodynamics EOS used in Section~\ref{ssec.BePlate} and 
\ref{ssec.twcol} is the volumetric part of the Neo-Hookean hyperelastic EOS that reads
\begin{equation}\label{eqn.NH}
	\varepsilon(\rho,p) = \frac{\Gs}{4\rho_0} \left((J-1)^2 +(\log(J))^2 \right) ,
	\qquad
	p = -\frac{\Gs}{2}\left (J-1 + \frac{\log(J)}{J}\right ),
	\qquad
	J = \frac{\rho}{\rho_0},
\end{equation}
where $ \Gs =\rho_0 \csh^2 $ is the shear modulus, which also corresponds to the second Lam{\'e} 
coefficient. Material mechanical properties are often described in terms of Young modulus $Y$ and 
Poisson ratio $\nu$. These parameters are linked as follows:
\begin{equation}
	\label{eq:Ynu_relation}
	\Gs = \Frac{Y}{2 \, (1+\nu)}.
\end{equation}
In this case the adiabatic sound speed can be computed as the square root of the bulk modulus $K$, hence
\begin{equation}
	\label{eqn.c0NH}
	c_0^2 = K = \frac{Y\, \nu}{(1+\nu)(1-2\nu)} + \frac{2\Gs}{3}.
\end{equation}

\subsection{Closure for inelastic deformations and fluid flows}	\label{ssec.irrevers}

The source terms in \eqref{eqn.cl4} and \eqref{eqn.cl5} describe the dissipative dynamics of the 
system and must be consistent with the second law of thermodynamics (entropy must be non-decreasing).
It can be shown that the specification of the source terms used in this paper guarantees that the 
entropy is indeed not decreasing, see \cite{DPRZ2016,SHTC-GENERIC-CMAT}.
Without providing further details (which can be found in \cite{DPRZ2016}) we specify them as 
follows. 

The source term in \eqref{eqn.cl4} handles inelastic deformations of the continuum with the 
relaxation function $ \Theta $ given by
\begin{equation}\label{eqn.theta}
\Theta = \tau_1 \frac{\csh^2}{3} |\Ge|^{-5/6}, 
\end{equation}
where $ |\Ge| = \det(\Ge)$ and $ \tau_1 $ is the strain relaxation time which, in 
general, is a function of state variables $ \tau_1 = \tau_1(\rho,T,\Ge) $. In particular, 
$ \tau_1 =const $ for Newtonian fluids, while $ \tau_1 \to \infty $ for pure elastic solids. 
The Navier-Stokes stress tensor is recovered in the relaxation limit of system \eqref{eqn.cl}, that is $ 
\epsilon = \tau_1 /t_{macro} \ll 1 $ at first order in $ \epsilon $, with $ t_{macro} $ being the macroscopic characteristic time 
of the physical process. The \emph{effective} shear 
viscosity 
is then given by \cite{DPRZ2016}
\begin{equation}\label{eqn.visc}
	\visc = \frac{1}{6}\rho_0 \tau_1 \csh^2.
\end{equation}

For elasto-plastic solids $ \tau_1 $ can be taken as \cite{BartonRom2010,HyperHypo2019,LagrangeHPR}
\begin{equation}\label{eqn.tau.plast}
	\tau_1 = \tau_{10} \left( \frac{\sigma_Y}{\sigma}\right)^n, 
	\qquad 
	\sigma = \sqrt{\frac{3}{2}\text{tr}(\devS^2)},
	\qquad
	\devS=\tensor{\sigma} - \frac13\text{tr}(\tensor{\sigma})\tensor{I},
\end{equation}
where $ \sigma_Y $ is the material yield strength under quasi-static loading conditions and $\tau_{10}=const$ is a material-specific constant.
The parameter $ n $ controls the rate-dependence of the plasticity mode, e.g. the greater $ n $ the 
more the material is rate-independent.
For other examples of $ \tau_1 $ suitable for the modeling of non-Newtonian fluids see 
\cite{nonNewtonian2021,Jackson2019a}.

\subsection{Heat conduction}	\label{ssec.heat}

The closure of the heat conducting part of system \eqref{eqn.cl} is also conditioned by the 
consistency with the second law of thermodynamics \cite{DPRZ2016,SHTC-GENERIC-CMAT} as well as with 
the Fourier law of heat conduction. Thus, the relaxation function in the right hand-side of 
\eqref{eqn.cl4} is given by
\begin{equation}
\Psi = \alpha^2 \tau_{20} \tau_2, 
\qquad
\tau_{20} = \frac{\rho}{\rho_0}\frac{T_0}{T}.
\end{equation}
The consistency with the first law of thermodynamics (total energy conservation) requires that the 
heat flux $ \q $ is defined as (see \cite{DPRZ2016})
\begin{equation}\label{eqn.q}
	\q = \alpha^2 T \J.
\end{equation}
At the mechanical and thermodynamical equilibrium, the characteristic velocity of the thermal 
perturbation propagation is 
\begin{equation}
	\ch^2 = \frac{\alpha^2}{\rho_0^2} \frac{T}{c_v}.
\end{equation}
The Fourier law of heat conduction $ \q = -\kappa \nabla T $ is recovered in the relaxation limit 
(local equilibrium) with $ \tau_2 \ll t_{macro}$,
at first order in $ \tau_2/t_{macro} $. The \emph{effective} heat conductivity is given by \cite{DPRZ2016}
\begin{equation}\label{eqn.kappa}
	\kappa = \tau_2 \alpha^2 \frac{T_0}{\rho_0}.
\end{equation}

\section{Lagrangian finite volume scheme on unstructured meshes}	\label{sec.numscheme}

\subsection{Discretization of the space-time computational domain}
We start by providing the details about the discretization of the time interval and the spatial computational domain.

\paragraph{Time computational domain} The time coordinate is defined in the time interval 
$t\in[0,t_f]$, where $t_f\in\mathbb{R}_0^+$ represents the final time. A sequence of discrete points 
$t^n$ approximates the temporal computational domain such that $t \in [t^n;t^{n+1}]$,
\begin{equation}
	t = t^n + \alpha \dt^n, \qquad  \alpha \in [0,1], 
	\label{eqn:time}
\end{equation} 
with $t^n$ and $\dt^n$ denoting the current time and time step, respectively. As done 
in \cite{phm109,CCL2020}, the time step is limited by a classical CFL stability condition 
combined with a criterion on the growth rate of the time step, hence
\begin{equation}
	\dt^n = \min \left( \textnormal{CFL} \cdot \min \limits_{T_i} \frac{h_i}{a_i}, \, C_m \dt^{n-1} 
	\right). 
	\label{eqn.timestep}
\end{equation}
The CFL number is chosen to be $\textnormal{CFL} \leq 1/d$ as usual on unstructured meshes, with $ 
d $ being the number of space dimensions, while $h_i=|T_i|^{1/d}$ is the 
characteristic mesh size with $|T_i|$ being the volume of the spatial element $T_i$. Finally, $a_i$ 
is the maximum eigenvalue estimate (equilibrium estimate) of the governing PDE given 
by
\begin{equation}
	a_i = \left. \sqrt{ c_0^2 + \frac{4}{3}\csh^2 + \ch^2} \, \right|_i .
	\label{eqn.a}
\end{equation}
The criterion for controlling the growth rate of the time step is set with $C_m=1.1$, so that 
the increase of the time step size is mild without sharp discontinuities.

\paragraph{Space computational domain} In the updated Lagrangian framework, the time-dependent 
computational domain $\Omega(t)$ is discretized at the current time $t^n$ by a set of 
non-overlapping unstructured control volumes $T^n_i$ with boundary $\partial T_i^n$, that are given 
by triangles ($d=2$) or tetrahedra ($d=3$). $N_E$ denotes the total number of elements contained in 
the domain and the union of all elements is called the \textit{current tessellation} 
$\mathcal{D}^n_{\Omega}$ of the domain 
\begin{equation}
	\mathcal{D}^n_{\Omega} = \bigcup \limits_{i=1}^{N_E}{T^n_i}. 
	\label{trian}
\end{equation}	
To ease notation, we omit the superscript $n$ in the following description, meaning that the spatial domain is discretized as follows at each time level $t^n$. The index $i$ refers to the element $T_i$, $f$ addresses a face (which is of dimension $d-1$) and $r$ denotes a node. The double index $ri$ is used to express a quantity defined at node $r$ from element $T_i$; likewise, $fi$ is used to express a quantity defined at the barycenter of face $f$ from element $T_i$. The Neumann neighbor of element $T_i$ which shares face $f$ is denoted with $T_j$, and the outward pointing unit normal vector is $\n_{fi}$. Finally, the cell volume is $|T_i|$ and the length (2D) or surface (3D) of a face $f$ is addressed with $s_f$.

\begin{figure}[!bp]
    \centering
	\includegraphics[height=6.5cm]{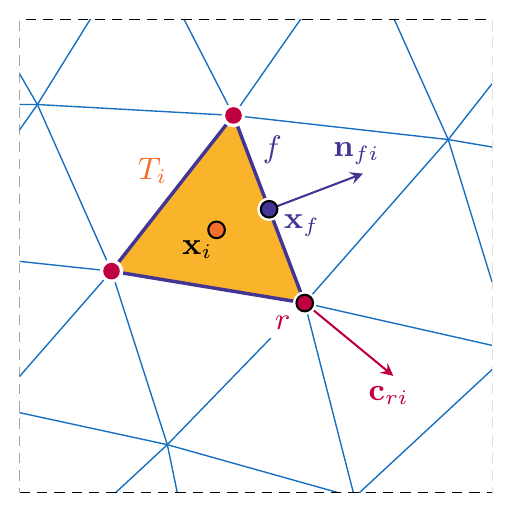}%
	\qquad\qquad
	\includegraphics[height=6.5cm]{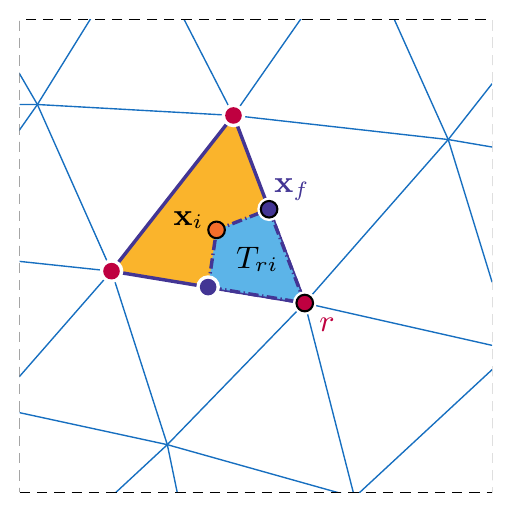}%
	\caption{Notation for the mesh. Left: cell barycenter $\x_i$ of element $T_i$, 
	    face barycenter $\x_f$ and normal vector $\n_{fi}$ relative to face $f$, 
	    corner vector $\cv_{ri}$ relative to node $r$. Right: subcell $T_{ri}$ associated with corner $r$ of element $T_i$.}
	\label{fig.mesh_notation}
\end{figure}

At the aid of {Figure \ref{fig.mesh_notation}}, the following notation is adopted for the definition of different topological sets:
\begin{itemize}
	\item $\mathcal{R}_i$ is the set of nodes $r$ of element $T_i$;
	\item $\mathcal{R}_f$ is the set of nodes $r$ of face $f$;
	\item $\mathcal{F}_i$ is the set of faces $f$ of element $T_i$;
	\item $\mathcal{F}_r$ is the set of faces $f$ sharing node $r$;
	\item $\mathcal{F}_{ri}$ is the set of faces $f$ sharing node $r$ and belonging to element $T_i$;
	\item $\mathcal{T}_{r}$ is the set of elements $T$ sharing node $r$, i.e. the Voronoi neighborhood of node $r$.
\end{itemize} 
Each cell can be divided into subcells $T_{ri}$, that are obtained by connecting the cell centroid with the barycenter of the faces belonging to $\mathcal{F}_{ri}$. The cell centroid $\x_i$ is defined as
\begin{equation}
	\x_i = \frac{1}{|T_i|} \int \limits_{T_i} \x \, \rmd \x,
\end{equation}
while the barycenter of a face  $\x_f$ is given by
\begin{equation}
	\x_f = \frac{1}{N_{rf}} \sum \limits_{r \in \mathcal{R}_f} \x_r,
\end{equation}
with $N_{rf}$ representing the total number of nodes which share face $f$, i.e. all nodes belonging to $\mathcal{R}_f$. Let us now introduce the corner vector $\cv_{ri}$ \cite{Despres2005,Maire2007}, which is a geometric object defined at node $r$ within element $T_i$:
\begin{equation}
	\cv_{ri} = \frac{1}{d} \sum \limits_{f \in \mathcal{F}_{ri}} s_f \n_{fi}.
\end{equation}
These corner vectors are a linear combination of the oriented surfaces $s_f \n_{fi}$, which assume a linear velocity field over the face so that they remain planar. Therefore, the corner vectors $\cv_{ri}$ provide a consistent discretization of the cell boundary $\partial T_i$ so that
\begin{equation}
	\sum \limits_{r \in \mathcal{R}_{i}} \cv_{ri} = 0,
\end{equation}
which is a second order discretization of Gauss theorem over the control volume $T_i$.

\begin{figure}[!bp]
    \centering
	\includegraphics[height=6.5cm]{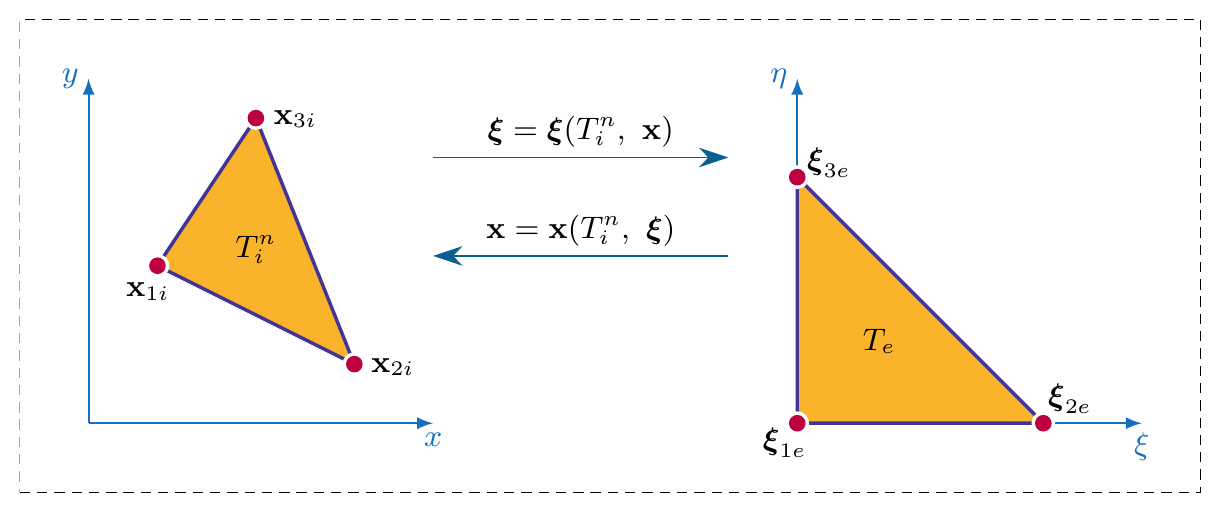}%
	\caption{Mapping from the physical element $T_i^n$ in the global coordinate system $\x = (x,\ y)$ to the reference element $T_e$ 
	   in the reference coordinate system $\boldsymbol{\xi} = (\xi,\ \eta)$.}
	\label{fig.refSystem}
\end{figure}

Each control volume defined in the physical space $\x$ can be mapped onto a reference element $T_e$ in the reference coordinate system $\boldsymbol{\xi}=(\xi,\eta,\zeta)$, see {Figure \ref{fig.refSystem}}. The reference element is the unit triangle in 2D with vertices $\boldsymbol{\xi}_{1e}=(\xi_{1e},\eta_{1e})=(0,0)$, $\boldsymbol{\xi}_{2e}=(\xi_{2e},\eta_{2e})=(1,0)$ and $\boldsymbol{\xi}_{3e}=(\xi_{3e},\eta_{3e})=(0,1)$, or the unit tetrahedron in 3D with vertices $\boldsymbol{\xi}_{1e}=(\xi_{1e},\eta_{1e},\zeta_{1e})=(0,0,0)$, $\boldsymbol{\xi}_{2e}=(\xi_{2e},\eta_{2e},\zeta_{2e})=(1,0,0)$, $\boldsymbol{\xi}_{3e}=(\xi_{3e},\eta_{3e},\zeta_{3e})=(0,1,0)$ and $\boldsymbol{\xi}_{4e}=(\xi_{4e},\eta_{4e},\zeta_{4e})=(0,0,1)$. The spatial mapping in 3D reads
\begin{equation} 
	\mathbf{x} = \mathbf{x}_{1i} + 
	\left( \mathbf{x}_{2i} - \mathbf{x}_{1i} \right) \xi + 
	\left( \mathbf{x}_{3i} - \mathbf{x}_{1i} \right) \eta + 
	\left( \mathbf{x}_{4i} - \mathbf{x}_{1i} \right) \zeta,
	\label{eqn.xi} 
\end{equation}
where $\mathbf{x}_{ki} = (x_{ki},y_{ki},z_{ki})$ represents the vector of physical spatial 
coordinates of the $k$-th vertex of element $T_i$ for $k=\{1,2,3,4\}$. Notice that the mapping 
\eqref{eqn.xi} corresponds to the linear isoparametric approximation of the Lagrange-Euler 
transformation between the deformed (physical) and undeformed (reference) configuration.

%

\subsection{Integral form of the governing equations}
Let $\Omega(t) \subset \mathbb{R}^d$ denote the time-dependent control volume in $d=\{2,3\}$ space dimensions, and $\partial \Omega(t)\subset \mathbb{R}^{d-1}$ its surface with $\n$ being the outward pointing unit normal vector. The time coordinate is $t \in \mathbb{R}_0^+$ and the physical space is defined by the position vector $\x=(x,y,z)$. To derive a conservative finite volume discretization, the governing equations \eqref{eqn.icl} must be integrated over the control volume $\Omega(t)$. Mass conservation directly follows from the Lagrangian description of continuum mechanics and holds true by construction, that is
\begin{equation}
	\mddt{} \int \limits_{\Omega(t)} \rho \, \rmd \x= 0.
	\label{eqn.mass_cons}
\end{equation}
After addressing a generic physical variable with $\phi=\phi(\x,t)$, let us introduce the transport relation
\begin{equation}
	\mddt{} \int \limits_{\Omega(t)} \rho \phi \, \rmd \x = \int \limits_{\Omega(t)} \rho 
	\mddt{\phi} \, \rmd \x,
\end{equation}
which is derived from the Reynolds transport theorem \cite{ContMechBook} using the mass conservation property \eqref{eqn.mass_cons}. The control volume formulation of the PDE system \eqref{eqn.icl} can then be expressed as follows:
\begin{subequations}
	\label{eqn.icl}
	\begin{align}
		&\mddt{} \int \limits_{\Omega(t)} \, \rmd \x - \int \limits_{\partial \Omega(t)} \vv \cdot 
		\n 
		\, \rmd s = 0,\label{eqn.icl1}\\[2mm]
		&\mddt{} \int \limits_{\Omega(t)} \rho \vv \, \rmd \x -\int \limits_{\partial \Omega(t)} 
		\tensor{T} \n \, \rmd s = \mathbf{0}, \label{eqn.icl2}\\[2mm]
		&\mddt{} \int \limits_{\Omega(t)} \rho E \, \rmd \x -\int \limits_{\partial \Omega(t)} 
		\left( 
		\tensor{T}\vv + \q \right) \cdot \n \, \rmd s=0, 
		\label{eqn.icl3}\\[2mm]
		&\mddt{} \int \limits_{\Omega(t)} \rho \J \, \rmd \x + \int \limits_{\partial \Omega(t)} T 
		\tensor{I} \, \n \, \rmd s = - \int \limits_{\Omega(t)} \frac{\rho \H}{\Psi} \, \rmd \x,
		\label{eqn.icl4}\\[2mm]
		&\int \limits_{\Omega(t)} \rho \mddt{\Ge} \, \rmd \x + \int \limits_{\Omega(t)} \rho \left( 
		\Ge \gradv(\vv) + \gradv(\vv)^\transpose \Ge \right) \, \rmd \x= \int 
		\limits_{\Omega(t)}
		\rho \frac{2}{\rho \sG} \tensor{\sigma} \, \rmd \x. \label{eqn.icl5}
	\end{align}
\end{subequations}
Equation \eqref{eqn.icl1} is also referred to as \textit{Geometric Conservation Law} (GCL) and it is equivalent to the local kinematic equation which governs the mesh motion:
\begin{equation}
	\mddt{\x} = \vv, \qquad \x(t=0)=\x_0,
	\label{eqn.trajODE}
\end{equation}
where $\x_0$ denotes the initial position at time $t=0$ of the material particle. The velocity 
gradient operator $\gradv(\vv)$ used in \eqref{eqn.icl5} is defined by applying the divergence 
theorem as
\begin{equation}
\gradv(\vv)=\del\vv=\frac{1}{|\Omega|} \int_{\partial \Omega} \vv \otimes \n\,\d s,
\label{eqn.Lop}
\end{equation}
which is the cornerstone for a compatible space discretization with the geometric conservation law 
\eqref{eqn.icl1}. 

The above integral formulation \eqref{eqn.icl} of the GPR model is suitable for the design of cell-centered finite volume methods.

\subsection{Fully discrete first order scheme}
Equation \eqref{eqn.mass_cons} implies that the mass of a generic control volume remains constant in time, that is
\begin{equation}
	m_i:=\int \limits_{T_i(t)} \rho \, \rmd \x,
\end{equation}
with $m_i$ denoting the mass of cell $T_i(t)$ over time. As usual for finite volume schemes, data are stored and evolved in time as piecewise constant cell averages. Therefore, for a generic variable $\phi(\x,t)$, its \textit{mass averaged value} $\phi_i$ over the element $T_i(t)$ is defined as
\begin{equation}
	\phi_i = \frac{1}{m_i} \int \limits_{T_i(t)} \rho \phi \, \rmd \x.
\end{equation}
Following \cite{Maire2011}, let us introduce the subcell force $\bm{f}_{ri}$, which is the traction force exerted on the outer boundary of the subcell $T_{ri}$ (see Figure \ref{fig.mesh_notation}). Specifically, the surface integral of the momentum equation \eqref{eqn.icl2} is split over the subcell boundaries, hence obtaining
\bea
\int \limits_{\partial T_i} \tensor{T} \n \, \rmd s = \sum \limits_{r \in \mathcal{R}_{i}} \, \int 
\limits_{\partial T_{ri} \cap \partial T_i} \tensor{T} \n \, \rmd s,
\eea
which yields the definition of the subcell force \cite{Maire2011,Maire_elasto} as 
\bea \label{eqn.fri}
\bm{f}_{ri} := \int \limits_{\partial T_{ri} \cap \partial T_i}   \tensor{T} \n \, \rmd s.
\eea
Using the definition of the subcell force and the mass average quantity introduced above, a fully discrete finite volume scheme for system \eqref{eqn.icl} writes
\begin{subequations}
	\label{eqn.fv}
	\begin{align}
		&\omega^{n+1}_i = \omega^{n}_i + \frac{\dt}{m_i} \sum \limits_{r \in \mathcal{R}_i} \tilde{\vv}^*_r \cdot \frac{1}{6} \left( \cv_{ri}^{n} + 4\cv_{ri}^{n+1/2} + \cv_{ri}^{n+1}\right),\label{eqn.fvcl1}\\[2mm]
		&\vv_i^{n+1} = \vv_i^{n} + \frac{\dt}{m_i} \sum \limits_{r \in \mathcal{R}_i} \tilde{\bm{f}}_{ri}^{*}, \label{eqn.fvcl2}\\[2mm]
		&E_i^{n+1} = E_i^n + \frac{\dt}{m_i} \left[ \sum \limits_{r \in \mathcal{R}_i} \tilde{\bm{f}}_{ri}^{*} \cdot \tilde{\vv}_r^* + \sum \limits_{f \in \mathcal{F}_i} \widehat{\q_{fi} \cdot \n_{fi}}^n \, s_{f}^n \right] =0, 
		\label{eqn.fvcl3}\\[2mm]
		& \J_i^{n+1} = \J_i^{n} - \frac{\dt}{m_i} \sum \limits_{f \in \mathcal{F}_i} 
		\widehat{T_{fi} \Id  \cdot \n_{fi}}^n \, s_{f}^n - \dt 
		\frac{\mathbf{H}_i^{n+1}}{\sH_i^{n+1}},
		\label{eqn.fvcl4}\\[2mm]
		&\Gei^{n+1} = \Gei^{n} - \dt \left( \Gei^n \, \gradv_i(\tilde{\vv}^*) + 
		\gradv_i(\tilde{\vv}^*)^{\transpose} \, \Gei^n \right) + \dt \frac{2 
		\tensor{\sigma}_i^{n+1}}{\rho \sG_i^{n+1}}, \label{eqn.fvcl5}
	\end{align}
\end{subequations}
which must be coupled with the discrete trajectory equation to move the grid nodes:
\begin{equation}
	\x_r^{n+1} = \x_r^n + \dt \, \tilde{\vv}_r^{*}.
	\label{eqn.trajR}
\end{equation} 
The vector of state variables is $\Q=\{\omega,\vv,E,\J,\Ge\}$, and these quantities are evolved in 
time with the discrete equations written above. Let us now examine more in details the finite 
volume scheme \eqref{eqn.fv}. 

\paragraph{Time discretization} All terms in \eqref{eqn.fv} are discretized explicitly, apart from 
the sources in the thermal impulse equation \eqref{eqn.fvcl4} and in the equation for the metric 
tensor \eqref{eqn.fvcl5}, because they might become stiff when approaching the limit cases of 
Navier-Stokes and Fourier. The stiff source terms are handled with an exponential integrator for 
ODE that will be presented in Section \ref{ssec.ExpInt}, which is applied to both the thermal 
impulse \eqref{eqn.fvcl4} and the metric tensor equations \eqref{eqn.fvcl5}. Consequently, an 
implicit discretization of $\tensor{G}_{e_i}^{n+1}$ implies that the stress tensor 
$\tensor{\sigma}_i^{n+1}$ is also taken implicitly. Therefore a formally \emph{implicit} treatment 
is 
considered for the total stresses in the momentum and energy equations 
\eqref{eqn.fvcl2} and \eqref{eqn.fvcl3}, embedded in the subcell force $\bm{f}_{ri}$ according to 
\eqref{eqn.fri}, as well as for the node velocity. These terms are marked with the asterisk 
superscript, i.e. $\tilde{\bm{f}}_{ri}^{*}$ and $\tilde{\vv}_r^*$. To overcome the resulting 
non-linearity in the stress tensor, i.e. 
$\tensor{\sigma}_i^{n+1}=\tensor{\sigma}_i^{n+1}(\tensor{G}_{e_i}^{n+1})$ according to 
\eqref{eqn.T}, an iterative Picard technique will be introduced later for the computation of the 
node velocity and the subcell force in the nodal solver. The tilde `$\tilde{\phantom{\vv}}$' and 
hat `$ \widehat{\phantom{\vv\vv}} $' in \eqref{eqn.fv} reference to two different types of 
numerical fluxes explained below.


According to \eqref{eqn.trajR}, the coordinate position $\x_r$ is a linear function of time. The corner vector $\cv_{ri}$ in \eqref{eqn.fvcl1} exhibits a linear or quadratic time dependency, in 2D or 3D, respectively. Therefore, the flux in the equation for the specific volume $\omega_i^{n+1}$ must be exactly integrated in time in order to ensure the satisfaction of the GCL at the discrete level. In \eqref{eqn.fvcl1}, the fourth order accurate Kepler quadrature rule is used, which might be replaced by a simple midpoint rule for 2D meshes. This ensures that the new density $\rho_i^{n+1}$ computed from the GCL \eqref{eqn.fvcl1} is equivalent to the density deduced from the mesh motion, i.e. $\rho_i^{n+1}=m_i/|T_i^{n+1}|$, see \cite{phm109,phmbn09} for further details.
	
\paragraph{Numerical fluxes} The finite volume scheme \eqref{eqn.fv} involves two different types of 
numerical fluxes for the discretization of the divergence operators, namely vertex-based and 
face-based fluxes. The computation of the nodal velocity is based on the nodal solver presented in 
Section \ref{ssec.NodalSol}, where multiple one-dimensional Riemann problems are simultaneously 
solved across all faces impinging to a node. Numerical dissipation is embedded into the nodal 
solver, thus the resulting node velocity $\tilde{\vv}_r^*$ is numerically stable. This velocity is 
subsequently used to compute the subcell force $\tilde{\bm{f}}_{ri}^{*}$, and the tilde symbol in 
\eqref{eqn.fv} indicates these vertex-based numerical fluxes. 

The second type of fluxes are face-based and those terms are marked with the hat symbol. A Rusanov-type numerical flux function \cite{Rusanov:1961a} is adopted, thus the heat flux in the energy equation \eqref{eqn.fvcl3} and the temperature gradient in the thermal impulse equation \eqref{eqn.fvcl4} are explicitly given by
\begin{eqnarray}
\widehat{\q_{fi} \cdot \n_{fi}} &=& \frac{1}{2} \left( (\alpha^2 T \J)_{fi} + (\alpha^2 T \J)_{fj} \right) \cdot \n_{fi} - \frac{1}{2} |\lambda_{f}| \left( E_{fj} - E_{fi} \right), \\
\widehat{T_{fi} \Id \cdot \n_{fi}} &=& \frac{1}{2} \left( (T \Id )_{fi} + (T \Id )_{fj} 
\right) \cdot \n_{fi} - \frac{1}{2} |\lambda_{f}| \left( \J_{fj} - \J_{fi} \right),
\end{eqnarray}	
where the subscripts ${fi}$ and ${fj}$ denote the value of the physical quantity computed at face $f$ from element $T_i$ and from the neighbor element $T_j$, respectively. The numerical dissipation is given in terms of the maximum eigenvalue estimate $\lambda_{f}$ in spatial normal direction between the elements sharing face $f$, which is evaluated as
\begin{equation}
	\lambda_{f} = \max (a_i,a_j),
	\label{eqn.lambda_max}
\end{equation}
with the wave speeds $(a_i,a_j)$ defined by \eqref{eqn.a}.
	
\paragraph{Non-conservative products} The last equation \eqref{eqn.fvcl5} contains non-conservative 
products between the metric tensor $\Ge$ and the velocity gradient $\gradv = \del\vv$. This 
equation 
can be also seen as a pure geometric 
relation and thus, according to 
\cite{CCL2020,Boscheri2021}, we do not consider the jump contribution across the cell boundaries of 
the 
non-conservative products as done in \cite{Lagrange3D}, but only the smooth part in the cell volume 
is taken into account. The discrete version of the velocity gradient \eqref{eqn.Lop} writes
\begin{equation}
	\gradv_{i}(\vv)=\frac{1}{|T_i|}\sum \limits_{r\in \mathcal{R}(i)} \vv_r \otimes \cv_{ri},
	\label{eqn.Loph}
\end{equation}
which is used to evaluate the products $\tensor{G}_{e_i}^n \, \gradv_i(\tilde{\vv}^*)$ and 
$\gradv_i(\tilde{\vv}^*)^{\transpose} \tensor{G}_{e_i}^n$ in \eqref{eqn.fvcl5}. Let us notice that 
the discrete gradient operator \eqref{eqn.Loph} is applied to the velocity field $\tilde{\vv}^*$, 
which results from the nodal solver.
	
\subsection{Nonlinear nodal solver}	\label{ssec.NodalSol}
The nodal fluxes are computed starting from the nodal solver algorithm proposed in a series of works \cite{Maire2007,Maire2011,Maire2011a,phmbn09,Despres2009}. The nodal solver aims at evaluating a unique velocity vector at each node $\tilde{\vv}_r^*$ by satisfying global conservation of momentum as well as compatibility with the second law of thermodynamics, meaning that the entropy variation in the cell is non-negative. 

The starting point is the knowledge of the state variables $\Q^n$ and the geometry $\x^n$ at time $t^{n}$. Next, the nodal velocity $\tilde{\vv}_r^*$ is obtained by solving the system
\bea
\label{eq:nsn}
\Mp \tilde{\vv}_{r}^* =  \sum \limits_{i \in \mathcal{T}_{r}} \Mcp \vv_i^n -\tensor{T}_i^* 
\cv_{ri}^n,
\eea
with the discrete subcell matrix $\Mcp$ and nodal matrix $\Mp$ given by
\bea
\label{eq:Mcp^n}
\Mcp = \sum \limits_{f \in \mathcal{F}_{ri}}  z_{i}^n \, s_f^n \, \n_{f}^n \otimes \n_{f}^n,  \qquad
\Mp =  \sum \limits_{i \in \mathcal{T}_{r}} \Mcp.
\eea
The numerical dissipation is introduced by the swept mass flux $z_{i}^n = \rho_i^n \, a_i^n$, with $a_i^n$ given by \eqref{eqn.a}. Notice that the above matrices are symmetric positive definite by construction, thus $\Mp$ in \eqref{eq:nsn} is always invertible. 
Once the node velocity is computed, the subcell force can easily be evaluated as
\bea
\label{eq:subforce^n}
\tilde{\bm{f}}_{ri}^* = \cv^n_{ri} \tensor{T}_i^* + \Mcp (\tilde{\vv}_{r}^*  -\vv_i^n ),
\eea
whereas the node coordinates are updated to the new position with the trajectory equation \eqref{eqn.trajR}.

The geometric quantities $\cv_{ri}^n$ and $\n_{f}^n$ are considered at the current time level $t^n$ 
as well as the cell quantities $\vv_i^n$ and $z_i^n$, which are evaluated at the cell centroid, 
i.e. $\vv_{r}(\x_i)^n$ and $z_i(\x_i)^n$. However, the total Cauchy stress $\tensor{T}_i^*$, and 
more 
specifically the tangential stress tensor $\tensor{\sigma}_i^*$, should take into account elastic, 
viscous 
or plasticity effects that might occur in the material, yielding a significant change in the node 
velocity. Therefore, a formally implicit discretization is used for the tangential stress tensor, 
that 
according to \eqref{eqn.T} leads to (for the energy potential \eqref{eqn.E})
\begin{equation}
	\tensor{\sigma}_i^{n+1} = -\rho_i^{n+1} \csh^2 \Gei^{n+1} \, \devGei^{n+1},
	\label{eqn.sigmah}
\end{equation}
so that the nodal solver \eqref{eq:nsn} must be coupled with the trajectory equation \eqref{eqn.trajR}, the GCL \eqref{eqn.fvcl1} and the equation for $\Gei$ \eqref{eqn.fvcl5}, in order to obtain $\x_r^{n+1}$, $\rho_i^{n+1}$ and $\Gei^{n+1}$, respectively. This choice wold lead to a strongly nonlinear system to be solved, whose convergence would become difficult to control. Consequently, we rely on a Picard iterative technique that has been already used in the context of all Mach flow solvers \cite{BDLTV2020,BDT_cns}. More precisely, let $l$ be the index for the Picard iteration, and let the viscous stress tensor be initialized at time $t^n$, i.e. $\tensor{\sigma}_i^{l,n+1}=\tensor{\sigma}_i^{n}$ for $l=0$. The nodal solver algorithm proceeds then iteratively as follows for $l=1,\ldots,\mathcal{L}$:
\begin{subequations}
	\label{eqn.nsPicard}
	\begin{align}
& \tilde{\vv}_{r}^{l+1,n+1} =  \left( \sum \limits_{i \in \mathcal{T}_{r}} \Mcp  \vv_i^n 
-\tensor{T}_i^{l,n+1} \cv_{ri}^n	\right) \Mp^{-1}, \label{eqn.nsP1}\\
& \x_r^{l+1,n+1} = \x_r^n + \dt \, \tilde{\vv}_{r}^{l+1,n+1}, \label{eqn.nsP2}\\
& \rho_i^{l+1,n+1} = \frac{m_i}{|T_i|^{l+1,n+1}}, \label{eqn.nsP3}\\
&\Gei^{l+1,n+1} = \Gei^{n} - \dt \left( \Gei^n \, \gradv_i(\tilde{\vv}^{l+1,n+1}) + 
\gradv_i(\tilde{\vv}^{l+1,n+1})^{\transpose} \, \Gei^n \right) + \dt \frac{2 
\tensor{\sigma}_i^{l+1,n+1}}{\rho^{l+1,n+1} \, \sG_i^{l+1,n+1}}, \label{eqn.nsP4}\\
& \tensor{\sigma}_i^{l+1,n+1} = -\rho_i^{l+1,n+1} \csh^2 \Gei^{l+1,n+1} \, \devGei^{l+1,n+1}. \label{eqn.nsP5}
    \end{align}
\end{subequations}	
The only fully implicit discretization is concerned with the computation of $\Gei^{l+1,n+1}$, where 
the non-linearity contained in the source term is solved by an exponential integrator discussed in 
the next section. To improve the efficiency of the iterative scheme \eqref{eqn.nsPicard}, the new 
density $\rho_i^{l+1,n+1}$ is deduced from the new volume $|T_i|^{l+1,n+1}$ computed with the new 
coordinates $\x_r^{l+1,n+1}$, thus the GCL \eqref{eqn.fvcl1} is not directly solved. However, due 
to the compatible time update of the fluxes in \eqref{eqn.fvcl1}, this approach is equivalent to 
explicitly solve the GCL. The multi-index $\{l,n+1\}$ is shorten to the asterisk index $\{*\}$, 
hence the result of the nonlinear nodal solver is the nodal velocity $\tilde{\vv}_{r}^*$ and the 
subcell force $\tilde{\bm{f}}_{ri}^*$, evaluated according to \eqref{eq:subforce^n}, with the total 
Cauchy stress $\tensor{T}_i^*$ (which is a function of $\tensor{\sigma}_i^*$) computed within the 
Picard loop \eqref{eqn.nsPicard}. The iterative scheme stops when one of the following exit 
conditions is satisfied, based on a tolerance set to $\delta=10^{-12}$.
\begin{itemize}
	\item The material is an ideal gas, thus the hydrodynamics limit of the model is reached:
	\bea
	\epsilon_h^{l+1} := \left|\Gei^{l+1,n+1}-\left(\frac{\rho^{l+1,n+1}}{\rho_0}\right)^{2/3} \, 
	\Id \right| \leq \delta,
	\label{eqn.rH}
	\eea
	which corresponds to the stiff relaxation limit for $\Ge$.
	\item The material is a purely elastic solid, therefore
	\bea
	\epsilon_e^{l+1}:=\left|\Gei^{l+1,n+1}-\tensor{G}_{i,*}^{l+1,n+1}\right|\leq \delta,
	\label{eqn.rE}
	\eea
	where $\tensor{G}_{i,*}^{l+1,n+1}$ is the solution of the homogeneous equation 
	related to \eqref{eqn.nsP4}. This condition means that the source term vanishes and no relaxation is 
	needed.
	\item Convergence is achieved between two consecutive iterations for any of the following residuals:
	\bea
	\left| \epsilon_h^{l+1} - \epsilon_h^{l}  \right| \leq \delta, \qquad \left| \epsilon_e^{l+1} - \epsilon_e^{l}  \right| \leq \delta.
	\label{eqn.rG}
	\eea 
\end{itemize} 
The maximum number of iterations is set to $\mathcal{L}=10$. In the numerical tests shown in 
Section \ref{sec.test}, the nonlinear nodal solver has always achieved convergence before reaching 
$\mathcal{L}$ iterations. Furthermore, let us point out that if the material is either an ideal 
fluid or an elastic solid, the nonlinear solver converges in one Picard iteration, namely one of 
the exit conditions between \eqref{eqn.rH} and \eqref{eqn.rE} is fulfilled. In these particular 
situations, the discretization of the tangential stresses $ \tensor{\sigma} $ becomes fully 
explicit and the nodal solver \textit{exactly} reduces to the EUCCLHYD scheme \cite{phm109} for 
ideal hydrodynamics, or the corresponding version for hyperelasticity materials presented in 
\cite{CCL2020,Boscheri2021}.	

\paragraph{Remark} Boundary conditions are imposed in a compatible way
with our numerical scheme. More precisely, either the Cauchy stress or the velocity normal to a 
boundary face must be prescribed in the nodal solver in order to ensure global conservation of momentum. 
All the details can be found in \cite{Maire2007,Maire2011a,LAM2018,Boscheri2021}. The face-based 
fluxes are computed with a standard finite volume scheme that requires ghost states to properly set 
the boundary condition. The reader is referred to \cite{BosARCME} and references therein for an 
exhaustive explanation on the treatment of boundary conditions in cell-centered finite volume 
schemes on moving unstructured meshes.
 	
\subsection{Exponential integrator for stiff source terms}	\label{ssec.ExpInt}

Due to the wide span of possible timescales characterizing the evolution of the elastic metric tensor
$\Ge$,  
the strain relaxation source on the right 
hand side of Eq. \eqref{eqn.cl5} oftentimes is of stiff nature and 
cannot be easily treated with an explicit method. The same holds for the simpler thermal 
impulse equation \eqref{eqn.cl4}.
For this reason we employ an efficient semi-implicit, semi-analytical integration scheme that
can handle such arbitrarily stiff source terms, and in particular can accurately recover the 
asymptotic equilibrium state associated with viscous fluxes in the Navier-Stokes equations and 
also relax to the Fourier limit of thermal conduction.

\subsubsection{Thermal impulse relaxation solver}
The thermal impulse equation \eqref{eqn.cl4} can be written in extended form as 
\begin{equation}
   \dfrac{\de{\vec{J}}}{\de{t}} = -\dfrac{\nabla T}{\rho} - \dfrac{\alpha^2\,T}{\rho\,\kappa}\,\vec{J} = 
   \Pstar - \dfrac{1}{\tau_2}\,\vec{J}, 
\end{equation}
which, once the discretization of $-{\nabla T}/{\rho}$ is fixed (to the constant value $\Pstar$), 
can be seen as a 
system of three uncoupled first order linear ordinary differential equations (ODEs) and 
an exact solution is indeed found thanks to the linearity and independence of
the three equations.
Explicitly, the solution is 
\begin{equation} \label{eqn.heatsolution}
    \vec{J}^{n+1} = (\vec{J}^n - \tau_2\,\Pstar)\,\exp(- \Delta t/\tau_2) + \tau_2\,\Pstar, 
\end{equation}
with $\Pstar = - \sum _{f \in \mathcal{F}_i} \widehat{T_{fi} \Id  \cdot \n_{fi}}^n \, 
s_{f}^n/{m_i}$ 
being approximated as the \textit{constant} discrete time-derivative of $\vec{J}$ as given by 
the update formula \eqref{eqn.fvcl4} when 
the source term is neglected. The only degenerate case to be considered is that
if $\Delta t/\tau_2$ is very small (of the order of $10^{-8}$),
i.e. if the source term is not stiff at all, then \eqref{eqn.heatsolution} might yield 
inaccurate results, due to floating point representation issues.
In this case, the solution algorithm simply opts to switch to explicit Euler integration, 
which for such mild (vanishing) sources yields perfectly valid solutions.

Finally the result can be combined with the implicit-explicit Runge--Kutta framework
by separating the (formally implicit) contribution due to the relaxation source 
$\Delta\vec{J}_{rel} = \vec{J}^{n+1} - (\vec{J}^n + \Delta t\,\Pstar)$ from the explicit
operator associated with the convective 
discretization $\Delta\vec{J}_{conv} := \Delta t\,\Pstar$.


\subsubsection{Strain relaxation solver}
The solver for the strain  relaxation source is based on the exponential integrator developed in \cite{tavellicrack} for the computation
of diffuse interface fractures and material failure, but exploits in a deeper manner the particular structure 
of the equation being solved, following the technique presented in \cite{chiocchettimueller} for 
the integration of stiff finite-rate pressure  relaxation sources.

An important aspect of the scheme is that it avoids fractional-step-type splitting, so that
the Navier-Stokes stress tensor and the Fourier heat flux can be recovered regardless 
of the ratio between the computational time step 
size and the relaxation timescales. This means that the global time step size need not to be adjusted 
to accommodate for the fast dynamics of the relaxation sources.

Moreover, we recall that the solver employed in \cite{tavellicrack} required, in general, the 
solution of a sequence of a non-homogeneous nine-by-nine systems of linear ordinary differential 
equations for the nine independent components of the distortion field $\Ae$, 
which involves the numerical computation of matrix exponentials and the inversion of 
the Jacobian matrix of the ODE system. Both these operations constitute delicate tasks in linear algebra
that require special care to be carried out in an efficient and accurate manner.

The approach used in this work entirely foregoes the solution of such nine-by-nine systems (six-by-six, in 
the case of the symmetric tensor $\Ge$) and the associated linear algebra intricacies. Thus, we 
compute the analytical solution to one of the several different linearized equations (automatically 
chosen 
by the solver) that 
approximate the nonlinear ODE
\begin{equation} \label{eqn.strainode}
\mddt{\Ge} =  \Lstar - \frac{6}{\tau_1}\,{\det{(\Ge)}}^{5/6}\,\Ge\,\devGe,
\end{equation}
while admitting simple solutions that can be evaluated in a robust fashion. Here,
with $\Lstar$ we denote a \textit{constant convective forcing term} to be 
given in the following paragraphs, in analogy to the previously defined $\Pstar$
discrete time derivative of thermal impulse.


This is achieved by first computing, cell by cell, 
the update to $\Ge$ 
associated with the left hand side of \eqref{eqn.cl5},
and then including its effects 
in \eqref{eqn.strainode}, in the form of the constant forcing term $\Lstar$.
Formally, this first step amounts to computing the 
solution ${\Ge}_\ast = \Ge(t_{n+1})$ to the initial value problem

\begin{equation} \label{eqn.ivplhs}
\left\{
   \begin{aligned}
   & \mddt{\Ge} + \Ge\,\nabla\vec{v} + \nabla\vec{v}^\transpose \, \Ge = \tensor{0},\\
   & \Ge(t^n) = {\Ge}^n,
   \end{aligned}
\right.
\end{equation}
which, in our case, is solved with the compatible
method applied to the equation for the metric tensor $\Ge$, yielding a point-wise update
to the cell averages that allows to define the constant convective forcing term
\begin{equation}
 \Lstar = \dfrac{\Ge_\ast - {\Ge}^n}{\Delta t}.
\end{equation}

Then a sub-timestepping loop with adaptive step size $\delta t^m$ is entered in order to approximate the solution of \eqref{eqn.strainode} 
with a sequence of solutions of linearized ODEs. Such a sub-timestepping loop is 
useful for ensuring 
the robustness and accuracy of the solver in complex flow configurations, 
but generally in our computations
the solver achieves convergence in one single sub-time-step $\delta t^m = \Delta t$.
We remand to \cite{chiocchettimueller,tavellicrack} for more details on the sub-timestepping approach, 
and we carry on our presentation of the method by listing three possible approximation choices for the solution 
of \eqref{eqn.strainode}.

\subsubsection{Approximate analytical solution for strain (1)}
When dealing with fluid flows (i.e., when the source term acts on fast timescales), 
rather often one may assume that $\Ge$ is a perturbation of a spherical 
tensor, that is, $\devGe$ can be assumed small.
Then, it is advantageous to rewrite the evolution equation for
the elastic metric tensor \eqref{eqn.strainode} as 
\begin{equation} \label{eqn.smalldevg}
\mddt{\Ge} = \Lstar - k\,\devGe\,\devGe + 
   k\,{\left(\dfrac{\tr \Ge}{3}\right)}^2\,\tensor{I} - k\,\dfrac{\tr \Ge}{3}\,\Ge,
\end{equation}
with $k = {6}\,{\det(\Ge)}^{5/6}/\tau_1$ taken constant for the sub-time-step. This splits the 
source in four
pieces. The first is the constant $\Lstar$, associated with convection which, by definition, cannot be stiff as its size
is limited by the CFL constraint of the global timestepping scheme. The second is a (small by hypothesis) quadratic term
in $\devGe$ which can be safely approximated as constant. 
The third is a function of $\tr \Ge$ only, again formally taken constant. This assumption 
can be justified by writing the evolution equation for the trace of the metric tensor $\Ge$
\begin{equation}
   \mddt{}\left(\tr{\Ge}\right) = \Lstar - k\, \tr{\left(\devGe\,\devGe\right)},
\end{equation}
which shows that either $\tr{\Ge}$ varies on a timescale associated with 
convection (by definition, slow), or as a quadratic function of $\devGe$ (small by assumption).
The approximate equation \eqref{eqn.smalldevg} then admits the simple exact solution

\begin{equation}
   \Ge^{m+1} = \Ge(t^m + \delta t^m) = \exp\left(-k\,\dfrac{\tr \Ge}{3}\,\delta t^m\right)\,\left(\Ge^m + \tensor{F}_0\right) - \tensor{F}_0,
\end{equation}
with 
\begin{equation}
    \tensor{F}_0 = -\dfrac{3}{k\,\tr{\Ge}}\,\left(\Lstar + 
        k\,{\left(\dfrac{\tr{\Ge}}{3}\right)}^2\,\tensor{I} - k\,\devGe\,\devGe\right).
\end{equation}

We should remark that \textit{nowhere} in this approximate solution we neglected the contributions due to 
$\devGe$, they only have taken to be constant for a sub-time-step, so that, in principle
there are no hypotheses restricting the use of such an approximation, beside increased requirements
imposed on the sub-timestepping scheme when $\devGe$ is not small. Specifically, in this work
our constant approximations 
are initially set to $\tr\Ge = \tr\Ge^m$ and $\devGe = {\devGe}^m$ and then updated within a fixed 
point iteration
as $\tr\Ge = \tr(\Ge^m + \Ge^{m+1})/2$ and $\devGe = \dev(\Ge^m + \Ge^{m+1})/2$.

\subsubsection{Approximate analytical solution for strain (2)}
Whenever the deviatoric part of $\Ge$ cannot be assumed small, i.e. in practice when 
\begin{equation}
   \sqrt{\tr{\left(\devGe^m\,\devGe^m \right)}} > \epsilon_1\,{\det{(\Ge^m)}}^{1/3},
\end{equation}
better accuracy in the approximation of \eqref{eqn.strainode} can be obtained by observing 
that it is possible to switch the order of the operands of the matrix product $\Ge\,\devGe$ 
appearing in \eqref{eqn.strainode}. In this work the coefficient 
$\epsilon_1$ is permanently set to $\epsilon_1 = 0.2$, which means that this approximation mode
has been used also when $\devGe$ is not \textit{strictly} small and indeed one may in principle 
choose higher values for $\epsilon_1$ without compromising the behavior of the solver.
Then we can rewrite \eqref{eqn.strainode} as
\begin{equation} \label{eqn.largedevg}
\dfrac{\de \Ge}{\de t} = \Lstar - k\,\Ge \devGe, 
\end{equation}
where we will take $\Lstar$, $k = {6}\,{\det{\Ge}}^{5/6}/\tau$, and $\devGe$ to be constant at 
each sub-time-step.
At the implementation level, in order to simplify 
the solution of \eqref{eqn.largedevg}, we 
work in the principal reference frame, which diagonalizes $\Ge$ and $\devGe$, 
i.e. we compute the orthonormal matrix
$\tensor{E}$ such that $\hat{\Ge} = \tensor{E}^{-1}\,{\Ge}\,\tensor{E}$ and 
$\hat{\devGe} = \tensor{E}^{-1}\,{\devGe}\,\tensor{E}$
are diagonal matrices and apply the associated change of basis to all vectors and tensors in our 
equation. 
In this way the exact solution to Eq. \eqref{eqn.largedevg} is 
\begin{equation}
   \Ge^{m+1} = \Ge(t^m + \delta t^m) = \tensor{E}\,\left[\exp\left(-k\,\hat{\devGe}\,\delta 
   t^m\vphantom{l^{l^l}}\right)\,\left(
   \tensor{E}^{-1}\,\Ge^n + \dfrac{1}{k}\,{\hat{\devGe}}^{-1}\,\tensor{E}^{-1}\,\Lstar
   \right)                 - 
   \dfrac{1}{k}\,{\hat{\devGe}}^{-1}\,\tensor{E}^{-1}\,\Lstar\right].
\end{equation}
The three-by-three matrix $\tensor{E}$ having for columns the eigenvectors of $\Ge$ can be quickly and robustly computed 
to arbitrary precision by means of Jacobi method for the eigenstructure of symmetric matrices, 
and its inverse is simply given by $\tensor{E}^{-1} = \tensor{E}^\transpose$.
Furthermore, $\hat{\devGe}$ can be inverted trivially in the principal reference frame by just 
taking 
the reciprocal of each diagonal entry. Like for the previous solution we iteratively update the constant 
estimate $\hat{\devGe} = \dev(\hat{\Ge}^m + \hat{\Ge}^{m+1})/2$ so to gain higher accuracy while
maintaining the linearity of the ODE being solved.


\subsubsection{Determinant constraint}
In the solution of the equation for the metric tensor $\Ge$, specifically when the computation
involves fluid-type behavior, special care must be paid to preserve the nonlinear algebraic constraint
$\det{\Ge}(t,\ \vec{x}) = (\rho(t,\ \vec{x})/\rho_0)^2$. For the purpose of notation compactness, in the following
we will denote the target determinant as $D(t,\ \vec{x})= (\rho(t,\ \vec{x})/\rho_0)^2$.
Thanks to the compatible discretization of the scheme presented in 
this work, the GCL is fulfilled, and the solution for density is evolved in a way that 
is compatible with the mesh motion, thus the convective part of the equation satisfies the 
determinant constraint up to the order of accuracy of the time discretization 
(see \cite{Boscheri2021} for further details). However, the same
constraint must be actively enforced when source terms are present. A simple approach to the problem
consists in uniformly multiplying all components of $\Ge$ by $(D/\det\Ge)^{1/3}$ so that the 
resulting determinant will be $D$.

The specific numerical value of the target determinant $D$ is clearly known (as a function of density) 
at the time levels $t^n$ and $t^{n+1}$, 
however it must be somehow approximated for all
the in-between times during which we operate our sub-timestepping procedure. 
In this work we impose that for a given sub-time-step indexed by $m$, connecting $t^m$ and 
$t^{m+1}$, the 
determinant $D$ is computed as 
\begin{equation} \label{eqn.determinant}
   D = \beta_s\,D_s + (1 - \beta_s)\,D_f,
\end{equation}
where we define $D_s = \det\left(\Ge + (t^m + \delta t^m - t^n)\,\Lstar\right)$ to be the value that
the determinant would have following a linear segment path connecting the two 
states $\Ge$ and $\Ge_\ast$, that is, the value that would allow
exact integration of the (zero) source term in the solid limit to be preserved. 
In the fluid limit instead, we take $D_f = \det{\Ge} + (t^m + \delta t^m - t^n)\,(\det{\Ge_\ast} - 
\det\Ge^{n})$ 
as a second order approximation of the determinant.
The mixing ratio $\beta_s$ for the two approximations $D_s$ and $D_f$ is 
a heuristic measure of how close to a solid can the material be considered and in particular the
expression we adopt is
\begin{equation}
   \beta_s = \min{\left(1,\ 
   \dfrac{||\Lstar||_2^2}{||6/\tau_1\,\det{(\Ge^m)}^{5/6}\,\Ge^m\,\devGe^m||_2^2 + 
   10^{-14}}\right)}^{4}, 
\end{equation}
with $||\bullet||_2^2$ denoting the square of the Frobenius norm of a given tensor.

\subsubsection{Approximate analytical solution for strain (3): fixed point iteration for the Navier-Stokes equilibrium state}
Oftentimes the timescale $\tau_1$ of strain relaxation is so fast that
one may decide, out of computational efficiency concerns, 
to just compute the strain state for which the forcing term due to
convection $\Lstar$ and the relaxation source are balanced yielding a \emph{local
equilibrium} state corresponding to the Navier-Stokes limit of the GPR model.
Such an equilibrium state can be easily computed by means of a fixed point iteration in the form
\begin{equation} \label{eqn.fpi}
    \begin{aligned}
       &\widetilde{\Ge}_l = 
       \dfrac{\tau_1}{6\,\det(\Ge_l)^{5/6}}\,\dev\left(\Ge^{-1}_l\,\Lstar\right) + 
       \dfrac{\tr\left(\Lstar - \Ge_l\,\dev{\left(\Ge^{-1}_l\,\Lstar\right)}\right)}{\tr\left(\Ge^{-1}_l\,\Lstar\right)}\,\tensor{I},\\
       &\Ge_{l+1} = \widetilde{\Ge}_l\,{\left(\dfrac{D}{\det\widetilde{\Ge}_l}\right)}^{1/3},
    \end{aligned}
 \end{equation} 
 with $D$ the target determinant as defined in \eqref{eqn.determinant} and $l$ the iteration index in the fixed point procedure.
 We found that the fixed point iteration \eqref{eqn.fpi} is always convergent regardless of the initial guess, 
 but nonetheless we care to provide a simple and efficient choice in the form
 \begin{equation}
    \Ge_{1} = \widetilde{\Ge}_0\,{\left(\dfrac{D}{\det\widetilde{\Ge}_0}\right)}^{1/3} , 
        \text{ with } \widetilde{\Ge}_0 = \tensor{I} + 
        \dfrac{\tau_1}{6\,\det(\Ge^m)^{5/6}}\,\dev{\Lstar}.
 \end{equation}
We remark that, when the relaxation timescale is much faster than the time step size, 
this approximate solution is not only very efficient and robust, 
but will also converge to the Navier-Stokes equilibrium state up to machine accuracy.

 \subsubsection{Summary of the selection procedure for the approximation method}
 At each sub-timestep between $t^m$ and $t^{m+1}$, our solver for the equation of the elastic metric tensor $\Ge$ has to select
 the optimal approximation method for the specific flow configuration at hand.
 The selection procedure is carried out as follows:
 \begin{enumerate}
    \item If the source is not stiff, i.e. in practice if $\beta_s > 1 - 10^{-14}$, then we use explicit Euler integration and compute
    the solution at the next time sub-level as $\Ge^{m+1} = \Ge^m + \Delta t\,\left(\Lstar - 
    \dfrac{6}{\tau_1\,\det{(\Ge^m)}^{5/6}\,\Ge\,\devGe}\right)$.
   \item Else, we define the indicator matrix $\tensor{\Lambda} = \abs{\left({\Ge^m}^{-1}\,\Lstar - 
   k\,\devGe^m\right)}$ and 
       if the sum of the off-diagonal components of $\tensor{\Lambda}$ is less than $\tr{\tensor{\Lambda}}/5$ while also $\delta t^m > \tau_1$ holds, 
       then
       the scheme opts for the fixed point iteration \eqref{eqn.fpi}. Here, $ \abs(\bullet) $ is 
       applied component-wise to a given matrix.
   \item Else, if $\sqrt{\tr{(\devGe^m\,\devGe^m)}} < \epsilon_1\,\det{(\Ge^m)}^{1/3}$ or if any of the 
   diagonal 
       entries of ${\hat{\devGe}}^m$ has magnitude smaller than $\epsilon_2\,\tr{\hat{\Ge}^m}$, then the 
       scheme selects \eqref{eqn.smalldevg}. In this work we permanently set $\epsilon_2 = 10^{-3}$ in order
       to prevent division by small numbers in \eqref{eqn.largedevg}.
   \item If none of the above, then we apply approximation \eqref{eqn.largedevg}.
\end{enumerate}
Regardless of the chosen approximation method, at the end of each sub-timestep, the result 
$\Ge^{m+1} = \Ge(t^{m+1})$
must be multiplied by ${(D/\det{(\Ge^{m+1})})}^{1/3}$ so that the determinant constraint is 
satisfied.

\subsection{Second order extension in space and time}
\label{ssec.2ndOrder}
The accuracy of the Godunov-type scheme \eqref{eqn.fv} is improved up to second order in space and time by performing a spatial TVD linear reconstruction of the state variables and an IMEX Runge-Kutta time stepping scheme, respectively.

\paragraph{TVD piecewise linear reconstruction} To achieve second order of accuracy in space, a piecewise linear reconstruction is performed for all variables in the state vector $\Q$. The starting point is given by the known cell averages $\Q_i^n$ and the mesh configuration at time $t^n$. As a result, piecewise linear polynomials $\mathbf{w}^n_i(\x^n)$ are obtained for each cell and each variable of the governing system \eqref{eqn.cl}. This reconstruction polynomial is expressed in terms of a set of piecewise linear spatial basis functions that form a modal basis, that is
\begin{equation}
	\mathbf{w}^n_i(\x^n) = \sum \limits_{l=1}^\mathcal{M} \psi_l(\boldsymbol{\xi}) \, \widehat{\w}^{n}_{l,i}:=\psi_l(\boldsymbol{\xi}) \, \widehat{\w}^{n}_{l,i},
	\label{eqn.wh}
\end{equation}
where the classical tensor index notation based on the Einstein summation convention is adopted, which 
implies summation over two equal indices. The modal basis functions $\psi_l(\boldsymbol{\xi})$ are 
defined in the reference system and count a total number $\mathcal{M}=d+1$ of unknown degrees of 
freedom. They explicitly write
\begin{equation}
	\psi_l(\boldsymbol{\xi}) = \left\{ \begin{array}{lll}
		\left( 1, \, \xi - 1/3, \, \eta - 1/3 \right)^\transpose & & \textnormal{in 2D} \\
		\left( 1, \, \xi - 1/4, \, \eta - 1/4, \, \zeta - 1/4 \right)^\transpose & & \textnormal{in 
		3D} \\
	\end{array} \right. .
\end{equation}
A so-called reconstruction stencil $\mathcal{S}_i = \bigcup \limits_{j=1}^{n_e} T^n_{m(j)}$ is needed,
where $1\leq j \leq n_e$ is a local index that counts the elements belonging to the stencil, while 
$m(j)$ maps the local index to the global element numbers used in the mesh configuration 
\eqref{trian}. Due to the unstructured mesh, neither the stencil nor the element configuration is 
symmetric, thus the stencil contains a total number of $n_e=d \cdot \mathcal{M}$ elements to avoid 
ill-conditioned reconstruction matrices, as suggested in \cite{StencilRec1990,kaeserjcp}. The 
reconstruction stencil is filled considering the Voronoi neighbors of $T_i^n$ (i.e. the neighbor elements 
sharing at least one vertex with element $T_i^n$), 
that are recursively added until $n_e$ is 
reached. The reconstruction is based on integral conservation for each element $T^n_j \in 
\mathcal{S}_i$, that is  
	\begin{equation}
		\frac{1}{|T^n_j|} \int \limits_{T^n_j} \psi_l(\boldsymbol{\xi}) \widehat{\w}^{c}_{l,i} \, 
		\rmd \x = \Q^n_j,  
		\qquad \forall T^n_j \in \mathcal{S}_i,     
		\label{eqn.intRec}
	\end{equation}
where the degrees of freedom $\widehat{\w}^{c}_{l,i}$ refer to the high order \textit{unlimited} polynomial obtained from a central reconstruction. To enforce conservation of the reconstruction polynomial we must at least require that the above expression holds 
exactly for cell $T_i^n$, that is $\frac{1}{|T^n_i|} \int \limits_{T^n_i} \psi_l(\boldsymbol{\xi}) 
\widehat{\w}^{c}_{l,i} \, \rmd \x = \Q^n_i$. This linear constraint is added to the 
overdetermined system \eqref{eqn.intRec}, which is solved only for the unknown expansion 
coefficients $\widehat{\w}^{c}_{l,i}$ relying on a classical Lagrangian multiplier approach, see 
\cite{Lagrange3D,Dumbser2007693}. The integrals appearing in \eqref{eqn.intRec} are evaluated with 
second order accurate quadrature rules and they are defined in the reference element $T_e$.

The above procedure generates a second order reconstruction polynomial within each cell, that is linear in the sense of Godunov \cite{godunov}, thus it is oscillatory and non-monotone across discontinuities of the numerical solution. To overcome this problem, the higher modes $\widehat{\w}^{c}_{l,i}$ of the expansion \eqref{eqn.wh} must be limited so that the cell gradient respects the following monotonicity condition, which also ensures the TVD (Total Variation Diminishing) property:
\begin{equation}
 \Q^{n,\min}_i \leq \mathbf{w}^n_i(\x^n) \leq \Q^{n,\max}_i, \qquad \Q^{n,\min}_i = \min \left\{\Q^n_i, \, \min_{j \in \mathcal{T}_{r(i)}} \Q_j^n \right\}, \qquad \Q^{n,\max}_i = \max \left\{\Q^n_i, \, \max_{j \in \mathcal{T}_{r(i)}} \Q_j^n \right\},
 \label{eqn.monotone}
\end{equation}
with $\mathcal{T}_{r(i)}$ denoting the set of elements which share at least one node with cell $T_i$. To respect the above relation \eqref{eqn.monotone}, let us define a scalar $b_i \in [0;1]$, that is used to limit the higher modes of the reconstruction polynomial. For each state variable and for $l = 2, \ldots, \mathcal{M}$ we compute
\begin{equation}
	\widehat{\w}^{n}_{l,i} = b_i \, \widehat{\w}^{c}_{l,i} \qquad \textnormal{with} \qquad b_i = \min_{r \in \mathcal{R}_i} b_{i,r}, 
	\label{eqn.lim}
\end{equation}
while the cell average $\widehat{\w}^{n}_{1,i}=\widehat{\w}^{c}_{1,i}$ is ensured by construction thanks to the constrained 
least square approach used to solve the linear system \eqref{eqn.intRec}. Once the limiting procedure \eqref{eqn.lim} has 
been carried out, the final degrees of freedom $\widehat{\w}^{n}_{l,i}$ are known and the limited reconstruction 
polynomial \eqref{eqn.wh} is fully defined for each cell. The coefficient $b_i$ in \eqref{eqn.lim} depends on 
the element node and is computed according to \cite{BarthJespersen} as
\begin{equation}
	b_{i,r} = \left\{ \begin{array}{lll}
		\min \left(1, \, \frac{\Q^{n,\max}_i - \Q^{n}_i}{\mathbf{w}^n_i(\x_r^n) - \Q^{n}_i} \right) & \textnormal{if} & \mathbf{w}^n_i(\x_r^n) > \Q_i^n \\
		\min \left(1, \, \frac{\Q^{n,\min}_i - \Q^{n}_i}{\mathbf{w}^n_i(\x_r^n) - \Q^{n}_i} \right) & \textnormal{if} & \mathbf{w}^n_i(\x_r^n) < \Q_i^n \\
		1 & \textnormal{if} & \mathbf{w}^n_i(\x_r^n) = \Q_i^n 
	\end{array} \right. .
\end{equation}  

\paragraph{Implicit-Explicit Runge-Kutta time stepping} The second order in time extension is based 
on an implicit-explicit (IMEX) Runge-Kutta approach \cite{AscRuuSpi,PR_IMEX}, that is particularly 
suitable for discretizing hyperbolic systems of PDE with relaxation source terms. For the sake of 
clarity and compact notation in the description of the second order in time method, we call 
$\mathcal{L}_{im}$ the implicit operator and $\mathcal{L}_{ex}$ the explicit one, hence obtaining 
the semi-discrete form of the governing equations \eqref{eqn.cl}:
\begin{equation}
	\mddt{\Q} = \mathcal{L}_{ex}(t,\Q,\nabla \Q) + \mathcal{L}_{im}(t,\Q),
\end{equation}
with
\begin{equation}
	\mathcal{L}_{ex}(t,\Q,\nabla \Q)  = \left[ \begin{array}{c}
		\rho^{-1} \, \nabla \cdot \vv \\ \rho^{-1} \, \nabla \cdot \tensor{T} \\ \rho^{-1} \, 
		\nabla \cdot \left( \tensor{T}\vv + \q \right) \\ \rho^{-1} \,\nabla \cdot T \, \Id  \\ 
		-(\Ge \del\vv + \del\vv^\transpose \Ge)
	\end{array} \right], \qquad \mathcal{L}_{im}(t,\Q) = \left[ \begin{array}{c}
		0 \\ 0 \\ 0 \\ -\H/\sH\\ 2\tensor{\sigma} /(\rho \sG) 
	\end{array} \right] .
\end{equation}
The vertex-based fluxes are considered explicitly as well as the flux in the trajectory equation \eqref{eqn.trajODE}, 
because the nonlinear nodal solver described in Section \ref{ssec.NodalSol} reduces to a fully explicit scheme in the 
stiff relaxation limits of the model \eqref{eqn.cl}. In this way, the first order finite volume scheme \eqref{eqn.fv} is compactly written in semi-discrete form as
\begin{equation}
	\label{eqn.fvt1}
	\Q^{n+1} = \Q^n - \dt \, \mathcal{L}_{ex}(t^n,\Q^n,\nabla \Q^n) + \dt \,  \mathcal{L}_{im}(t^{n+1},\Q^{n+1}).
\end{equation}
To improve the time accuracy of the first order method \eqref{eqn.fvt1}, we use the 
second-order ARS(2,2,2) IMEX scheme \cite{AscRuuSpi} detailed in Table \ref{tab:Butcher_tableaux_ARS222} by
means of its two tableaux, where $\beta = 1 - \sqrt{2}/ 2$ and $\alpha=1-1/(2\beta)$.
\begin{table}[!ht]
	\begin{center}
		\begin{tabular}{cc|ccc}
			\multirow{4}{*}{\rotatebox{90}{Explicit}}
			&0       & 0           & 0           & 0 \\
			&$\beta$ & $\beta$     & 0           & 0 \\
			&1       & $\alpha$  & 1 - $\alpha$ & 0 \\
			\cline{2-5}
			&& $\alpha$  & 1 - $\alpha$ & 0
		\end{tabular} \hspace{1cm}
		\begin{tabular}{cc|ccc}
			\multirow{4}{*}{\rotatebox{90}{Implicit}}
			&0       & 0 & 0           & 0       \\
			&$\beta$ & 0 & $\beta$     & 0       \\
			&1       & 0 & 1 - $\beta$ & $\beta$ \\ 
			\cline{2-5}
			&& 0 & 1 - $\beta$ & $\beta$
		\end{tabular} 
		\caption{Butcher tableaux for the ARS(2,2,2) time discretization.
			Left panel: explicit tableau.
			Right panel: implicit tableau.
			$\beta = 1 - \Frac{\sqrt{2}}{2}$, $\alpha=\beta-1$.}
		\label{tab:Butcher_tableaux_ARS222}
	\end{center}
\end{table}
Remarking that $\alpha=\beta-1$ and $1-\alpha=2-\beta$, the following two-step scheme for the Lagrangian GPR model \eqref{eqn.cl} is obtained:
\begin{eqnarray}
	\frac{\Q^{(1)}- \Q^n}{\dt} &=& \beta \, \mathcal{L}_{ex}(t^n,\Q^n,\nabla \Q^n) + \beta \,  \mathcal{L}_{im}(t^{(1)},\Q^{(1)}), \label{eqn.IMEX1} \\
	\frac{\Q^{n+1}- \Q^n}{\dt} &=& (\beta-1) \, \mathcal{L}_{ex}(t^n,\Q^n,\nabla \Q^n) + (2-\beta) \, \mathcal{L}_{ex}(t^{(1)},\Q^{(1)},\nabla \Q^{(1)}) \nonumber \\
	&+& \beta \mathcal{L}_{im}(\Q^{(1)}) + (1-\beta) \, \mathcal{L}_{im}(t^{(1)},\Q^{(1)}) + \beta \,  \mathcal{L}_{im}(t^{n+1},\Q^{n+1}), \label{eqn.IMEX2}
\end{eqnarray}
where the superscript $(1)$ denotes the intermediate Runge-Kutta stage. The trajectory ODE also follows 
the explicit tableau, since it contains only vertex-based fluxes, thus one gets
\begin{eqnarray}
	\frac{\x^{(1)} - \x^n}{\dt} &=& \beta \, \tilde{\vv}^{*,n}, \\
	\frac{\x^{n+1} - \x^n}{\dt} &=& (\beta-1) \, \tilde{\vv}^{*,n} + (2-\beta) \, \vv^{*,(1)}.
\end{eqnarray}
Let us notice that the ARS(2,2,2) IMEX scheme in Table \ref{tab:Butcher_tableaux_ARS222} has the 
property of being stiffly accurate (SA), meaning that the last Runge-Kutta stage coincides with 
the new solution at time $t^{n+1}$. This is a crucial feature for developing asymptotic preserving 
time discretizations, as detailed in \cite{PR_IMEX}. Furthermore, the scheme is only diagonally 
implicit, therefore it requires the solution of one implicit equation for the heat source and for 
the source of the metric tensor equation. This is also confirmed by looking at the time stepping 
scheme \eqref{eqn.IMEX1}-\eqref{eqn.IMEX2}, where the implicit contribution is only taken into 
account by the terms $\beta \, \mathcal{L}_{im}(\Q^{(1)})$ in the first stage \eqref{eqn.IMEX1} 
and $\beta \, \mathcal{L}_{im}(\Q^{n+1})$ in the second stage \eqref{eqn.IMEX2}.

\section{Asymptotic analysis of the scheme}	\label{sec.AP}
In this section we analyze the stiff relaxation limits of the fully discrete first order finite 
volume scheme \eqref{eqn.fv}, which are retrieved for $\tau_1 \to 0$ and $\tau_2 \to 0$. The 
numerical method is proven to be \emph{asymptotic preserving} for the viscous stress tensor 
$\tensor{\sigma}$ and for the thermal impulse $\J$. To that aim, let us introduce the $k$-th order 
Chapman-Enskog expansion of a generic variable $\phi$ in powers of the stiffness parameter $\tau$ 
(i.e. $ \tau=\tau_1 $ or  $ \tau=\tau_2 $), 
that reads
\bea
\phi = \phi_{(0)} + \tau \phi_{(1)} + \tau^2 \phi_{(2)} + \ldots + \mathcal{O}(\tau^k).
\label{eqn.exp}
\eea
Application of the expansion \eqref{eqn.exp} up to the first order in $ \tau $ to $\Gei$ and $\J_i$ 
in the 
cell-centered evolution 
equations \eqref{eqn.fvcl4} and \eqref{eqn.fvcl5} yields
\begin{eqnarray}
	\frac{\Geio^{n+1} + \tau_1\Geioo^{n+1} - \Geio^{n} - \tau_1\Geioo^{n}}{\dt} &=& -\left( 
	\Geio^{n} + \tau_1\Geioo^{n} \right) \, \gradv_i(\tilde{\vv}^{*}) - 
	\gradv_i(\tilde{\vv}^{*})^{\transpose} \, \left( \Geio^{n} + \tau_1\Geioo^{n} \right)\nonumber 
	\\
	&+& \frac{6}{\tau_1} \left|  \Geio^{n+1} + \tau_1\Geioo^{n+1} \right|^{5/6} \, \left(\Geio^{n+1} + \tau_1\Geioo^{n+1}\right) \, \left(\mathring{\tensor{G}}_{e_i(0)}^{n+1} + \tau_1 \mathring{\tensor{G}}_{e_i(1)}^{n+1}\right) \nonumber \\
	&+& \mathcal{O}(\tau_1^2), \label{eqn.Gexp} \\
	\frac{\J_{i(0)}^{n+1}+\tau_2 \J_{i(1)}^{n+1}-\J_{i(0)}^{n}-\tau_2 \J_{i(1)}^{n}}{\dt} &=& - 
	\frac{1}{m_i} \sum \limits_{f \in \mathcal{F}_i} \widehat{T_{fi} \Id  \cdot \n_{fi}}^n \, 
	s_{f}^n -  \frac{1}{\tau_2} \frac{T_i^n}{T_0} \frac{\rho_0}{\rho_i^{n+1}} \, 
	\left(\J_{i(0)}^{n+1}+\tau_2 \J_{i(1)}^{n+1}\right)  + \mathcal{O}(\tau_2^2) \label{eqn.Jexp}.
\end{eqnarray}	

\paragraph{Asymptotic limit of the viscous stress tensor at zeroth order approximation} Multiplying by $\tau_1$ equation \eqref{eqn.Gexp}, letting $\tau_1 \to 0$ and retaining leading order terms for $\Gei$ leads to
\begin{equation}
	6 \left| \Geio^{n+1} \right|^{5/6} \, \Geio^{n+1} \mathring{\tensor{G}}_{e_i(0)}^{n+1} = \mathcal{O}(\tau_1^2).
	\label{eqn.zeroG}
\end{equation} 
According to the GPR model at the continuous level \cite{HPR2016} and due to the compatible 
discretization of the GCL \eqref{eqn.fvcl1} along the lines of \cite{Boscheri2021}, the finite 
volume scheme \eqref{eqn.fv} ensures the following compatibility condition
\begin{equation}
	\sqrt{|\Gei^{n+1}|}=\frac{\rho_i^{n+1}}{\rho_0} + \mathcal{O}(\dt) > 0.
	\label{eqn.Grho} 
\end{equation} 	
Therefore, $\left| \Geio^{n+1} \right|^{5/6}>0$ in \eqref{eqn.zeroG} and at zeroth leading order the discrete metric tensor becomes trace-free for $\tau_1\to 0$, i.e. $\mathring{\tensor{G}}_{e_i(0)}^{n+1} = \mathcal{O}(\tau_1^2)$. Inserting this result in the discrete viscous stress tensor \eqref{eqn.sigmah}, we obtain
\begin{equation}
	\tensor{\sigma}_i^{n+1} = -\rho_i^{n+1} \csh^2 \Geio^{n+1} \, 
	\mathring{\tensor{G}}_{e_i(0)}^{n+1} = \mathcal{O}(\tau_1^2).
\end{equation}
Hence, viscous stresses vanish and we retrieve the inviscid case which corresponds to the 
compressible Euler equations. Let us notice that, in the absence of heat conduction ($\alpha=\kappa=0$), 
the finite volume scheme \eqref{eqn.fvcl1}-\eqref{eqn.fvcl3} exactly reduces to the EUCCLHYD scheme 
presented in \cite{phm109}. Furthermore, we can write
\begin{equation}
	\Geio^{n+1} = \frac13\text{tr}\Geio^{n+1} \, \Id  + \mathcal{O}(\tau_1^2):=g_i^{n+1} \Id , 
	\label{eqn.G0g}
\end{equation}
thus allowing us to write the Chapman-Enskog expansion of $\Gei^{n+1}$ as
\begin{equation}
	\Gei^{n+1} = g_i^{n+1} \Id  + \tau_1 \Geioo^{n+1} + \mathcal{O}(\tau_1^2).
	\label{eqn.Gei1}
\end{equation}
The unknown coefficient $g_i^{n+1}$ can be determined by computing the determinant of $\Gei^{n+1}$ from the above relation and neglecting terms of the order $\mathcal{O}(\tau_1)$, that is
\begin{equation}
	|\Gei^{n+1}| = (g_i^{n+1})^3 + \mathcal{O}(\tau_1), \quad \Rightarrow \quad g_i^{n+1} = \left(\frac{\rho_i^{n+1}}{\rho_0}\right)^{2/3} + \mathcal{O}(\tau_1).
	\label{eqn.gconst}
\end{equation} 

\paragraph{Asymptotic limit of the viscous stress tensor at first order approximation} Here, a first order approximation of the viscous stress tensor is analyzed in the stiff limit $\tau_1\to 0$. Let us start from the previous expansion \eqref{eqn.Gei1}, which also implies $\devGei^{n+1}=\tau_1 \mathring{\tensor{G}}_{e_i(1)}^{n+1}+ \mathcal{O}(\tau_1^2)$. The compatibility condition \eqref{eqn.Grho} allows the density $\rho_i^{n+1}$ to be computed as
\begin{equation}
	\rho_i^{n+1} = \rho_0 \left| g_i^{n+1} +  \tau_1 \Geioo^{n+1} \right|^{1/2} + \mathcal{O}(\tau_1^2) = \rho_0 \left( (g_i^{n+1})^{3/2} + \frac{\tau_1}{2} (g_i^{n+1})^{1/2} \text{tr}\Geioo^{n+1} \right) + \mathcal{O}(\tau_1^2),
\end{equation}
thus a first order expansion of the discrete viscous stress tensor \eqref{eqn.sigmah} writes
\begin{equation}
\tensor{\sigma}_i^{n+1} = -\rho_0 \csh^2 \left( (g_i^{n+1})^{3/2} + \frac{\tau_1}{2} 
(g_i^{n+1})^{1/2} \text{tr}\Geioo^{n+1} \right) \left( g_i^{n+1} \Id  + \tau_1 \Geioo^{n+1} \right) 
\tau_1 \mathring{\tensor{G}}_{e_i(1)}^{n+1} + \mathcal{O}(\tau_1^4).
\end{equation}
By neglecting high order terms in $\tau_1$, the above expression leads to
\begin{equation}
	\tensor{\sigma}_i^{n+1} = - \tau_1 \rho_0 \csh^2 (g_i^{n+1})^{5/2} \mathring{\tensor{G}}_{e_i(1)}^{n+1} + \mathcal{O}(\tau_1^2).
	\label{eqn.sigma1}
\end{equation}	
Following \cite{DPRZ2016} and recalling that $\mathring{\tensor{G}}_{e_i(0)}^{n+1}=\mathcal{O}(\tau_1^2)$ from 
the zeroth approximation asymptotic analysis, the deviatoric operator is applied to the discrete 
equation \eqref{eqn.fvcl5}, then only first order terms for $\tau_1$ are retained, hence obtaining
\begin{equation}
	\Geio^{n+1} \, \gradv_i(\tilde{\vv}^{*}) + \gradv_i(\tilde{\vv}^{*})^{\transpose} \, 
	\Geio^{n+1} - \frac{2}{3} \text{tr}\left(\Geio^{n+1} \, \gradv_i(\tilde{\vv}^{*}) \right) 
	\Id  = -6 \left| \Geio^{n+1} \right|^{7/6} \mathring{\tensor{G}}_{e_i(1)}^{n+1} + 
	\mathcal{O}(\tau_1^2).
\end{equation}	
Using \eqref{eqn.G0g}, the above equation simplifies to
\begin{eqnarray}
	&& g_i^{n+1} \left( \gradv_i(\tilde{\vv}^{*}) +  \gradv_i(\tilde{\vv}^{*})^{\transpose} - 
	\frac{2}{3} \text{tr} \left(\gradv_i(\tilde{\vv}^{*}) \right) \Id  \right) = -6 
	\left(g_i^{n+1}\right)^{7/2} \mathring{\tensor{G}}_{e_i(1)}^{n+1} + \mathcal{O}(\tau_1^2), 
	\nonumber \\
	&& \left(g_i^{n+1}\right)^{5/2} \mathring{\tensor{G}}_{e_i(1)}^{n+1} = -\frac{1}{6} \left( 
	\gradv_i(\tilde{\vv}^{*}) +  \gradv_i(\tilde{\vv}^{*})^{\transpose} - \frac{2}{3} \text{tr} 
	\left(\gradv_i(\tilde{\vv}^{*}) \right) \Id  \right) + \mathcal{O}(\tau_1^2).
\end{eqnarray}
By substituting now the left hand-side of the above expression into the first order approximation 
of the viscous stress tensor \eqref{eqn.sigma1}, we eventually obtain
\begin{equation}
	\tensor{\sigma}_i^{n+1} = \frac{\tau_1}{6} \rho_0 \csh^2 \left( \gradv_i(\tilde{\vv}^{*}) +  
	\gradv_i(\tilde{\vv}^{*})^{\transpose} - \frac{2}{3} \text{tr} \left(\gradv_i(\tilde{\vv}^{*}) 
	\right) \Id  \right) + \mathcal{O}(\tau_1^2),
\end{equation} 	
which is the discrete version of the stress tensor of the compressible Navier-Stokes equations 
based on Stokes hypothesis, in which the effective dynamic viscosity coefficient is defined by
\begin{equation}
	\visc = \frac{1}{6} \rho_0 \tau_1 \csh^2,
	\label{eqn.visc-tau}
\end{equation} 
that retrieves the analytical formulation \eqref{eqn.visc}.
	
\paragraph{Asymptotic limit of the heat flux at first order approximation} From the leading order term $\tau_2^{-1}$ in \eqref{eqn.Jexp}, it follows that 
\begin{equation}
	\J_{i(0)}^{n+1} = \mathcal{O}(\tau_2^2).
	\label{eqn.J0}
\end{equation}
Assuming well-prepared initial data, i.e $\J_{i(0)}^{n} = \mathcal{O}(\tau_2^2)$, and inserting the above result into \eqref{eqn.Jexp} while letting $\tau_2 \to 0$ yields
\begin{equation}
	\J_{i(1)}^{n+1} = - \frac{T_0}{T_i^n} \frac{\rho_i^{n+1}}{\rho_0} \, \frac{1}{m_i} \sum 
	\limits_{f \in \mathcal{F}_i} \widehat{T_{fi} \Id  \cdot \n_{fi}}^n \, s_{f}^n + 
	\mathcal{O}(\tau_2^2).
	\label{eqn.J1}
\end{equation}
Finally, using \eqref{eqn.J0} and \eqref{eqn.J1} in the Chapman-Enskog expansion \eqref{eqn.exp}, 
we obtain
\begin{equation}
	\J_{i}^{n+1} = - \tau_2 \frac{T_0}{T_i^n} \frac{1}{\rho_0} \, \frac{1}{|T_i^{n+1}|} \sum 
	\limits_{f \in \mathcal{F}_i} \widehat{T_{fi} \Id  \cdot \n_{fi}}^n \, s_{f}^n + 
	\mathcal{O}(\tau_2^2), \qquad |T_i^{n+1}|=m_i/\rho_i^{n+1}.
	\label{eqn.Jheat}
\end{equation}
Therefore, in the stiff relaxation limit, the discrete heat flux vector $\q_i^{n}=\alpha^2 
T_i^n\J_{i}^{n}$ in the energy equation \eqref{eqn.fvcl3} becomes 
\begin{equation}
	\q_i^{n} = - \tau_2 \alpha^2 \frac{T_0}{\rho_0}\, \frac{1}{|T_i^{n}|} \sum \limits_{f \in 
	\mathcal{F}_i} \widehat{T_{fi} \Id  \cdot \n_{fi}}^n \, s_{f}^n + \mathcal{O}(\tau_2^2).
	\label{eqn.qheatAP}
\end{equation}
Notice that the term $1/|T_i^{n}|\sum \limits_{f \in \mathcal{F}_i} \widehat{T_{fi} \Id  
\cdot \n_{fi}}^n \, s_{f}^n$ in \eqref{eqn.qheatAP} is nothing but the discrete gradient operator 
for the temperature, i.e. $\nabla T$. Hence, taking the heat conduction coefficient $\kappa$ 
as in \eqref{eqn.kappa}, equation \eqref{eqn.qheatAP}  is an asymptotically consistent
discretization of the Fourier 
law of heat conduction $\q=-\kappa \nabla T$, that is recovered
by the finite volume scheme \eqref{eqn.fv} in the stiff relaxation limit $\tau_2 \to 0$. 
	
	
\section{Numerical results}	\label{sec.test}
In this section, we perform some numerical applications which aim at demonstrating the accuracy and the robustness of the new Lagrangian finite volume scheme \eqref{eqn.fv} for the solution of the GPR model \eqref{eqn.cl}.	We refer to this numerical algorithm with the abbreviation LGPR. A wide range of test cases is proposed, that covers simulations for ideal and viscous heat conducting fluids, elasto-plastic solids as well as purely elastic solids. We emphasize that for each test case the \textit{full} GPR model is solved, without neglecting any evolution equation. 

The setup of the test problems shown hereafter requires some parameters that will be specified, the 
most important being the relaxation times $\tau_1$ and $\tau_2$ which determine the mechanical amd 
thermodynamical behavior of the material under consideration. The parameters of the GPR model are summarized in Table \ref{tab.GPRpar}, together with the equation of state used for each 
test problem. For solids, the shear modulus is always computed as $G=\rho_0\,\csh^2$, while for fluids either the viscosity $\visc$ or the relaxation time $\tau_1$ must be prescribed, which are linked by the relation \eqref{eqn.visc}. In the case of heat conducting fluids, the relaxation time $\tau_2$ is deduced from \eqref{eqn.kappa}. The adiabatic sound speed $c_0$ for ideal gas EOS is evaluated according to \eqref{eqn.IG}, for Mie-Grüneisen EOS it is a parameter of the material, whereas for Neo-Hookean materials relation \eqref{eqn.c0NH} holds true. 

\begin{table}[!htbp]  
	\caption{Parameters of the GPR model for the test problems 
	shown in Section \ref{sec.test}. The hydrodynamics equation of state (EOS) is chosen among the 
	ideal gas (IG), Mie-Grüneisen (MG) or Neo-Hookean (NH).}  
	\begin{center} 
  		\renewcommand{\arraystretch}{1.2}
		\begin{tabular}{l|ccccccccc|c}
		Test & $\rho_0$ & $\csh$ & $\visc$ & $c_v$ & $\alpha$ & $\kappa$ & $T_0$ &  $\tau_1$ & $\tau_2$ & EOS \\
	    \hline
	    Isentropic vortex & 1 & 0.5 & - & 1 & 0 & 0 & 1 & $10^{-14}$ & $10^{-14}$ & IG   \\
	    Swinging plate & 1100 & 73 & - & 1 & 0 & 0 & 1 & $10^{\phantom{-}14}$ & $10^{-14}$ & NH \\
	    Kidder & 1 & 0.5 & - & 1 & 0 & 0 & 1 & $10^{-14}$ & $10^{-14}$ & IG   \\
	    Saltzman & 1 & 0.5 & - & 1 & 0 & 0 & 1 & $10^{-14}$ & $10^{-14}$ & IG   \\
	    Sedov & 1 & 0.5 & - & 1 & 0 & 0 & 1 & $10^{-14}$ & $10^{-14}$ & IG   \\
	    Riemann problems & 1 & 10 & $[10^{-3};10^{-2}]$ & 2.5 & 0 & 0 & 1 & \eqref{eqn.visc} & $10^{-14}$ & IG   \\
	    Heat conduction & 1 & 1 & $10^{-2}$ & 2.5 & 2 & $10^{-2}$ & 1 & \eqref{eqn.visc} & \eqref{eqn.kappa} & IG   \\
	    Viscous shock & 1 & 50 & $2 \cdot 10^{-2}$ & 2.5 & 50 & $9.33 \cdot 10^{-2}$ & 1 & \eqref{eqn.visc} & \eqref{eqn.kappa} & IG   \\
	    Shell collapse & 1845 & 9073.62 & - & 1 & 0 & 0 & 1 & \eqref{eqn.tau.plast} & $10^{-14}$ & MG   \\
	    2D projectile & 2785 & 3150 & - & 1 & 0 & 0 & 1 & \eqref{eqn.tau.plast} & $10^{-14}$ & MG   \\
	    3D Taylor bar & 8930 & 2245 & - & 1 & 0 & 0 & 1 & \eqref{eqn.tau.plast} & $10^{-14}$ & MG   \\
	    Elastic plate & 1845 & 9046.59 & - & 1 & 0 & 0 & 1 & $10^{\phantom{-}14}$ & $10^{-14}$ & NH   \\
	    Twisting column & 1100 & 73 & - & 1 & 0 & 0 & 1 & $10^{\phantom{-}14}$ & $10^{-14}$ & NH   \\
		\end{tabular}
	\end{center}
   \label{tab.GPRpar}
\end{table}	
				
If not specified, the material is initially unloaded, thus we set $\Ge=(\rho/\rho_0)^{2/3}\,\Id$ 
according to \eqref{eqn.gconst}. Furthermore, the thermal impulse vector is initialized with $\J=\mathbf{0}$ and 
the initial density distribution is $\rho=\rho_0$. The CFL number in \eqref{eqn.timestep} is taken to 
be $\textnormal{CFL}=0.45$ in 2D and $\textnormal{CFL}=0.3$ in 3D. The simulations are run using the fully second order space-time scheme described in Section \ref{ssec.2ndOrder}.

Physical units are based on the $[m, kg, s]$ unit system, thus Young and shear moduli for solids are measured in $[Pa]$, as well as pressure and stresses. Energy is expressed in Joule $[J]$.

\subsection{Numerical convergence studies} \label{ssec.conv}
The convergence studies of the LGPR scheme are carried out by considering both limits of the model, 
namely for $\tau_1 \to 0$ and $\tau_1 \to \infty$. The first case corresponds to the simulation of 
an ideal inviscid fluid (see Section \ref{sec.AP}), while the latter is concerned with a purely 
elastic material. These particular limits of the model permit to retrieve at the discrete level 
the EUCCLHYD scheme \cite{phm109} for hydrodynamics, and the cell-centered Lagrangian algorithm 
for nearly incompressible hyperelastic solids presented in \cite{Boscheri2021}.


\paragraph{Convergence studies in the stiff limit $\tau_1 \to 0$} We consider the isentropic vortex problem initially introduced in a two-dimensional setting \cite{HuShuTri} for the compressible Euler equations. The initial computational domain is $\Omega^{2D}(0)=[0;10] \times [0;10]$ and $\Omega^{3D}(0)=[0;10] \times [0;10] \times [0;5]$ with periodic boundaries. 
The ideal gas EOS \eqref{eqn.IG} is used with $\gamma=1.4$, and the relaxation time is set to $\tau_1=10^{-14}$, so that the hydrodynamic limit is retrieved. The fluid is characterized by a homogeneous background field on the top of which some perturbations are added, thus
\begin{equation}
	\rho = \rho_0 + \delta \rho, \quad u = u_c + \delta u,  \quad  v = v_c + \delta v, \quad  w=w_c, \quad  T = T_0 + \delta T,
\end{equation}
where the perturbations for density and pressure read
\begin{equation}
	\label{rhopressDelta}
	\delta \rho = (1+\delta T)^{\frac{1}{\gamma-1}}-1, \quad \delta p = (1+\delta T)^{\frac{\gamma}{\gamma-1}}-1, 
\end{equation}
with the temperature fluctuation $\delta T = -\frac{(\gamma-1){\lambda}^2}{8\gamma\pi^2}e^{1-r^2}$. 
According to \cite{HuShuTri}, the vortex strength is ${\lambda}=5$ and it moves with a convective velocity $\vv_c=(u_c,v_c,w_c)=(1,1,0)$. The velocity field is affected by the following perturbations:
\begin{equation}
	\label{ShuVortDelta}
	\left(\begin{array}{c} \delta u \\ \delta v \\ \end{array}\right) = \frac{{\lambda}}{2\pi}e^{\frac{1-r^2}{2}} \left(\begin{array}{c} -y' \\ \phantom{-}x' \end{array}\right), \qquad x' = x - x_c, \quad y' = y - y_c,
\end{equation}
where the center of the vortex is $(x_c,y_c,z_c)=(5,5,2.5)$ and the generic radial position is $r = \sqrt{{x'}^2 + {y'}^2}$. The final time of the simulation is chosen to be $t_f=0.1$ and the exact solution is simply given by the time-advected initial condition with convective velocity $\vv_c$. The simulation is performed on a sequence of four successively refined unstructured grids in 2D and in 3D, and the numerical convergence results are reported in Table \ref{tab.conv_hydro}. The errors are measured in $L_2$ norm for the variables $\{\omega,u,E\}$, and $h(\Omega(t_f))$ represents the mesh size which is taken to
be the maximum diameter of the circumspheres or the circumcircles of the elements in the final domain configuration $\Omega(t_f)$. The expected first and second order of accuracy are achieved in all cases.
	
	
	\begin{table}[!htbp]  
		\caption{Numerical convergence results for the isentropic vortex test in 2D and in 3D 
		using the LGPR scheme with relaxation time $\tau_1=10^{-14}$ on a sequence of refined 
		unstructured meshes of size $h(\Omega(t_f))$ measured at the end of the simulation 
		$t_f=0.1$. The errors are measured in $L_2$ norm for specific volume $\svol$, horizontal 
		velocity $u$ and total energy $E$.}  
		\begin{center} 
			\begin{small}
				\renewcommand{\arraystretch}{1.0}
				\begin{tabular}{c|cccccc}
					\multicolumn{7}{c}{2D LGPR $\mathcal{O}\,1$ ($\tau_1=10^{-14}$)} \\
					\hline
					$h(\Omega(t_f))$ & $(\svol)_{L_2}$ & $\mathcal{O}(1/\rho)$ & $u_{L_2}$ & $\mathcal{O}(u)$ & $E_{L_2}$ & $\mathcal{O}(E)$ \\ 
					\hline
					3.26E-01  & 5.405E-02 & -    & 1.547E-01 & -    & 2.579E-01 & -    \\ 
					2.47E-01  & 4.164E-02 & 0.96 & 1.219E-01 & 0.88 & 2.044E-01 & 0.86 \\ 
					1.63E-01  & 3.053E-02 & 0.74 & 8.866E-02 & 0.76 & 1.471E-01 & 0.78 \\ 
					1.28E-01  & 2.286E-02 & 1.20 & 7.041E-02 & 0.96 & 1.164E-01 & 0.97 \\ 
					\multicolumn{7}{c}{}\\
					\multicolumn{7}{c}{2D LGPR $\mathcal{O}\,2$ ($\tau_1=10^{-14}$)} \\
					\hline
					$h(\Omega(t_f))$ & $(\svol)_{L_2}$ & $\mathcal{O}(1/\rho)$ & $u_{L_2}$ & $\mathcal{O}(u)$ & $E_{L_2}$ & $\mathcal{O}(E)$ \\ 
					\hline
					3.26E-01  & 4.996E-02 & -    & 4.895E-02 & -    & 9.281E-02 & -    \\ 
					2.47E-01  & 3.312E-02 & 1.49 & 3.020E-02 & 1.76 & 5.509E-02 & 1.90 \\ 
					1.63E-01  & 1.913E-02 & 1.32 & 1.534E-02 & 1.63 & 2.858E-02 & 1.58 \\ 
					1.28E-01  & 1.327E-02 & 1.51 & 9.153E-03 & 2.13 & 1.770E-02 & 1.98 \\ 
					\multicolumn{7}{c}{}\\
					\multicolumn{7}{c}{}\\
					\multicolumn{7}{c}{3D LGPR $\mathcal{O}\,1$ ($\tau_1=10^{-14}$)} \\
					\hline
					$h(\Omega(t_f))$ & $(\svol)_{L_2}$ & $\mathcal{O}(1/\rho)$ & $u_{L_2}$ & $\mathcal{O}(u)$ & $E_{L_2}$ & $\mathcal{O}(E)$ \\ 
					\hline
					5.29E-01  & 2.389E-01 & -    & 5.600E-01 & -    & 8.781E-01 & -    \\ 
					3.62E-01  & 2.013E-01 & 0.35 & 4.075E-01 & 0.65 & 6.660E-01 & 0.56 \\ 
					2.31E-01  & 1.752E-01 & 0.31 & 2.882E-01 & 0.77 & 4.877E-01 & 0.69 \\ 
					1.81E-01  & 1.454E-01 & 0.76 & 2.301E-01 & 0.91 & 3.974E-01 & 0.83 \\ 
					\multicolumn{7}{c}{}\\
					\multicolumn{7}{c}{3D LGPR $\mathcal{O}\,2$ ($\tau_1=10^{-14}$)} \\
					\hline
					$h(\Omega(t_f))$ & $(\svol)_{L_2}$ & $\mathcal{O}(1/\rho)$ & $u_{L_2}$ & $\mathcal{O}(u)$ & $E_{L_2}$ & $\mathcal{O}(E)$ \\ 
					\hline
					5.29E-01  & 2.899E-01 & -    & 2.946E-01 & -    & 5.185E-01 & -    \\ 
					3.62E-01  & 1.426E-01 & 1.44 & 1.188E-01 & 1.85 & 2.275E-01 & 1.67 \\ 
					2.31E-01  & 8.304E-02 & 1.20 & 5.829E-02 & 1.59 & 1.099E-01 & 1.62 \\ 
					1.81E-01  & 5.931E-02 & 1.37 & 3.600E-02 & 1.96 & 7.206E-02 & 1.72 \\
				\end{tabular}
			\end{small}
		\end{center}
		\label{tab.conv_hydro}
	\end{table}
	
\paragraph{Convergence studies in the limit $\tau_1 \to \infty$} Here, we study the convergence of the LGPR scheme when $\tau_1=10^{14}$, thus the source term in the metric tensor equation vanishes and purely elastic solids are simulated. The swinging plate/cube test problem is considered, according to the setup proposed in \cite{scovazzi2,scovazzi3}. The initial computational domain is $\Omega(0)=[0;2]^d$ and the analytical solution for the velocity is given by
\begin{eqnarray}
\vv^{2D}(t,\x) = \Lambda U_0 \cos ( \Lambda t)
\left( \begin{array}{l}
	-\sin \left( \Frac\pi2 x \right) \cos \left( \Frac\pi2 y \right) \\
	\phantom{-}\cos \left( \Frac\pi2 x \right) \sin \left( \Frac\pi2 y \right)
\end{array} \right), \qquad \Lambda=\frac{\pi}{2} \sqrt{\frac{2\Gs}{\rho_0}}, \\
\vv^{3D}(t,\x) = \Lambda U_0 \cos ( \Lambda t)
\left( \begin{array}{l}
	-2\sin \left( \Frac\pi2 x \right) \cos \left( \Frac\pi2 y \right) \cos \left( \Frac\pi2 z \right) \\
	\phantom{-2}\cos \left( \Frac\pi2 x \right) \sin \left( \Frac\pi2 y \right) \cos \left( \Frac\pi2 z \right) \\
	\phantom{-2}\cos \left( \Frac\pi2 x \right) \cos \left( \Frac\pi2 y \right) \sin \left( \Frac\pi2 z \right)
\end{array} \right), \qquad \Lambda=\pi \sqrt{\frac{3\Gs}{4\rho_0}},
\label{eqn.SwingIni}
\end{eqnarray}
with $U_0=5\cdot 10^{-4}$. 
The Neo-Hookean hydrodynamic EOS \eqref{eqn.NH} and shear elastic energy \eqref{eqn.EE} are 
adopted, with Poisson 
ratio $\nu=0.45$ and 
Young modulus $Y=1.7 \cdot 10^7$, thus a nearly incompressible solid is modeled. The velocity and 
displacement fields are divergence-free, leading to the exact pressure $p=0$. Space-time dependent 
boundary conditions are prescribed for the normal velocities, according to the exact solution 
\eqref{eqn.SwingIni}. The final time of the simulation is chosen to be $t_f=\pi/\Lambda$, thus the 
final displacement corresponds to the initial one. As done for the stiff limit case, the swinging 
plate/cube test is run on a sequence of triangular and tetrahedral meshes that become finer, hence 
allowing the convergence rates to be computed. The results are collected in Table 
\ref{tab.conv_elasticity}, demonstrating that the formal order of accuracy is achieved even when 
simulating elastic solids with the same set of equations and the same LGPR scheme \eqref{eqn.fv} 
used for the hydrodynamics limit.
	
	
	\begin{table}[!htbp]  
		\caption{Numerical convergence results for the swinging plate test in 2D and in 3D using 
		the LGPR scheme with relaxation time $\tau_1=10^{14}$ on a sequence of refined unstructured 
		meshes of size $h(\Omega(t_f))$ measured at the end of the simulation $t_f=\pi/\omega$. The 
		errors are measured in $L_2$ norm for horizontal velocity $u$, total energy $E$, metric 
		tensor component $\tensor{G}_{e_{11}}$ and Cauchy stress component $\tensor{T}_{11}$.}  
		\begin{center} 
			\begin{small}
				\renewcommand{\arraystretch}{1.0}
				\begin{tabular}{c|ccccccccc}
					\multicolumn{9}{c}{2D LGPR $\mathcal{O}\,1$ ($\tau_1=10^{14}$)} \\
					\hline
					$h(\Omega(t_f))$ & $u_{L_2}$ & $\mathcal{O}(u)$ & $E_{L_2}$ & $\mathcal{O}(E)$ & $(\tensor{G}_{e_{11}})_{L_2}$ & $\mathcal{O}(\tensor{G}_{e_{11}})$ & $(\tensor{T}_{11})_{L_2}$ & $\mathcal{O}(\tensor{T}_{11})$ \\ 
					\hline
					1.56E-01  & 6.928E-02 & -    & 3.926E-03 & -    & 1.815E-04 & -    & 1.126E+03 & -    \\ 
					7.78E-02  & 5.291E-02 & 0.39 & 3.100E-03 & 0.34 & 1.093E-04 & 0.73 & 6.707E+02 & 0.74 \\ 
					5.46E-02  & 4.134E-02 & 0.70 & 2.712E-03 & 0.38 & 6.766E-05 & 1.35 & 4.225E+02 & 1.30 \\ 
					3.92E-02  & 3.374E-02 & 0.61 & 2.338E-03 & 0.45 & 4.620E-05 & 1.15 & 2.958E+02 & 1.08 \\ 
					\multicolumn{9}{c}{}\\
					\multicolumn{9}{c}{2D LGPR $\mathcal{O}\,2$ ($\tau_1=10^{14}$)} \\
					\hline
					$h(\Omega(t_f))$ & $u_{L_2}$ & $\mathcal{O}(u)$ & $E_{L_2}$ & $\mathcal{O}(E)$ & $(\tensor{G}_{e_{11}})_{L_2}$ & $\mathcal{O}(\tensor{G}_{e_{11}})$ & $(\tensor{T}_{11})_{L_2}$ & $\mathcal{O}(\tensor{T}_{11})$ \\ 
					\hline
					1.56E-01  & 1.377E-02 & -    & 1.505E-03 & -    & 1.696E-04 & -    & 1.303E+03 & -    \\ 
					7.78E-02  & 2.888E-03 & 2.24 & 2.730E-04 & 2.45 & 4.224E-04 & 1.99 & 2.845E+02 & 2.18 \\ 
					5.46E-02  & 1.131E-03 & 2.65 & 9.735E-05 & 2.91 & 1.860E-05 & 2.31 & 1.228E+02 & 2.37 \\ 
					3.92E-02  & 6.081E-04 & 1.87 & 5.080E-05 & 1.96 & 1.088E-05 & 1.62 & 6.990E+01 & 1.70 \\
					\multicolumn{9}{c}{}\\
					\multicolumn{9}{c}{}\\
					\multicolumn{9}{c}{3D LGPR $\mathcal{O}\,1$ ($\tau_1=10^{14}$)} \\
					\hline
					$h(\Omega(t_f))$ & $u_{L_2}$ & $\mathcal{O}(u)$ & $E_{L_2}$ & $\mathcal{O}(E)$ & $(\tensor{G}_{e_{11}})_{L_2}$ & $\mathcal{O}(\tensor{G}_{e_{11}})$ & $(\tensor{T}_{11})_{L_2}$ & $\mathcal{O}(\tensor{T}_{11})$ \\ 
					\hline
					1.56E-01  & 6.179E-02 & -    & 5.892E-03 & -    & 1.859E-04 & -    & 5.111E+03 & -    \\ 
					7.78E-02  & 5.190E-02 & 0.77 & 5.286E-03 & 0.48 & 1.639E-04 & 0.56 & 4.676E+03 & 0.39 \\ 
					5.46E-02  & 4.845E-02 & 0.91 & 4.971E-03 & 0.81 & 1.603E-04 & 0.30 & 4.631E+03 & 0.15 \\ 
					3.92E-02  & 4.219E-02 & 1.06 & 4.433E-03 & 0.88 & 1.473E-04 & 0.65 & 4.279E+03 & 0.61 \\ 
					\multicolumn{9}{c}{}\\
					\multicolumn{9}{c}{3D LGPR $\mathcal{O}\,2$ ($\tau_1=10^{14}$)} \\
					\hline
					$h(\Omega(t_f))$ & $u_{L_2}$ & $\mathcal{O}(u)$ & $E_{L_2}$ & $\mathcal{O}(E)$ & $(\tensor{G}_{e_{11}})_{L_2}$ & $\mathcal{O}(\tensor{G}_{e_{11}})$ & $(\tensor{T}_{11})_{L_2}$ & $\mathcal{O}(\tensor{T}_{11})$ \\ 
					\hline
					1.56E-01  & 1.792E-02 & -    & 2.752E-03 & -    & 1.297E-04 & -    & 4.262E+03 & -    \\ 
					7.78E-02  & 1.171E-02 & 2.08 & 1.864E-03 & 1.91 & 1.080E-04 & 0.89 & 3.606E+03 & 0.82 \\ 
					5.46E-02  & 7.487E-03 & 1.69 & 1.223E-03 & 1.59 & 8.105E-05 & 1.09 & 2.723E+03 & 1.06 \\ 
					3.92E-02  & 6.235E-03 & 2.24 & 1.013E-03 & 2.31 & 7.160E-05 & 1.52 & 2.405E+03 & 1.52 \\
				\end{tabular}
			\end{small}
		\end{center}
		\label{tab.conv_elasticity}
	\end{table}

\subsection{Kidder problem} \label{ssec.kid}
The Kidder test case \cite{Kidder1976} describes the isentropic compression of a shell filled with perfect gas
that is initially bounded between the internal radius $r_{int}=0.9$ and the external radius 
$r_{ext}=1.0$. The initial 
condition is given in terms of the general radial coordinate $r=\sqrt{\mathbf{x}^2}$ and it reads
\begin{equation}
	\left( \begin{array}{c} \rho(0,r) \\ \mathbf{v}(0,r) \\ p(0,r) \end{array}  \right) = \left( 
	\begin{array}{c}  
	\left(\frac{r_{ext}^2-r^2}{r_{ext}^2-r_{int}^2}\rho_{int}^{\gamma-1}+\frac{r^2-r_{int}^2}{r_{ext}^2-r_{ext}^2}\rho_{ext}^{\gamma-1}\right)^{-\frac{1}{\gamma-1}}
	 \\ 0 \\ s\rho(0,r)^\gamma \end{array}  \right), 
	\label{eqn.KidderIC}
\end{equation}
where $\rho_{int}=1$ and $\rho_{ext}=2$ are the initial values of density at the corresponding 
frontier of 
the shell, while $\gamma=5/3$ is the ratio of specific heats for the ideal gas EOS. The initial 
entropy distribution is assumed to be uniform, hence $s=\frac{p}{\rho^\gamma}=1$. The final time 
$t_f=\frac{\sqrt{3}}{2}\tilde{t}$ is determined in such a way that the shell is bounded by $0.45 
\leq r \leq 0.5$, which provides a reference solution for the shell configuration at the end of the 
simulation. The focalisation time is given by $\tilde{t} = 
\sqrt{\frac{\gamma-1}{2}\frac{(r_{ext}^2-r_{int}^2)}{c_{ext}^2-c_{int}^2}}$,
with the internal and external sound speeds 
$c_{int}=\sqrt{\gamma\frac{p_{int}}{\rho_{int}}}$ and
$c_{ext}=\sqrt{\gamma\frac{p_{ext}}{\rho_{ext}}}$. Pressure boundary conditions are set on 
the internal and external frontiers of the shell according to the exact solution available in 
\cite{Kidder1976}. Since this is a smooth problem, the numerical convergence is studied again on a 
sequence of successively refined computational meshes. The errors are measured in $L_\infty$ norm 
and reported in Table \ref{tab.conv_kidder} for the second order version of the LGPR scheme, where 
the formal accuracy is reached in 2D/3D. The corresponding errors related to the internal and 
external frontier location, i.e. $R_{int}(t)$ and $R_{ext}(t)$, are estimated as the arithmetic 
average of the difference between 
the analytical and the numerical radial coordinate for each node lying on the internal and external frontier. 

\begin{table}[!htbp]  
	\caption{Numerical convergence results for the Kidder problem in 2D and in 3D using the LGPR 
	scheme with relaxation time $\tau_1=10^{-14}$ on a sequence of refined unstructured triangular 
	and tetrahedral meshes of size $h(\Omega(t_f))$ measured at the end of the simulation $t_f$. 
	The errors are measured in $L_\infty$ norm for specific volume $1/\rho$, total energy $E$, 
	final position of the internal $R_{int}$ and external $R_{ext}$ radius of the shell.}  
	\begin{center} 
		\begin{small}
			\renewcommand{\arraystretch}{1.0}
			\begin{tabular}{c|ccccccccc}
				\multicolumn{9}{c}{2D LGPR $\mathcal{O}\,2$} \\
				\hline
				$h(\Omega(t_f))$ & $(1/\rho)_{L_\infty}$ & $\mathcal{O}(1/\rho)$ & $E_{L_\infty}$ & $\mathcal{O}(E)$ & $(R_{int})_{L_\infty}$ & $\mathcal{O}(R_{int})$ & $(R_{ext})_{L_\infty}$ & $\mathcal{O}(R_{ext})$ \\ 
				\hline
				2.15E-02  & 7.648E-02 & -    & 5.302E+00 & -    & 5.538E-03 & -    & 5.001E-03 & -    \\ 
				1.42E-02  & 2.900E-02 & 2.36 & 2.316E+00 & 2.01 & 2.061E-03 & 2.40 & 2.122E-03 & 2.08 \\ 
				7.95E-03  & 1.356E-02 & 1.30 & 4.888E-01 & 2.67 & 5.991E-04 & 2.12 & 7.702E-04 & 1.74 \\ 
				4.14E-03  & 3.836E-03 & 1.94 & 1.670E-01 & 1.65 & 2.145E-04 & 1.57 & 2.528E-04 & 1.71 \\ 
				\multicolumn{9}{c}{}\\
				\multicolumn{9}{c}{3D LGPR $\mathcal{O}\,2$} \\
				\hline
				$h(\Omega(t_f))$ & $(1/\rho)_{L_\infty}$ & $\mathcal{O}(1/\rho)$ & $E_{L_\infty}$ & $\mathcal{O}(E)$ & $(R_{int})_{L_\infty}$ & $\mathcal{O}(R_{int})$ & $(R_{ext})_{L_\infty}$ & $\mathcal{O}(R_{ext})$ \\ 
				\hline
				2.15E-02  & 4.150E+00 & -    & 2.162E+01 & -    & 9.825E-02 & -    & 1.810E-02 & -    \\ 
				1.42E-02  & 2.982E+00 & 2.27 & 1.757E+01 & 1.42 & 5.620E-02 & 3.84 & 1.101E-02 & 3.41 \\ 
				7.95E-03  & 1.327E-01 & 7.22 & 9.911E+00 & 1.33 & 3.009E-02 & 1.45 & 7.078E-03 & 1.03 \\ 
				4.14E-03  & 5.634E-02 & 1.80 & 3.437E+00 & 2.22 & 1.161E-02 & 2.00 & 3.343E-03 & 1.58 \\
			\end{tabular}
		\end{small}
	\end{center}
	\label{tab.conv_kidder}
\end{table}

\begin{figure}[!htbp]
		\begin{center}
			\begin{tabular}{cc}
				\includegraphics[width=0.47\textwidth,draft=false,draft=false]{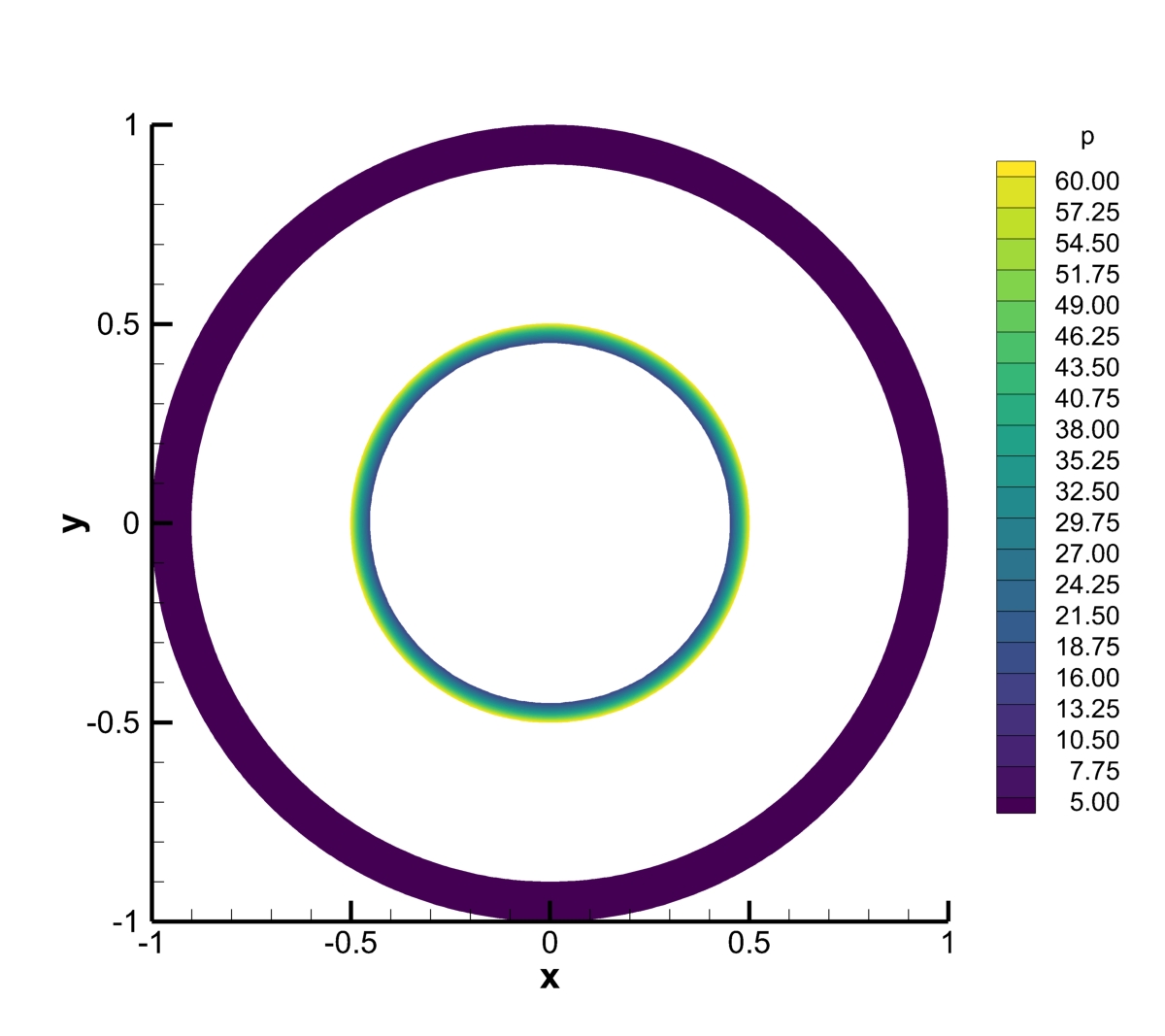}  
				&          
				\includegraphics[width=0.47\textwidth,draft=false,draft=false]{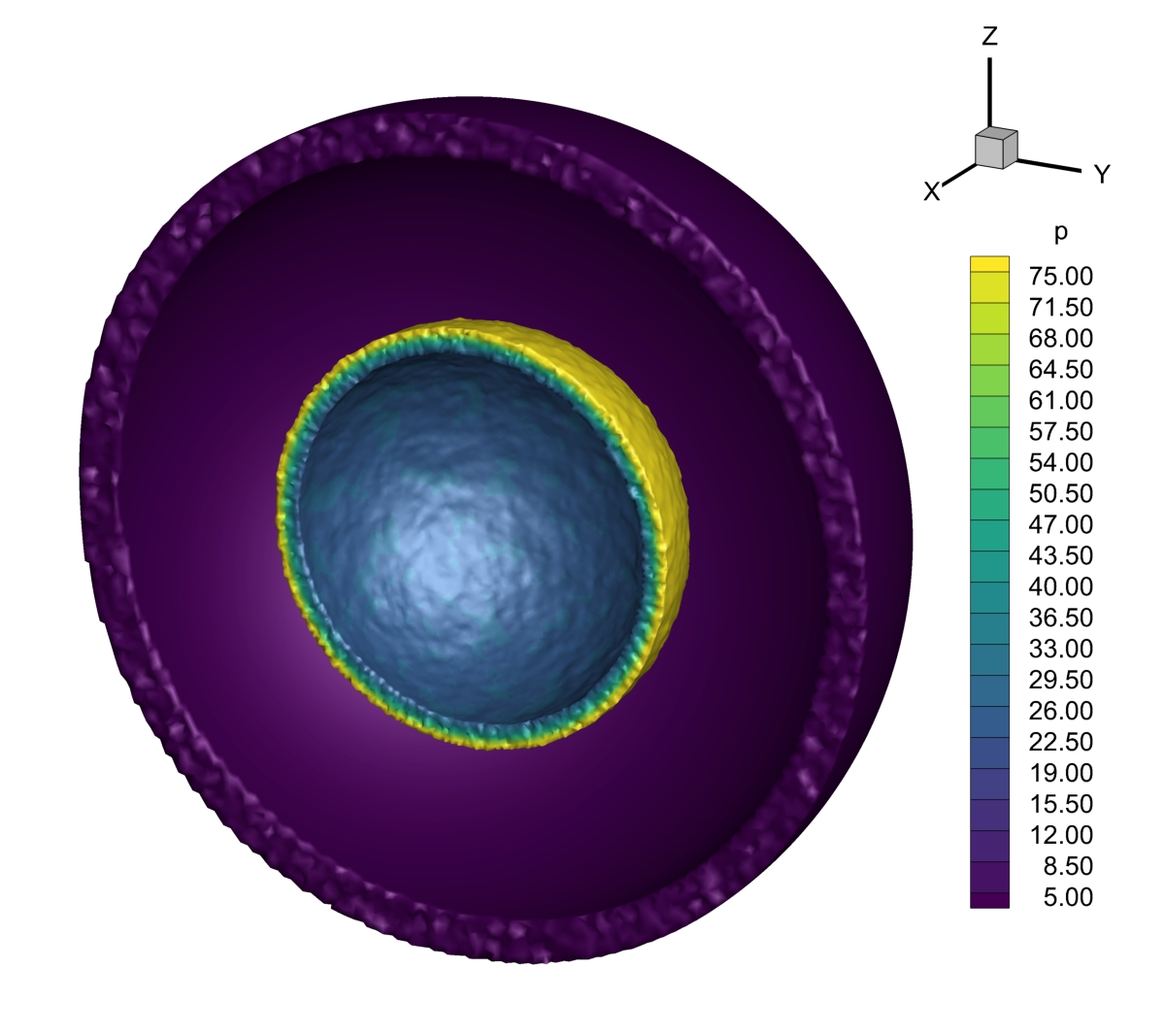} \\
				\includegraphics[width=0.47\textwidth,draft=false,draft=false]{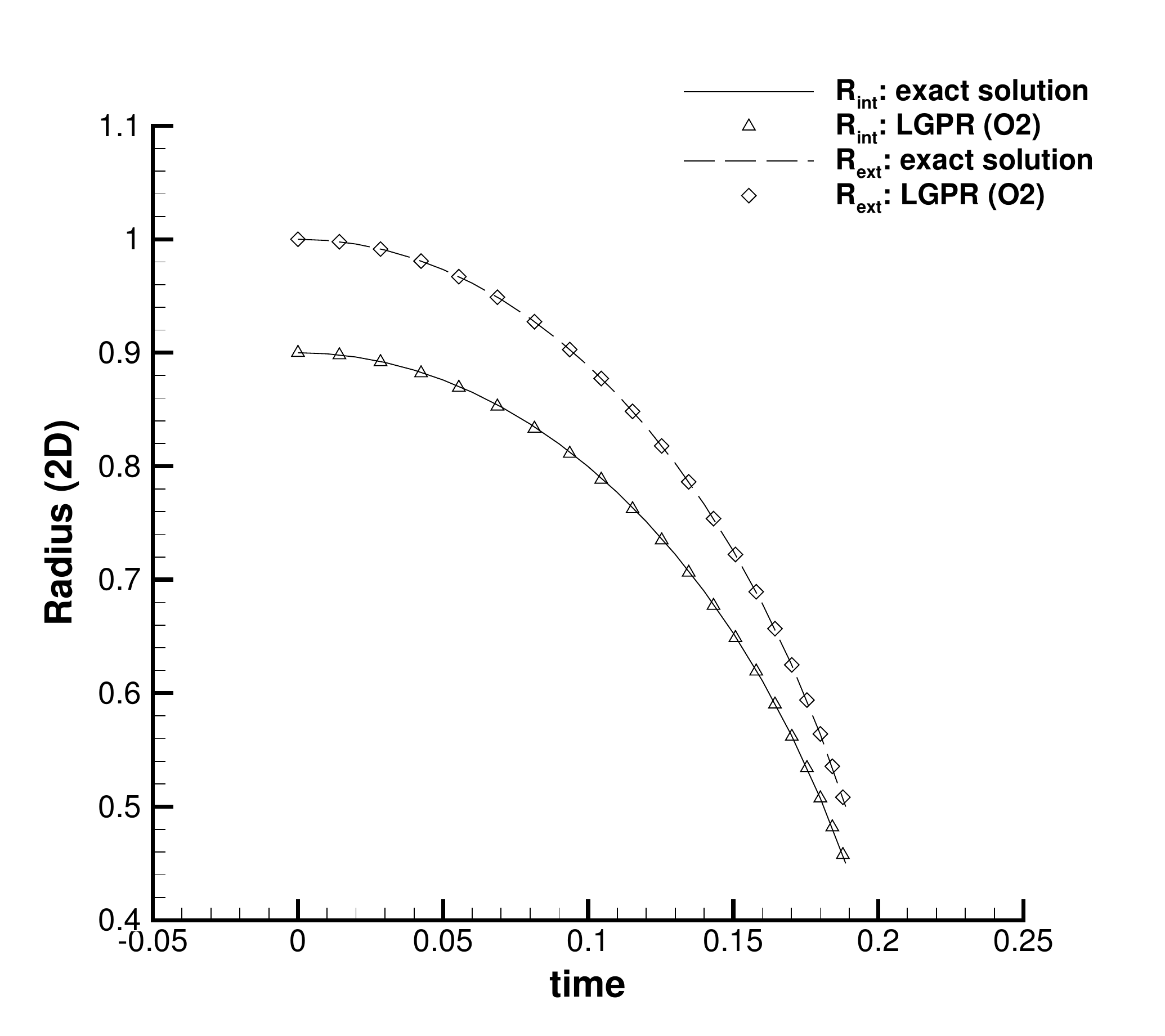}  
				&          
				\includegraphics[width=0.47\textwidth,draft=false,draft=false]{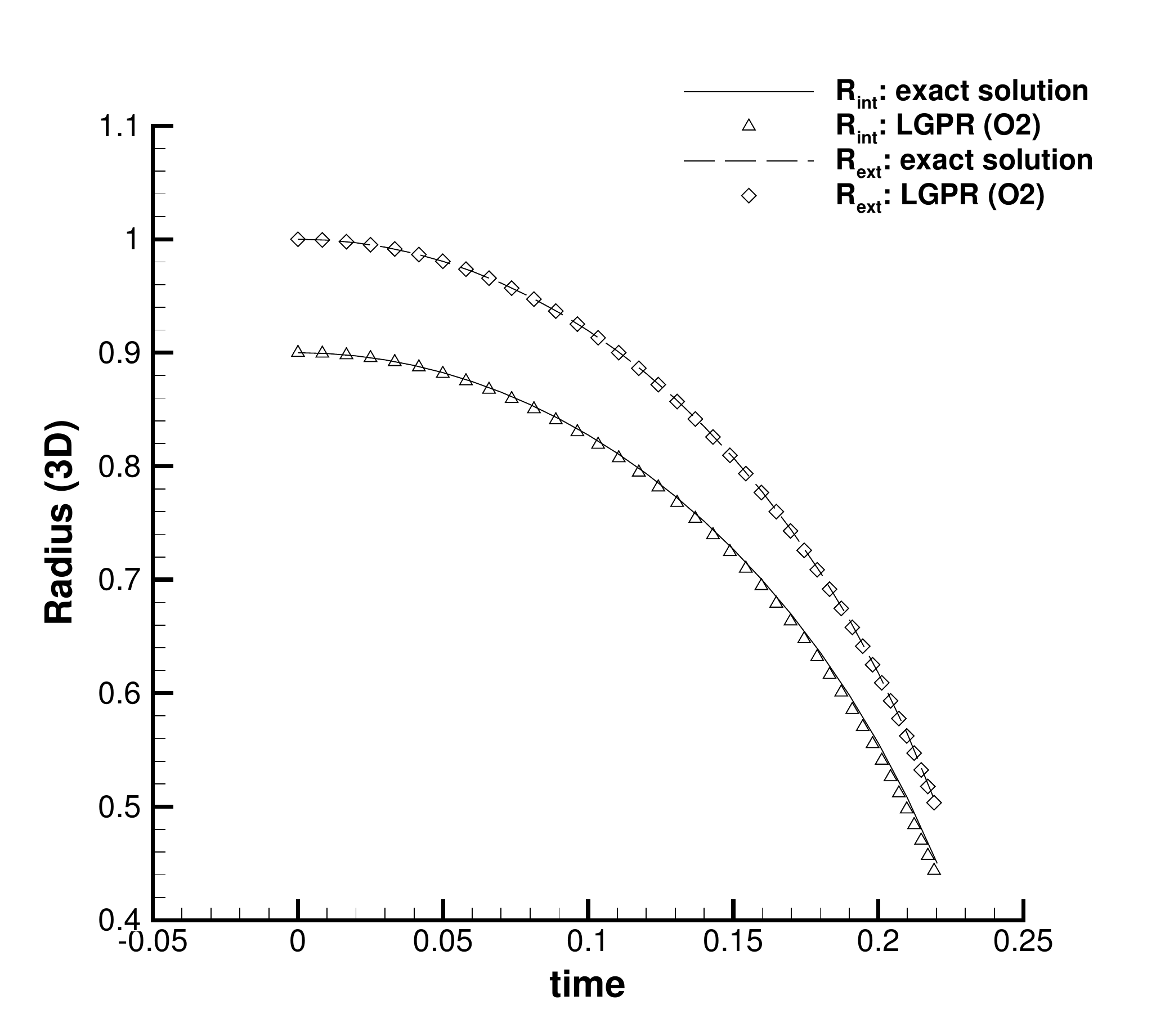} \\ 
			\end{tabular}
			\caption{Kidder problem in 2D (left) and in 3D (right). Top: pressure distribution and geometry configuration at the initial and final time of the simulation. Bottom: evolution of the internal $R_{int}$ and external $R_{ext}$ radius of the shell and comparison between analytical and numerical solution.}
			\label{fig.Kidder-radius}
		\end{center}
	\end{figure}

Figure \ref{fig.Kidder-radius} depicts the initial and final pressure distribution in 2D and 3D, which preserves a good symmetry despite the unstructured meshes. The time evolution of the internal and external frontiers of the shell is compared against the analytical solution during the entire simulation, exhibiting an excellent agreement.

\subsection{Saltzman problem} \label{ssec.Saltz}
\begin{figure}[!htbp]
	\begin{center}
		\begin{tabular}{cc}
			\includegraphics[width=0.5\textwidth,draft=false]{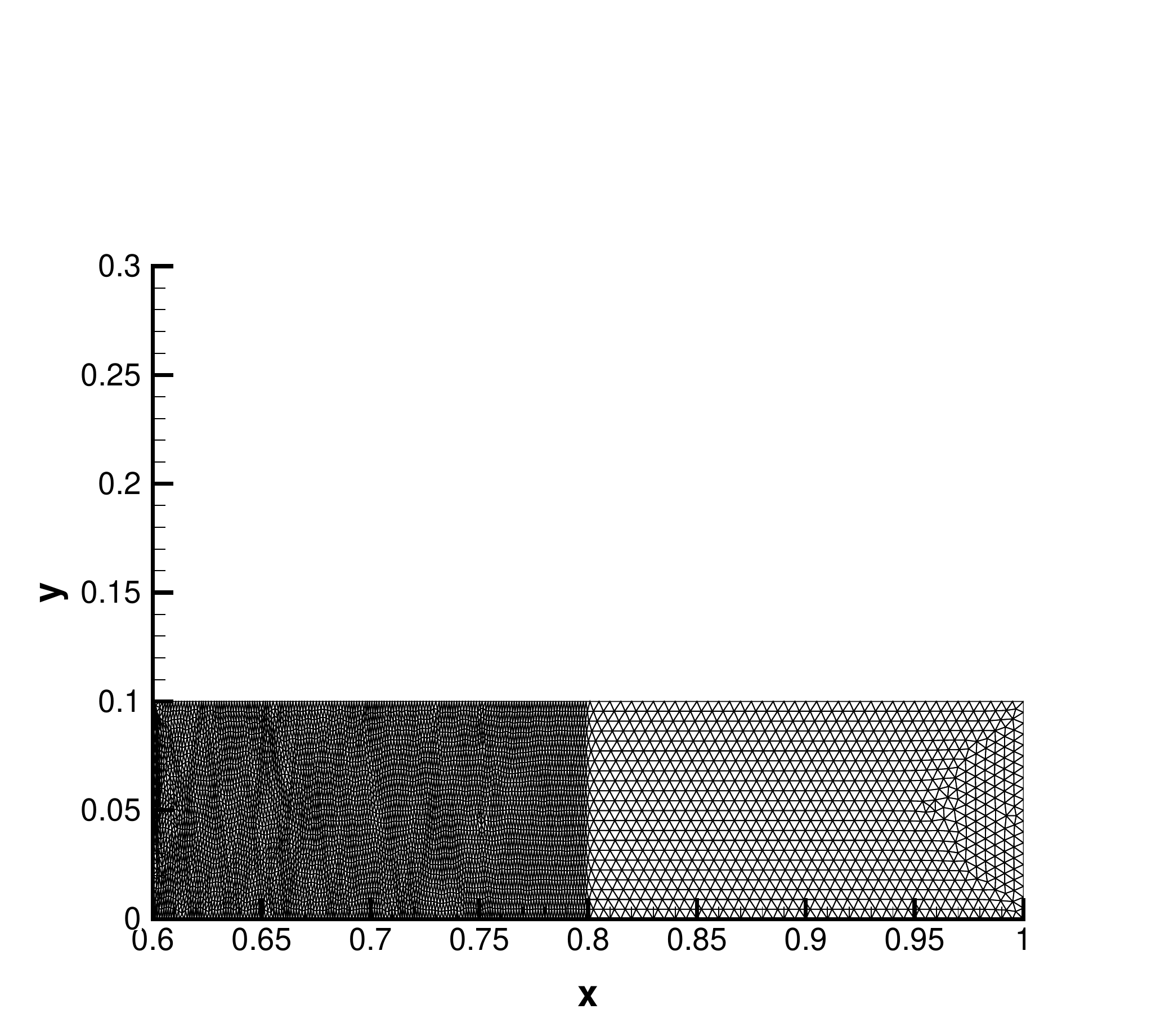}  &          
			\includegraphics[width=0.4\textwidth,draft=false]{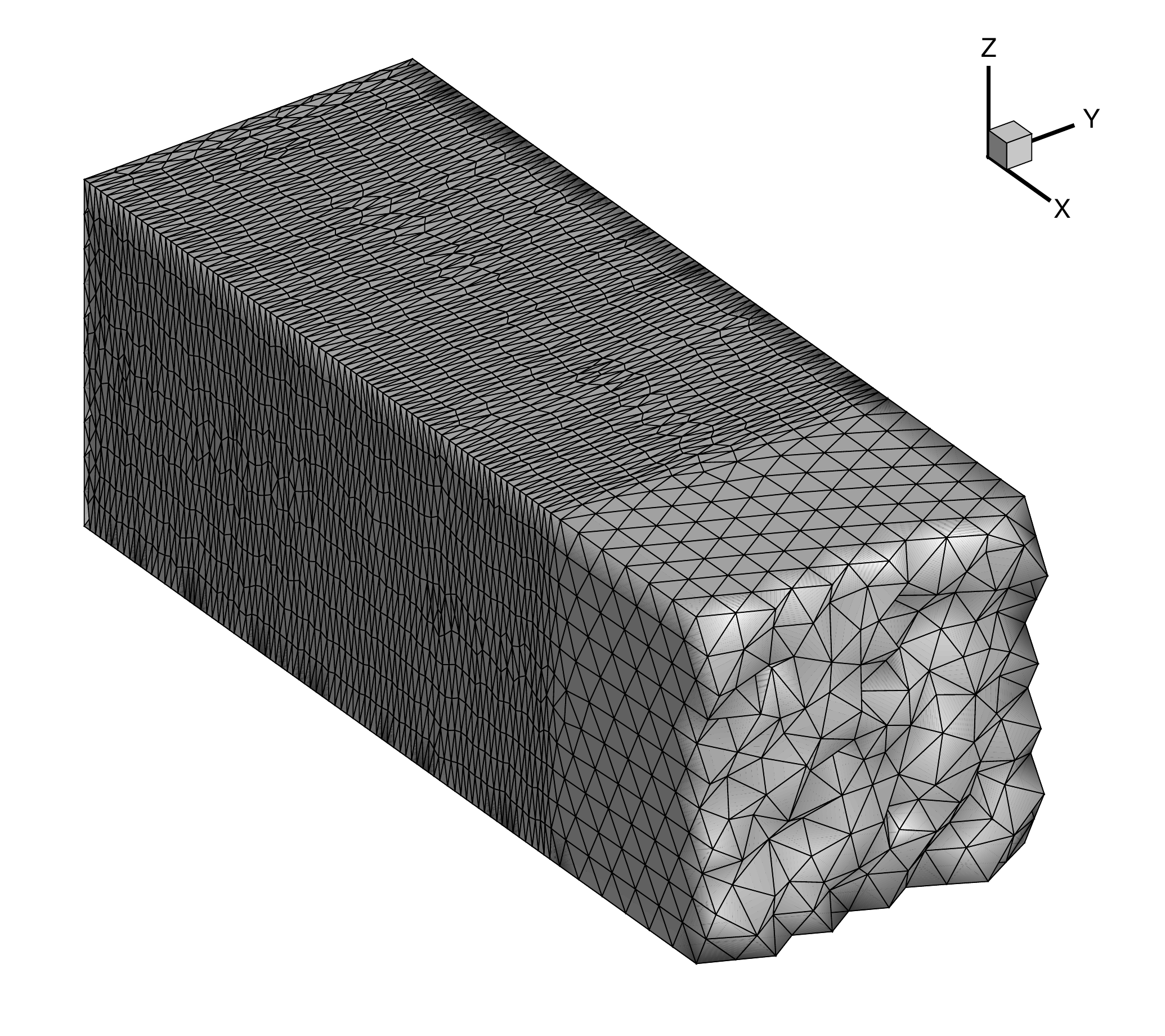} \\
		\end{tabular}
		\caption{Saltzman problem. Mesh configuration in 2D (left) and in 3D (right) at time $t=0.6$.}
		\label{fig.Saltz_mesh}
	\end{center}
\end{figure}

The Saltzman problem involves a strong shock wave that is generated by the motion of a piston traveling along the main direction of a rectangular box. This is a challenging test case first proposed in \cite{SaltzmanOrg,SaltzmanOrg3D} on a skewed quadrangular mesh to verify the robustness of Lagrangian schemes when the mesh is not aligned with the fluid flow. The initial computational domain is $\Omega^{2D}(0)=[0;1] \times [0;0.1]$ and $\Omega^{3D}(0)=[0;1] \times [0;0.1] \times [0;0.1]$, with zero velocity slip-wall boundary conditions everywhere, except for the piston, which is assigned with moving slip-wall boundary condition.	The usage of unstructured meshes allows any prescribed skewness to be ignored, since the computational mesh does not exhibit any face aligned with the fluid flow, as depicted in Figure \ref{fig.Saltz_mesh} at the final time of the simulation $t_f=0.6$. The characteristic mesh size is $h=1/100$ for both 2D and 3D grids.

\begin{figure}[!htbp]
	\begin{center}
		\begin{tabular}{cc}
			\includegraphics[width=0.47\textwidth,draft=false]{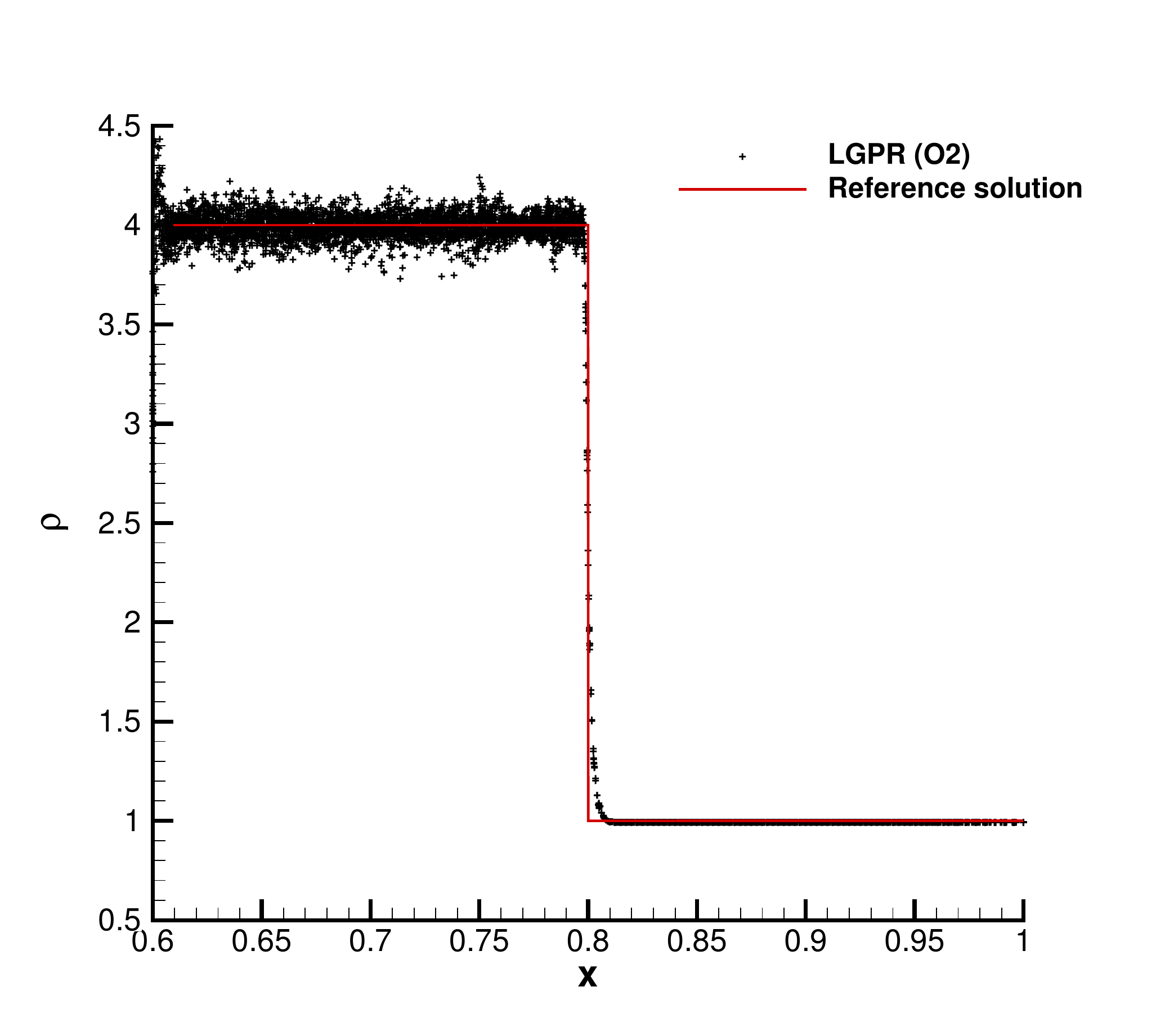}  &          
			\includegraphics[width=0.47\textwidth,draft=false]{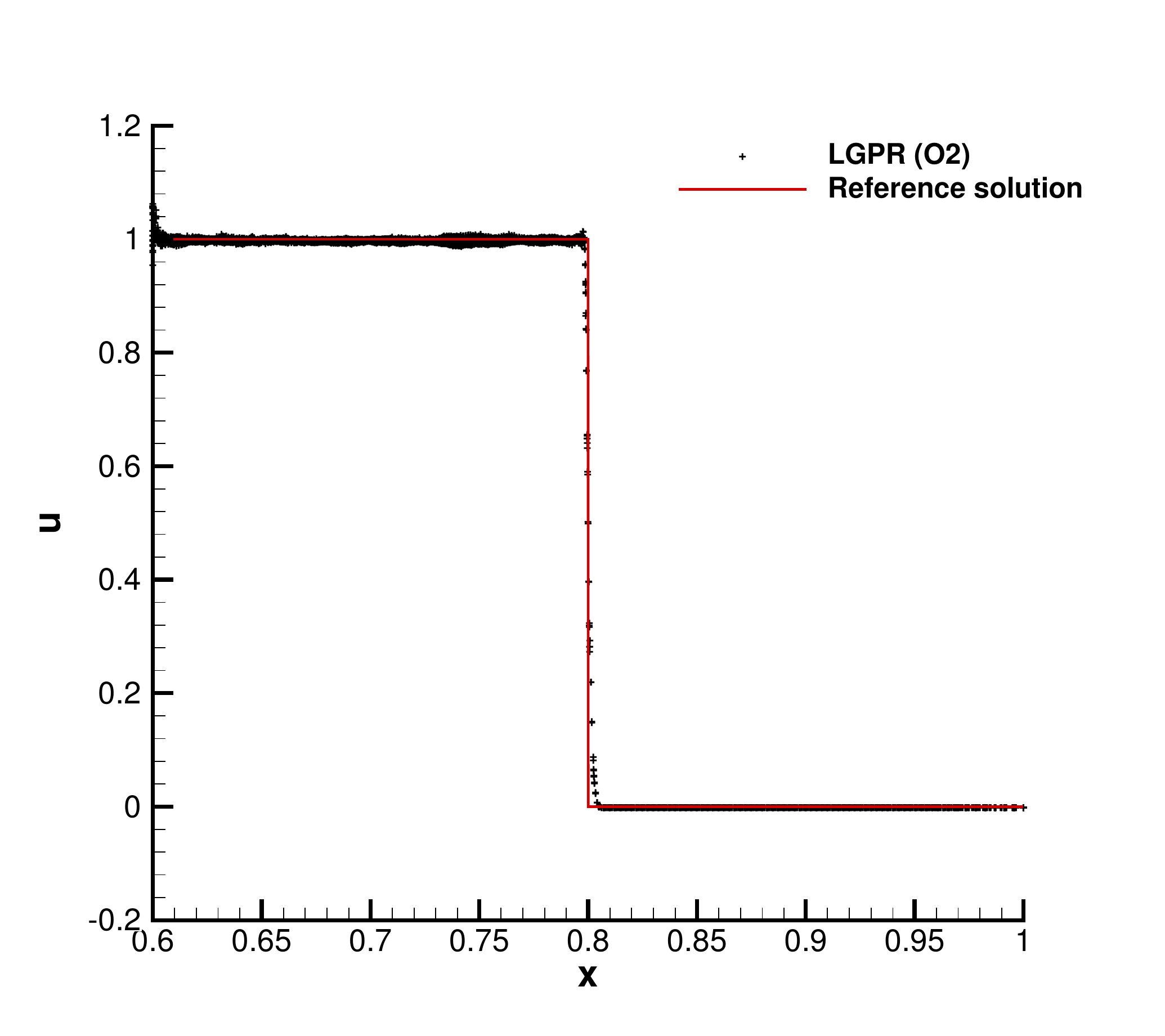} \\
			\includegraphics[width=0.47\textwidth,draft=false]{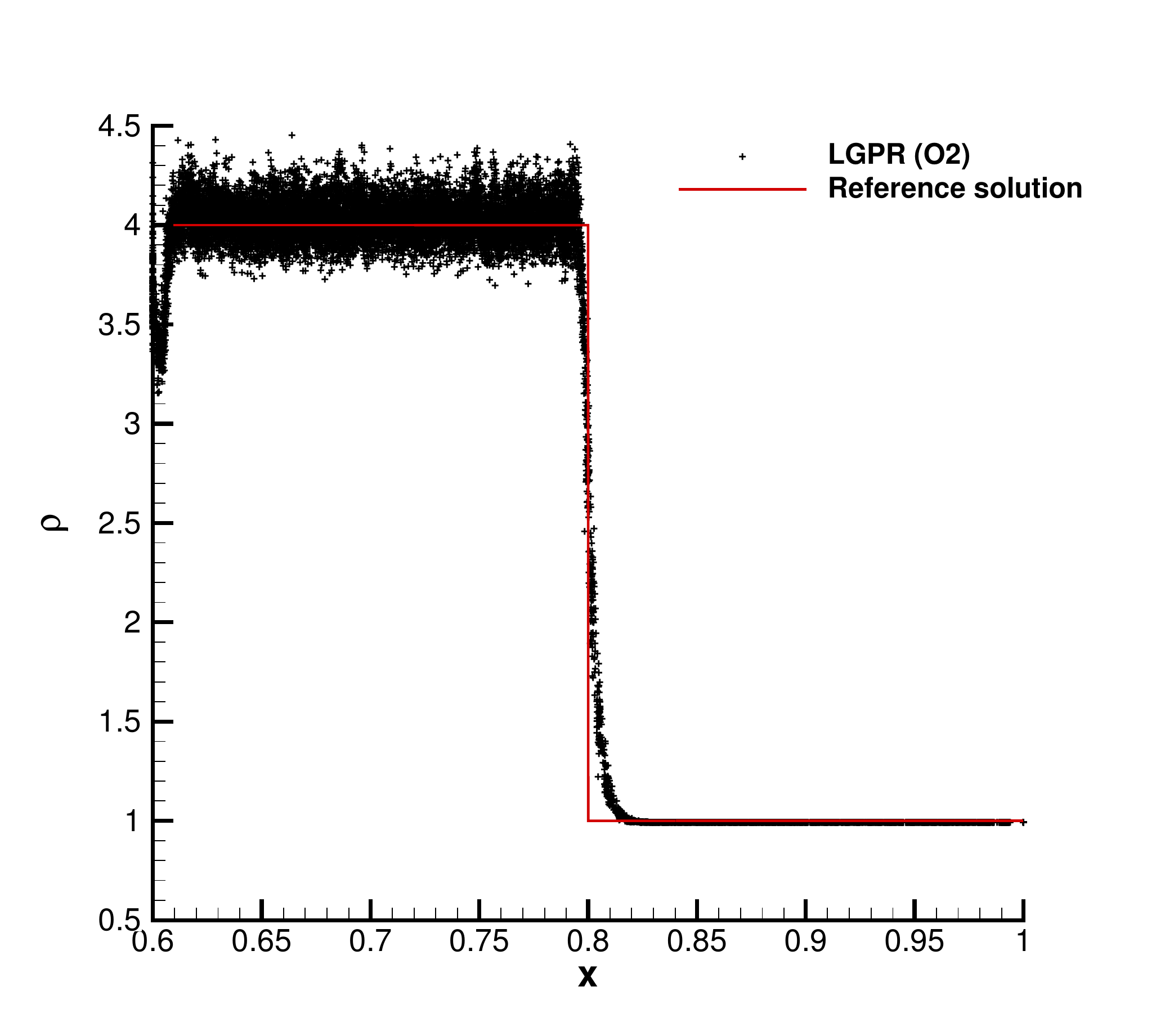}  &          
			\includegraphics[width=0.47\textwidth,draft=false]{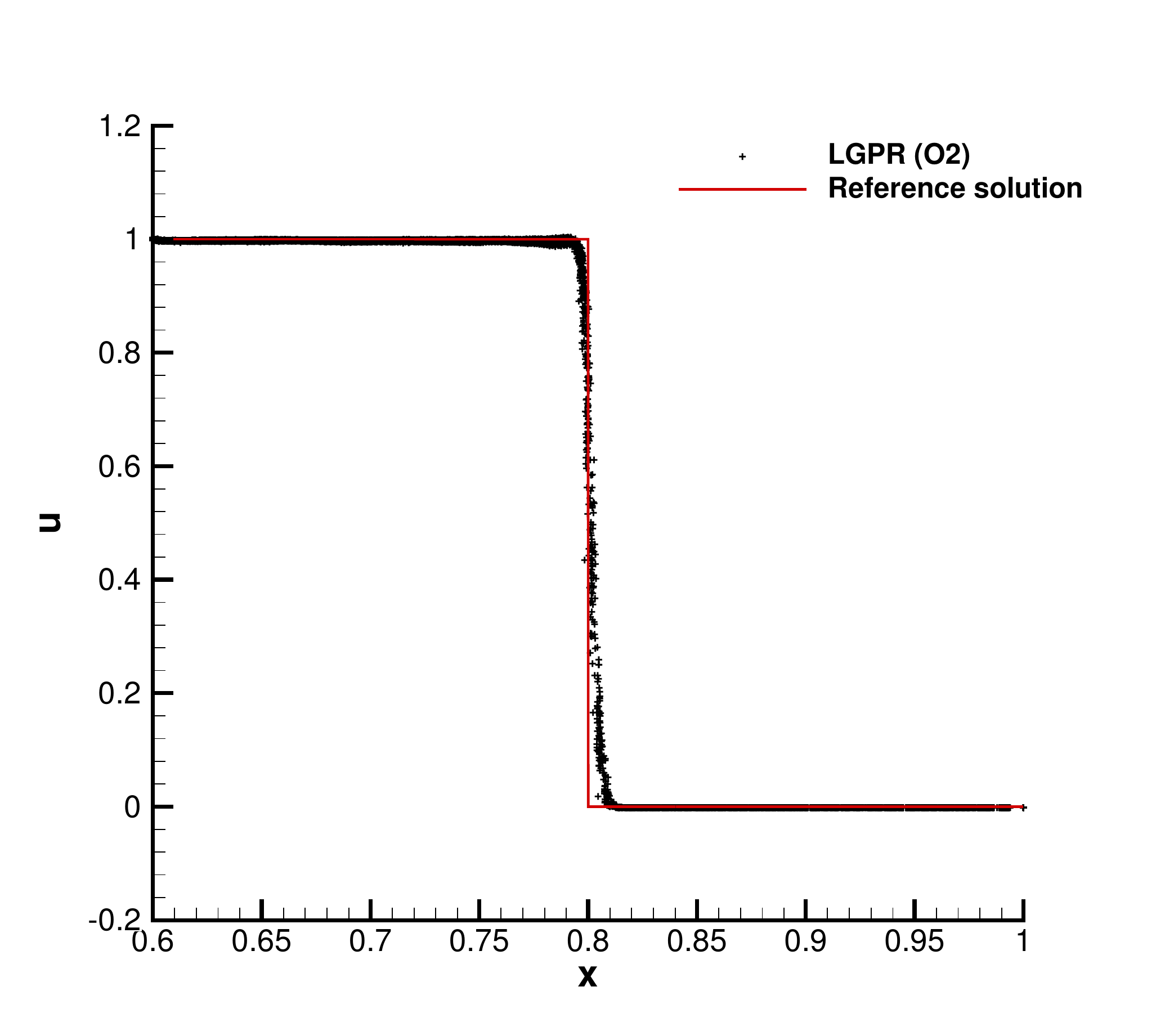} \\
		\end{tabular}
		\caption{Saltzman problem in 2D (top) and in 3D (bottom).  Scatter plot of cell density (left) and horizontal velocity (right) compared against the reference solution at time $t=0.6$.}
		\label{fig.Saltz_scatter}
	\end{center}
\end{figure}

The initial pressure of the fluid is $p=10^{-6}$ and the ideal gas EOS is used with $\gamma=5/3$. 
The piston lies on the left side of the domain and moves with velocity $\v_p=(1,0,0)$. The CFL 
number is set to $\textnormal{CFL}=0.01$ up to time $t=0.01$ in order to prevent the generation of 
invalid elements in those cells lying near the piston that are instantaneously highly compressed 
for $t>0$. The numerical results are shown in Figure \ref{fig.Saltz_scatter} as a scatter plot of 
all cell values. Compared to the analytical solution \cite{Lagrange2D,ToroBook}, a good 
approximation of the shock plateau and the shock wave location is observed. The decrease of density 
near the piston, especially in 3D, is due to the well known \textit{wall-heating problem}, see 
\cite{toro.anomalies.2002}.

\subsection{Sedov problem} \label{ssec.Sedov}
\begin{figure}[!htbp]
	\begin{center}
		\begin{tabular}{cc}
			\includegraphics[width=0.47\textwidth,draft=false]{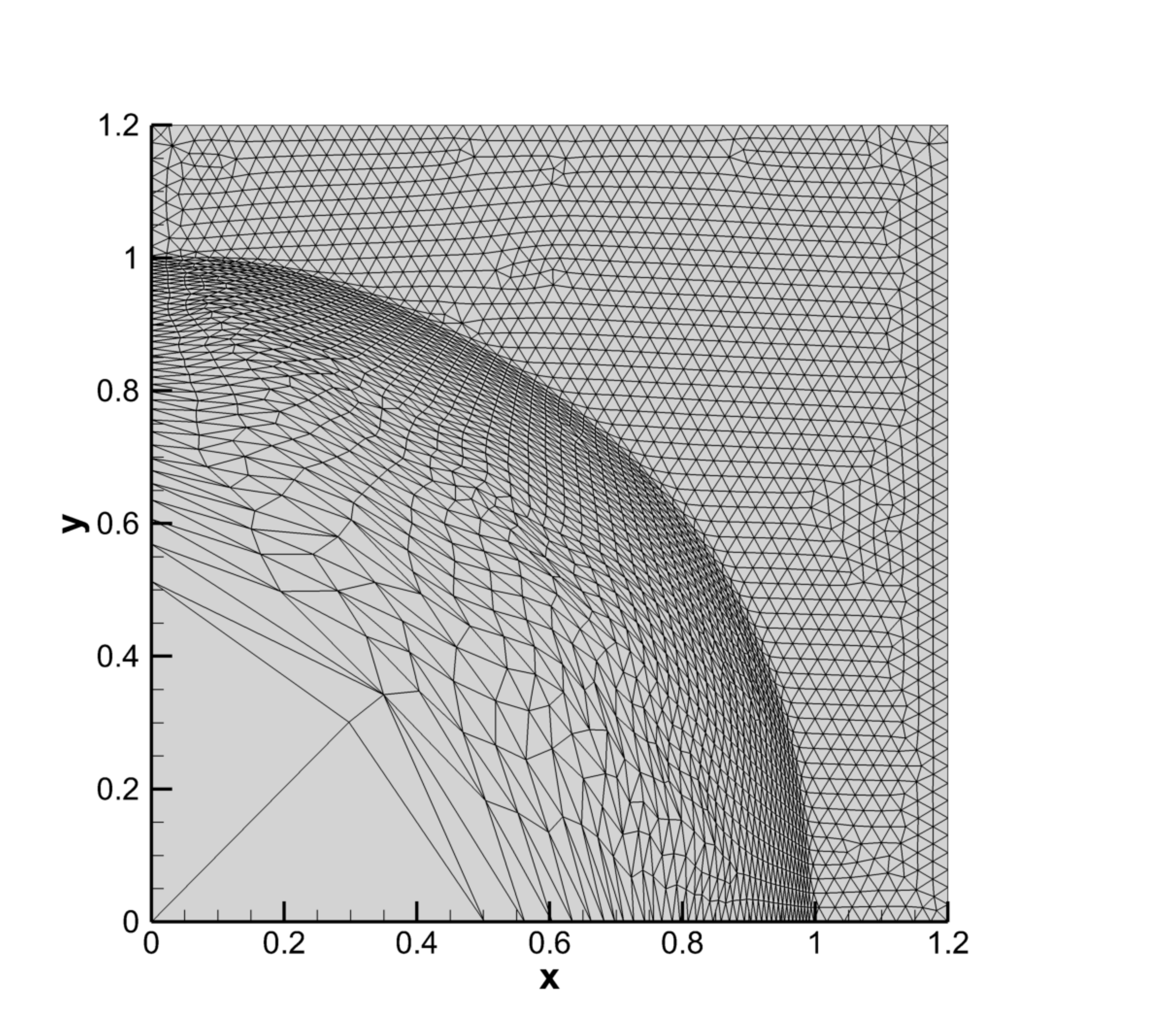}  &          
			\includegraphics[width=0.47\textwidth,draft=false]{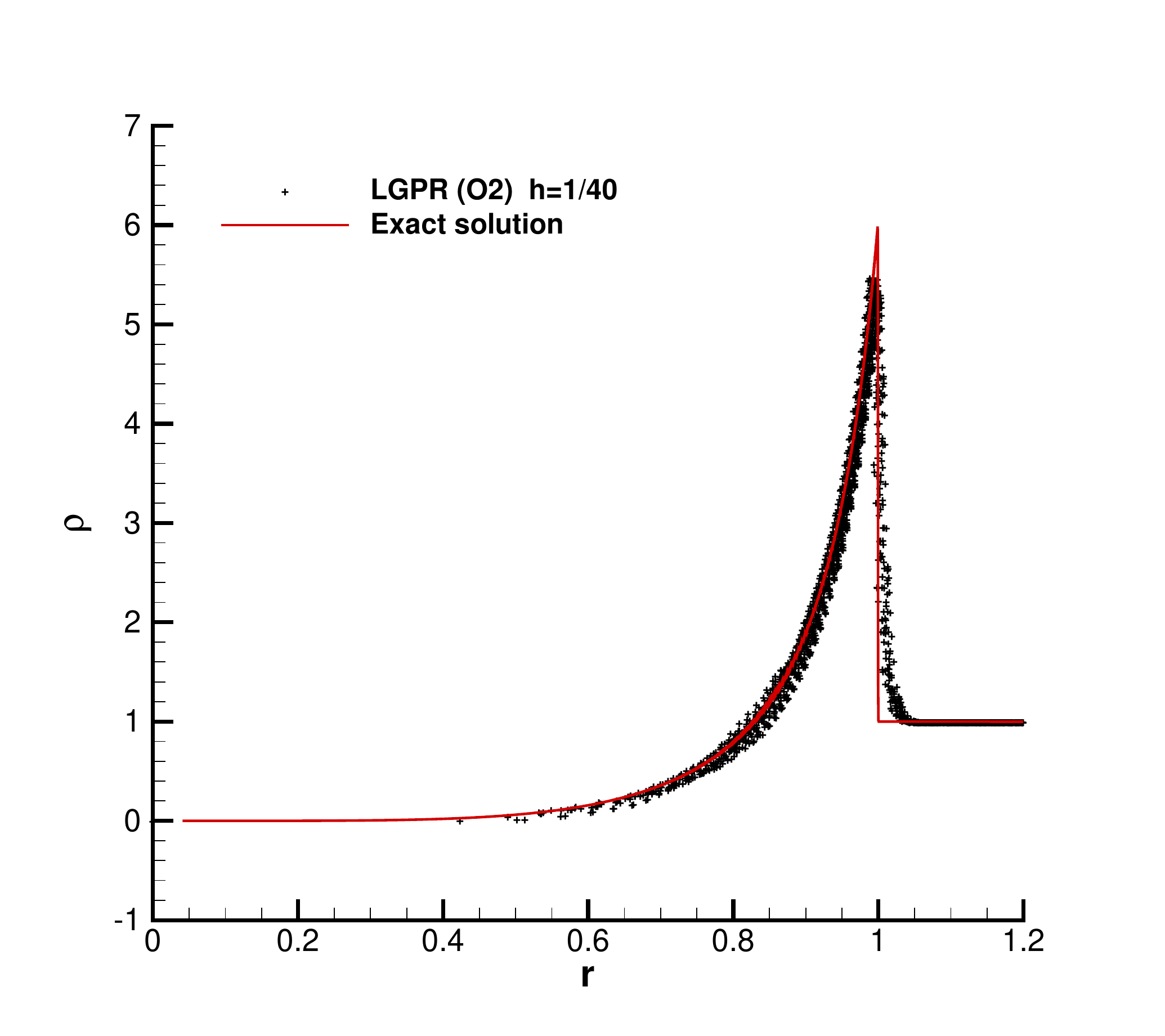} \\
			\includegraphics[width=0.47\textwidth,draft=false]{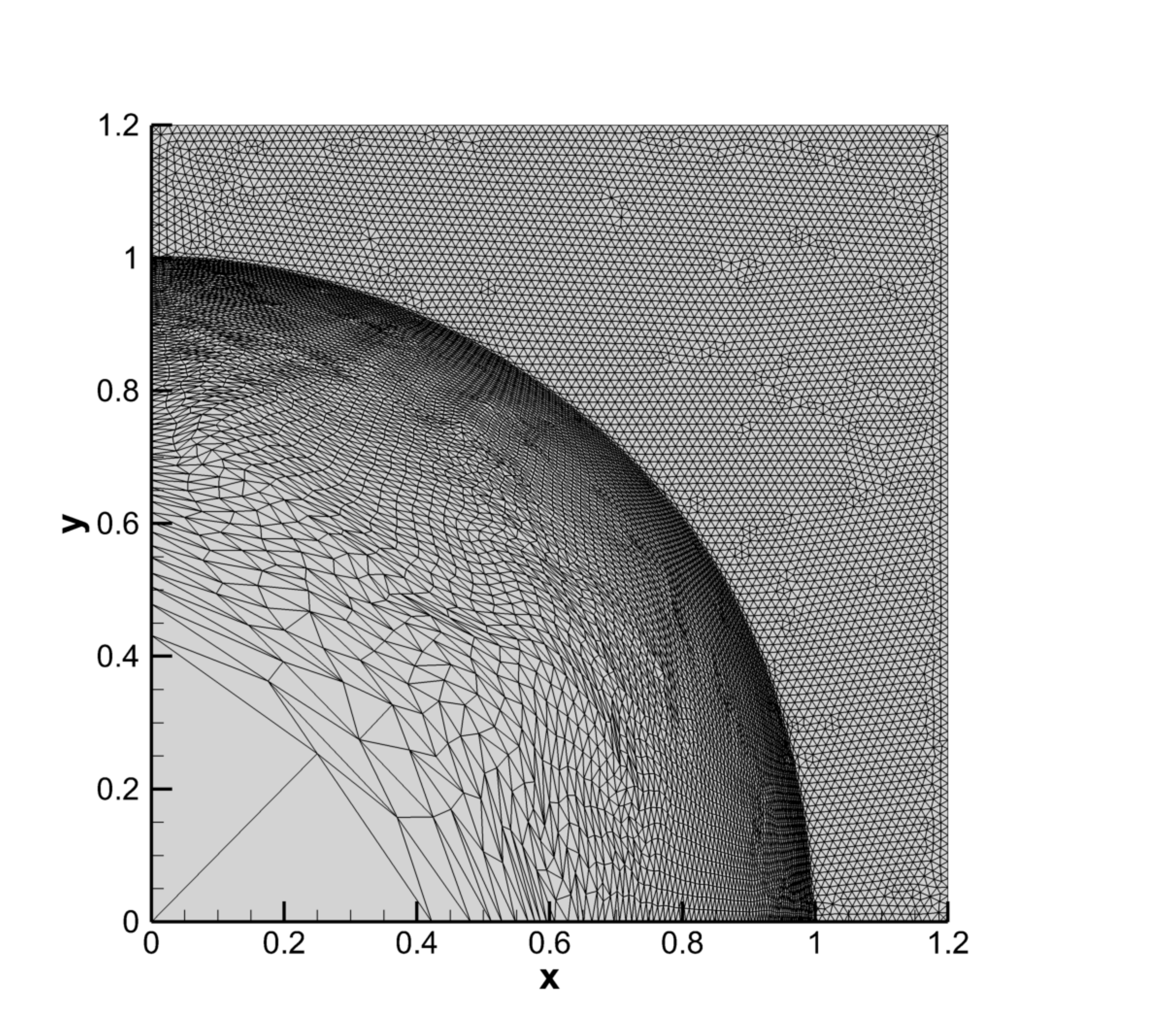}  &          
			\includegraphics[width=0.47\textwidth,draft=false]{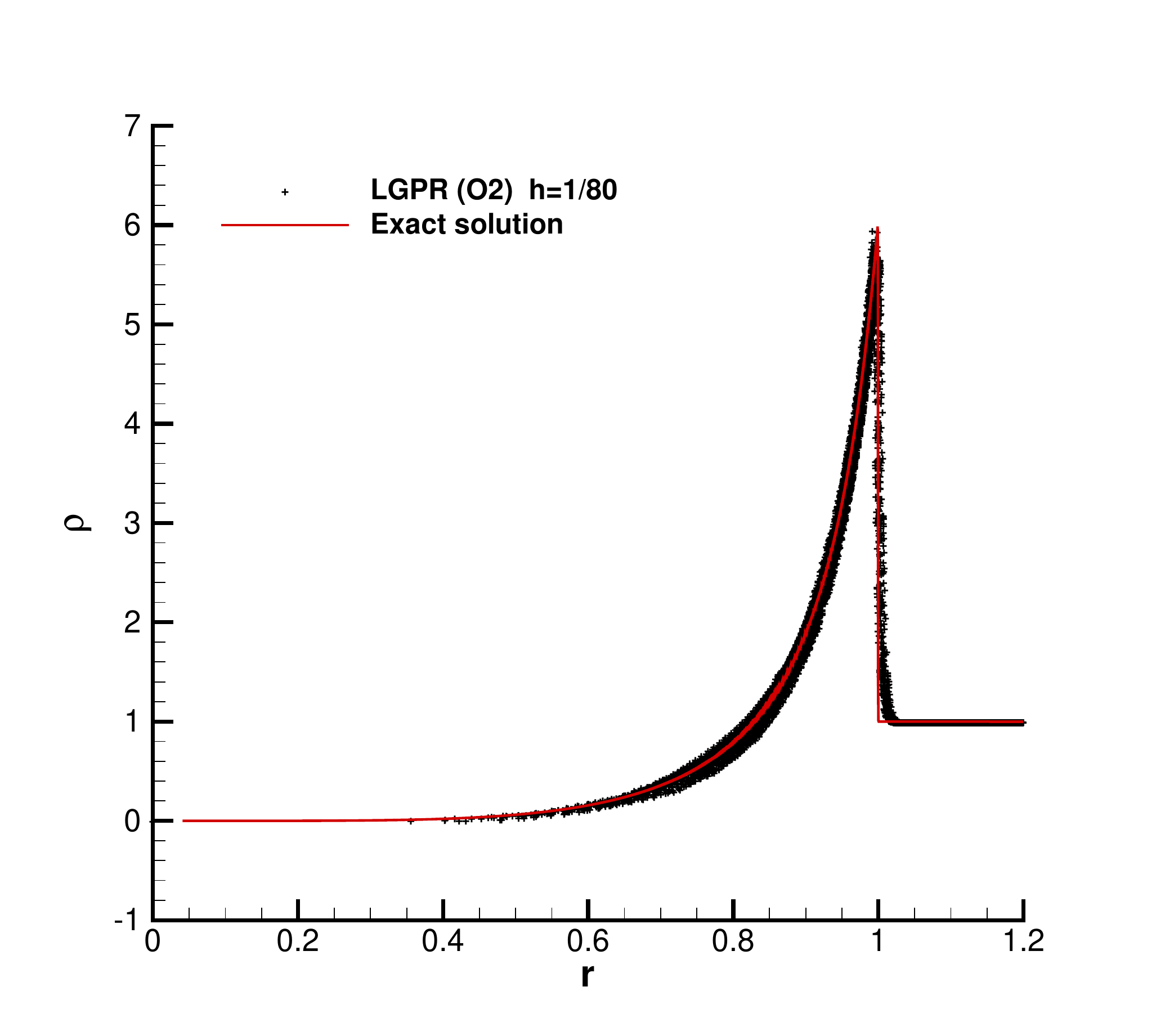} \\
			\includegraphics[width=0.47\textwidth,draft=false]{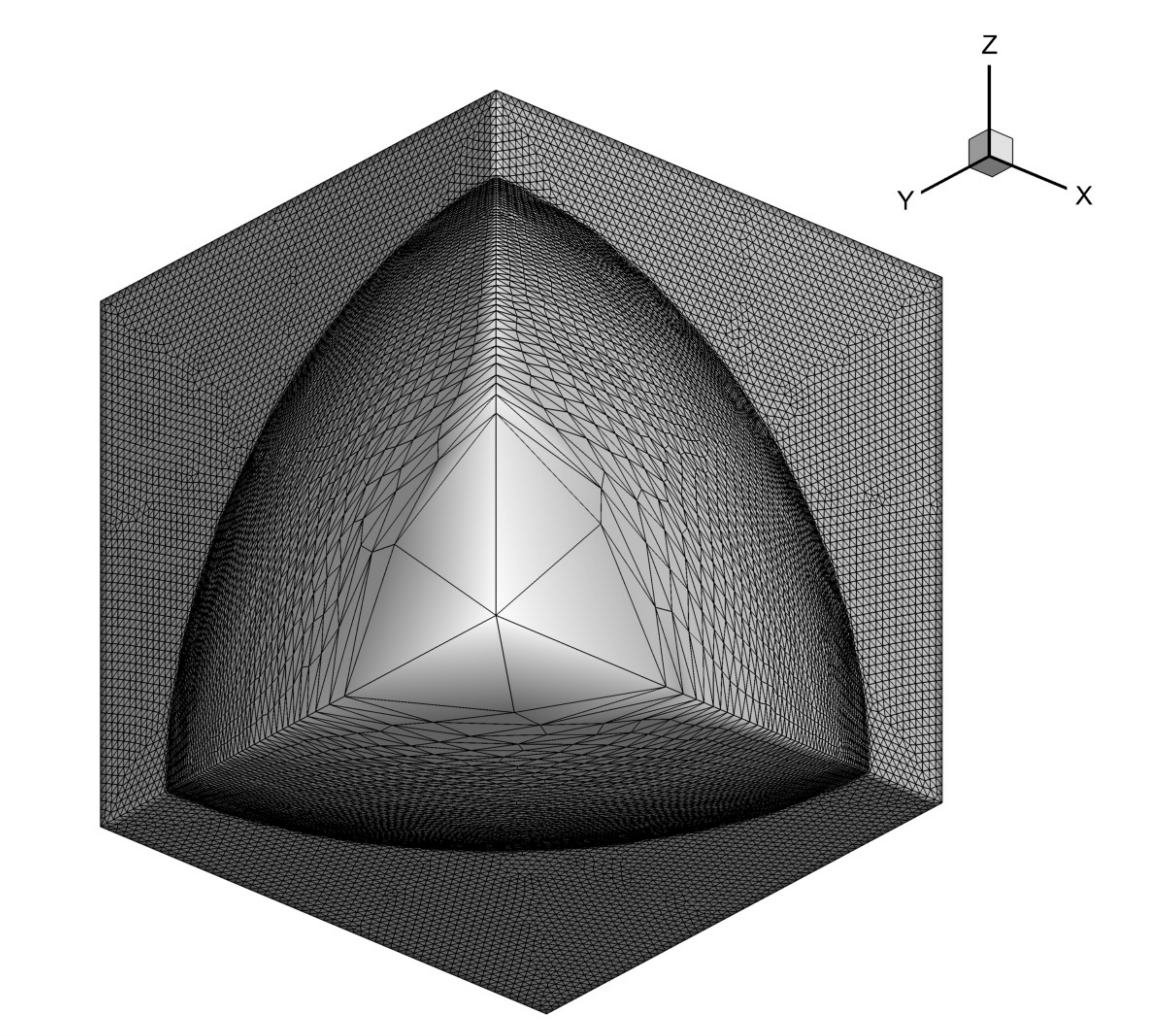}  &          
			\includegraphics[width=0.47\textwidth,draft=false]{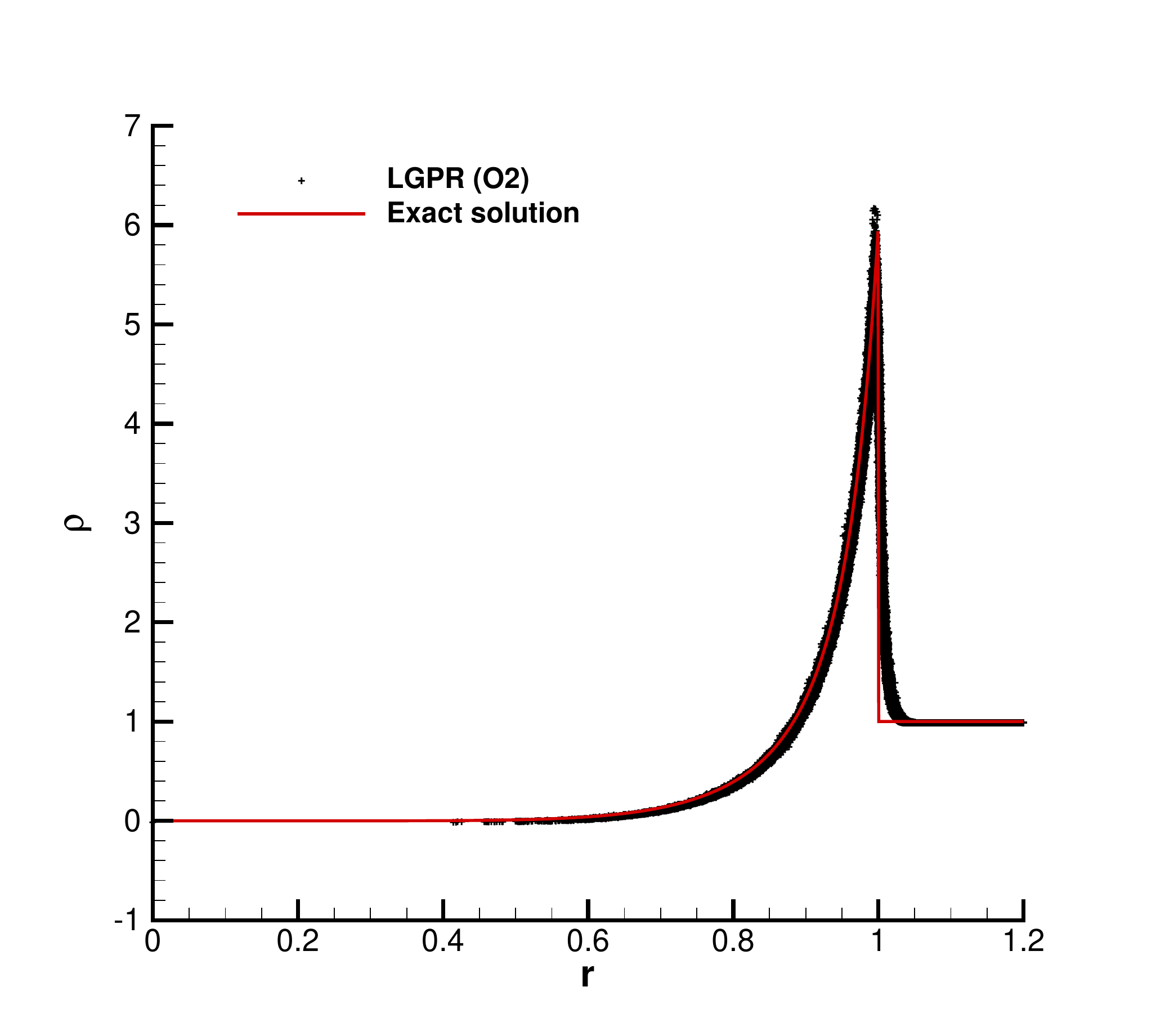} \\
		\end{tabular}
		\caption{Sedov problem. Mesh configuration (left) and scatter plot of cell density (right) at the final time $t=1$ obtained with mesh size $h=1/40$ in 2D (top),  $h=1/80$ in 2D (middle) and $h=1/60$ in 3D (bottom).}
		\label{fig.Sedov}
	\end{center}
\end{figure}

The Sedov problem is concerned with the evolution of a blast wave with cylindrical or spherical symmetry, generated at the origin $\mathbf{O}=(x,y,z)=(0,0,0)$ of the initial computational domain $\Omega(0)=[0;1.2]^d$. Symmetry boundary conditions are imposed on those faces which share the origin, while the remaining sides are treated as slip-walls. The characteristic mesh size is $h=1/60$ in 3D, whereas we use two different computational meshes in 2D, namely with $h=1/40$ and $h=1/80$. The ideal gas with $\gamma=1.4$ is initially at rest with an initial pressure of $p=10^{-6}$ in the entire computational domain, apart in those cells containing the origin $\mathbf{O}$ where the pressure is prescribed as
\begin{equation}
	p_{or} = (\gamma-1)\rho_0 \frac{E_{tot}}{\alpha \,  V_{or}} \quad \textnormal{ with } \quad 
	E_{tot} = \left\{ \begin{array}{l} 0.244816 \; \textnormal{ if } d=2 \\
		0.851072 \; \textnormal{ if } d=3 \end{array}  \right., 
	\label{eqn.p0.sedov}
\end{equation}
with $V_{or}$ denoting the volume of the elements attached to the origin. The factor $\alpha$ takes into account the cylindrical or spherical symmetry and is set to $\alpha=4$ and $\alpha=8$, respectively. An analytical solution can be derived from self-similarity arguments \cite{SedovExact}, making this test widely used in literature \cite{phm109}. Figure \ref{fig.Sedov} collects the results for all three simulations, plotting the final mesh configuration as well as a scatter plot of the cell density compared against the exact solution. The shock wave is correctly captured by the conservative LGPR scheme, and a very good symmetry of the density distribution can be appreciated. The density peak is well retrieved thanks to the Lagrangian nature of the scheme, which introduces much less numerical viscosity compared to Eulerian schemes on fixed grids at the same order of accuracy.

\subsection{Riemann problems with viscous fluids} \label{ssec.Sod}
\begin{figure}[!htbp]
	\begin{center}
		\begin{tabular}{ccc}
			\includegraphics[width=0.3\textwidth,draft=false]{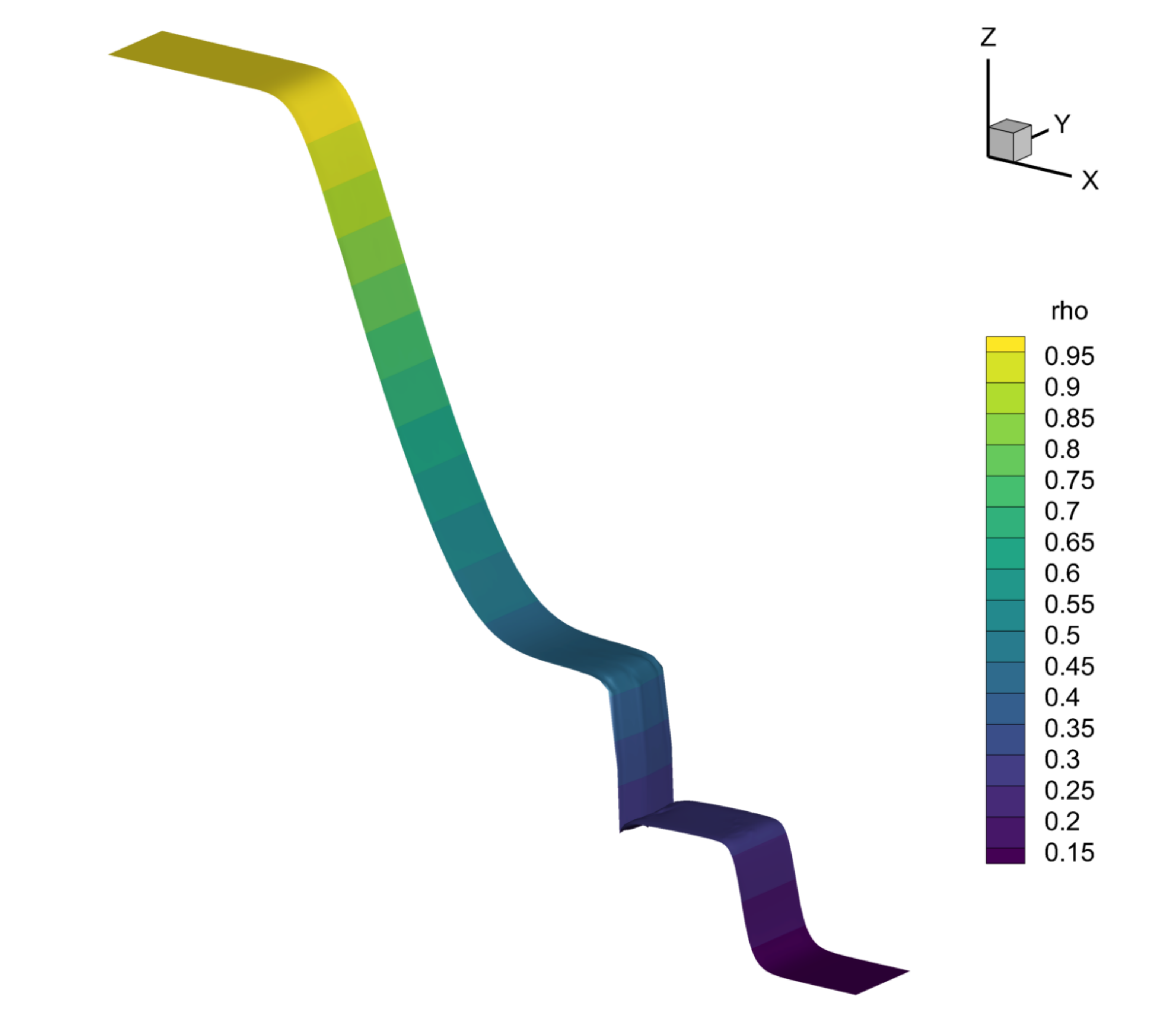}  &          
			\includegraphics[width=0.3\textwidth,draft=false]{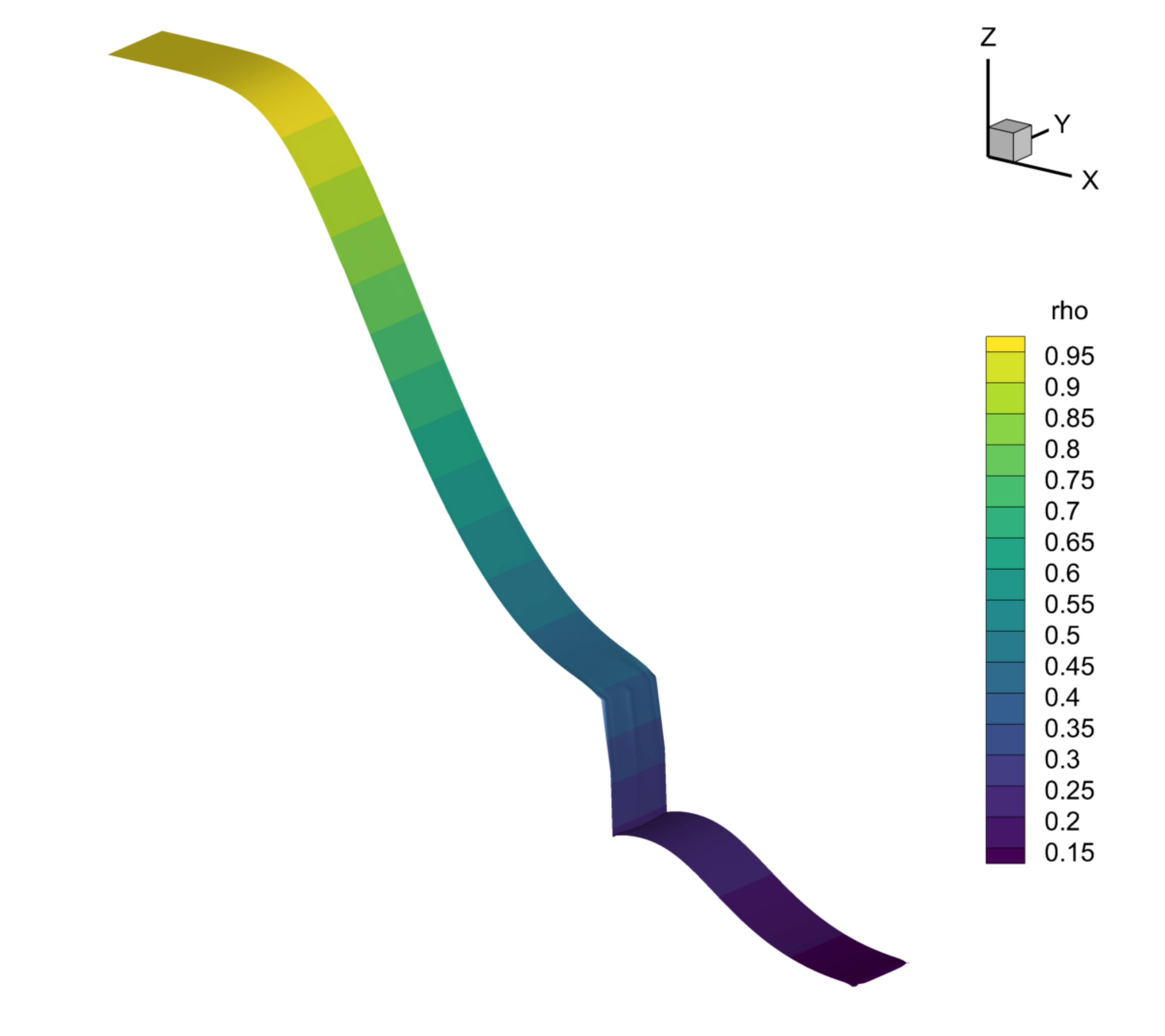} &
			\includegraphics[width=0.40\textwidth,draft=false]{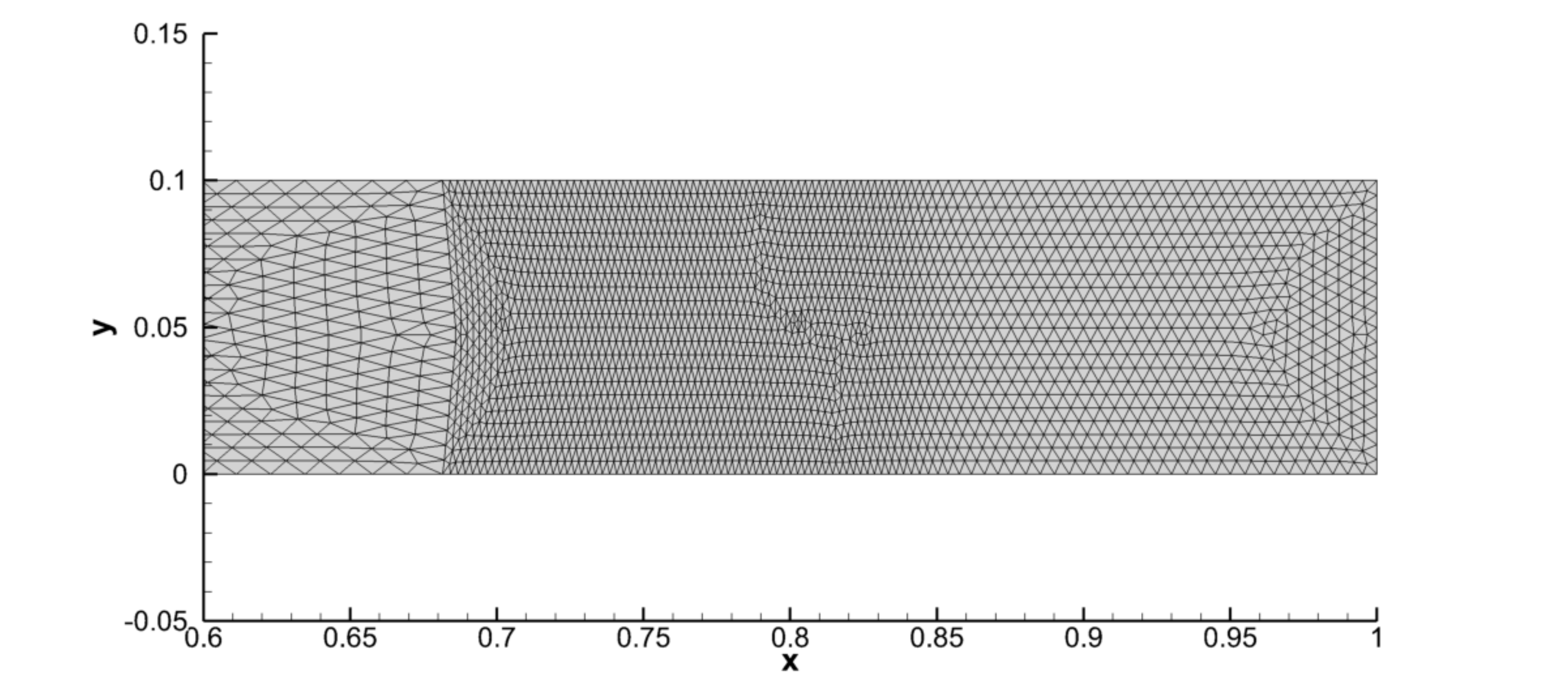} \\
		\end{tabular}
		\caption{Riemann problems with viscous fluids. Three-dimensional density distribution at the final 
		time $t=0.2$ with viscosity coefficient $\visc=10^{-3}$ (left) and $\visc=10^{-2}$ 
		(middle). A zoom in view of the mesh configuration across the contact discontinuity for 
		$\visc=10^{-3}$ (right).}
		\label{fig.Sodvisc3D}
	\end{center}
\end{figure}

\begin{figure}[!htbp]
	\begin{center}
		\begin{tabular}{ccc}
			\includegraphics[width=0.33\textwidth,draft=false]{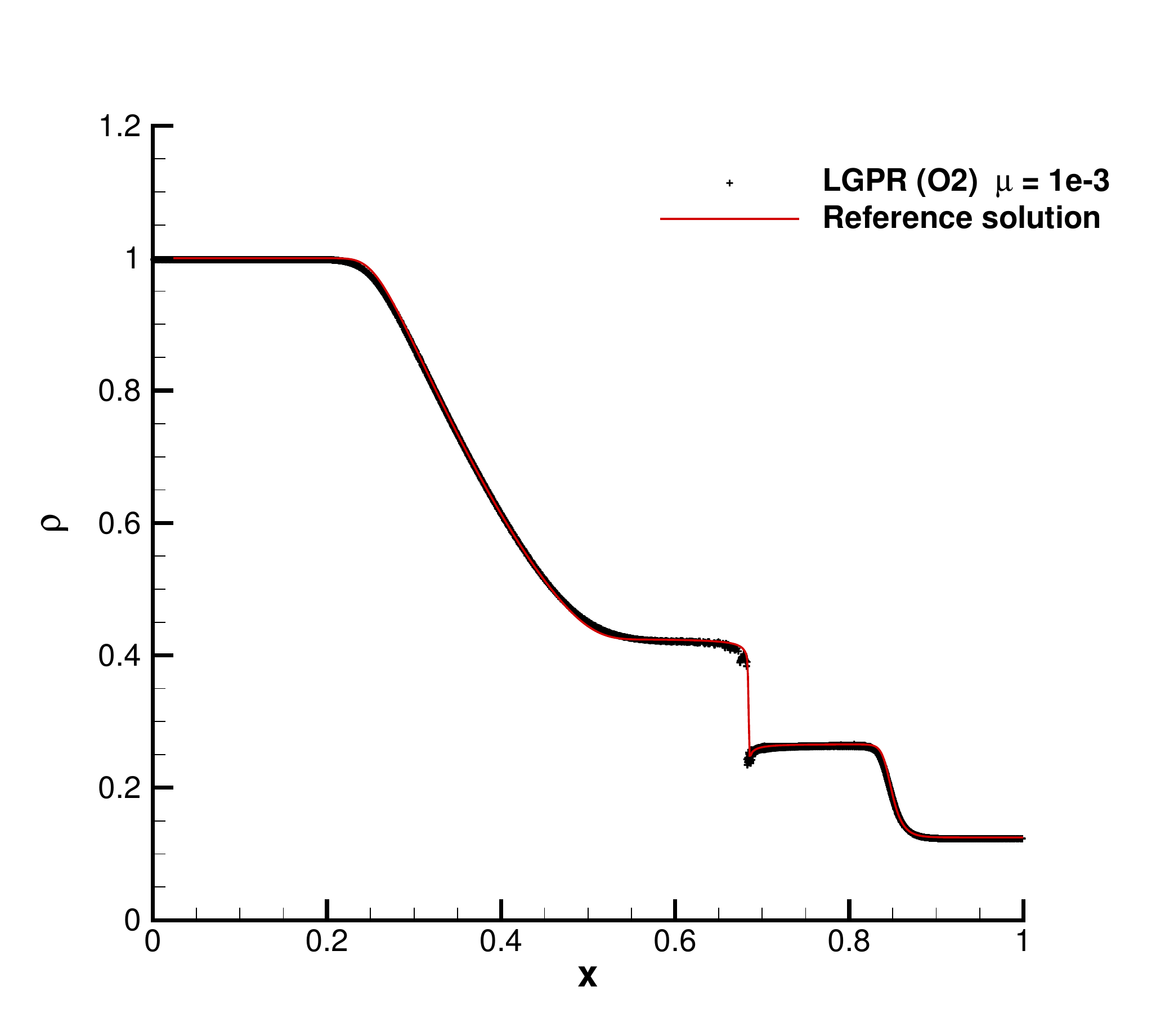}  &          
			\includegraphics[width=0.33\textwidth,draft=false]{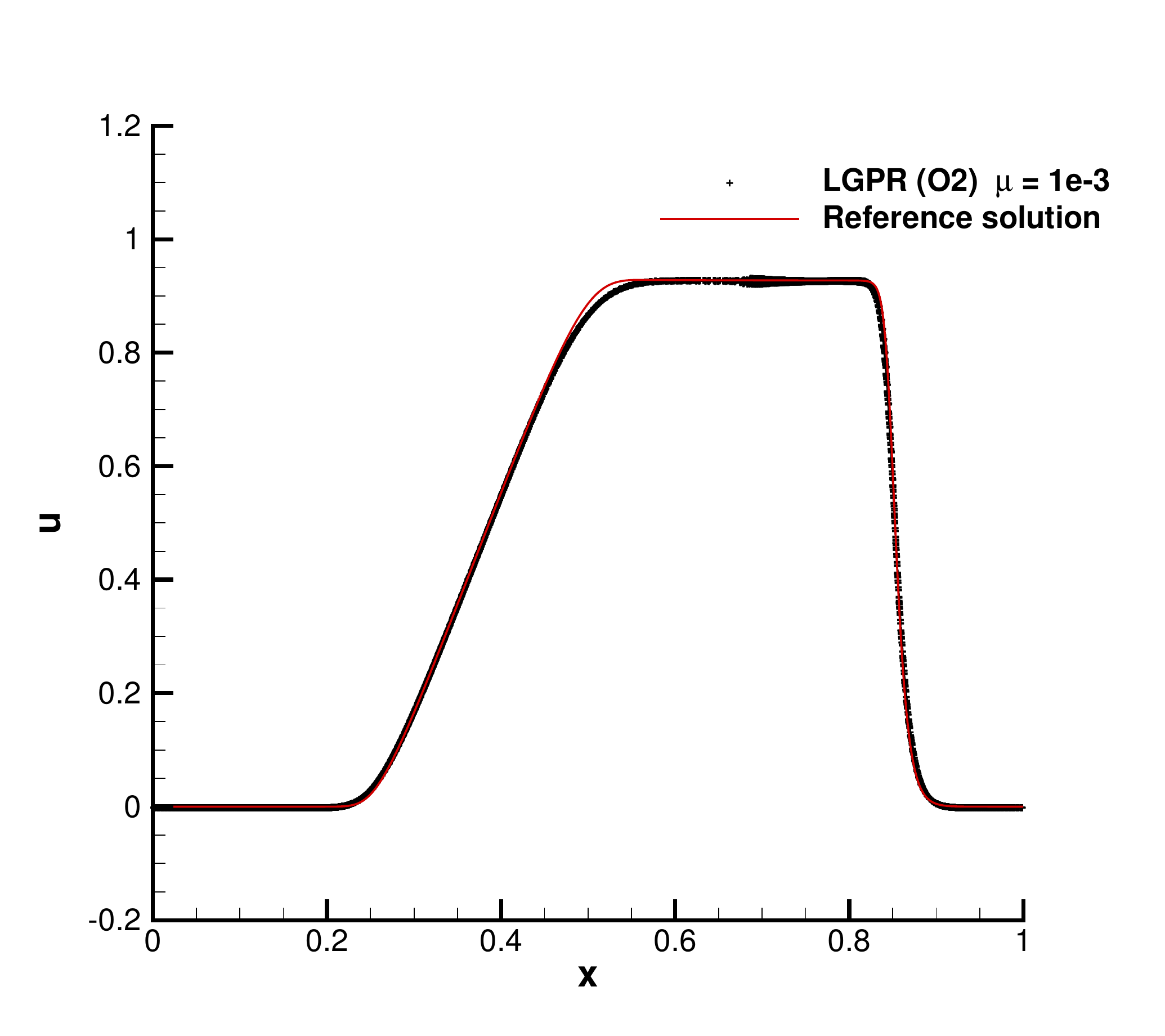} &
			\includegraphics[width=0.33\textwidth,draft=false]{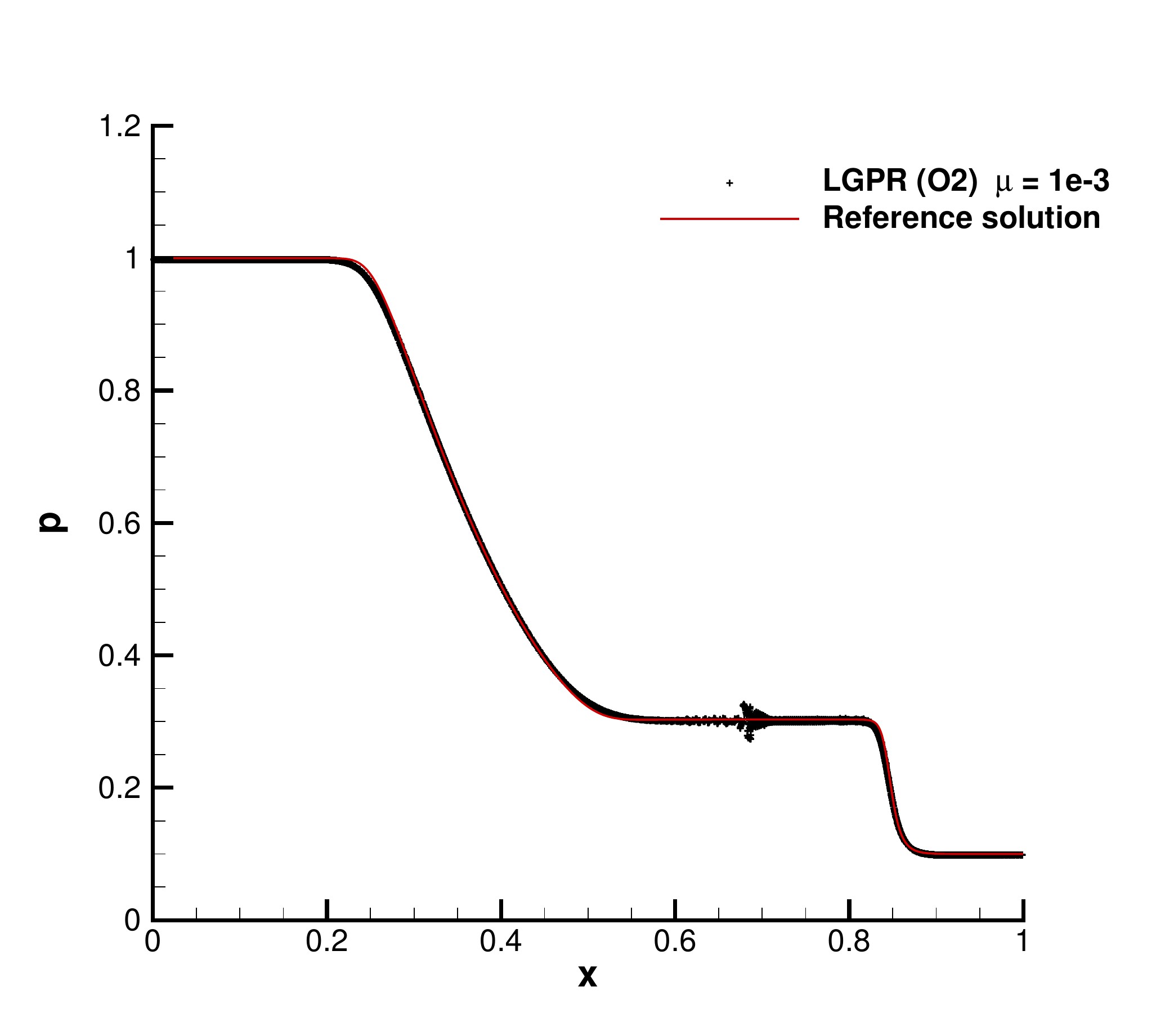} \\
			\includegraphics[width=0.33\textwidth,draft=false]{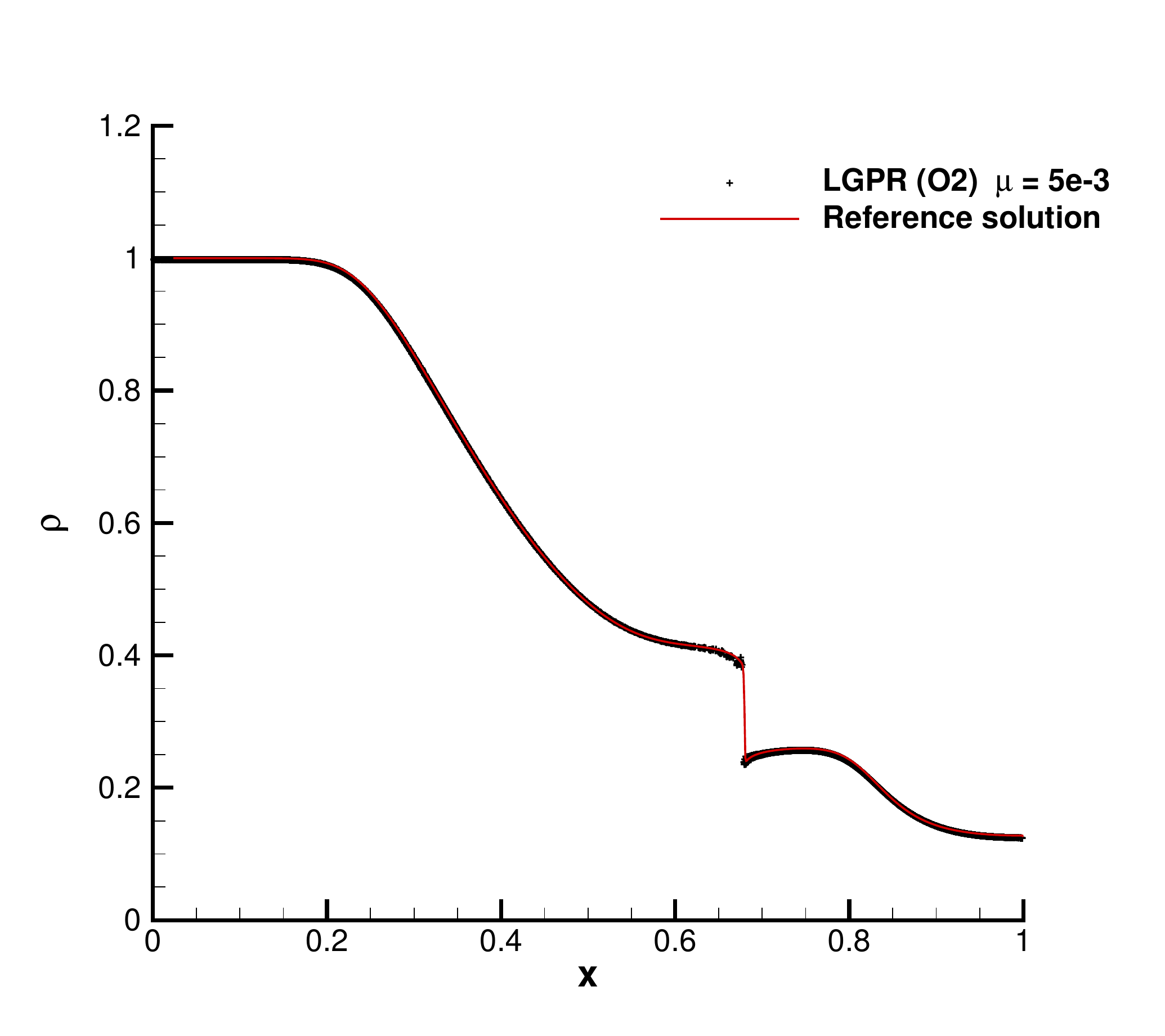}  &          
			\includegraphics[width=0.33\textwidth,draft=false]{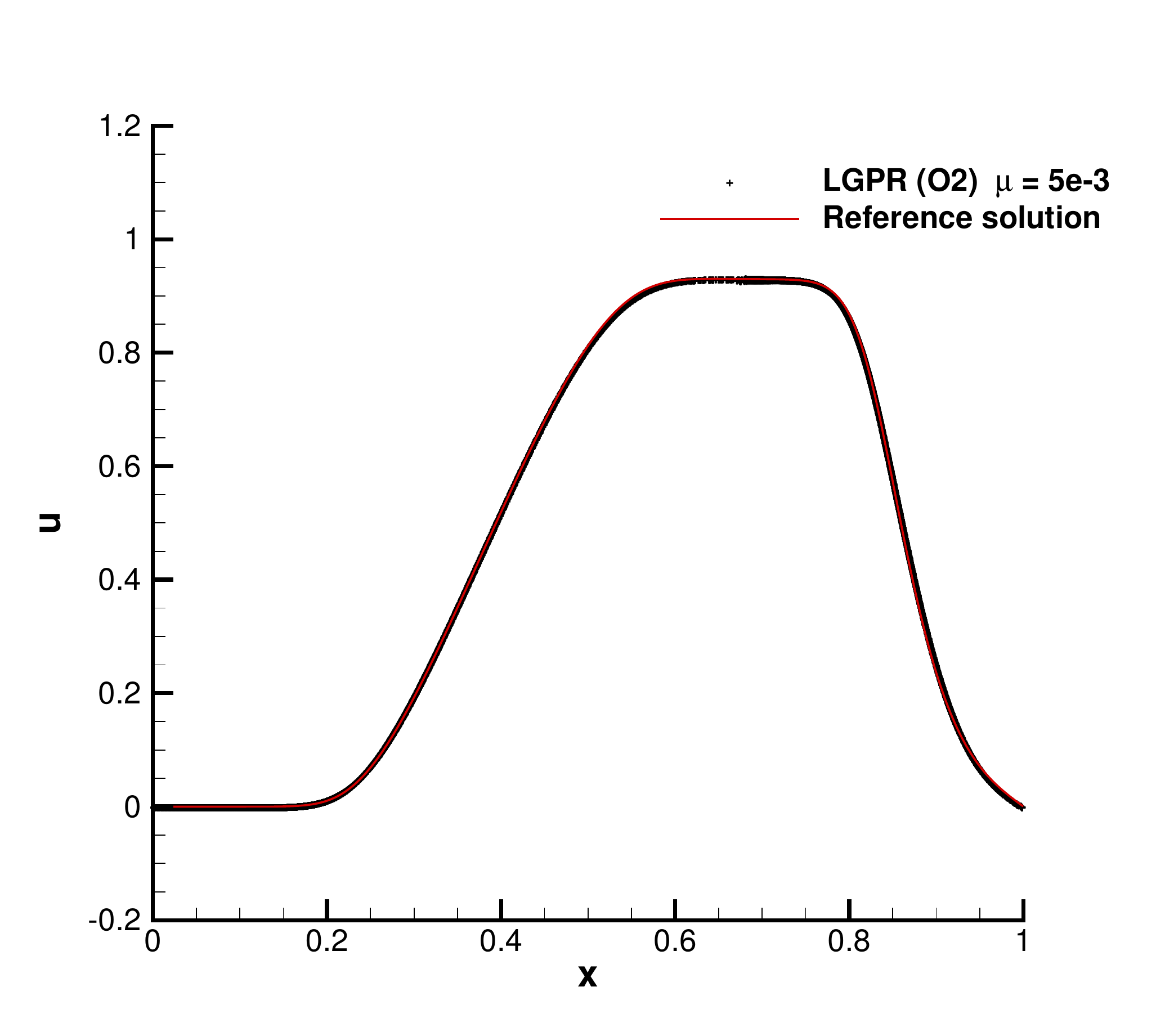} &
			\includegraphics[width=0.33\textwidth,draft=false]{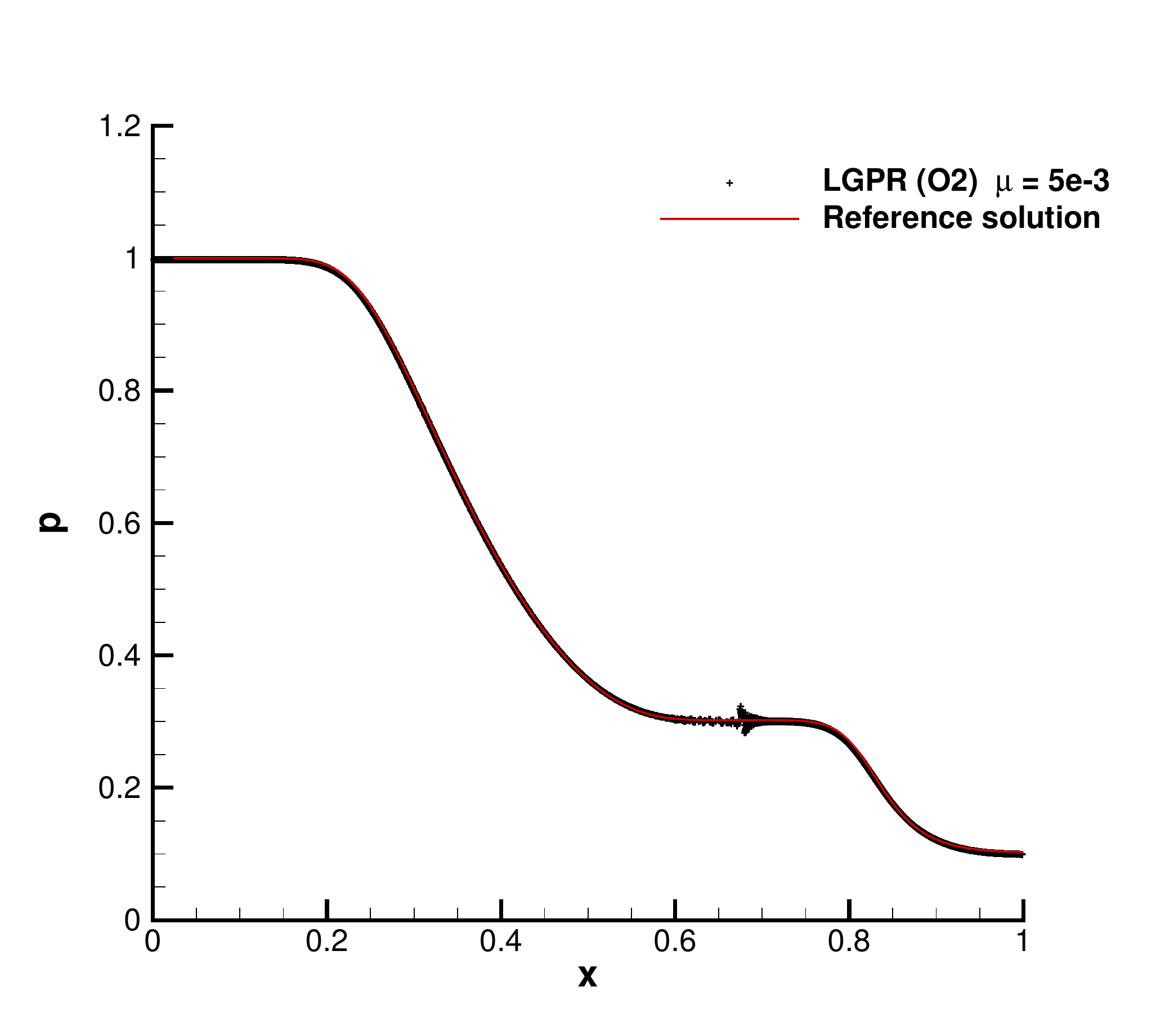} \\
			\includegraphics[width=0.33\textwidth,draft=false]{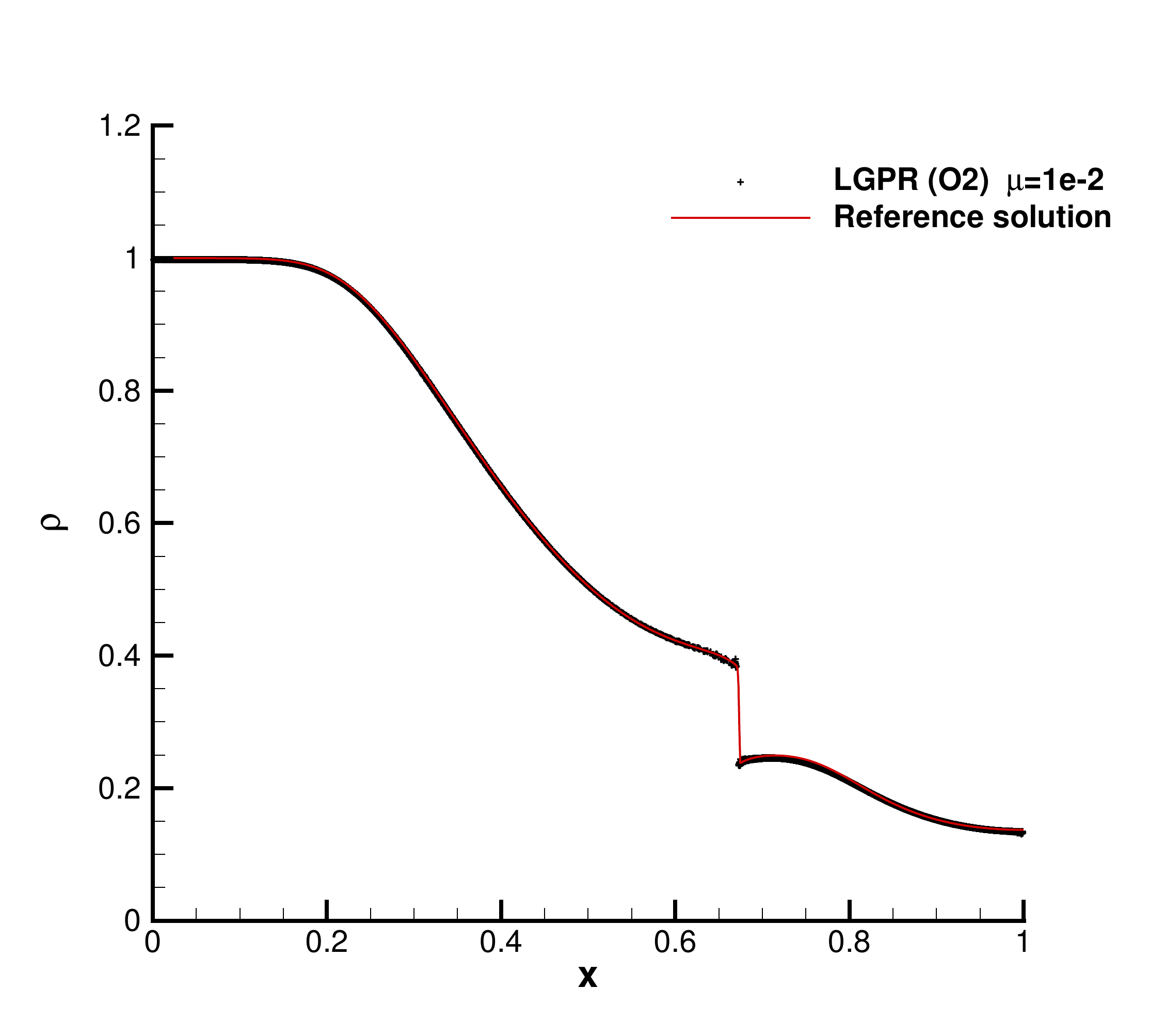}  &          
			\includegraphics[width=0.33\textwidth,draft=false]{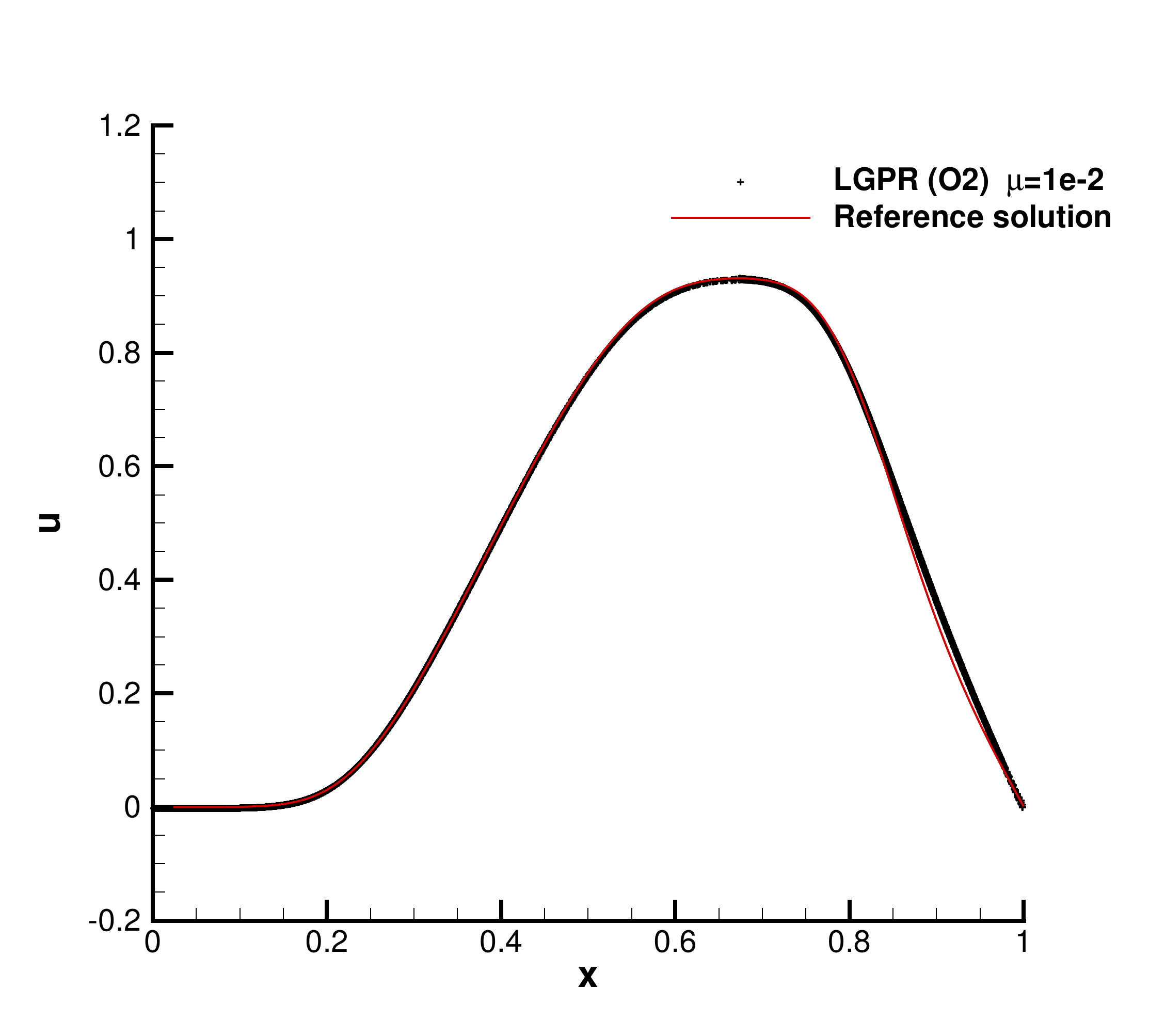} &
			\includegraphics[width=0.33\textwidth,draft=false]{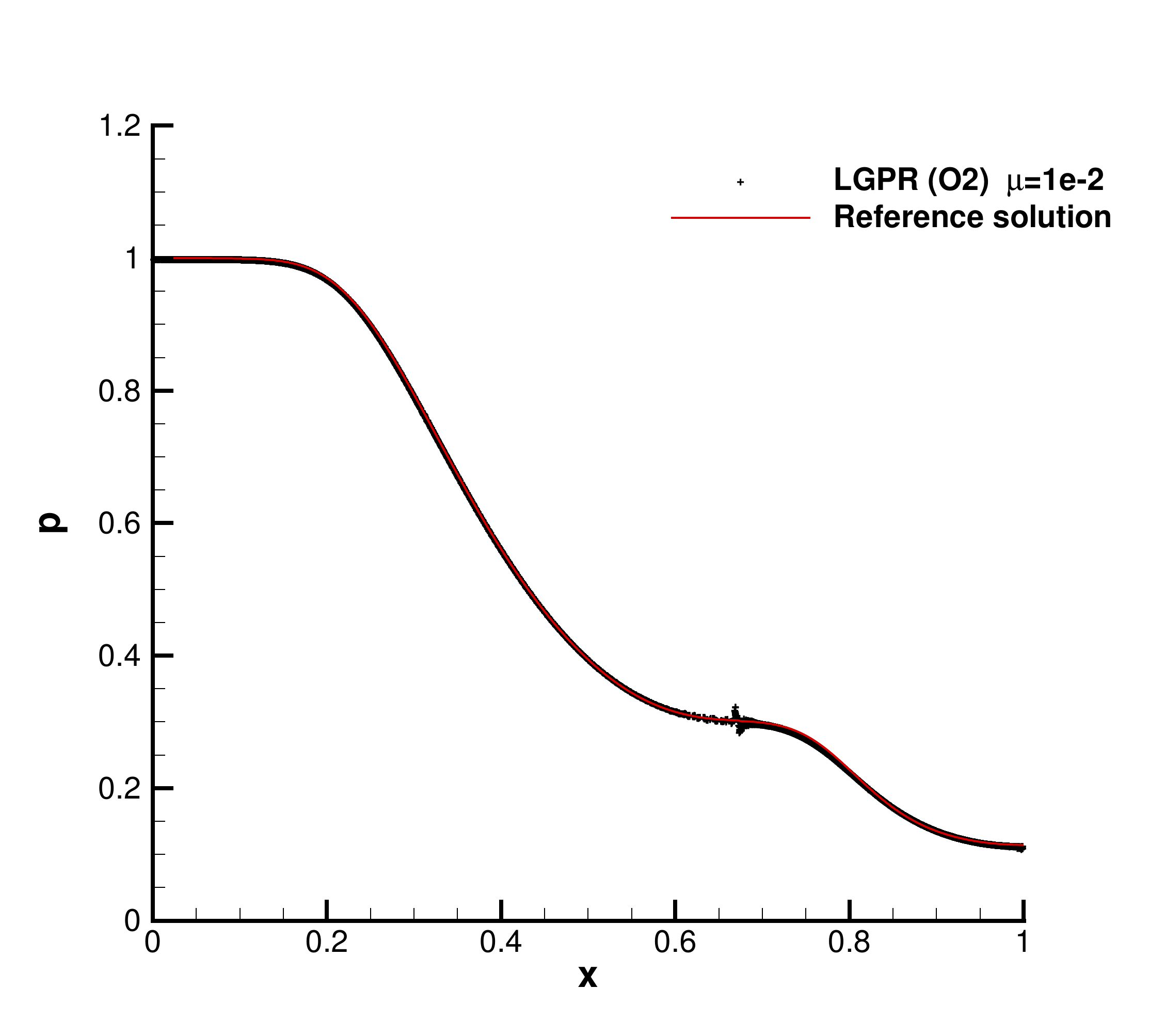} \\
		\end{tabular}
		\caption{Viscous shock tube problem. Scatter plot of cell density (left), horizontal 
		velocity (middle) and pressure (right) at the final time $t=0.2$ with viscosity coefficient 
		$\visc=10^{-3}$ (top), $\visc=5 \cdot 10^{-3}$ (middle) and $\visc=10^{-2}$ (bottom). 
		Comparison against the reference solution of the Navier-Stokes equations (solid red line).}
		\label{fig.Sodvisc}
	\end{center}
\end{figure}

The test cases presented so far deal with the stiff hydrodynamics limit of the GPR model \eqref{eqn.cl}. 
Here, we want to test the new LGPR scheme for the simulation of viscous fluids using the \textit{same} 
set of equations. To this purpose, we solve in a two-dimensional setting the well-known Sod shock tube 
problem, which is a classical one-dimensional test problem that involves a rarefaction wave traveling 
towards the left boundary as well as a right-moving contact discontinuity and a shock wave 
traveling to the right. However, instead of ideal fluids, viscous ones are now considered, with 
three different viscosity coefficients, namely $\visc=10^{-3}$, $\visc=5 \cdot 10^{-3}$ and $\visc=10^{-2}$. 
The initial computational domain is the rectangular box $\Omega(0)=[0;1]\times[0;0.1]$ which is 
discretized with a characteristic mesh size of $h=1/200$. Slip-wall boundaries are set everywhere, 
hence prescribing zero normal velocity in the nodal solver. The final time is $t_f=0.2$ so that all 
waves remain bounded in the computational domain. The initial condition consists in a discontinuity 
located at $x_0=0.5$ between the left and the right state: 
\begin{equation}
	(\rho(x),u(x),p(x)) = \left\{ \begin{array}{lll} \left( 1, 0, 1 \right), & \textnormal{ if } & x \leq x_0, \\ 
		\left( 0.125,0, 0.1 \right), & \textnormal{ if } & x > x_0. 
	\end{array}  \right. , \qquad t=0. 
	\label{eq:Sod_IC}
\end{equation}
The fluid is initially at rest and the ideal gas EOS is adopted with $\gamma=1.4$, while thermal 
conduction is neglected. According to the asymptotic analysis presented in Section \ref{sec.AP}, the 
viscous stress tensor of the Navier-Stokes equation should be retrieved, therefore the results of the 
LGPR scheme are compared with a reference solution of the Navier-Stokes equations without thermal 
conduction computed on a very fine one-dimensional mesh of $10\,000$ cells with a MUSCL-TVD finite 
volume method. Figure \ref{fig.Sodvisc3D} plots a three-dimensional view of the density distribution 
as well as the mesh configuration across the contact discontinuity at the final time, whereas 
Figure \ref{fig.Sodvisc} shows a scatter plot for density, horizontal velocity and pressure of the 
numerical solution compared against the reference solution of the Navier-Stokes model. An excellent 
agreement is obtained for all different values of the viscosity coefficient, demonstrating that the 
novel LGPR scheme is 
capable of simulating viscous fluids as well.
	
\subsection{Heat conduction in a gas} \label{ssec.Heat2D}
\begin{figure}[!htbp]
	\begin{center}
		\begin{tabular}{cc}
			\includegraphics[width=0.47\textwidth,draft=false]{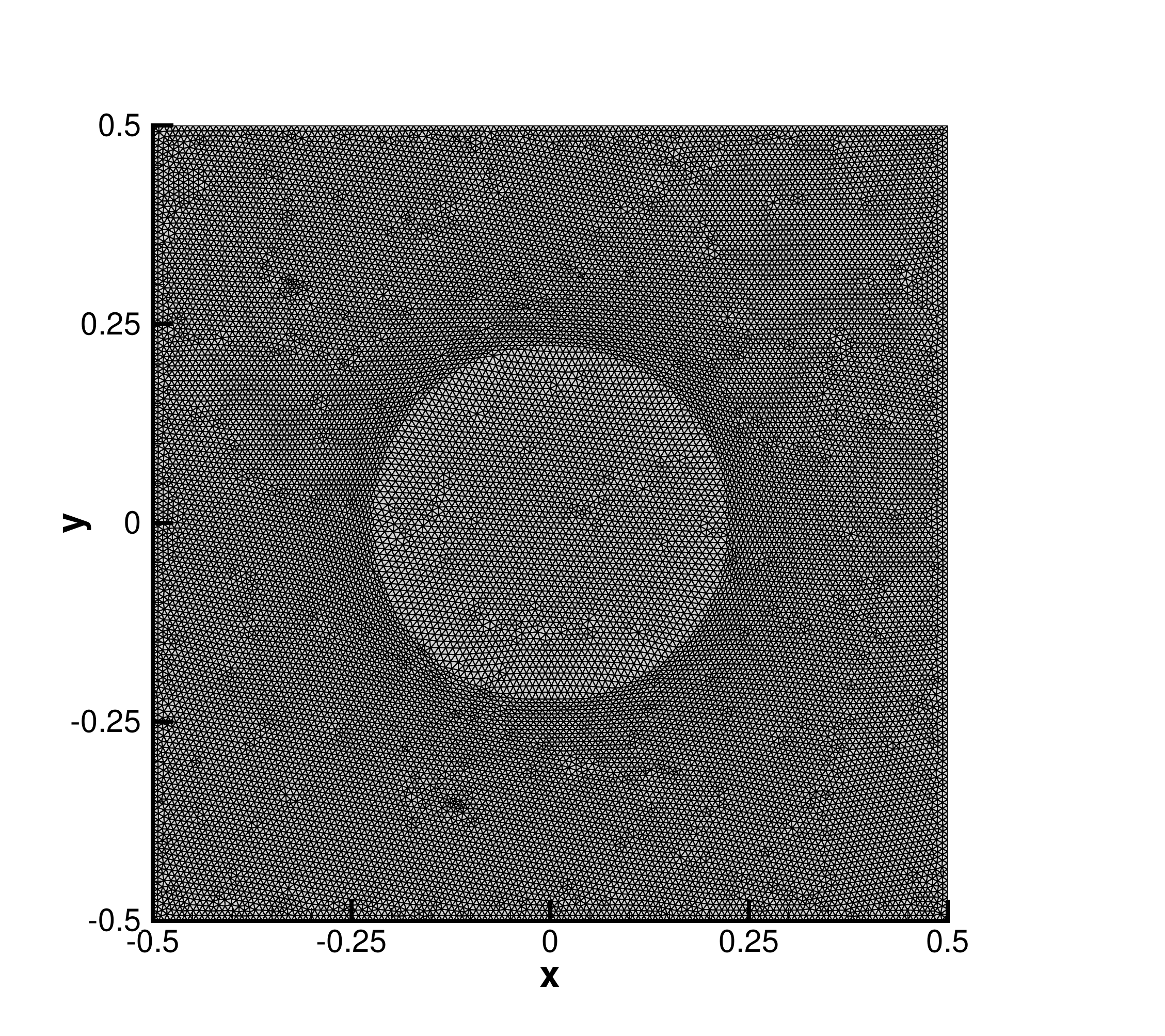} &
			\includegraphics[width=0.47\textwidth,draft=false]{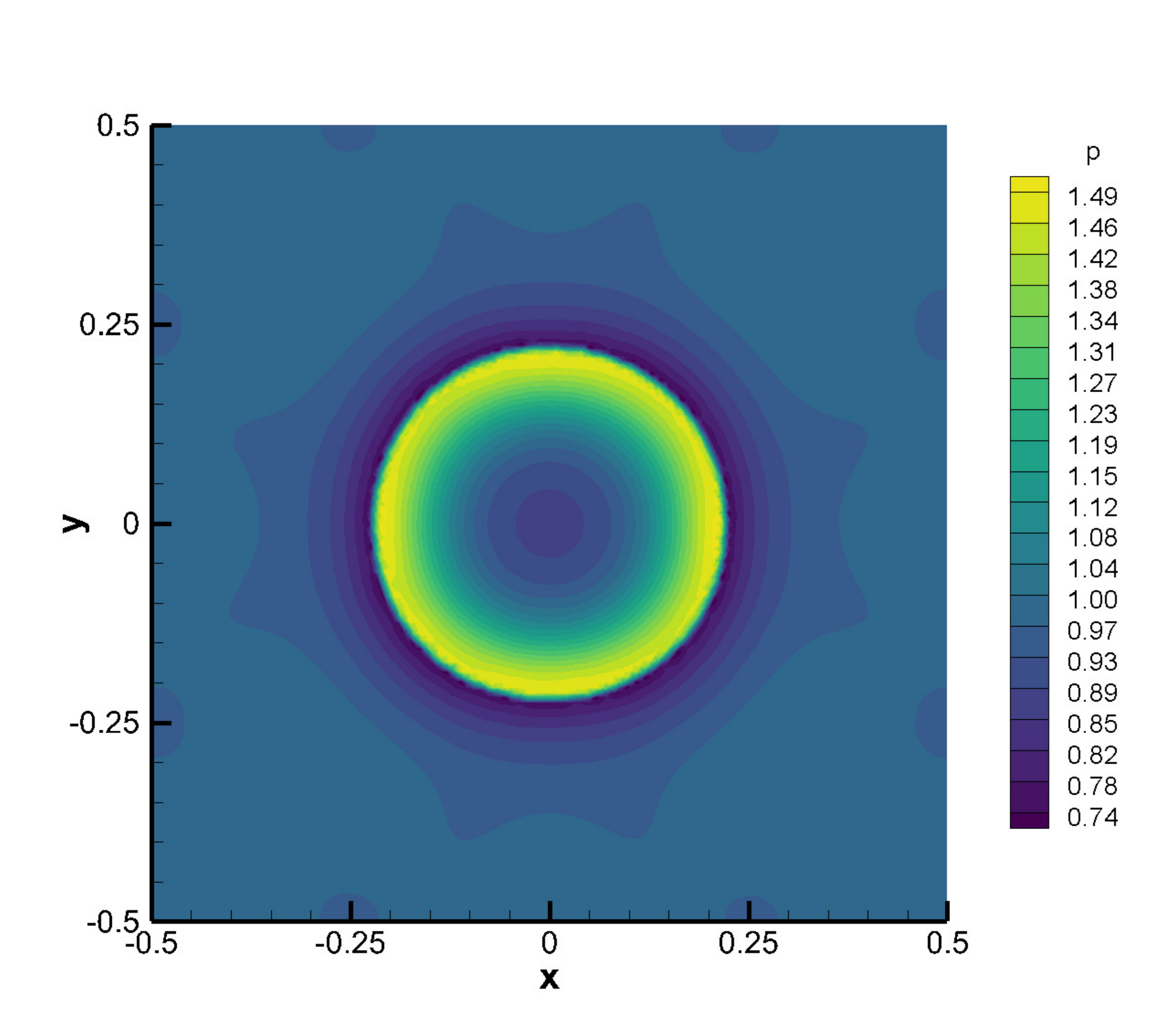} \\
			\includegraphics[width=0.47\textwidth,draft=false]{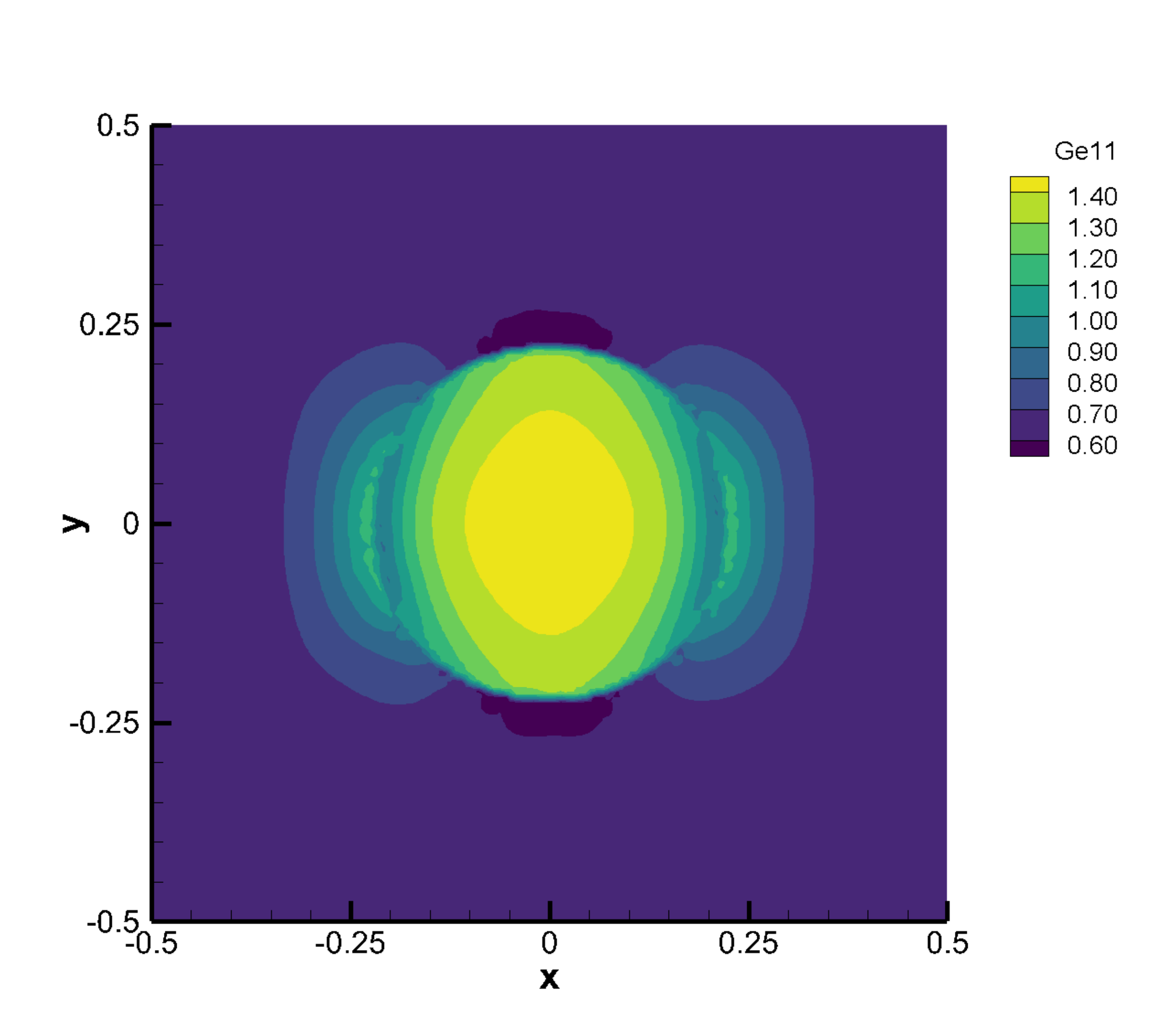} &
			\includegraphics[width=0.47\textwidth,draft=false]{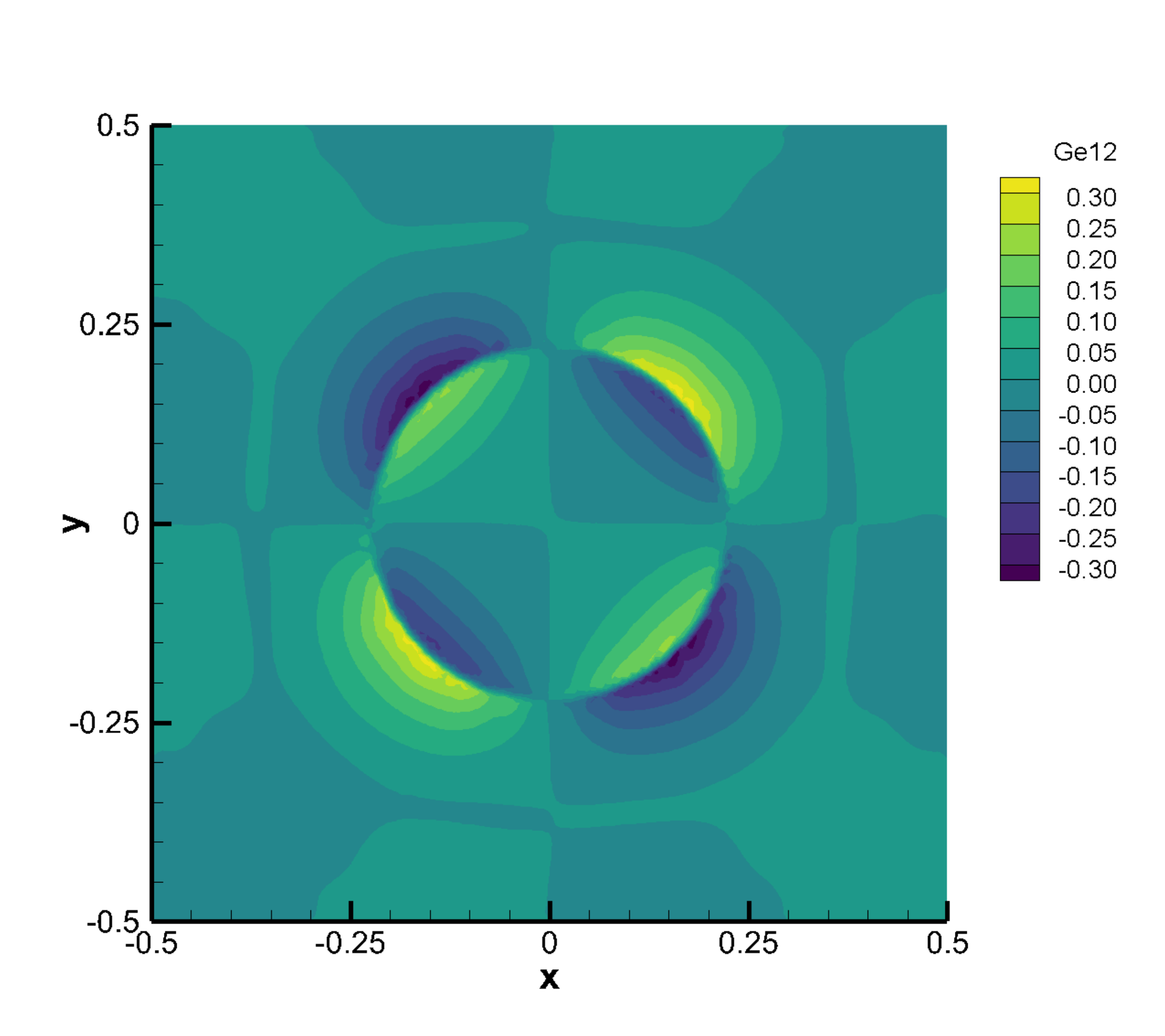} \\
			\includegraphics[width=0.47\textwidth,draft=false]{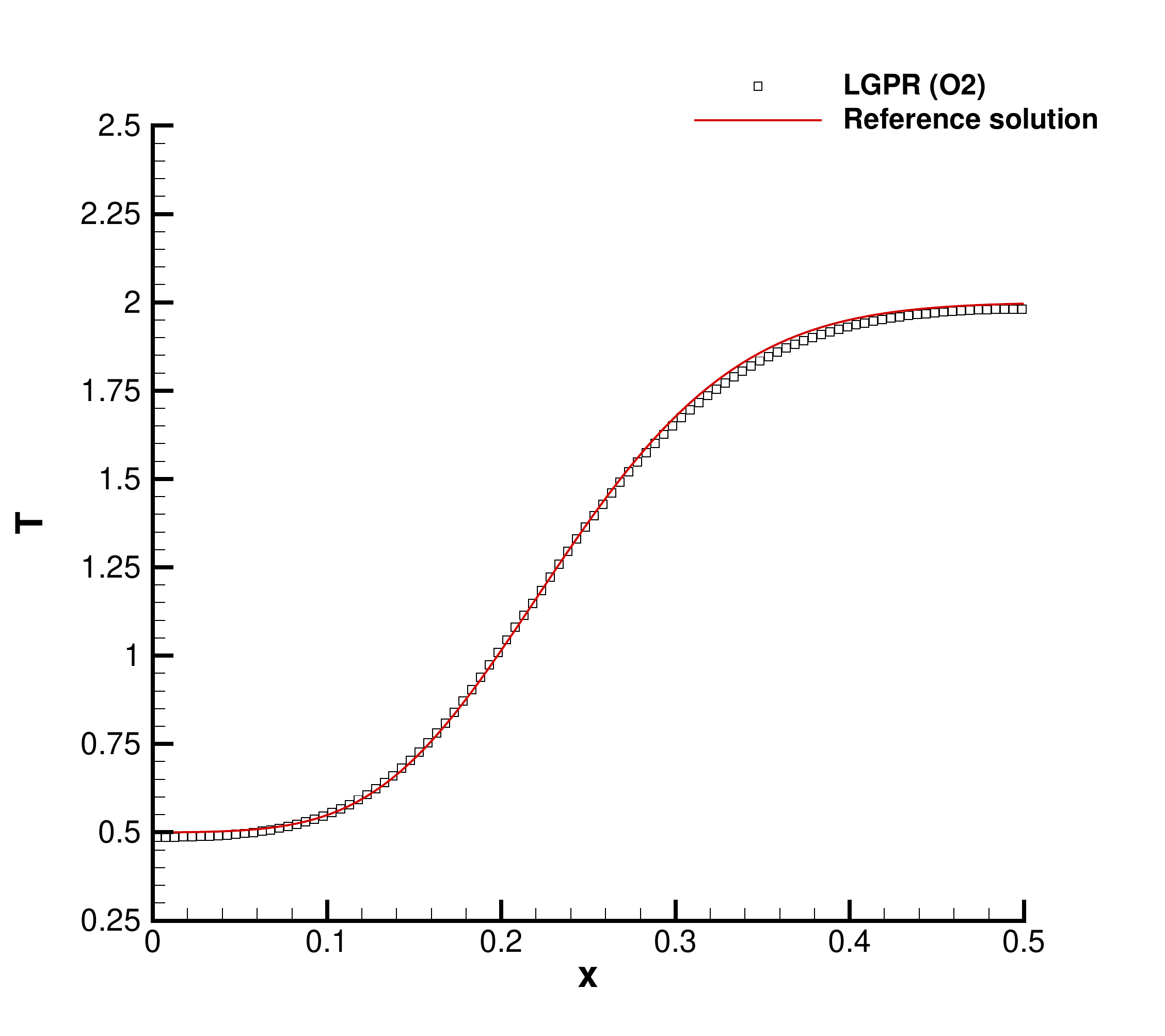} &
			\includegraphics[width=0.47\textwidth,draft=false]{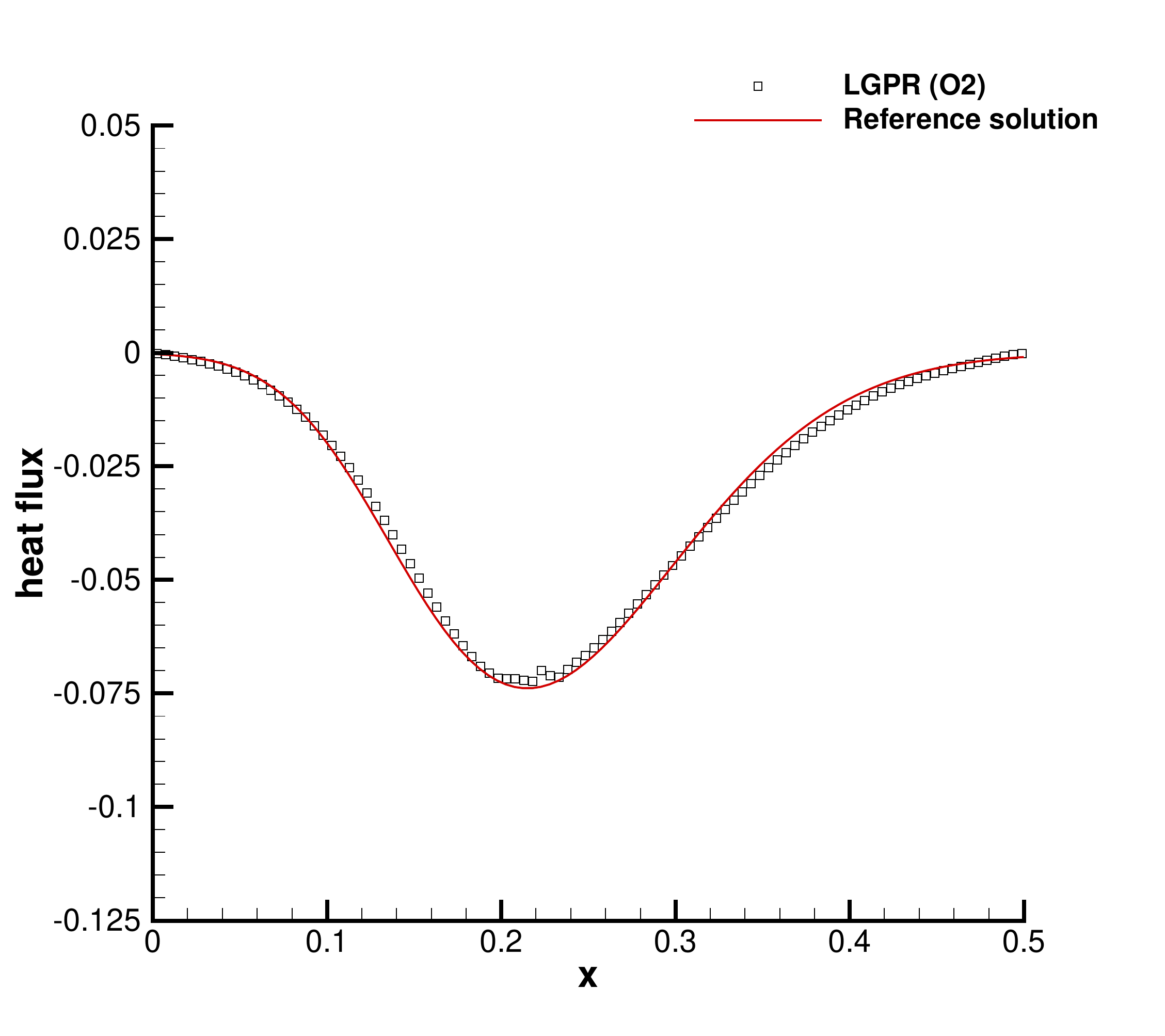} \\
		\end{tabular}
		\caption{Heat conduction in a gas. Top: mesh configuration  (left) and pressure distribution (right) Initial (left) at time $t=1$. Middle: tensor component $\tensor{G}_{e_{11}}$ (left) and $\tensor{G}_{e_{12}}$ (right). Bottom: temperature distribution (left) and heat flux (right). For the Navier–Stokes solution, the classical Fourier heat flux $q_1=-\kappa T_x$ is shown, while for the GPR model, we plot $q_1=\alpha^2 T J_1$.}
		\label{fig.heat2D}
	\end{center}
\end{figure}

The aim of this test case is to verify the correct behavior of the novel LGPR scheme in the case of a problem dominated by heat transfer via heat conduction. A high density circle of gas is initialized at the center of the computational domain $\Omega(0)=[-0.5,0.5]^2$, that is
\begin{equation}
	\rho(0,r) = \left\{
	\begin{array}{ll}
		2   & r \leq R_0 \\
		0.5 &  x > R_0
	\end{array}\right., 
\end{equation}
with $r=\sqrt{x^2+y^2}$ representing the generic radial coordinate and $R_0=0.2$ being the radius 
of the circle containing the high density gas. The fluid is initially at rest ($u=v=w=0$) with 
constant pressure $p=1$ and obeys an ideal gas law with $\gamma=1.4$. The heat wave velocity and 
heat conduction coefficient are specified in Table \ref{tab.GPRpar}. The computational domain is 
discretized with a characteristic mesh size of $h=1/100$ and slip-wall boundaries are imposed on 
all sides. Figure \ref{fig.heat2D} shows the pressure distribution as well as the metric tensor 
components $(\tensor{G}_{e_{11}},\tensor{G}_{e_{12}})$ at the final time. Despite the highly 
unstructured and non-symmetric mesh shown in Figure \ref{fig.heat2D}, the numerical solution 
exhibits excellent symmetry. Finally, the solution along a 1D cut with 200 points along the 
$x$-direction at $y=0$ is compared against the reference solution of the Navier-Stokes-Fourier 
equations for temperature and the heat flux, achieving an excellent agreement. This demonstrates 
that the stiff limit of the heat conduction equation \eqref{eqn.cl4} is properly retrieved by the 
LGPR scheme \eqref{eqn.fvcl4}, hence giving numerical evidence of the asymptotic preserving 
property of the scheme studied in Section \ref{sec.AP}.

\subsection{Viscous shock profile} \label{ssec.shockNS}
\begin{figure}[!htbp]
	\begin{center}
		\begin{tabular}{c}
			\includegraphics[width=0.85\textwidth,draft=false]{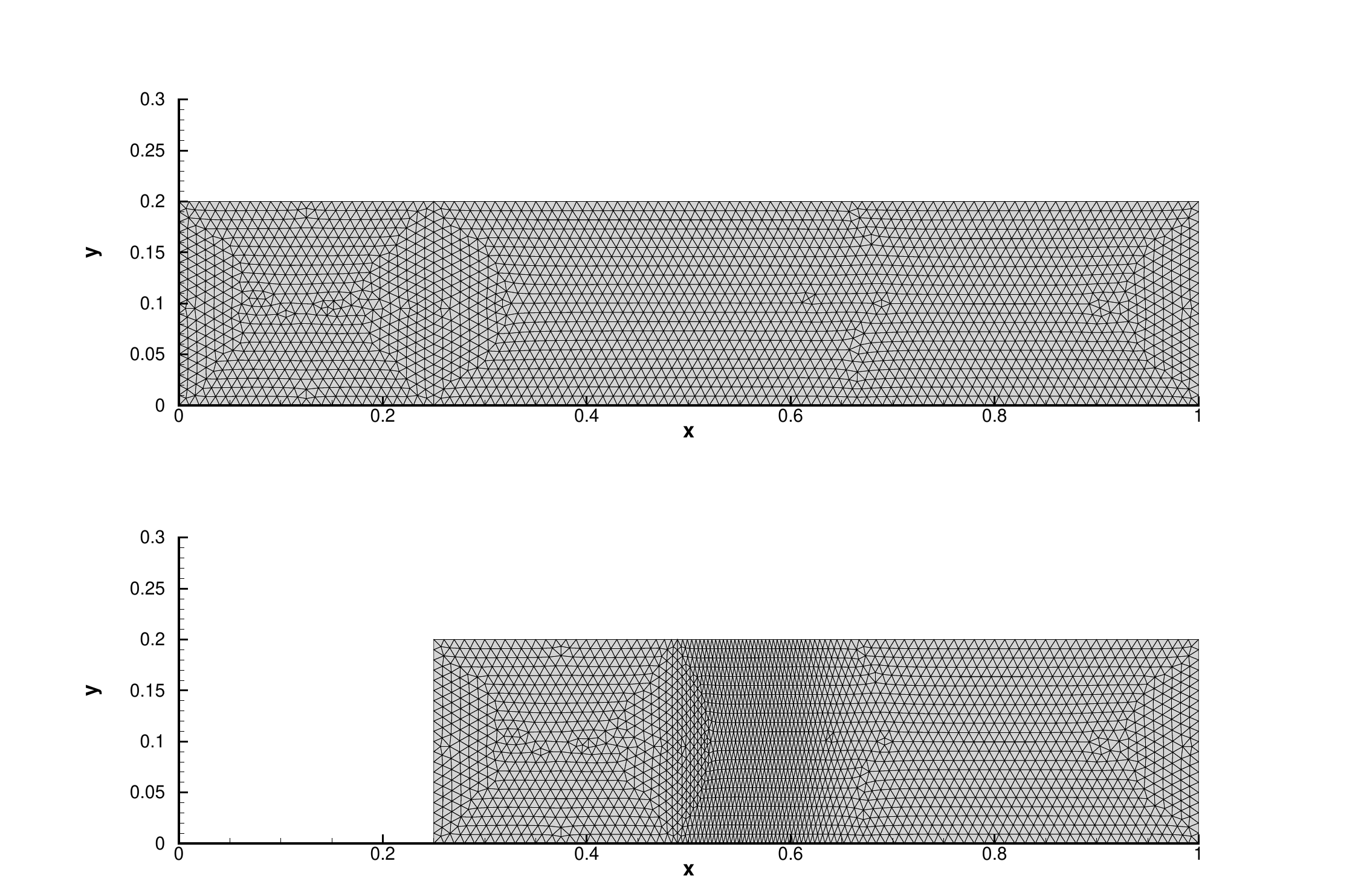}
		\end{tabular}
		\caption{Viscous shock profile. Initial (left) and final (right) mesh configuration with $h=1/100$.}
		\label{fig.shockNS_mesh}
	\end{center}
\end{figure}

In order to verify the numerical method against supersonic viscous flows, we propose to solve the 
problem of an isolated viscous shock wave which is traveling into a viscous heat conducting fluid 
at rest with a shock Mach number of $M_s=2$. The analytical solution can be found in 
\cite{Becker1923}, where the compressible Navier-Stokes-Fourier equations are solved for the 
special case 
of a stationary shock wave at Prandtl number $Pr= 0.75$ with constant viscosity. The exact 
solution for the 
dimensionless velocity $\bar u = \frac{u}{M_s \, c_0}$ of this stationary shock wave with shock Mach number $M_s$ is then given by the root of 
the following equation, see \cite{Becker1923}:
\begin{equation} 
	\label{eqn.alg.u} 
	\frac{|\bar u - 1|}{|\bar u - \kappa^2|^{\kappa^2}} = \left| \frac{1-\kappa^2}{2} \right|^{(1-\kappa^2)} 
	\exp{\left( \frac{3}{4} \textnormal{Re}_s \frac{M_s^2 - 1}{\gamma M_s^2} x \right)},
\end{equation}
with
\begin{equation}
	\kappa^2 = \frac{1+ \frac{\gamma-1}{2}M_s^2}{\frac{\gamma+1}{2}M_s^2}.
\end{equation}
Equation \eqref{eqn.alg.u} allows the dimensionless velocity $\bar u$ to be obtained as a function of $x$. The form of the viscous profile of the dimensionless pressure $\bar p = \frac{p-p_0}{\rho_0 c_0^2 M_s^2}$ is given by
the relation 
\begin{equation}
	\label{eqn.alg.p} 
	\bar p = 1 - \bar u +  \frac{1}{2 \gamma}
	\frac{\gamma+1}{\gamma-1} \frac{(\bar u - 1 )}{\bar u} (\bar u - \kappa^2).  
\end{equation}
Finally, the profile of the dimensionless density $\bar \rho = \frac{\rho}{\rho_0}$ is derived from the integrated continuity equation: $\bar \rho \bar u = 1$. To obtain an unsteady shock wave traveling into a medium at rest, it is sufficient to superimpose a constant velocity field $u = M_s c_0$ to the solution of the stationary shock wave found in the previous steps.
The initial computational domain is the rectangular channel $\Omega(0)=[0,1]\times[0,0.2]$ which is paved with two different triangular meshes of characteristic mesh size $h=1/100$ and $h=1/200$. On the left side of the domain ($x=0$) the constant inflow velocity is prescribed, whereas periodic boundaries are set along the $y$ direction and a constant pressure is imposed at $x=1$. The initial condition involves a shock wave centered at $x=0.25$ propagating at Mach $M_s=2$ from left to right with a Reynolds number of $Re=100$. The polytropic index of the ideal gas is $\gamma=1.4$ and the upstream shock state is defined by $\rho=1$, $u=v=0$, $p=1/\gamma$. Figure \ref{fig.shockNS_mesh} shows the initial and final mesh configuration at time $t_f=0.2$ with the shock front located at $x=0.65$. 

\begin{figure}[!htbp]
	\begin{center}
		\begin{tabular}{ccc}
			\includegraphics[width=0.33\textwidth,draft=false]{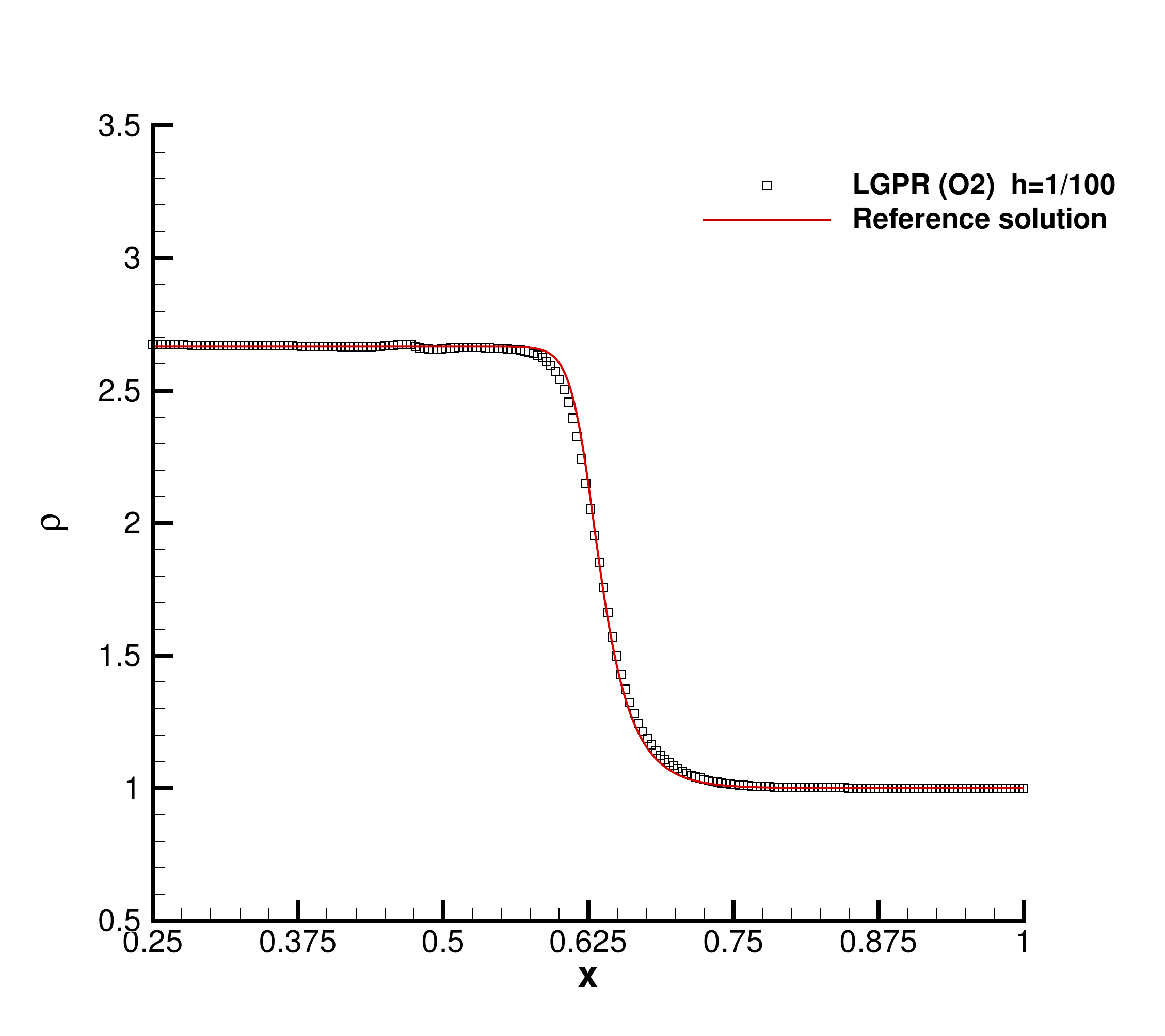}  &          
			\includegraphics[width=0.33\textwidth,draft=false]{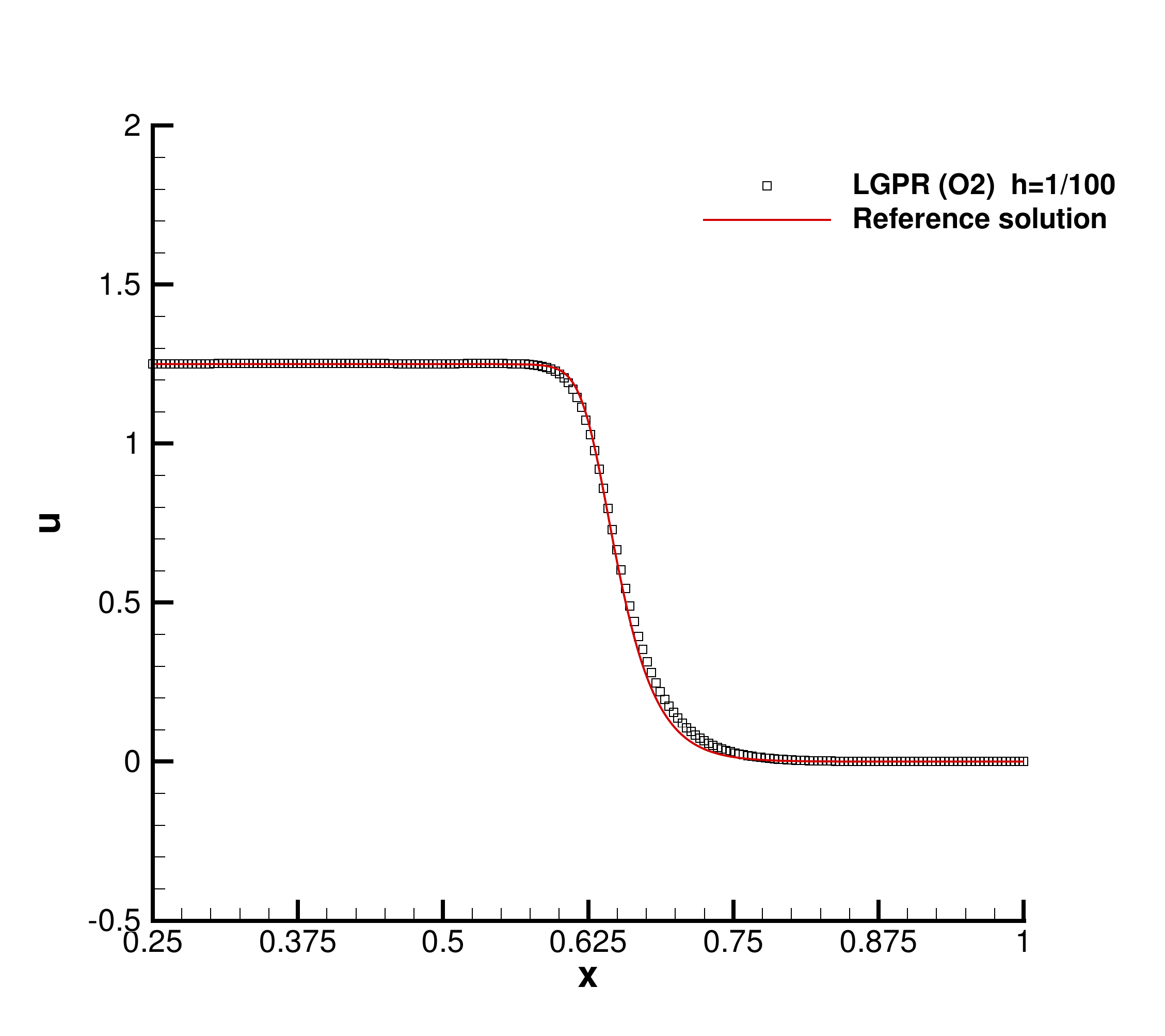} &
			\includegraphics[width=0.33\textwidth,draft=false]{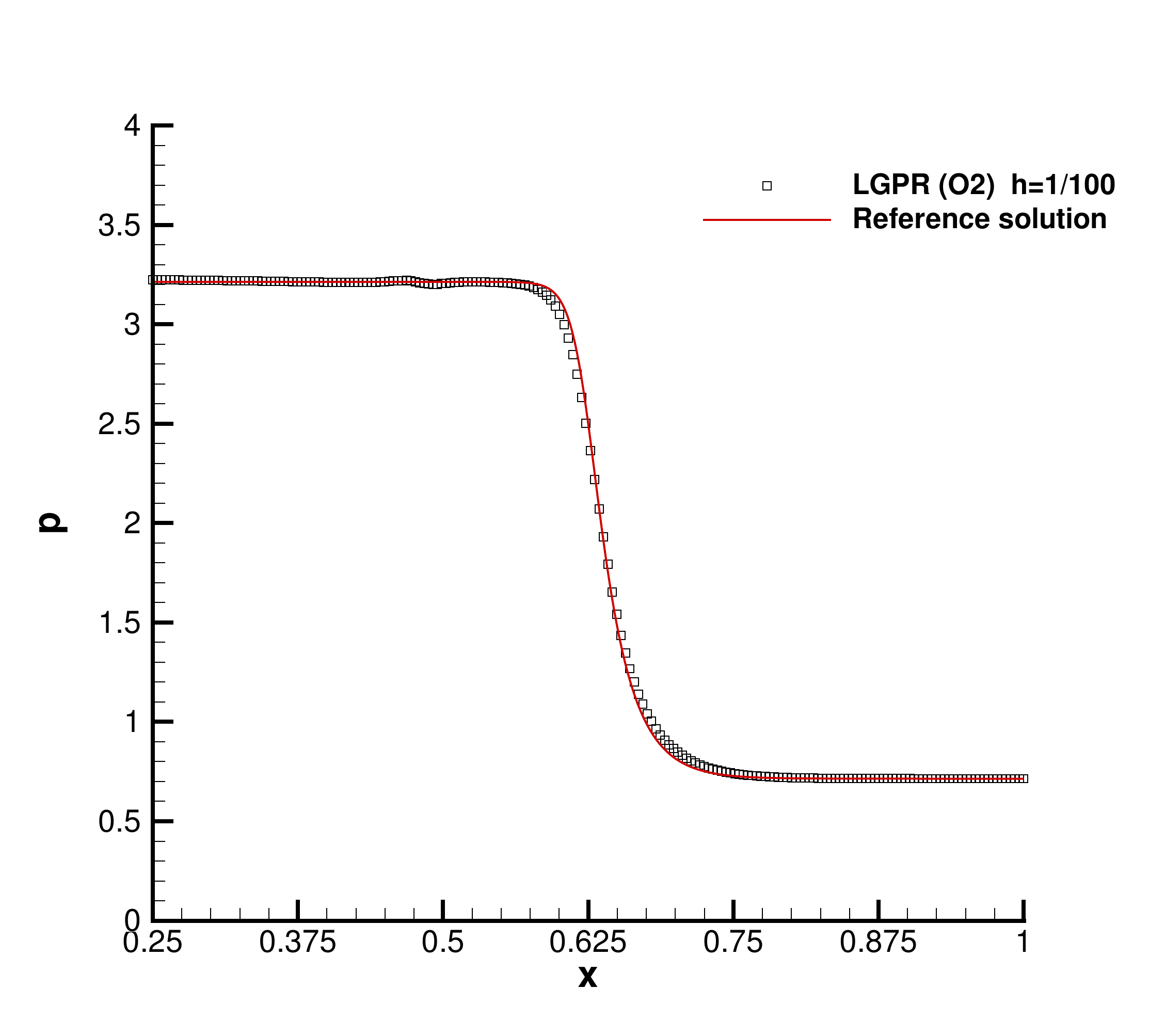} \\
			\includegraphics[width=0.33\textwidth,draft=false]{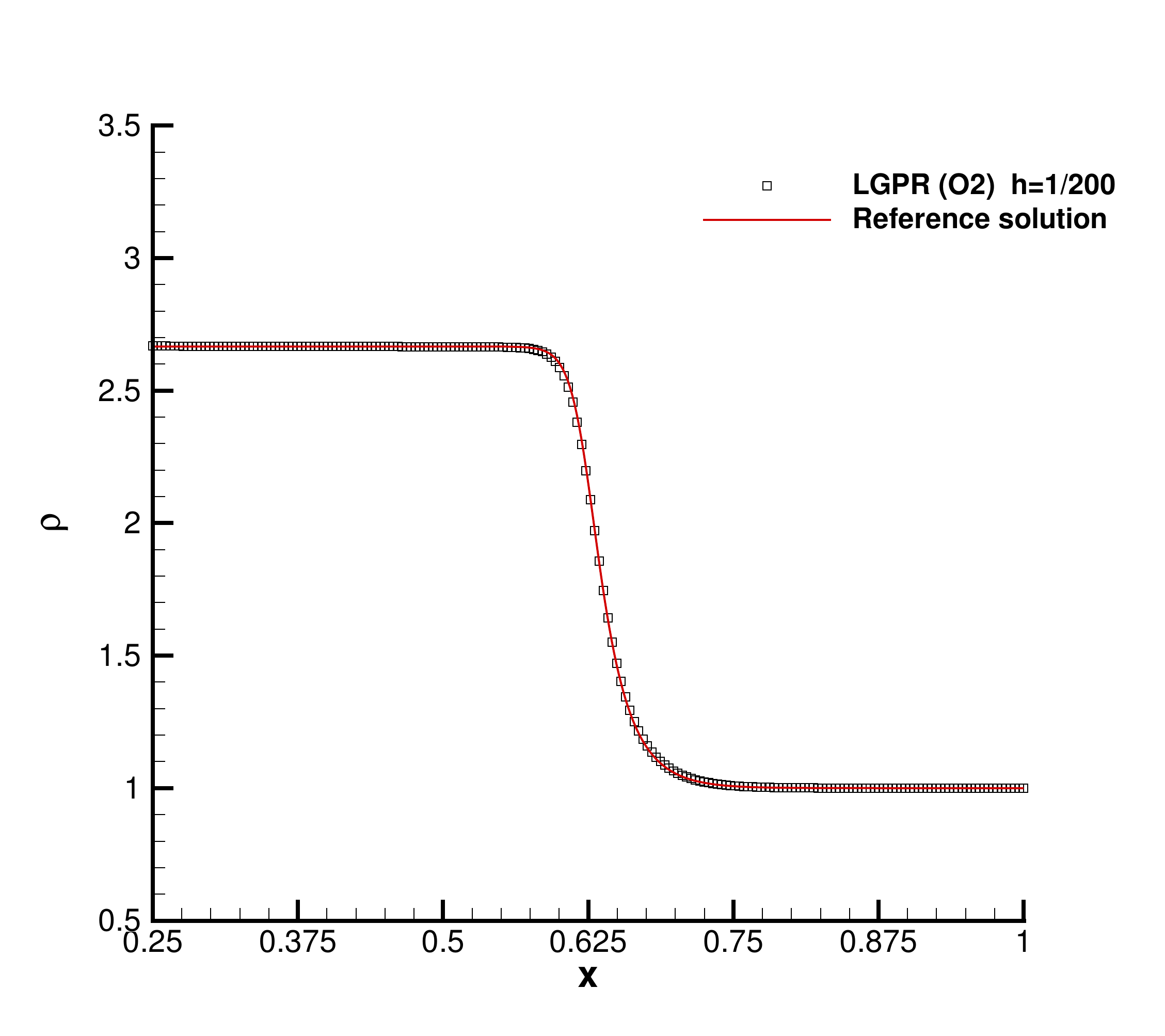}  &          
			\includegraphics[width=0.33\textwidth,draft=false]{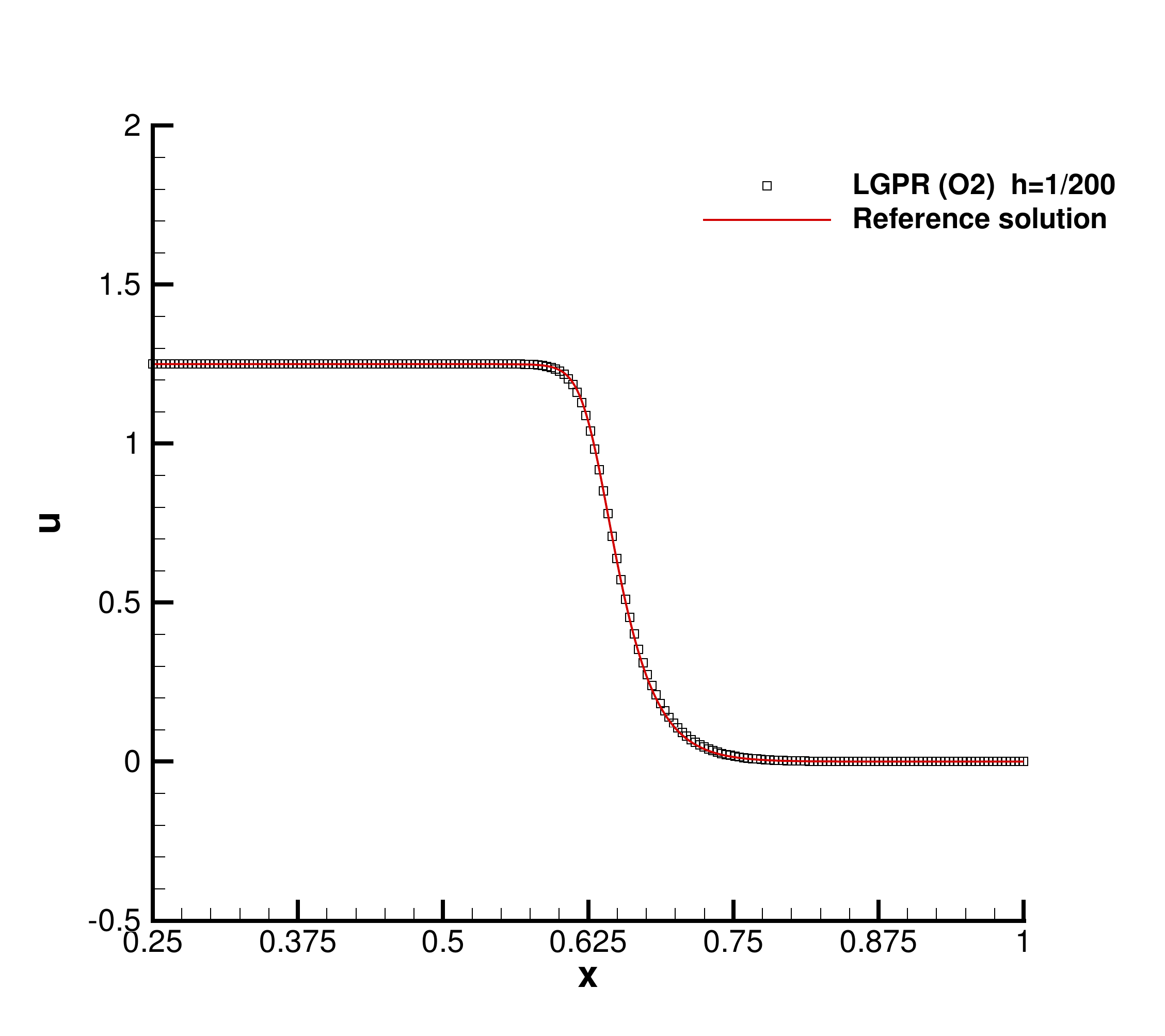} &
			\includegraphics[width=0.33\textwidth,draft=false]{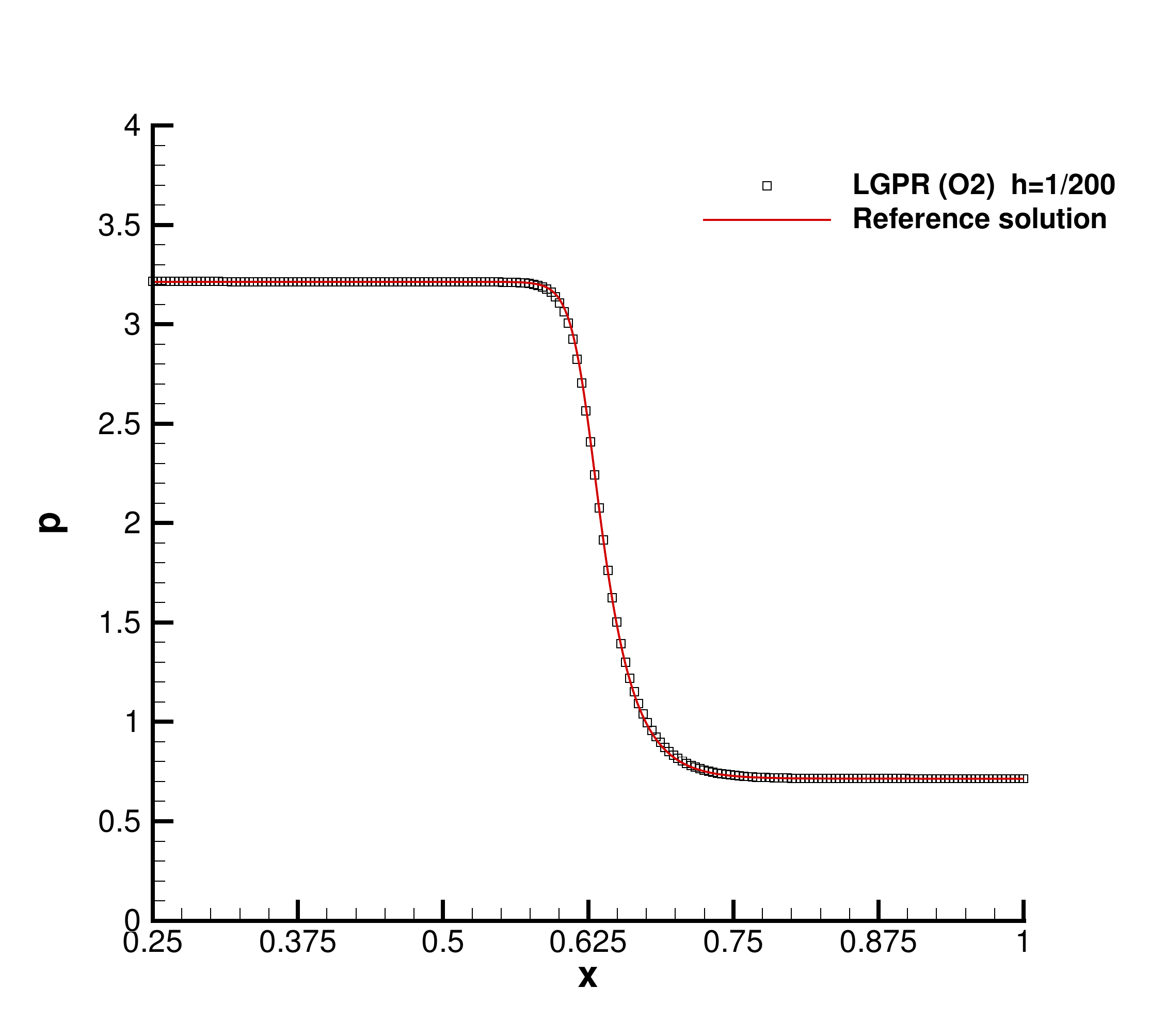} \\
		\end{tabular}
		\caption{Viscous shock profile for the Mach number $ M_s =2 $. Comparison of the exact 
		solution 
		of the compressible Navier–Stokes equations with the Lagrangian GPR model for density 
		(left), horizontal velocity (middle) and pressure (right) with mesh size $h=1/100$ (top) 
		and $h=1/200$ (bottom).}
		\label{fig.shockNS_1}
	\end{center}
\end{figure}

\begin{figure}[!htbp]
	\begin{center}
		\begin{tabular}{cc}
			\includegraphics[width=0.47\textwidth,draft=false]{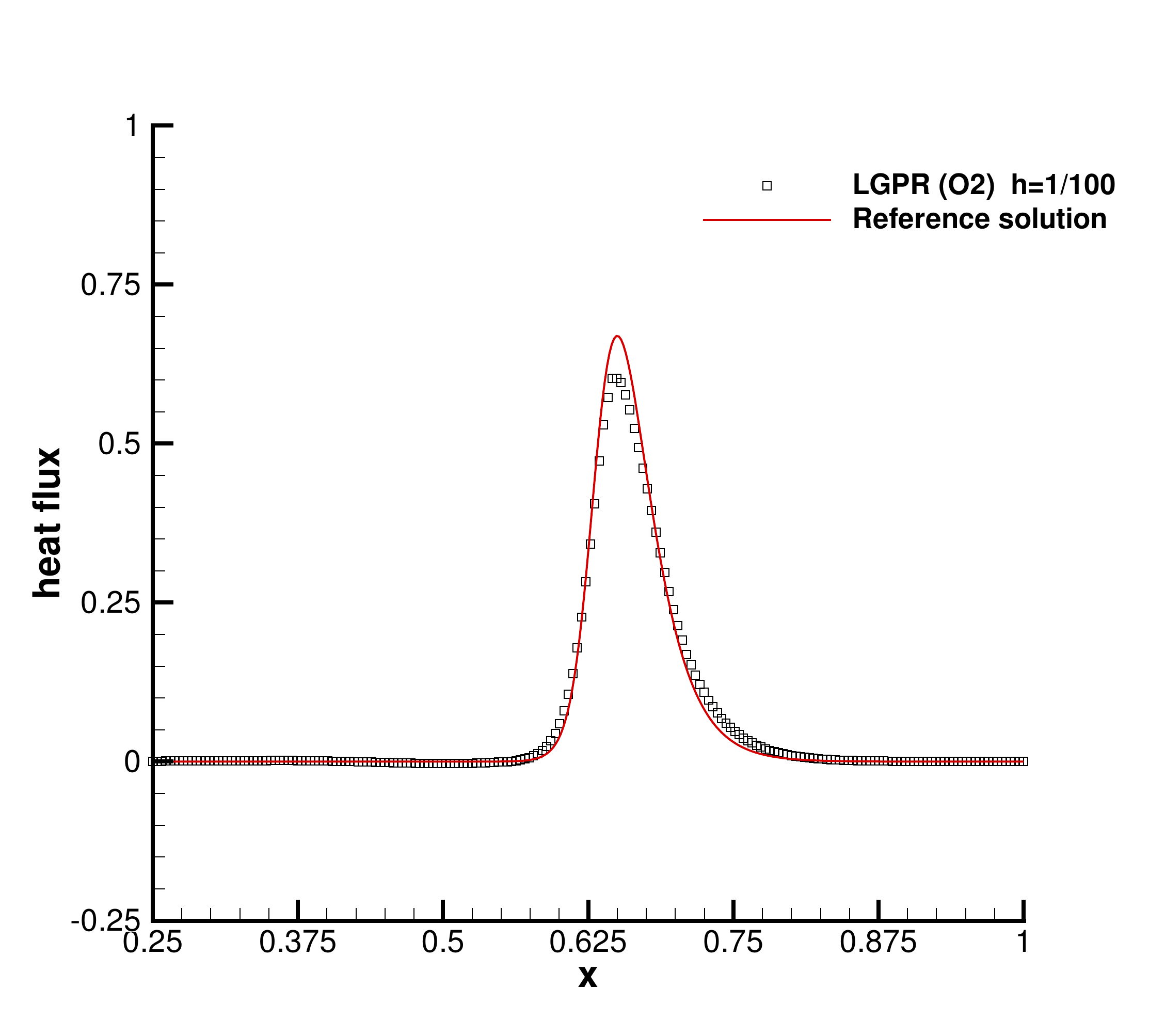}  &          
			\includegraphics[width=0.47\textwidth,draft=false]{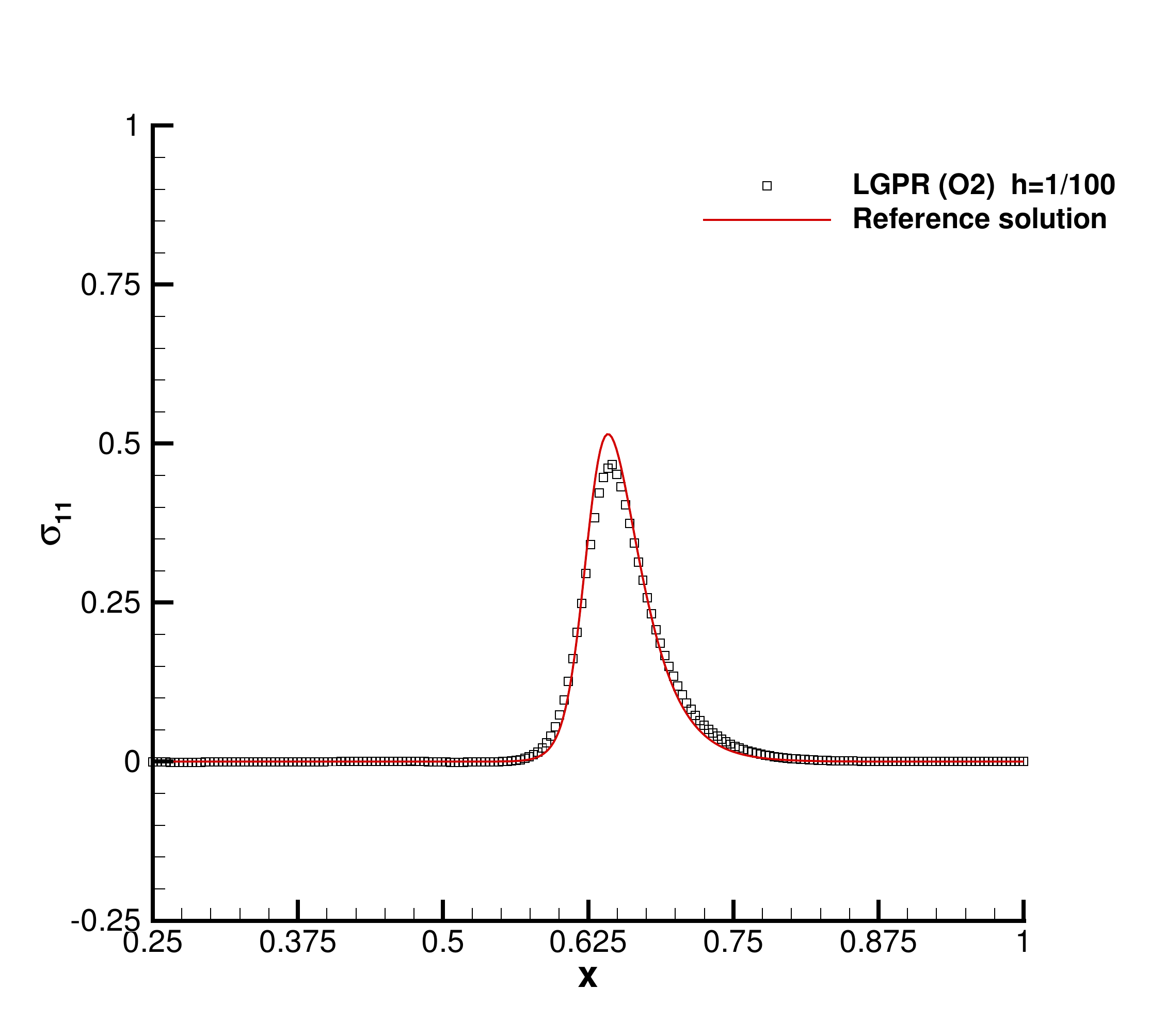} \\
			\includegraphics[width=0.47\textwidth,draft=false]{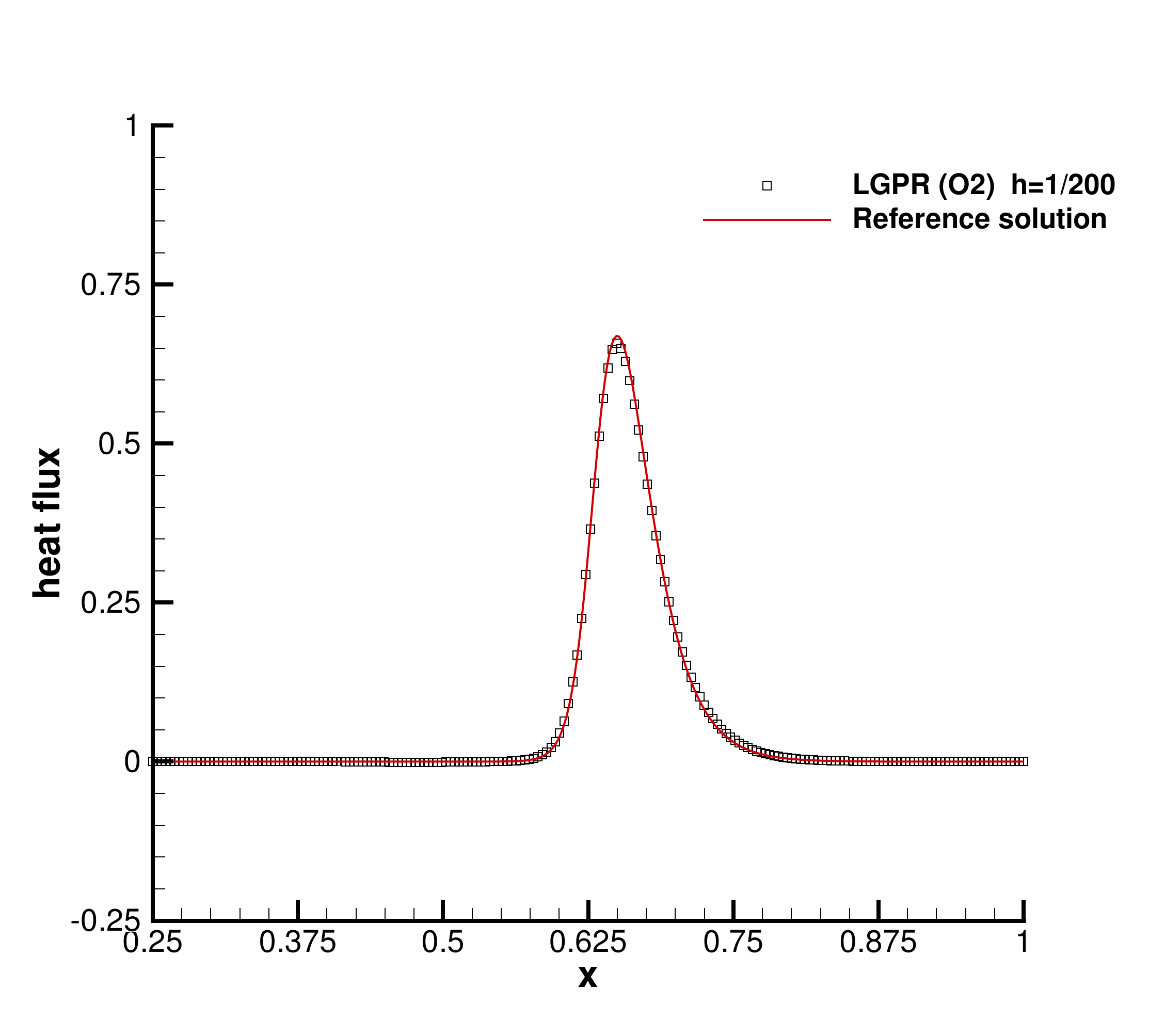}  &          
			\includegraphics[width=0.47\textwidth,draft=false]{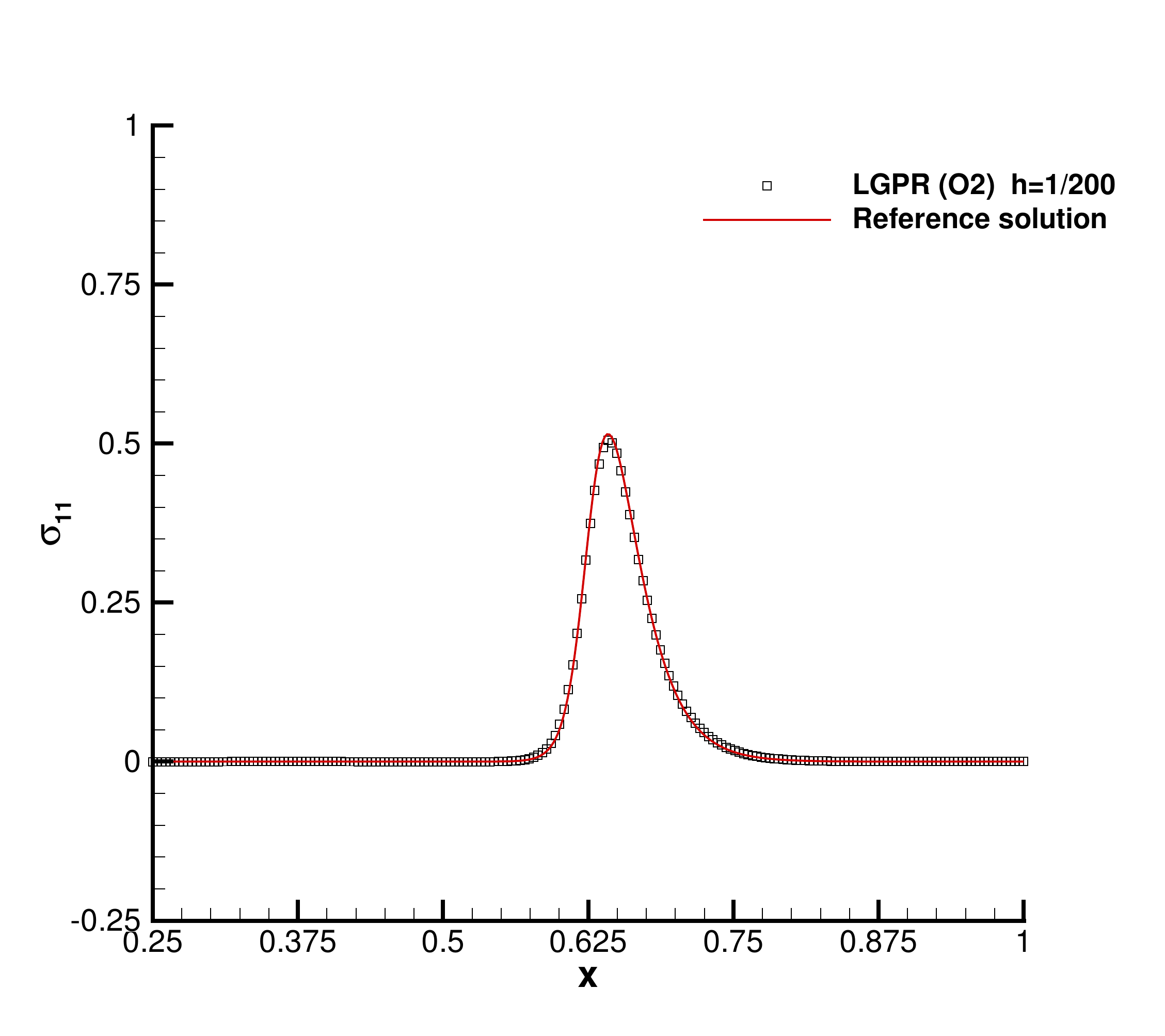} \\
		\end{tabular}
		\caption{Viscous shock profile. Comparison of the exact solution of the compressible Navier–Stokes equations with the Lagrangian GPR model for viscous stress tensor component $\sigma_{11}$ (left) and heat flux (right) with mesh size $h=1/100$ (top) and $h=1/200$ (bottom). For the Navier–Stokes solution, the classical Fourier heat flux $q_1=-\kappa T_x$ is shown, while for the GPR model, we plot $q_1=\alpha^2 T J_1$.}
		\label{fig.shockNS_2}
	\end{center}
\end{figure}

Figure \ref{fig.shockNS_1} illustrates a comparison against the analytical solution at the final 
time, where one can note an excellent matching. We compare the exact solution and the numerical 
solution, extracted as a 1D cut with 200 points along the $x$-direction at $y=0.1$, for density, 
horizontal velocity, pressure and heat flux. Mesh convergence is also qualitatively demonstrated by 
the numerical results obtained with $h=1/100$ and $h=1/200$. Finally, the heat flux and the viscous 
stress component $\tensor{\sigma}_{11}$ are depicted in Figure \ref{fig.shockNS_2} and compared 
against the Navier-Stokes-Fourier model, where the Navier-Stokes stress tensor is recovered in the 
stiff limit by the LGPR scheme, as proven in Section \ref{sec.AP}.
	
\subsection{Collapse of a thick-walled cylindrical beryllium shell} \label{ssec.Shell}
\begin{figure}[!htbp]
	\begin{center}
		\begin{tabular}{cc}
			\includegraphics[width=0.47\textwidth,draft=false]{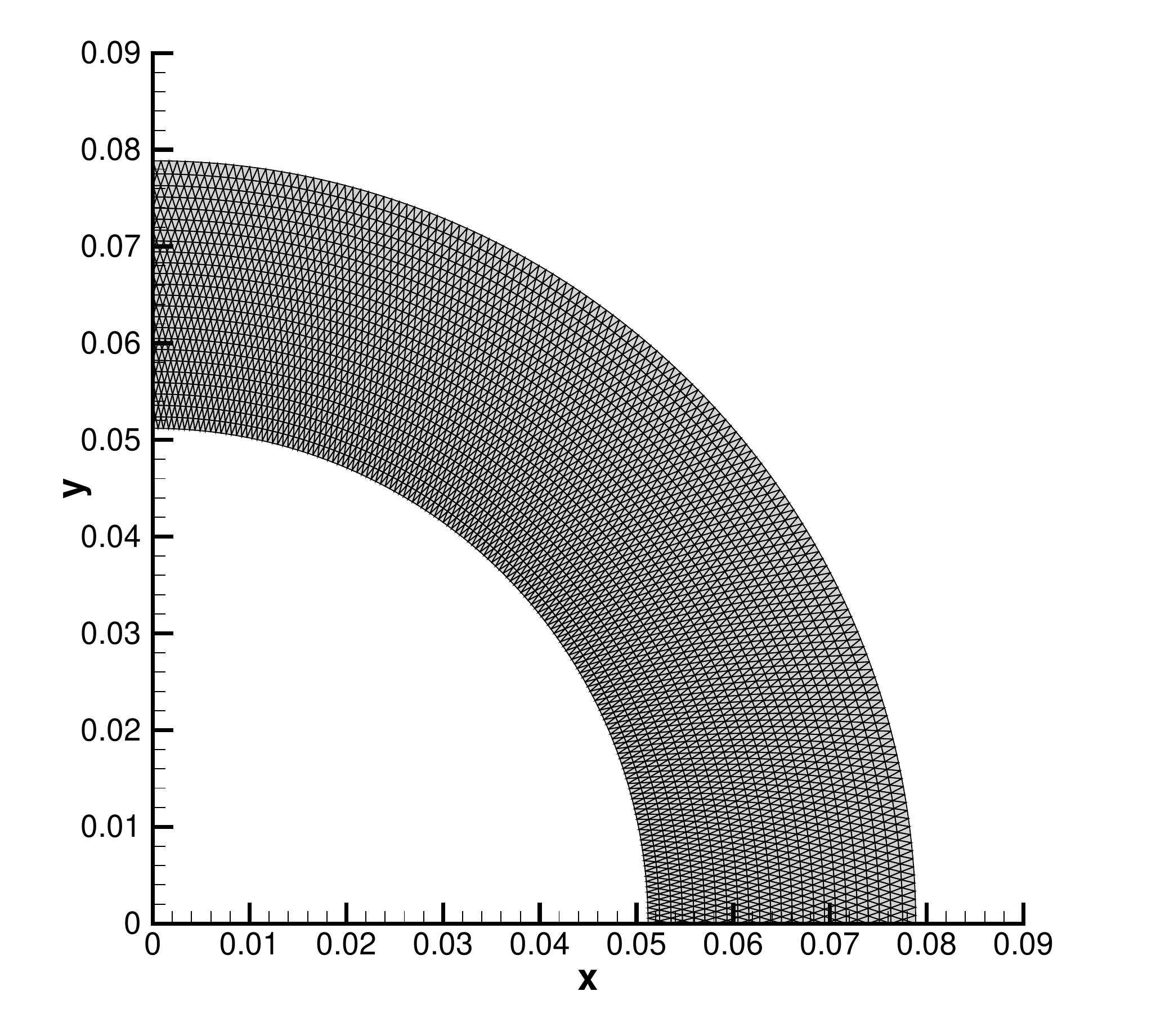}  &               
			\includegraphics[width=0.47\textwidth,draft=false]{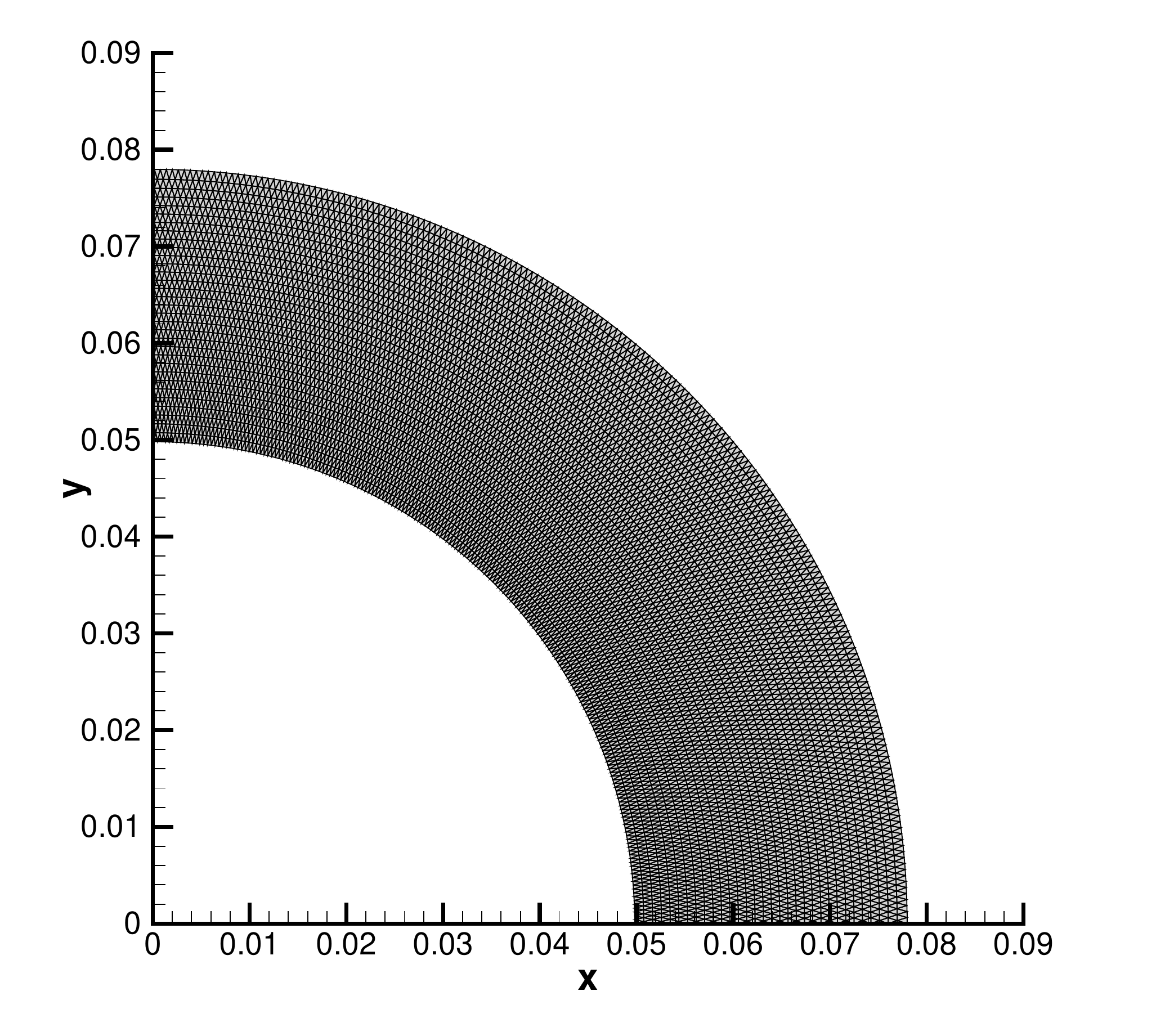} \\
		\end{tabular}
		\caption{Collapse of beryllium shell for the test cases with $V_0^{(1)}=417.1$ (left). Final mesh configuration at time $t_f=125.67 \cdot 10^{-6}$ with mesh size $h=1/100$ (left) and $h=1/140$ (right).}
		\label{fig.shell_mesh}
	\end{center}
\end{figure}

\begin{figure}[!htbp]
	\begin{center}
		\begin{tabular}{ccc}
			\includegraphics[width=0.33\textwidth,draft=false]{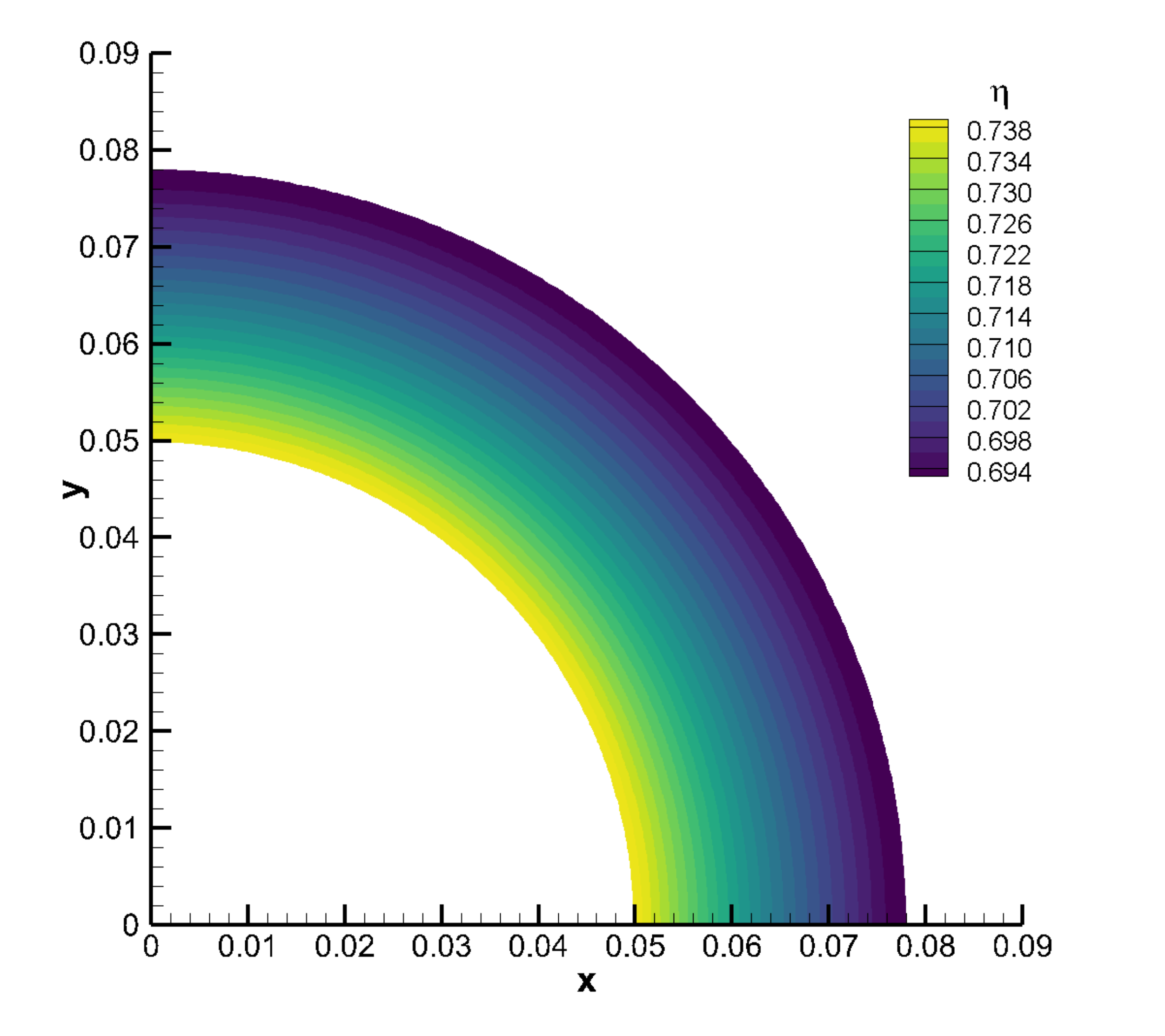}  &          
			\includegraphics[width=0.33\textwidth,draft=false]{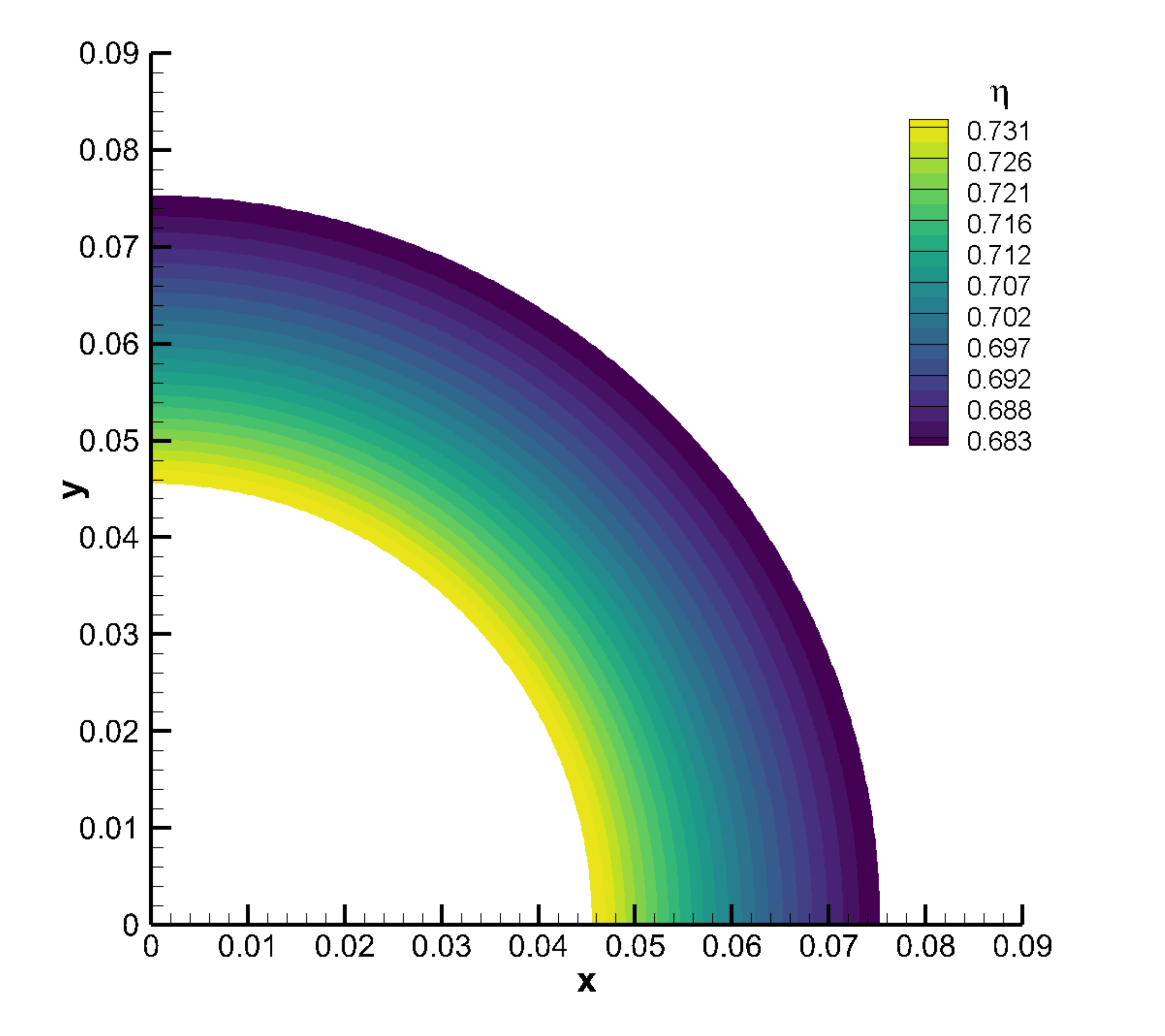} &        
			\includegraphics[width=0.33\textwidth,draft=false]{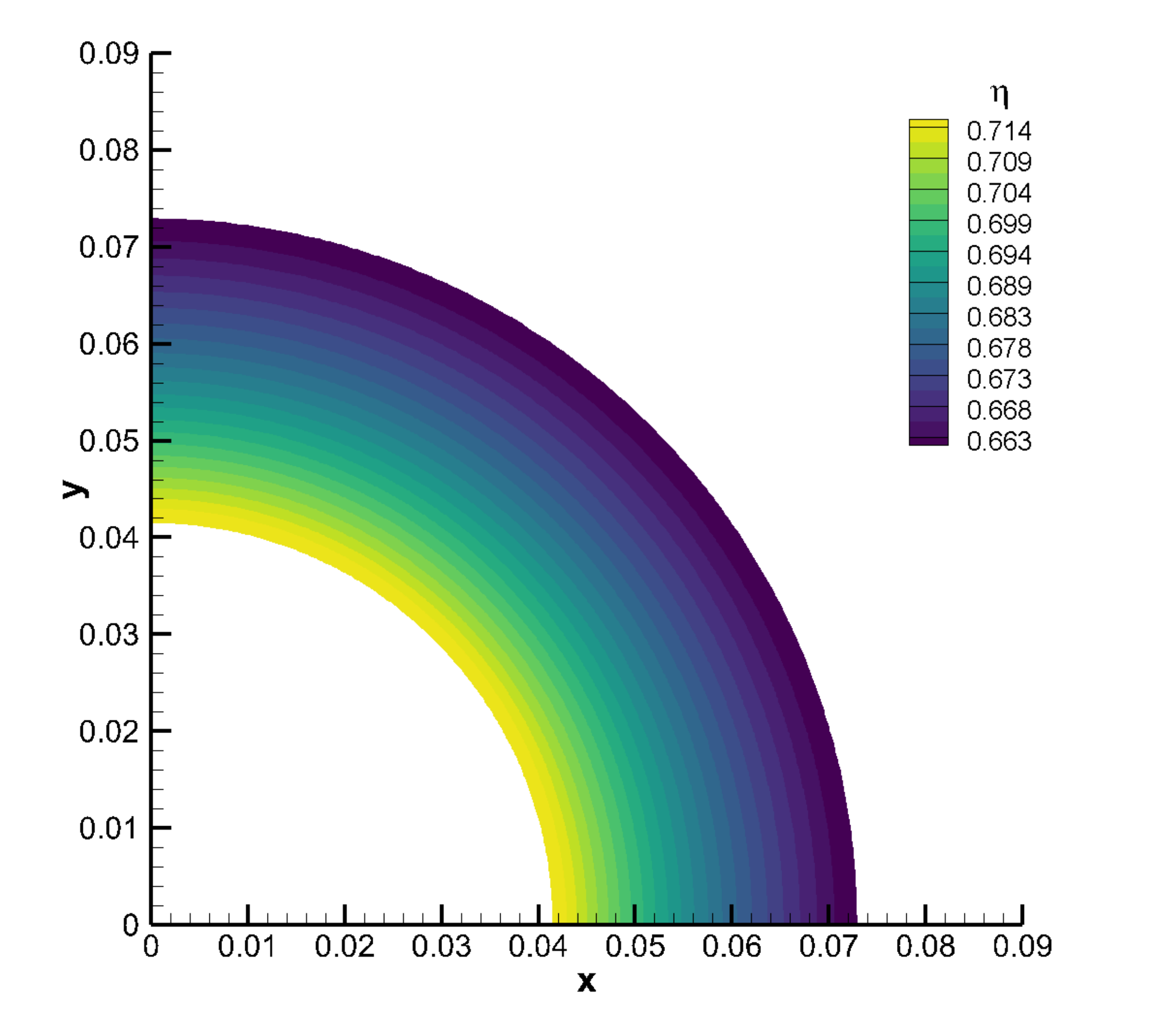} \\
			\includegraphics[width=0.33\textwidth,draft=false]{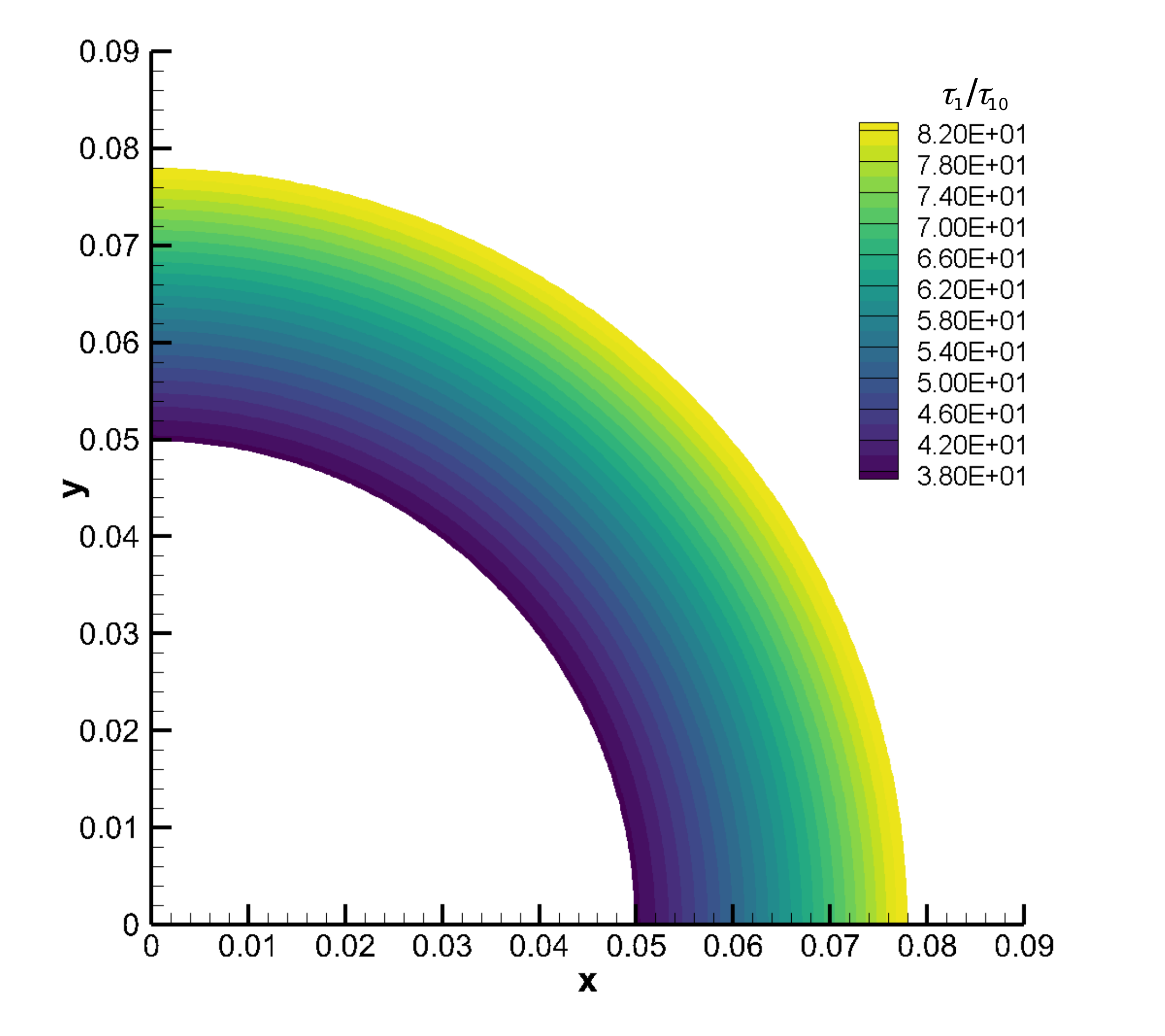}  &          
			\includegraphics[width=0.33\textwidth,draft=false]{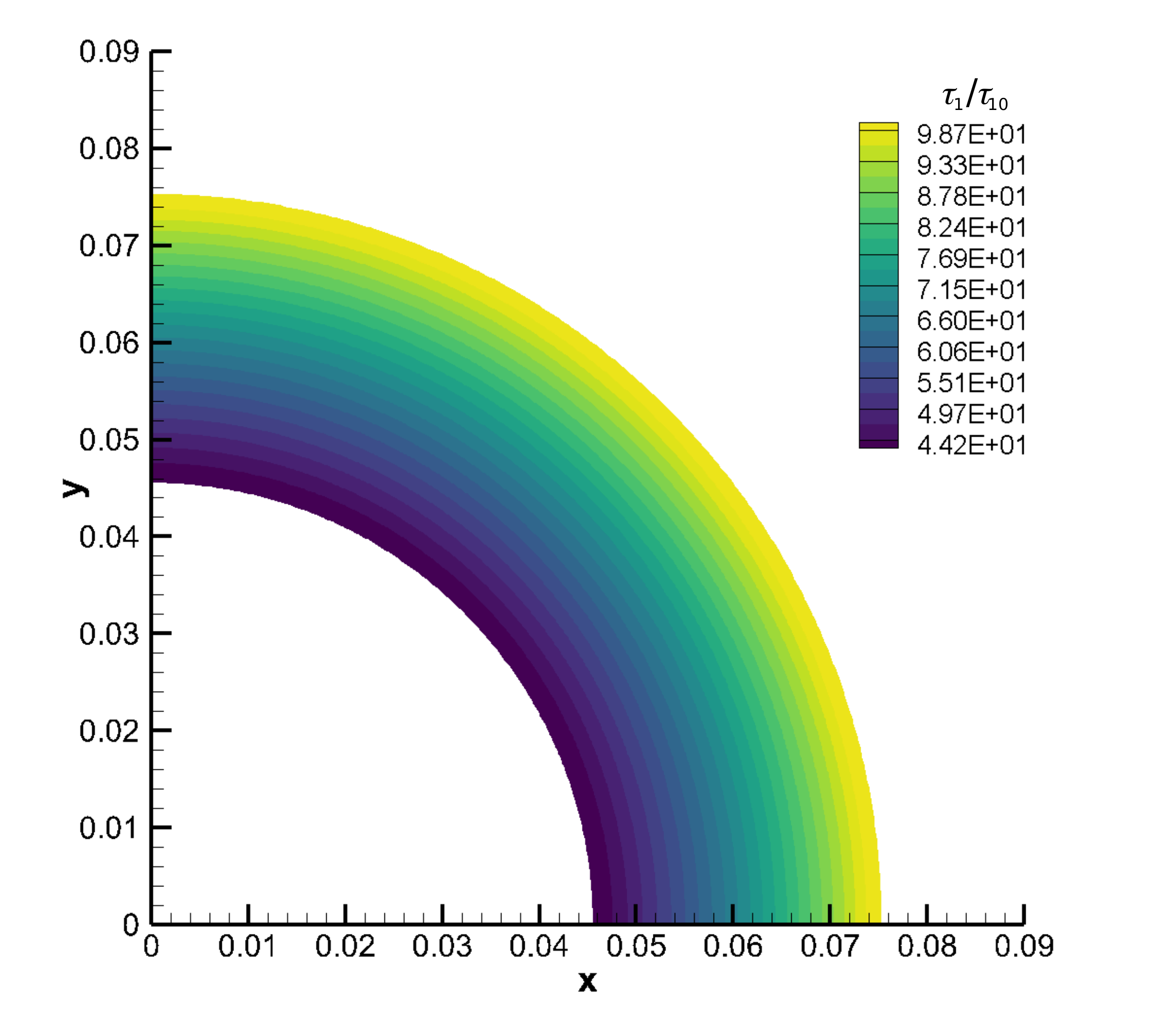} &        
			\includegraphics[width=0.33\textwidth,draft=false]{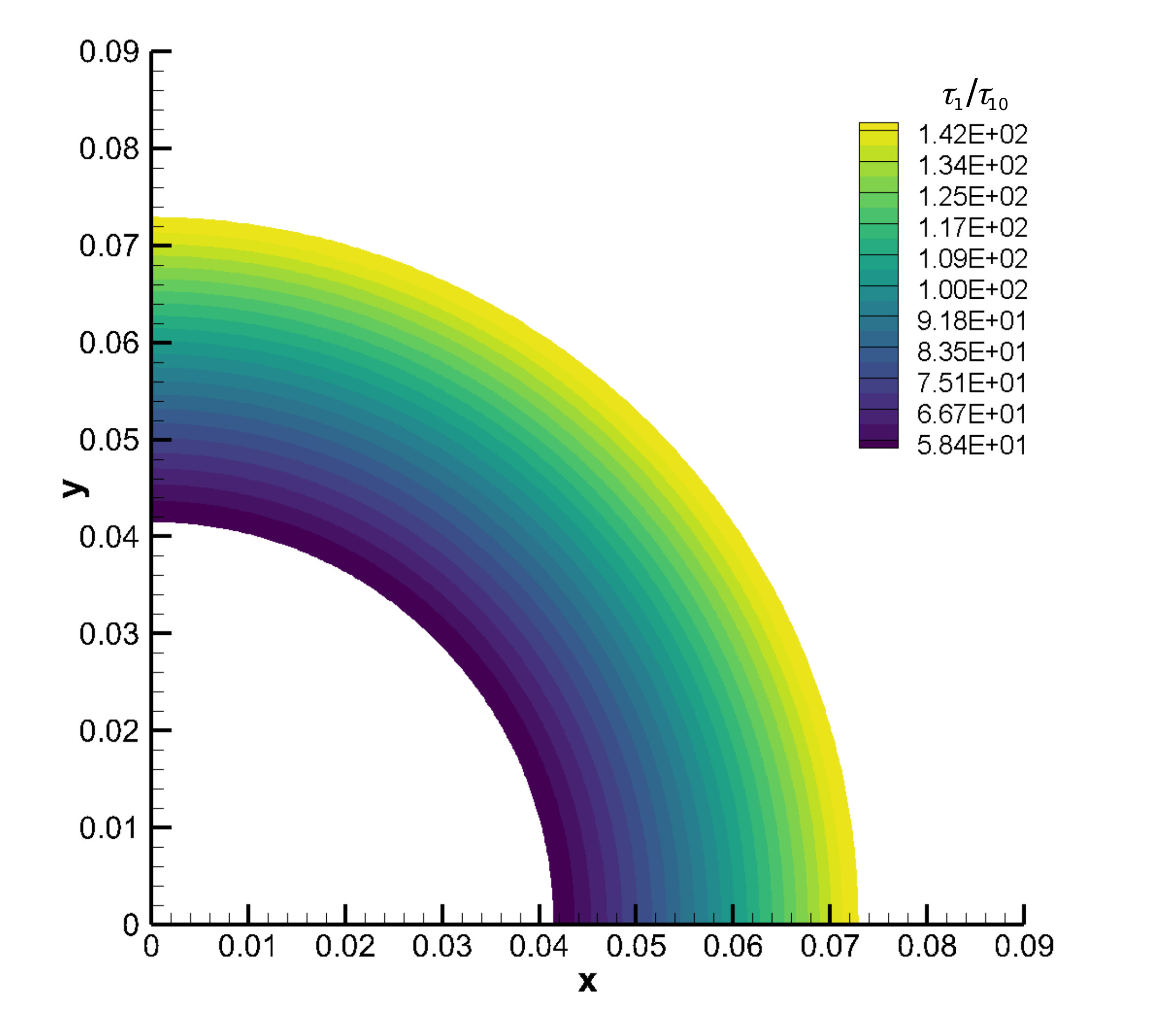} \\
		\end{tabular}
		\caption{Collapse of beryllium shell for the three test cases: $V_0^{(1)}=417.1$ (left), 
		$V_0^{(2)}=454.7$ (middle) and $V_0^{(3)}=490.2$ (right) with mesh size $h=1/140$. Numerical 
		distribution of plastic map ($\eta$) and normalized relaxation time ($\tau_1/\tau_{10}$) at 
		the final time of each test case.}
		\label{fig.shell1}
	\end{center}
\end{figure}

The next test case aims at exploiting the capability of the GPR model \eqref{eqn.cl} to simulate 
elasto-plastic solids. We consider a test problem firstly proposed in \cite{Howell2002}, which 
describes the collapse of a cylindrical beryllium shell responding to an initial radial velocity 
field directed towards its center. The initial setup is taken from \cite{KammLANL08}, thus the 
computational domain $\Omega(0)$ is the shell with inner and outer radii $r_{int}=8\times 10^{-2}$ 
and 
$r_{ext}=10\times10^{-2}$. Free-traction boundary conditions are considered everywhere. The 
material is 
beryllium and the Mie-Grüneisen equation of state \eqref{eqn.MG} is used with parameters 
$\Gamma_0=1.16$ and $s=1.124$. The adiabatic sound speed results to be $c_0=12870$. The initial 
pressure is $p=0$ and the radial velocity magnitude is given by  
$v_r(0,r)=-V_0\left(\frac{r_{int}}{r}\right)^2$, with $r=\sqrt{x^2+y^2}$. According to 
\cite{Aguirre2014}, three different 
values of $V_0$ are considered, namely $V_0^{(1)}=417.1$, $V_0^{(2)}=454.7$ and $V_0^{(3)}=490.2$. The corresponding 
final times are $t_f^{(1)}=125.67 \cdot 10^{-6}$, $t_f^{(2)}=131.6 \cdot 10^{-6}$ and $t_f^{(3)}=136.26 \cdot 10^{-6}$. The 
initial kinetic energy due to the velocity distribution is entirely dissipated by the plastic 
deformation of the material leading to a deceleration of the shell. Therefore, in \cite{Howell2002} 
a closed form solution at the stopping time is derived under the ideal plasticity assumption, which 
leads to a relationship between the initial velocity $V_0$ and the inner and outer stopping radii. 
In this test case, the material experiences elasto-plastic deformations, thus the closure relations 
\eqref{eqn.tau.plast} are adopted to account for plasticity effects. Specifically, the yield 
strength is $\sigma_Y=330 \cdot 10^{6}$ and the parameters of the power law for the computation of 
the relaxation time $\tau_1$ are the exponent $n=12$ and $\tau_{10}=10^{-7}$. The test case is run 
for all the three velocities $V_0$ on two different unstructured meshes with characteristic mesh 
size of $h=1/100$ and $h=1/140$, in order to show mesh convergence of the numerical solution. The 
final mesh configuration for $V_0=417.1$ is depicted in Figure \ref{fig.shell_mesh}, while the 
plastic map $\eta=\sigma/\sigma_Y$ and the normalized relaxation time $\tau_1/\tau_{10}$ at the 
final time of each simulation are shown in Figure \ref{fig.shell1}.

\begin{figure}[!htbp]
	\begin{center}
		\begin{tabular}{ccc}
			\includegraphics[width=0.33\textwidth,draft=false]{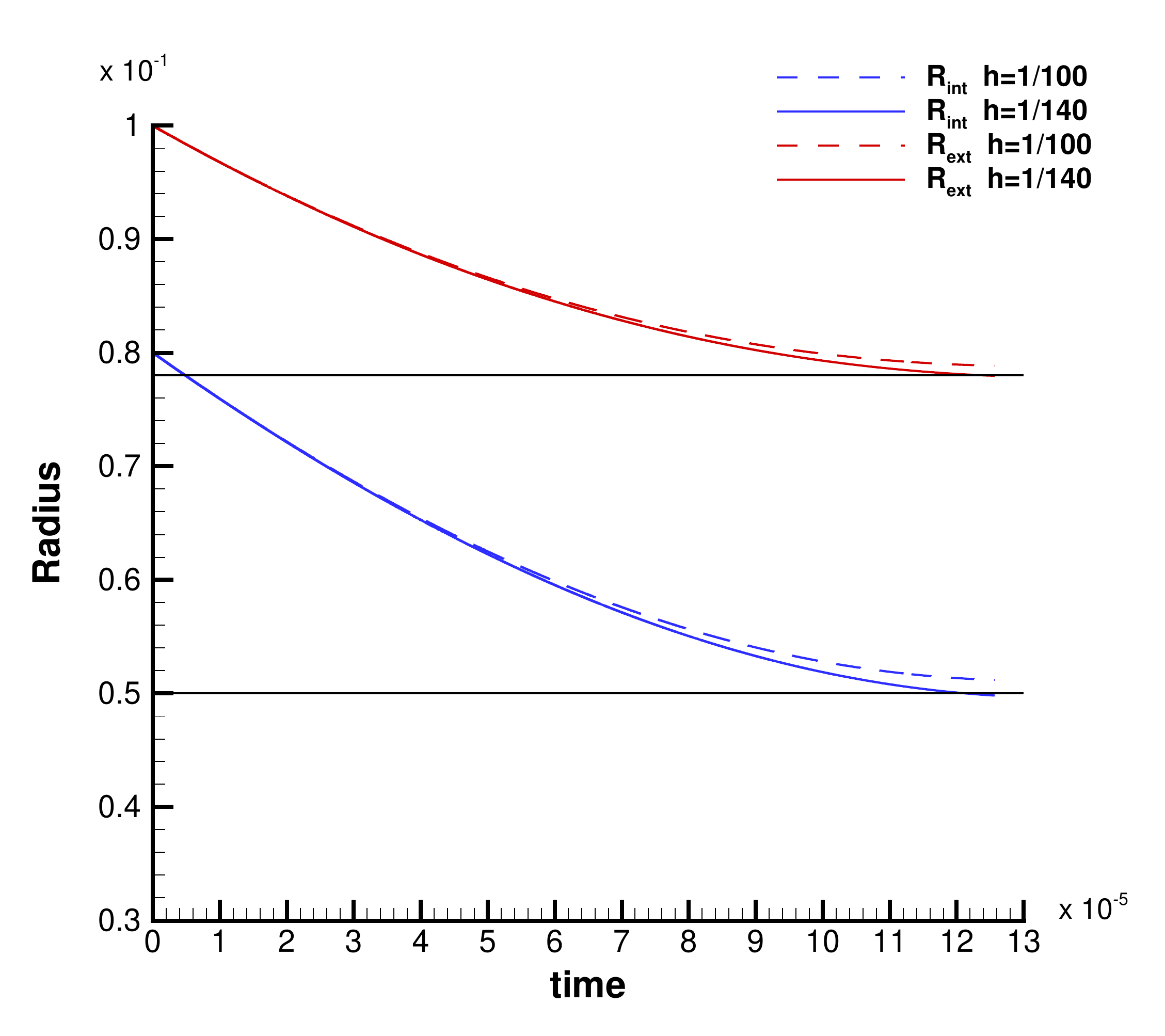}  &          
			\includegraphics[width=0.33\textwidth,draft=false]{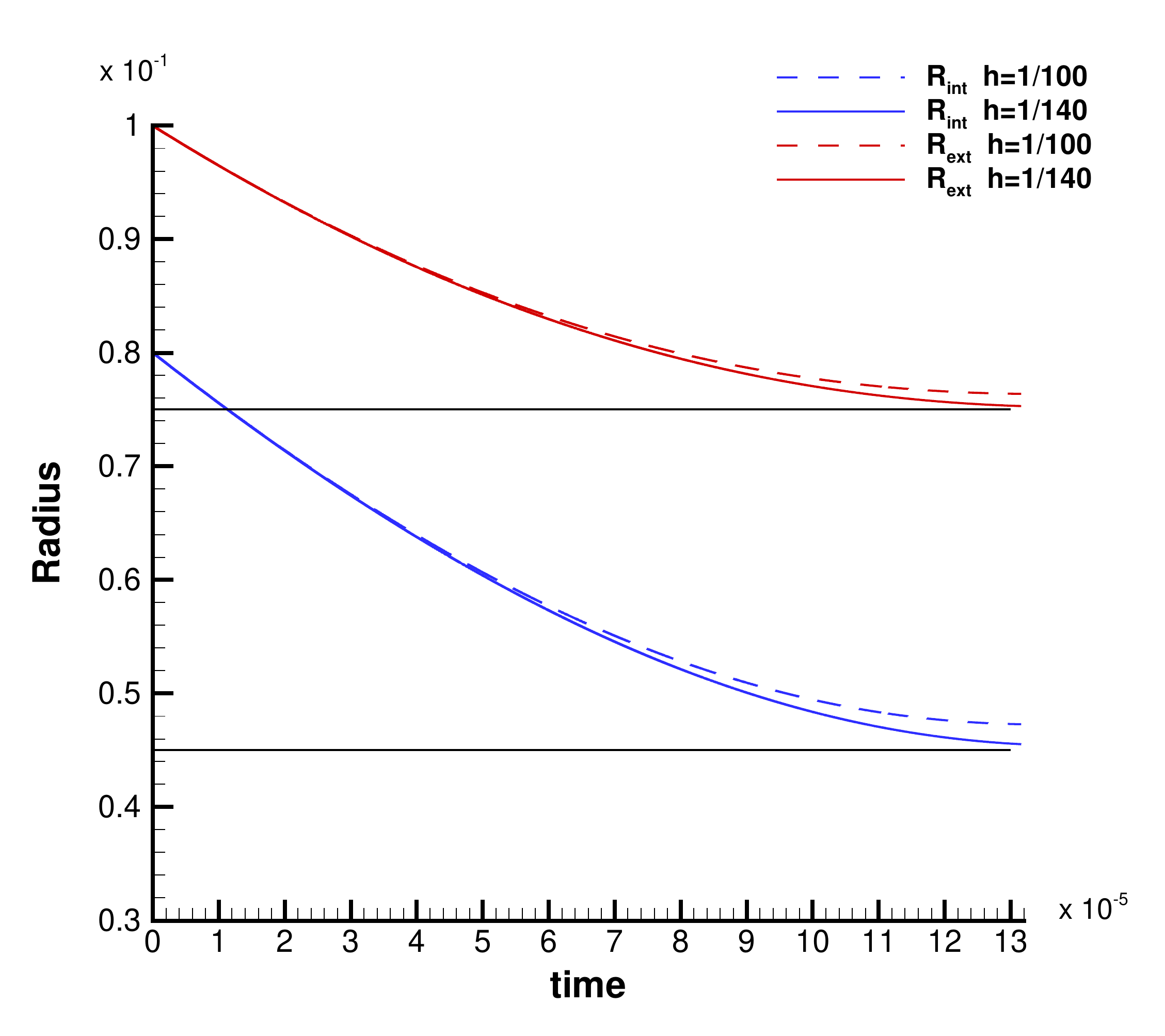} &        
			\includegraphics[width=0.33\textwidth,draft=false]{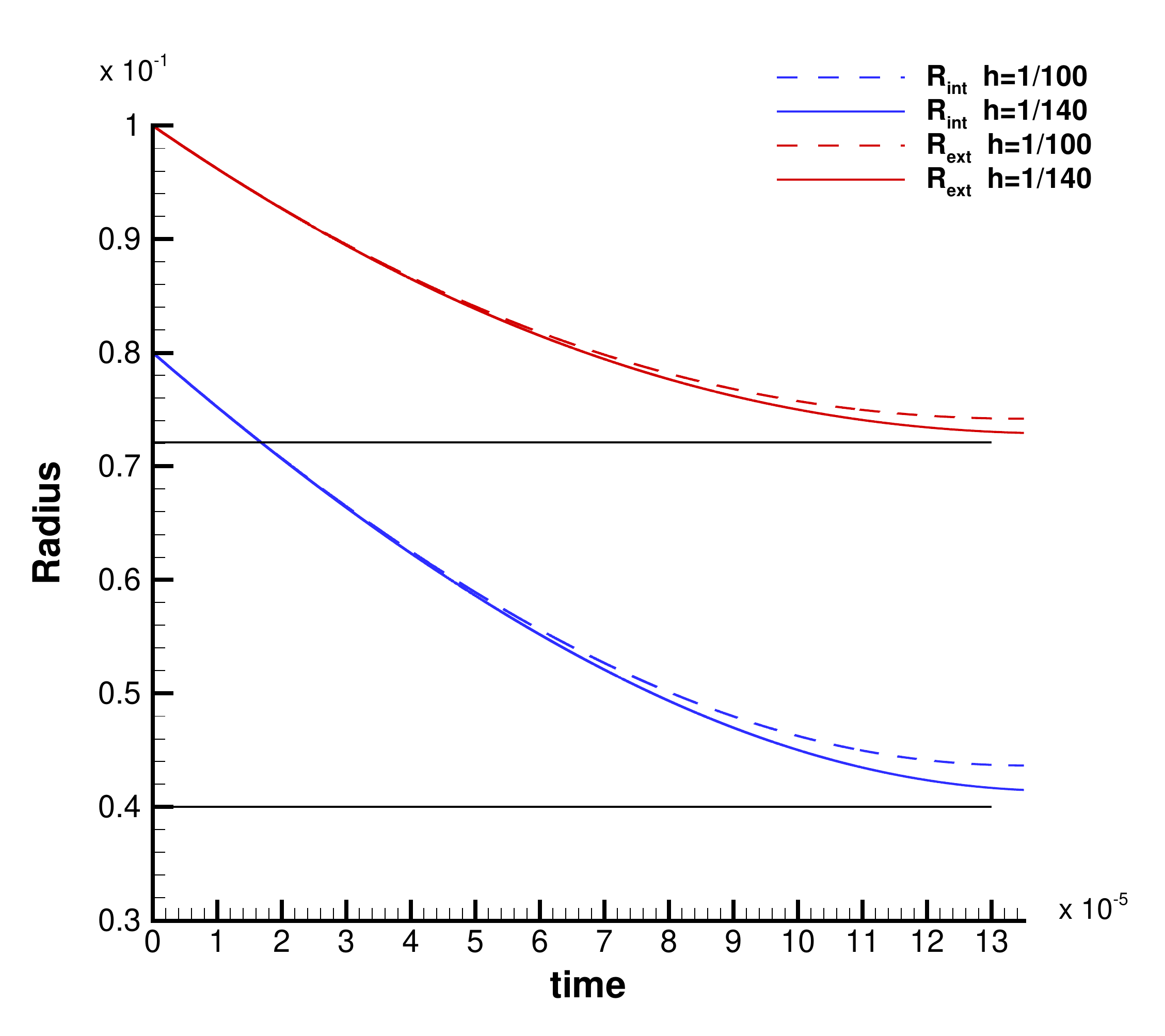} \\
			\includegraphics[width=0.33\textwidth,draft=false]{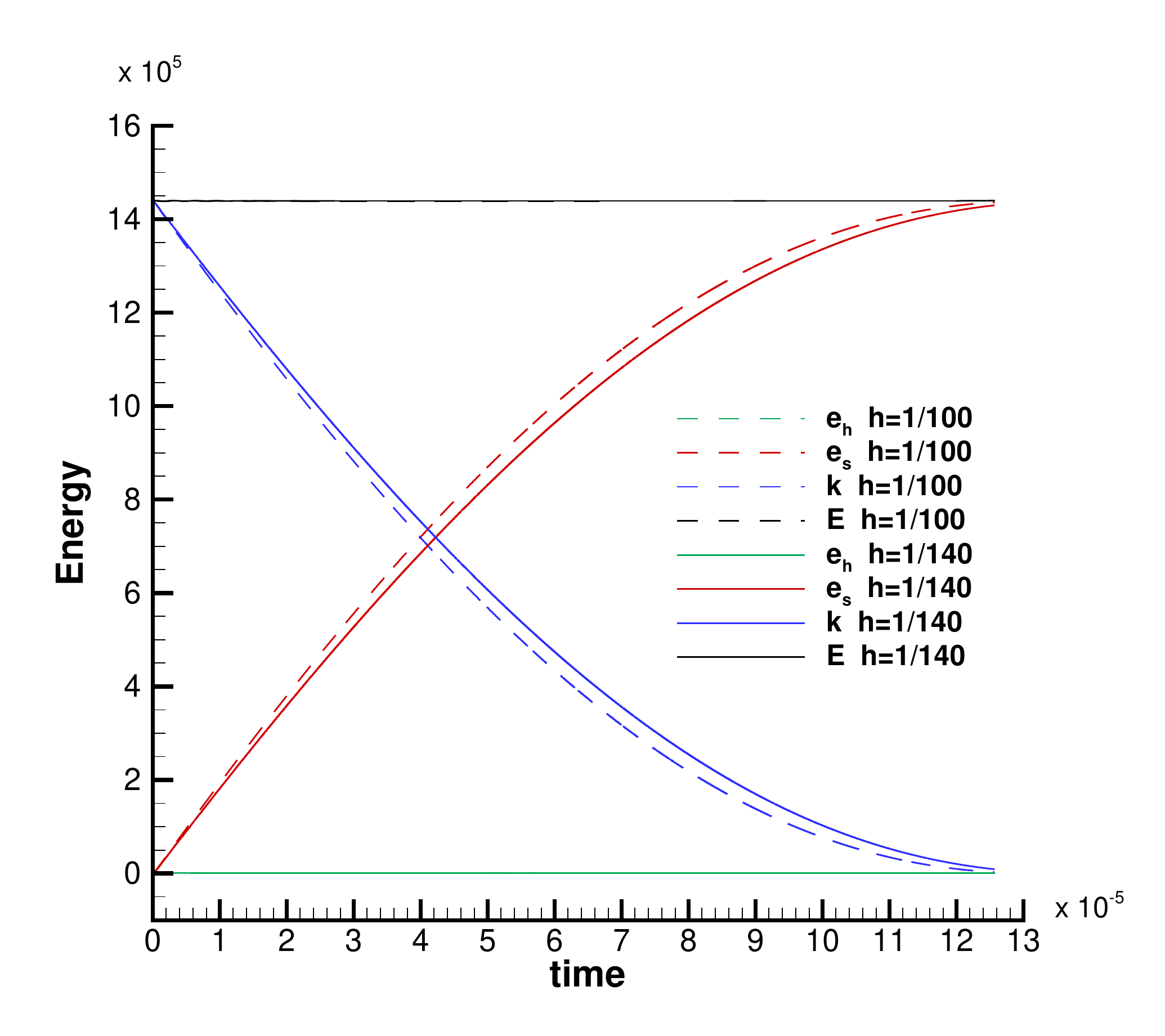}  &          
			\includegraphics[width=0.33\textwidth,draft=false]{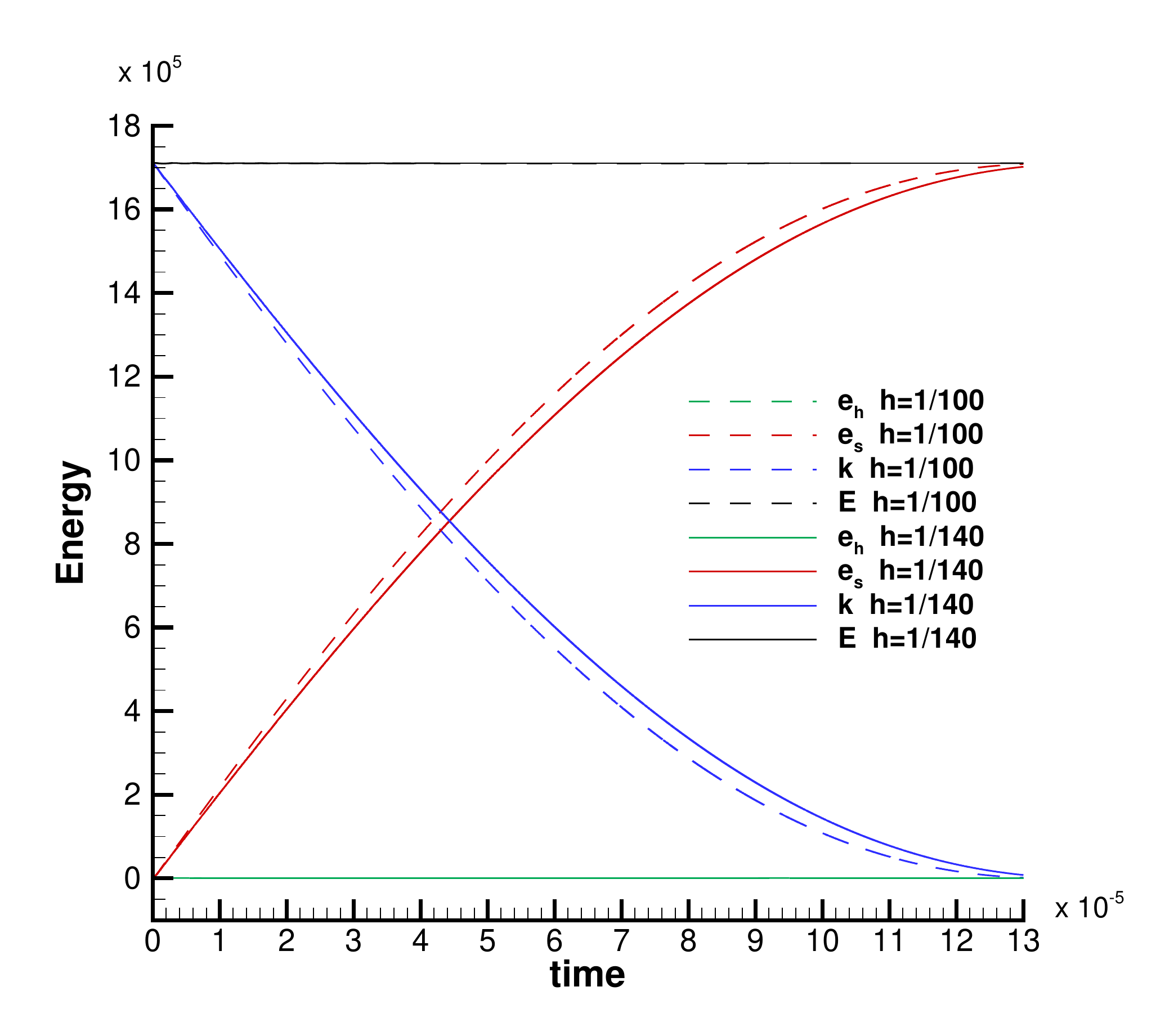} &        
			\includegraphics[width=0.33\textwidth,draft=false]{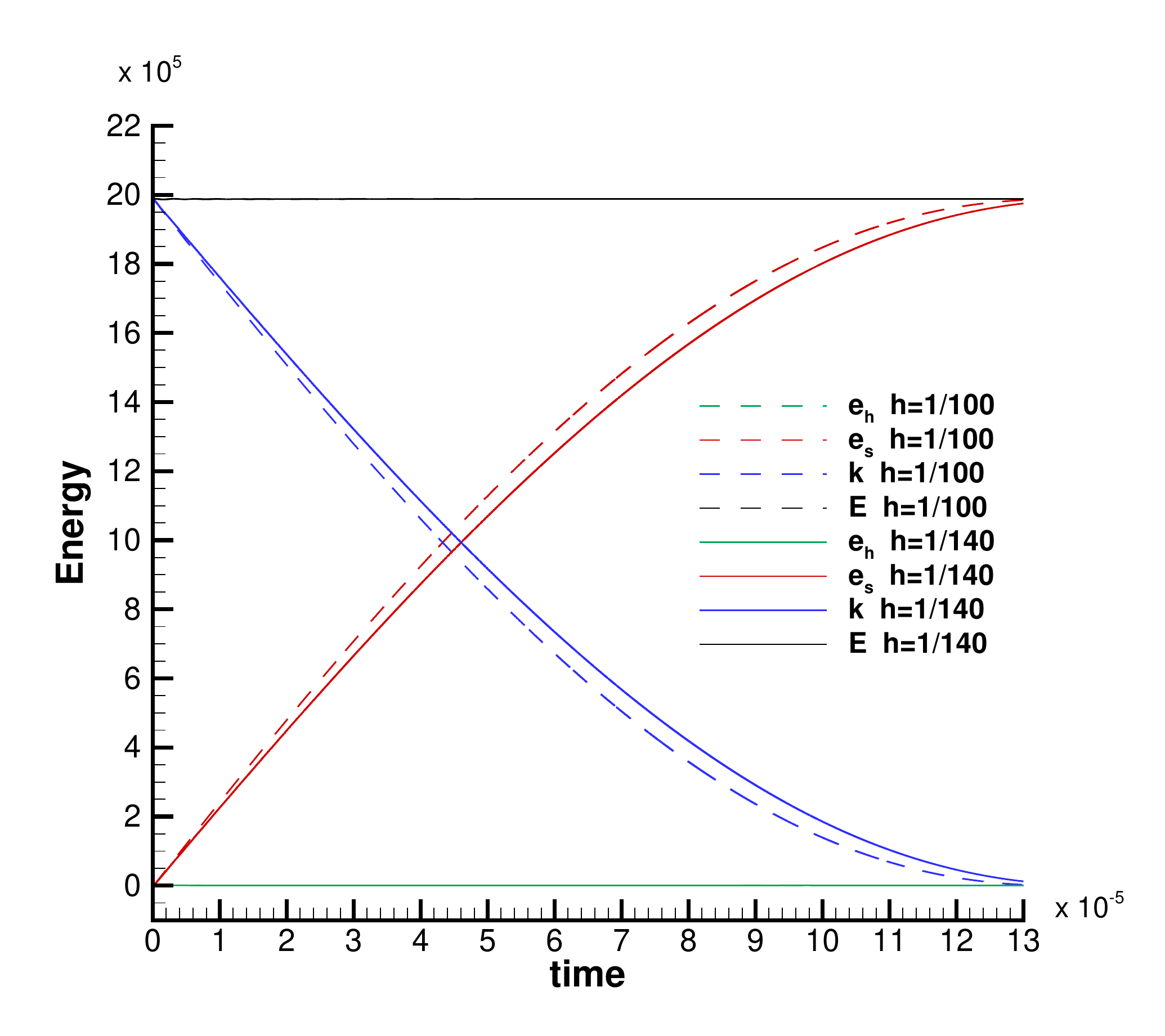} \\
		\end{tabular}
		\caption{Collapse of beryllium shell for the three test cases: $V_0^{(1)}=417.1$ (left), $V_0^{(2)}=454.7$ (middle) and $V_0^{(3)}=490.2$ (right). Top: time evolution of the internal $R_{int}$ and external $R_{ext}$ radius of the shell and comparison between analytical and numerical solution. Bottom: analysis of energy conservation in terms of volumetric and shear internal energy ($E_h$ and $E_e$, respectively),  kinetic ($E_k$) and total ($E$) energy. Results with mesh size $h=1/100$ are drawn with dashed lines, while solid lines refer to mesh size $h=1/140$.}
		\label{fig.shell2}
	\end{center}
\end{figure}

The time evolution of the internal and external frontiers of the shell, i.e. $R_{int}(t)$ and $R_{ext}(t)$, is plot in Figure \ref{fig.shell2}, which is compared against the exact displacement of the shell at the final time of each simulation. A good agreement can be observed and the finer mesh correctly retrieves a more accurate solution. We underline that the usage of a high quality, though unstructured, computational mesh is crucial for maintaining the symmetry of the numerical solution as demonstrated by Figure \ref{fig.shell1}. Finally, the analysis of energy conservation over time is plot in Figure \ref{fig.shell2}, where we report the volumetric and shear internal energy ($E_h$ and $E_e$, respectively),  kinetic ($E_k$) and total ($E$) energy. We clearly see that all initial kinetic energy is dissipated into internal energy due to elasto-plastic deformations that occur in the material. Indeed, the volumetric energy contribution is rather small compared to the shear internal energy $E_e$.

\subsection{2D projectile impact} \label{ssec.TaylorBar2D}
\begin{figure}[!htbp]
	\begin{center}
		\begin{tabular}{cc}
			\includegraphics[width=0.47\textwidth,draft=false]{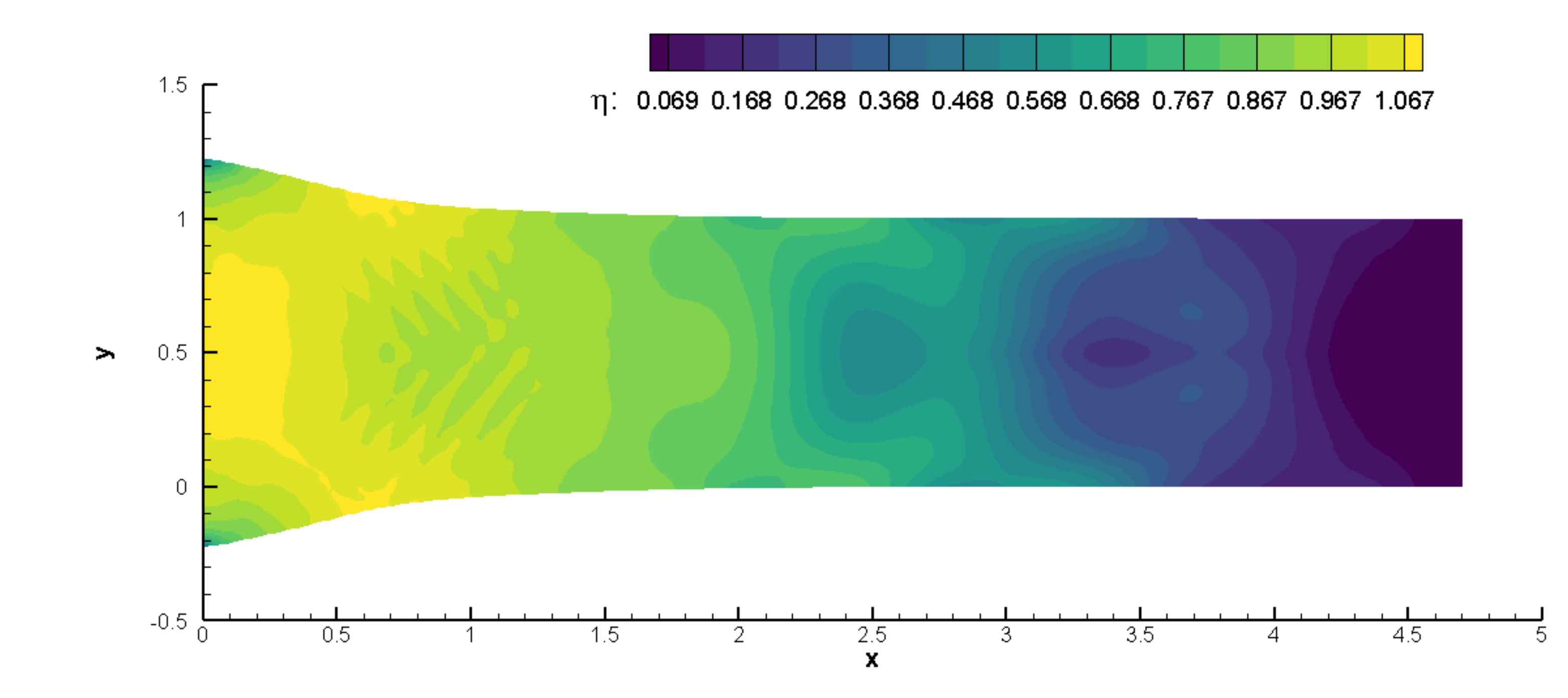}  &          
			\includegraphics[width=0.47\textwidth,draft=false]{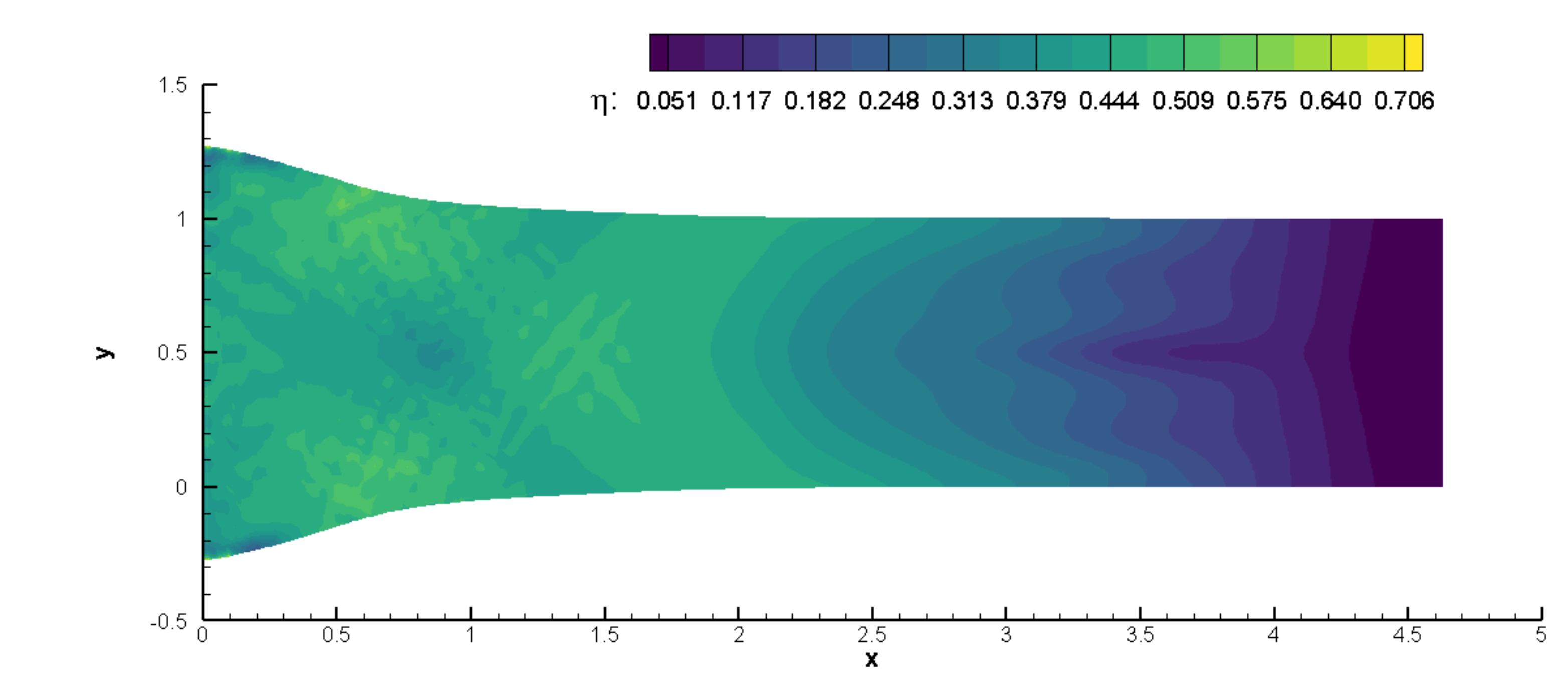} \\
			\includegraphics[width=0.47\textwidth,draft=false]{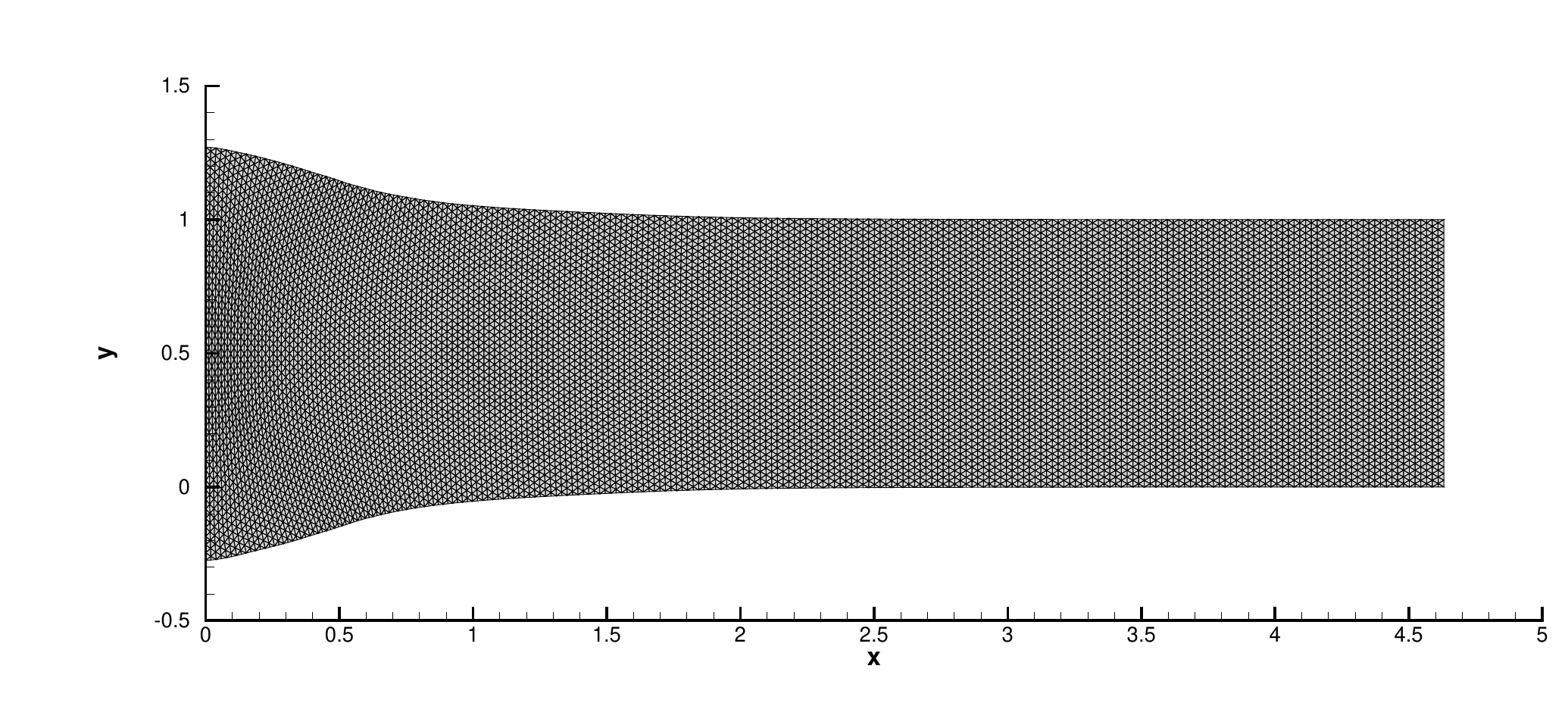}  
			&          
			\includegraphics[width=0.47\textwidth,draft=false]{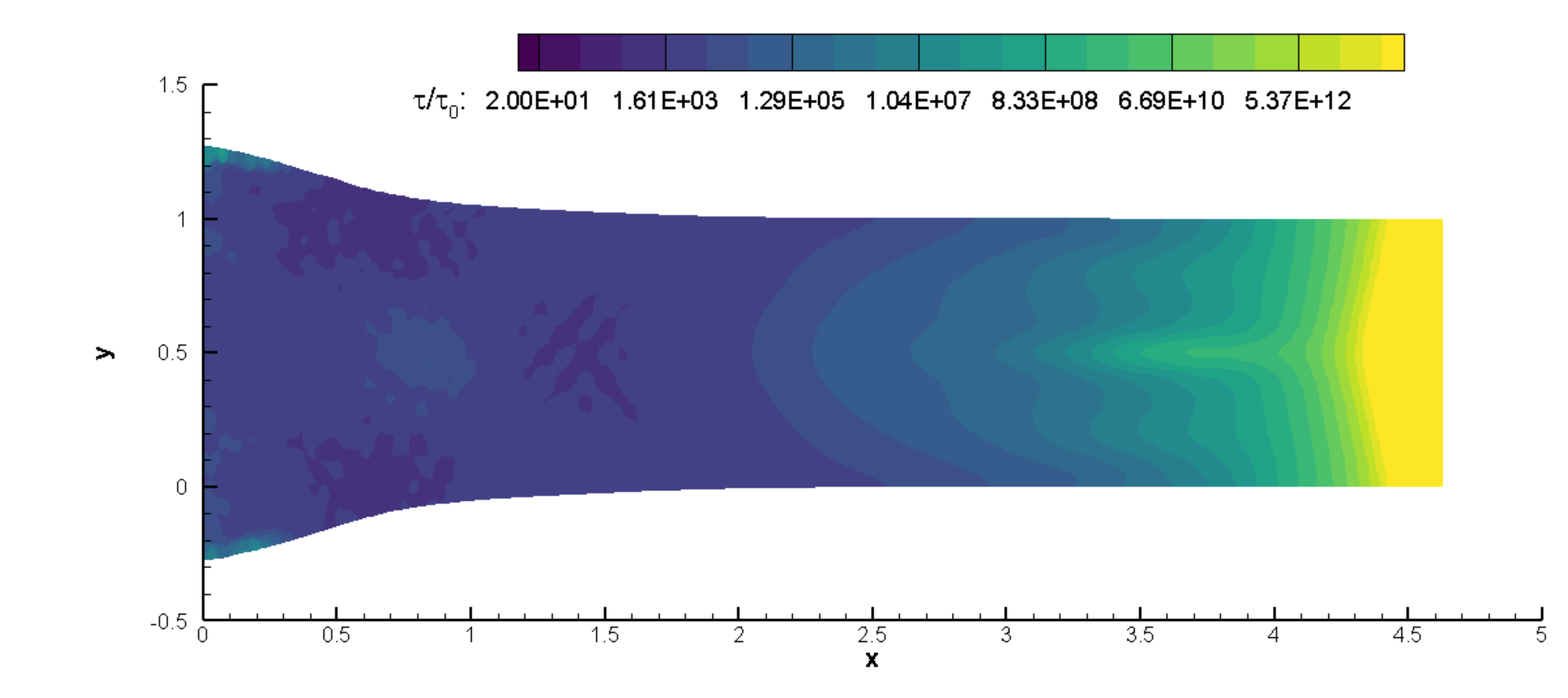} \\
		\end{tabular}
		\caption{2D projectile impact. Top: numerical distribution of plastic map ($\eta$) at time 
		$t=2.5\cdot 10^{-3}$ (left) and $t=5\cdot 10^{-3}$ (right). Bottom: numerical distribution 
		of normalized relaxation time ($\tau/\tau_{10}$) at the final time $t=5\cdot 10^{-3}$ 
		(right) and mesh configuration (left).}
		\label{fig.TaylorBar2D}
	\end{center}
\end{figure}

\begin{figure}[!htbp]
	\begin{center}
		\begin{tabular}{cc}
			\includegraphics[width=0.47\textwidth,draft=false]{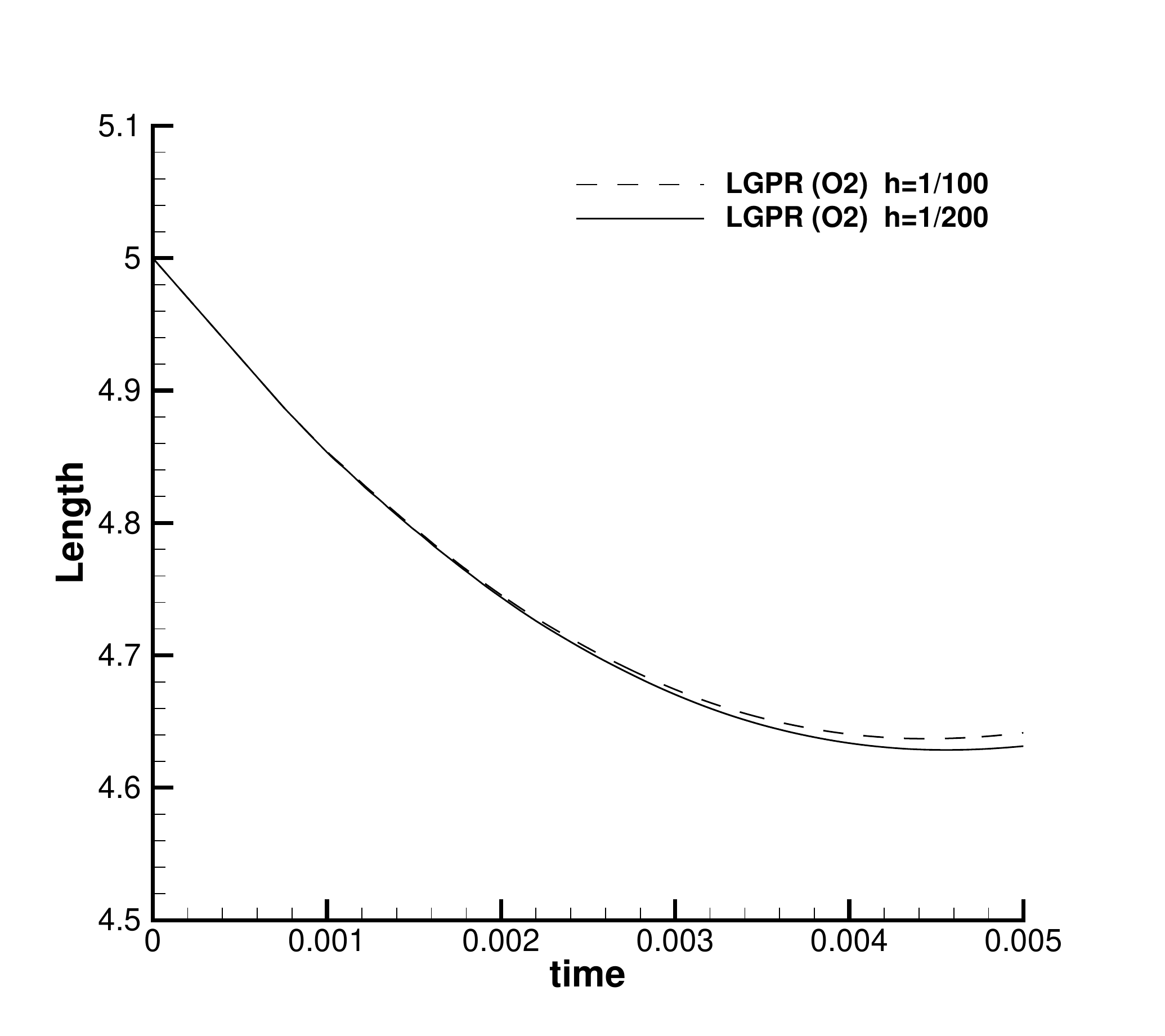}  &          
			\includegraphics[width=0.47\textwidth,draft=false]{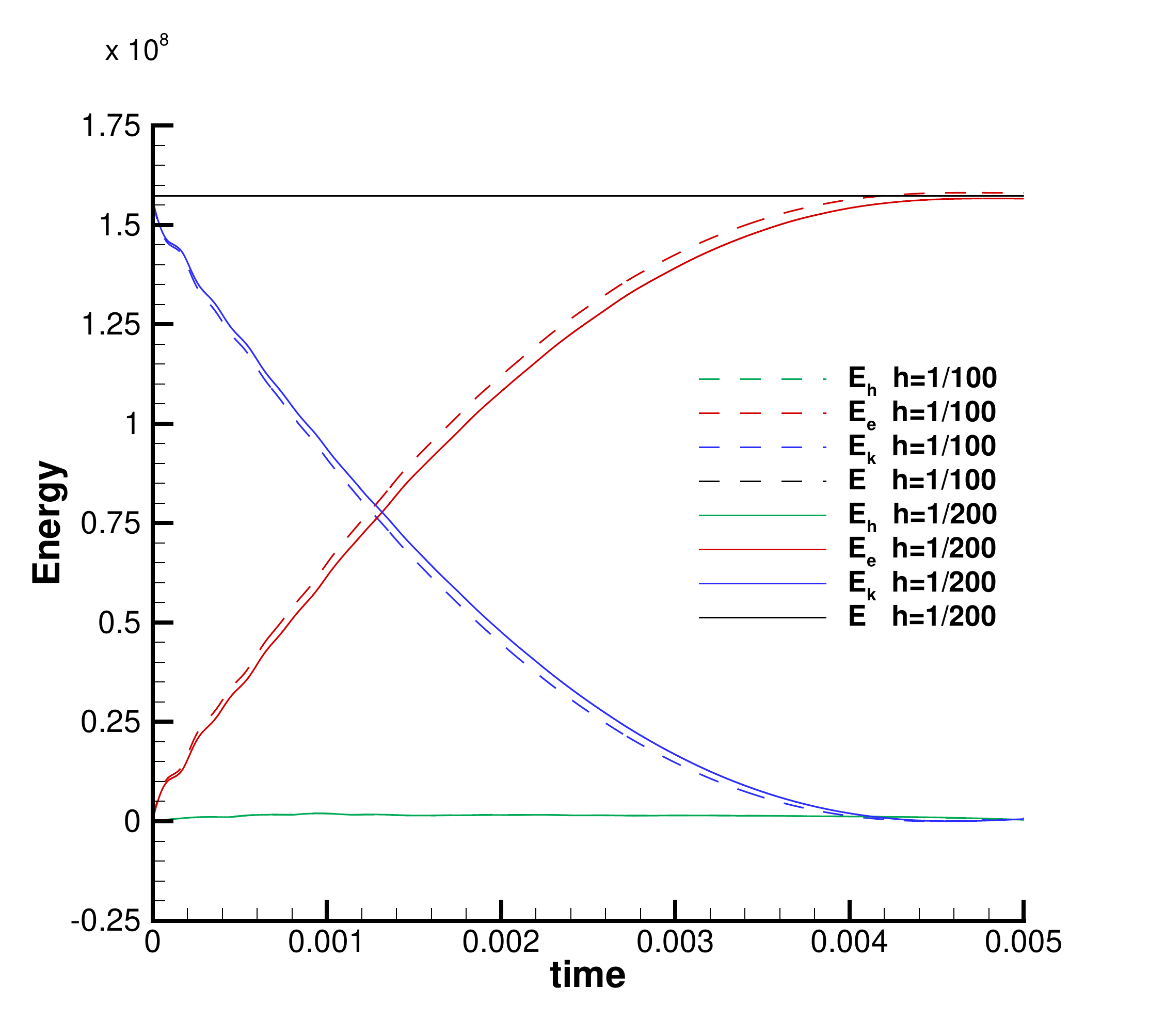} \\
		\end{tabular}
		\caption{2D projectile impact. Left: time evolution of the projectile length. Right: 
		analysis of energy conservation in terms of volumetric and elastic energy ($E_h$ and $E_e$, 
		respectively),  kinetic ($E_k$) and total ($E$) energy. Results with mesh size $h=1/100$ 
		are drawn with dashed lines, while solid lines refer to mesh size $h=1/200$. }
		\label{fig.TaylorBar2D-energy}
	\end{center}
\end{figure}

This problem consists of the impact of a two-dimensional aluminum bar impacting on a
rigid wall and the setup is taken from \cite{Maire_elasto}. The computational domain is the initial 
projectile $\Omega(0)=[0;5]\times [0;1]$ that is paved with two different triangular meshes of 
characteristic size $h=1/100$ and $h=1/200$. A slip-wall boundary is set on the left side of the 
domain, while free-traction boundary conditions are imposed elsewhere. The final time of the 
simulation is $t_f=0.005$ and the material undergoes plastic deformations that convert the initial 
kinetic energy into shear internal energy. The closure relation for the hydrodynamics part of the 
energy  is the Mie-Grüneisen EOS with 
$\Gamma_0=2$ and $s=1.338$, while the relaxation time $\tau_1$ is computed relying on the 
formulation \eqref{eqn.tau.plast} with $n=10$ and $\tau_{10}=5 \cdot 10^{-4}$. The yield stress for 
aluminum is $\sigma_Y=300 \cdot 10^{6}$ and the adiabatic sound speed is $c_0=5328$. Initially, the 
material is assigned $p=0$ with a velocity set to $\vv=(-150,0,0)$, thus the projectile hits the 
wall located at the left side of the domain. Although there exists no exact solution for this 
problem, it is nonetheless employed for verifying robustness and accuracy. Figure 
\ref{fig.TaylorBar2D} depicts the plastic map $\eta=\sigma/\sigma_Y$ at different output times as 
well as the final mesh configuration and the normalized relaxation time $\tau/\tau_{10}$. 
Plasticity effects are experienced by the material during the impact and the stresses at the end of 
the simulation visible in Figure \ref{fig.TaylorBar2D} are residual stresses. This observation is 
also 
confirmed by the energy conservation analysis presented in Figure \ref{fig.TaylorBar2D-energy}, 
where all the initial kinetic energy $E_k$ is converted into shear internal energy $E_e$. 
Finally, the time evolution of the maximum length of the projectile is depicted, obtaining a final 
length of approximately $L=4.62$ which is in good agreement with the literature 
\cite{Maire_elasto,HyperHypo2019}. Mesh convergence is shown as well by comparing the results with 
$h=1/100$ and $h=1/200$.

\subsection{3D Taylor bar impact on a wall} \label{ssec.TaylorBar3D}
\begin{figure}[!htbp]
	\begin{center}
		\begin{tabular}{cc}
			\includegraphics[width=0.35\textwidth,draft=false]{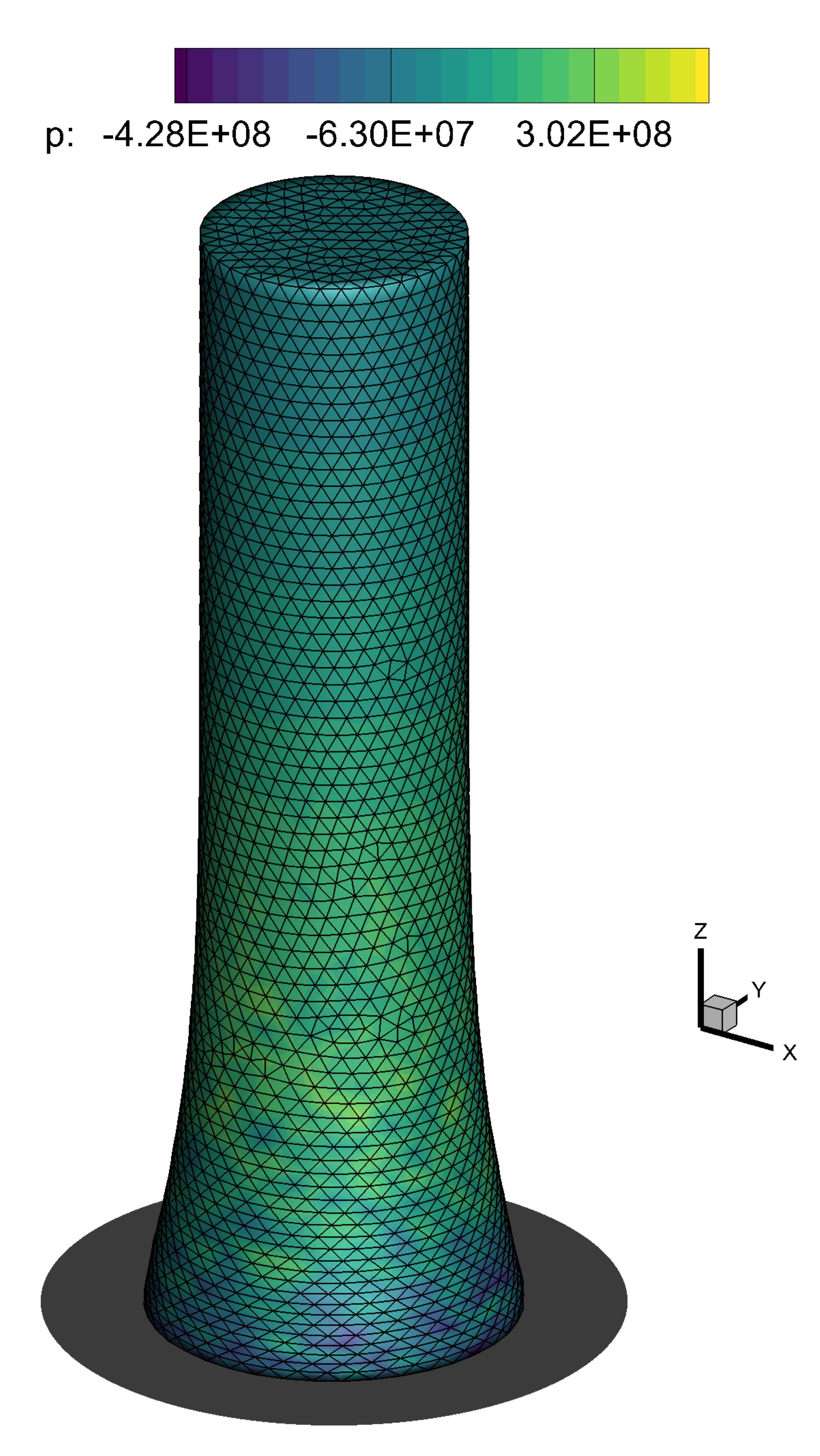}  &          
			\includegraphics[width=0.35\textwidth,draft=false]{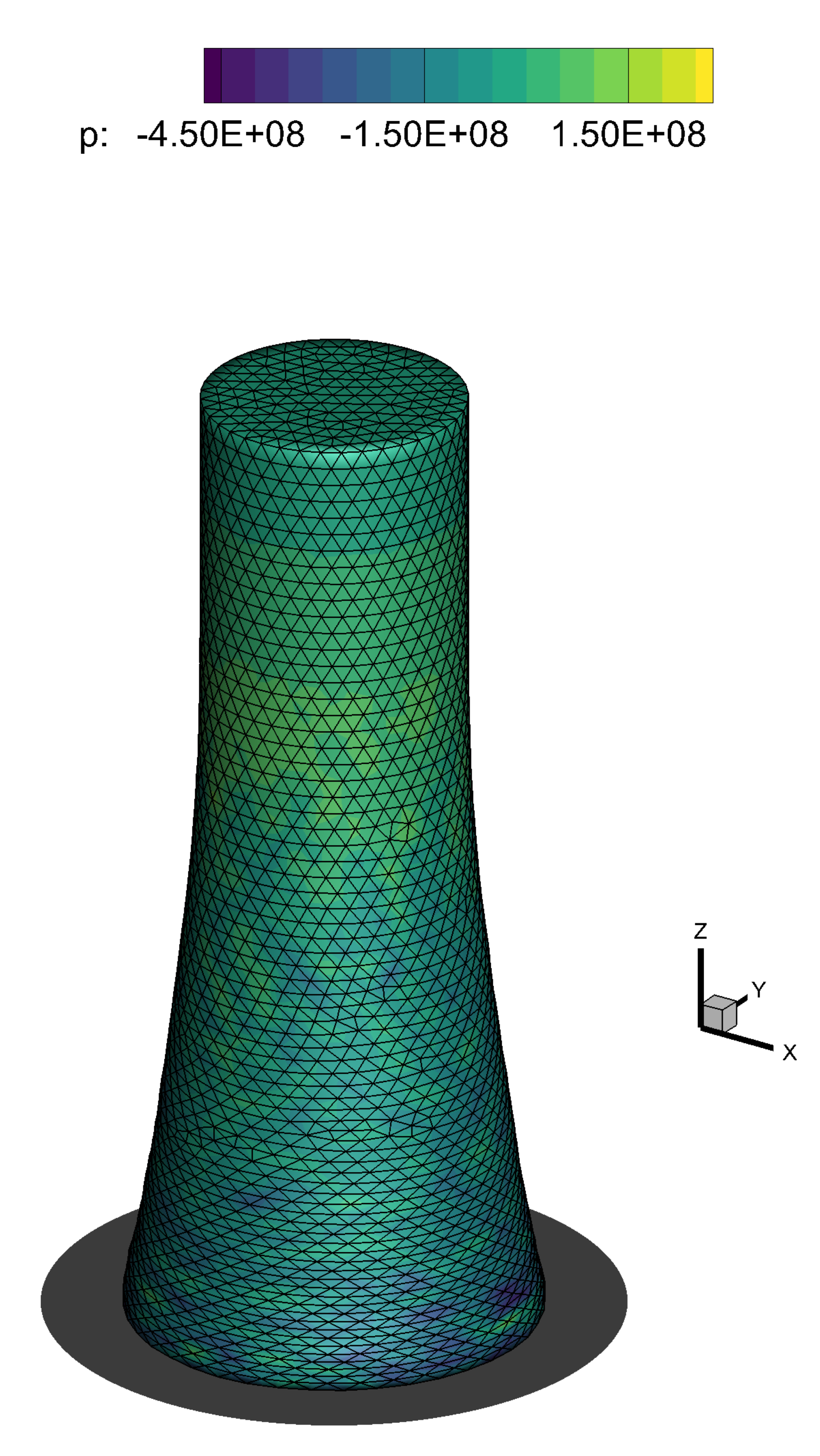} \\
			\includegraphics[width=0.35\textwidth,draft=false]{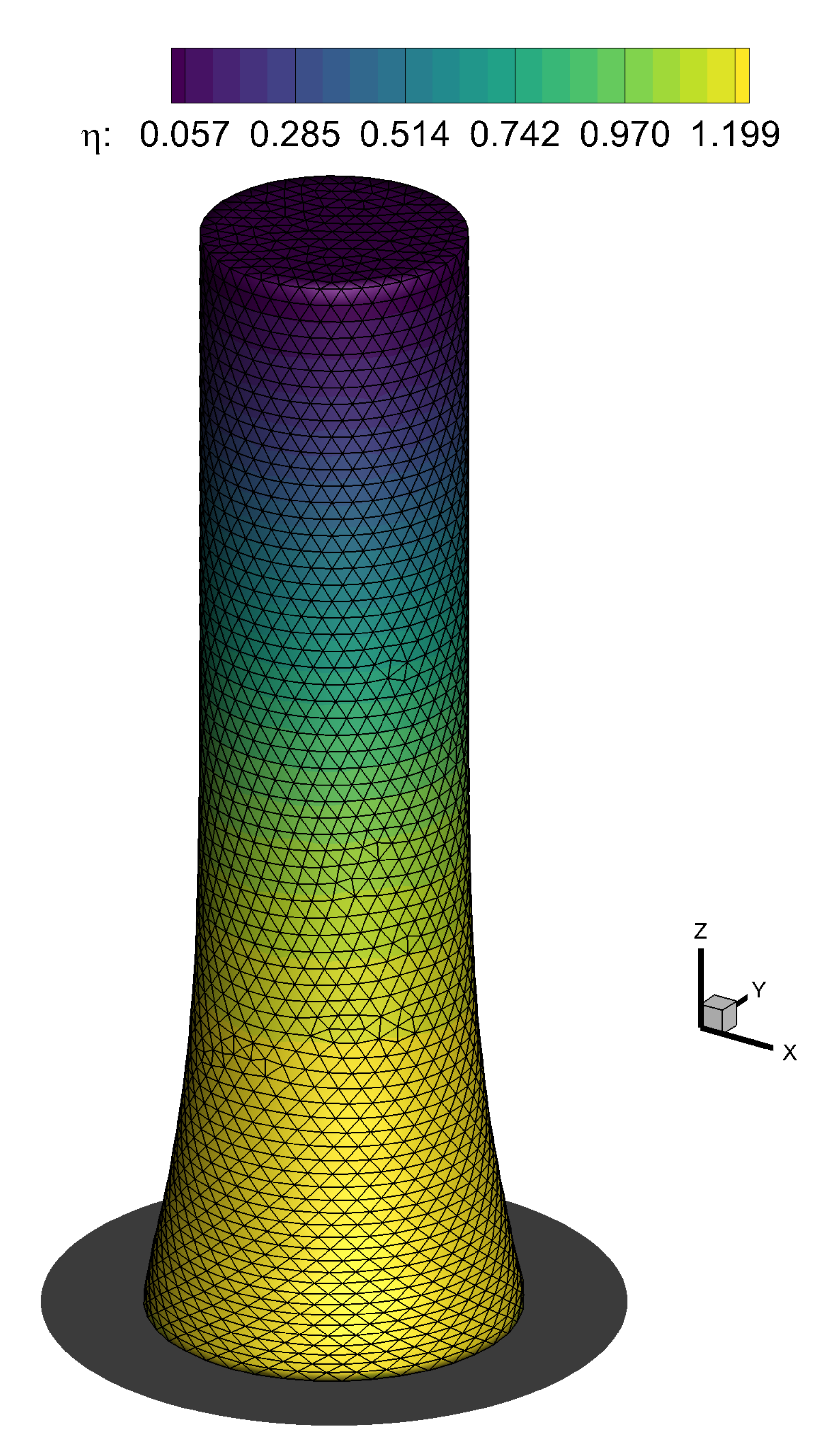}  &          
			\includegraphics[width=0.35\textwidth,draft=false]{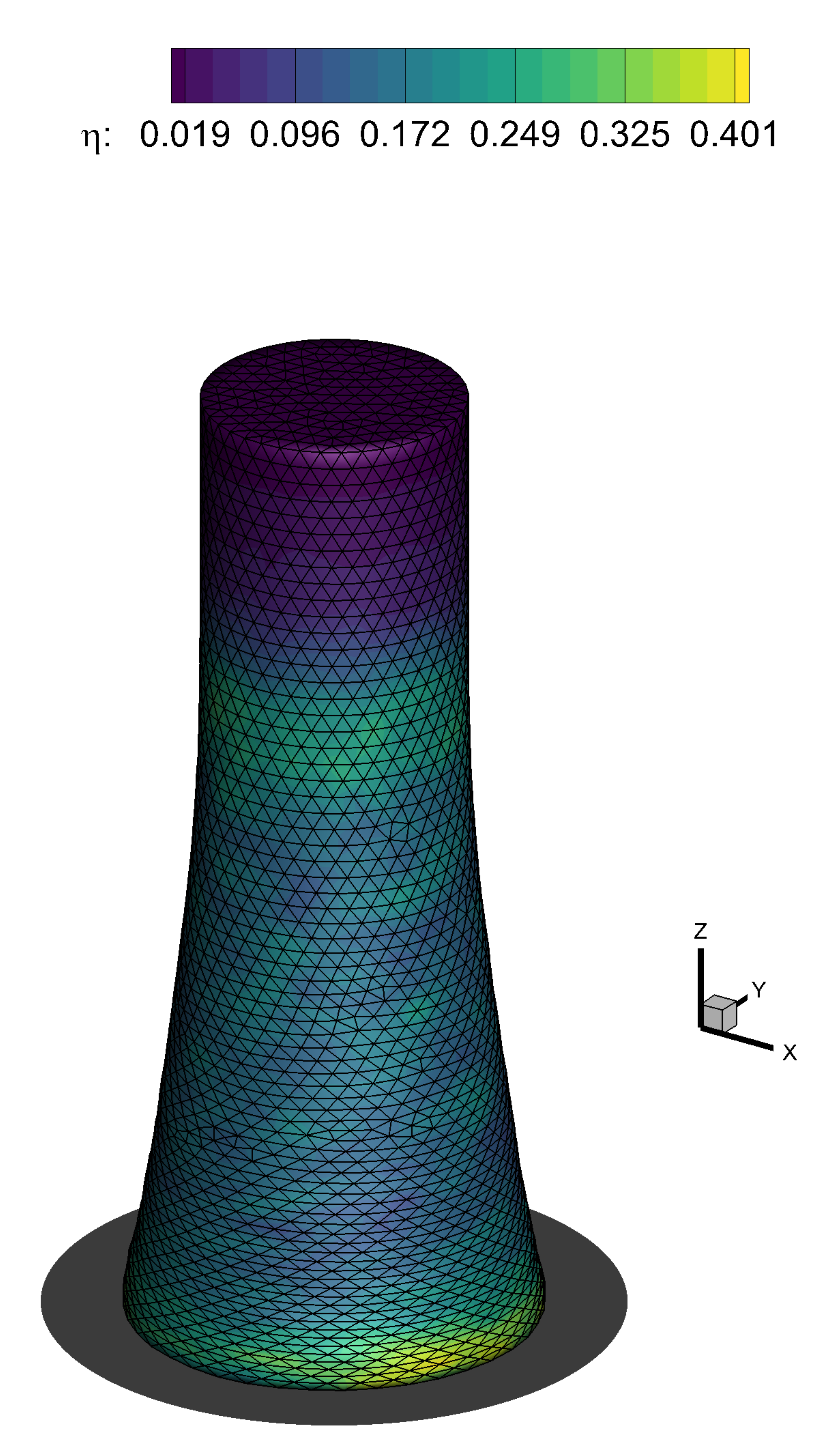} \\
		\end{tabular}
		\caption{3D Taylor bar impact on a wall. Top:  numerical distribution of pressure ($p$) at 
		time $t=2\cdot 10^{-5}$ (left) and $t=8\cdot 10^{-5}$ (right). Bottom: numerical 
		distribution of plastic map ($\eta$) at time $t=2\cdot 10^{-5}$ (left) and $t=8\cdot 
		10^{-5}$ (right).}
		\label{fig.TaylorBar3D}
	\end{center}
\end{figure}

\begin{figure}[!htbp]
	\begin{center}
		\begin{tabular}{cc}
			\includegraphics[width=0.47\textwidth,draft=false]{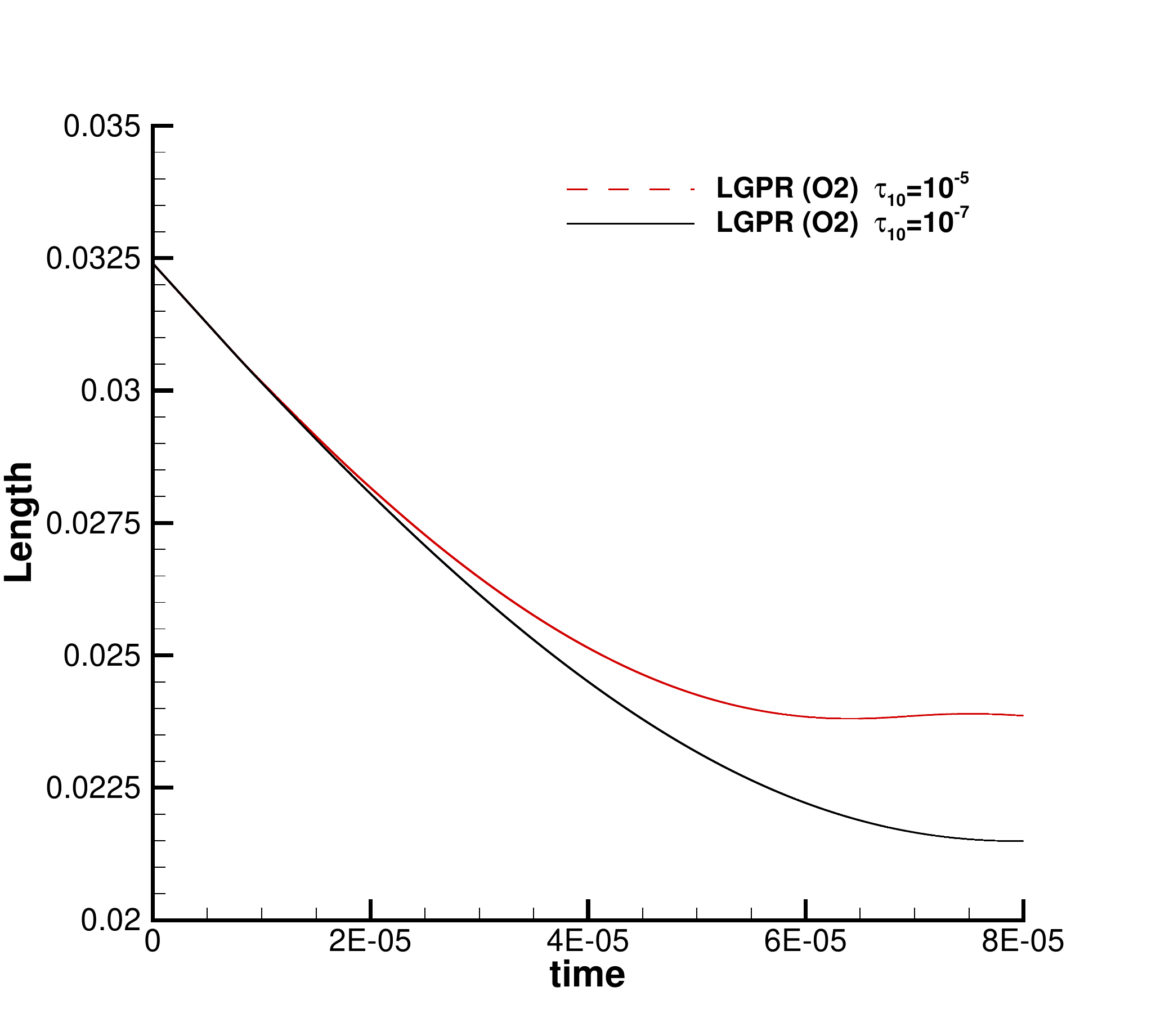}  &          
			\includegraphics[width=0.47\textwidth,draft=false]{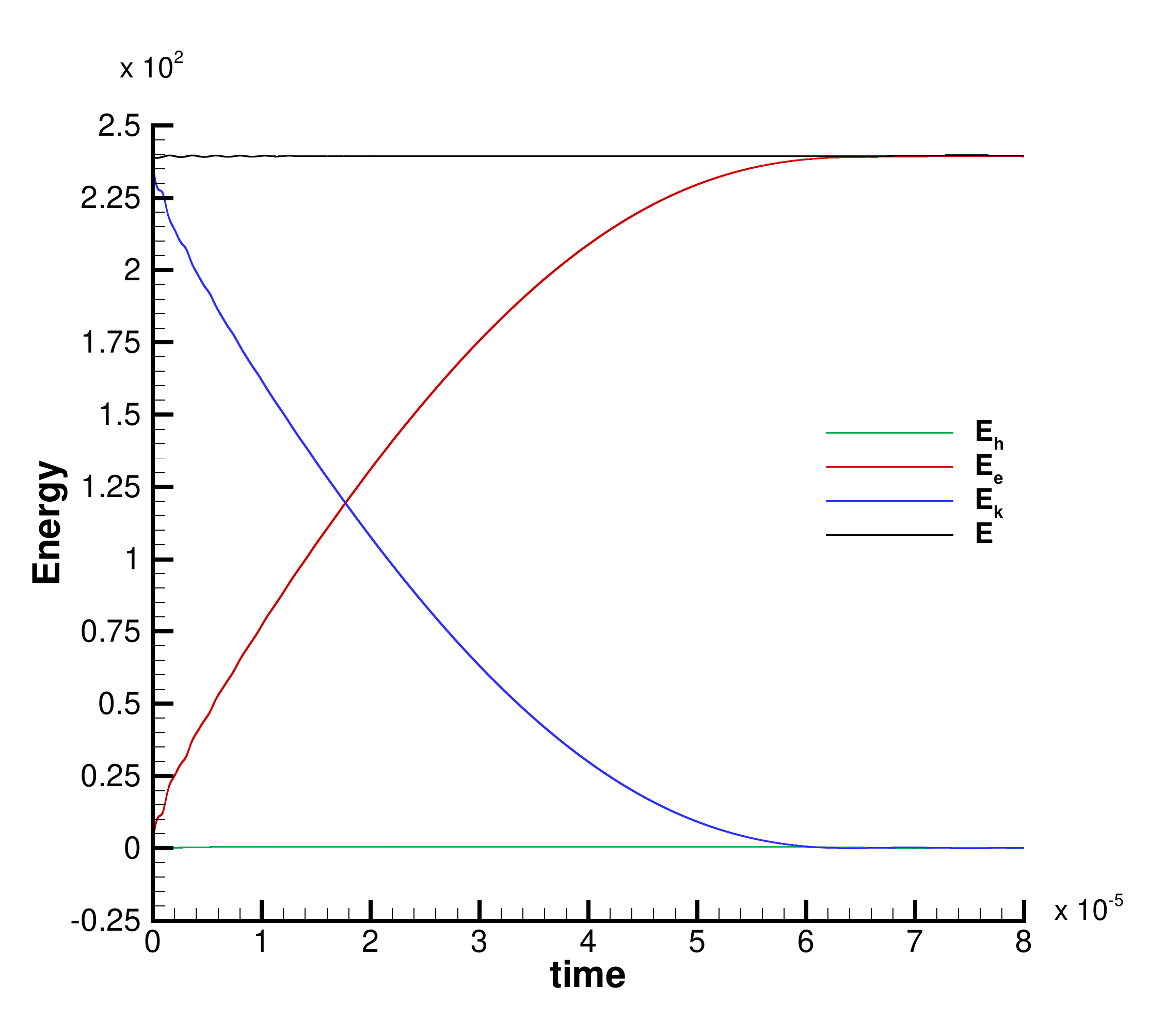} \\
		\end{tabular}
		\caption{3D Taylor bar impact on a wall. Left: time evolution of the bar length with 
		$\tau_{10}=10^{-7}$ (solid line) and $\tau_{10}=10^{-5}$ (dashed line). Right: analysis of energy 
		conservation in terms of volumetric and elastic energy ($E_h$ and $E_e$, respectively),  
		kinetic ($E_k$) and total ($E$) energy for $\tau_{10}=10^{-5}$.}
		\label{fig.TaylorBar3D-energy}
	\end{center}
\end{figure}

Next, we consider the fully three-dimensional simulation of a copper target that impacts a solid 
wall, according to the setup provided in \cite{Taylor}. A cylindrical rod made of copper has an 
initial length $L_0=0.0324$ and an initial radius $R_0=0.0032$. At $t=0$, it hits a rigid flat 
plate with velocity $\vv(0,\x)=(0,0,-227)$ and pressure $p=0$. The computational mesh is composed 
of $N_E=16464$ tetrahedra with characteristic mesh size $h=1/50$. Free-traction boundary conditions 
are set everywhere, apart on the impact surface where wall boundaries are imposed. The simulation is carried on until the final time $t_f=80 \cdot 10^{-6}$ and the 
Mie-Grüneisen EOS is used for copper with $\Gamma_0=2$, $s=1.48$ and yield strength 
$\sigma_Y=400\cdot 10^6$. To take into account plastic deformations, the relaxation time $\tau_1$ 
is dynamically computed with \eqref{eqn.tau.plast} using $n=10$, while we study the different 
behavior of the material by setting $\tau_{10}=10^{-7}$ and $\tau_{10}=10^{-5}$. The plastic map 
$\eta$ and the pressure distribution are shown in Figure \ref{fig.TaylorBar3D} for 
$ 
\tau_{10}=10^{-5} $ at time $t=2\cdot 
10^{-5}$ and $t=t_f$, highlighting that most of the plasticity effects take place close to the 
wall at the initial instants of the impact. Figure \ref{fig.TaylorBar3D-energy} presents the study of energy conservation over time as 
well as the time evolution of the length of the copper rod. As $\tau_{10}$ increases, the material 
behavior gets closer to elastic solids, while a more realistic setting is recovered using 
$\tau_{10}=10^{-7}$, which correctly accounts for plastic deformations. 

\subsection{Elastic vibration of a beryllium plate} \label{ssec.BePlate}
\begin{figure}[!htbp]
	\begin{center}
		\begin{tabular}{cc}
			\includegraphics[width=0.47\textwidth,draft=false]{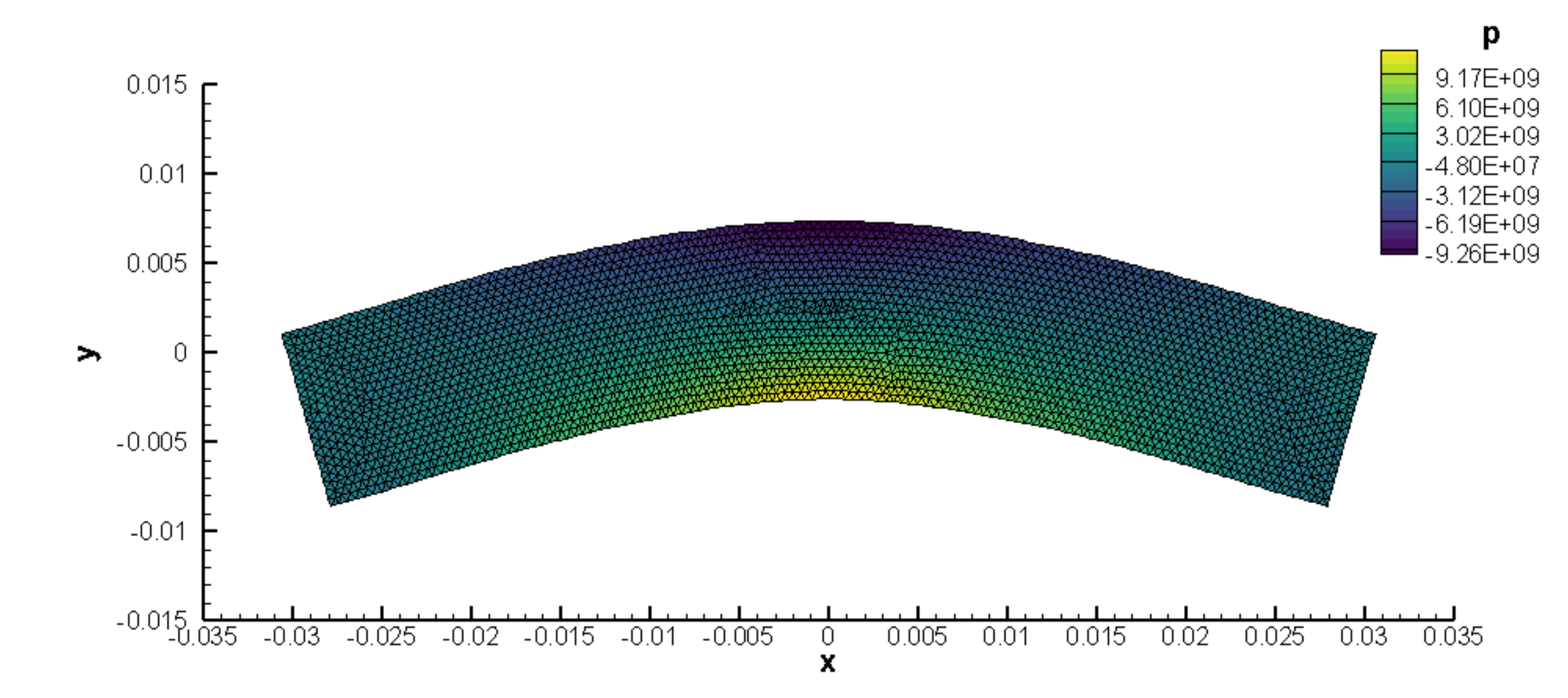}  &          
			\includegraphics[width=0.47\textwidth,draft=false]{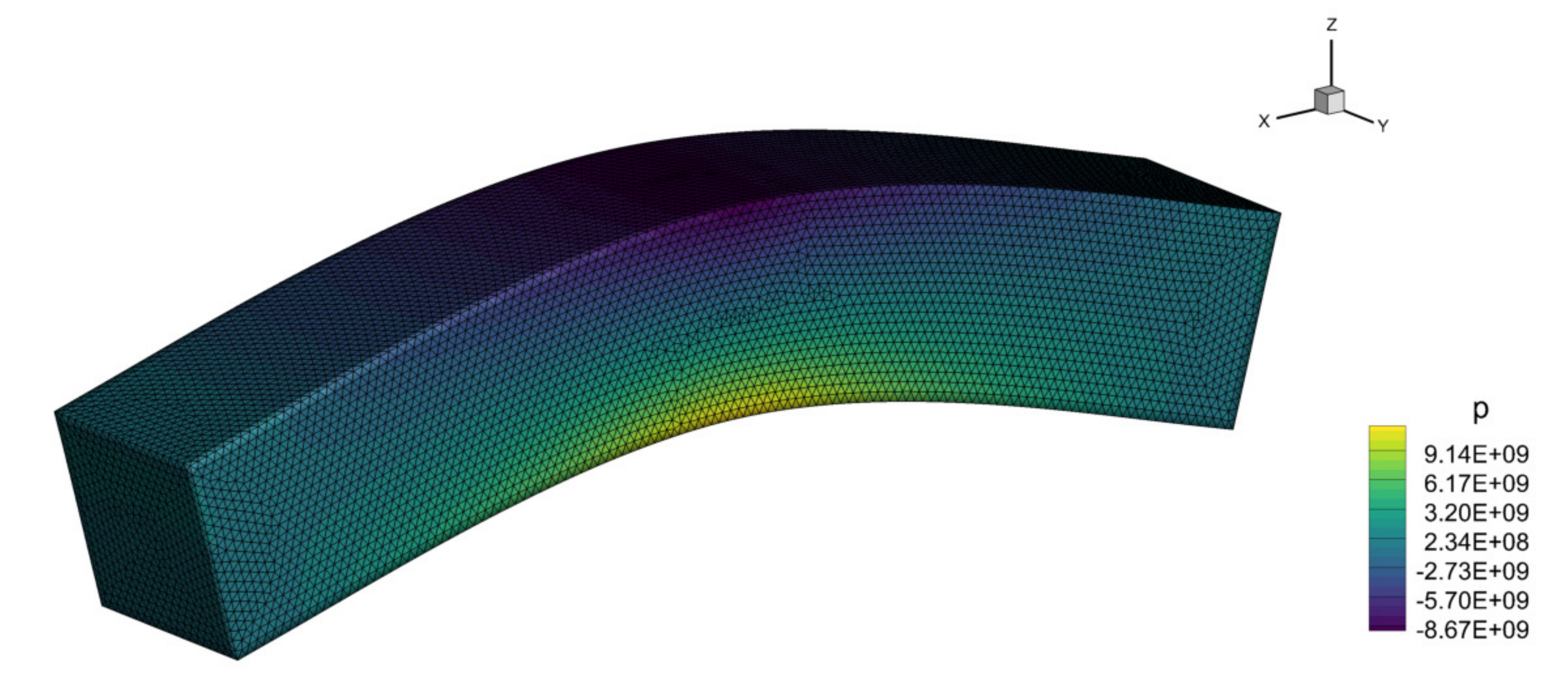} \\
			\includegraphics[width=0.47\textwidth,draft=false]{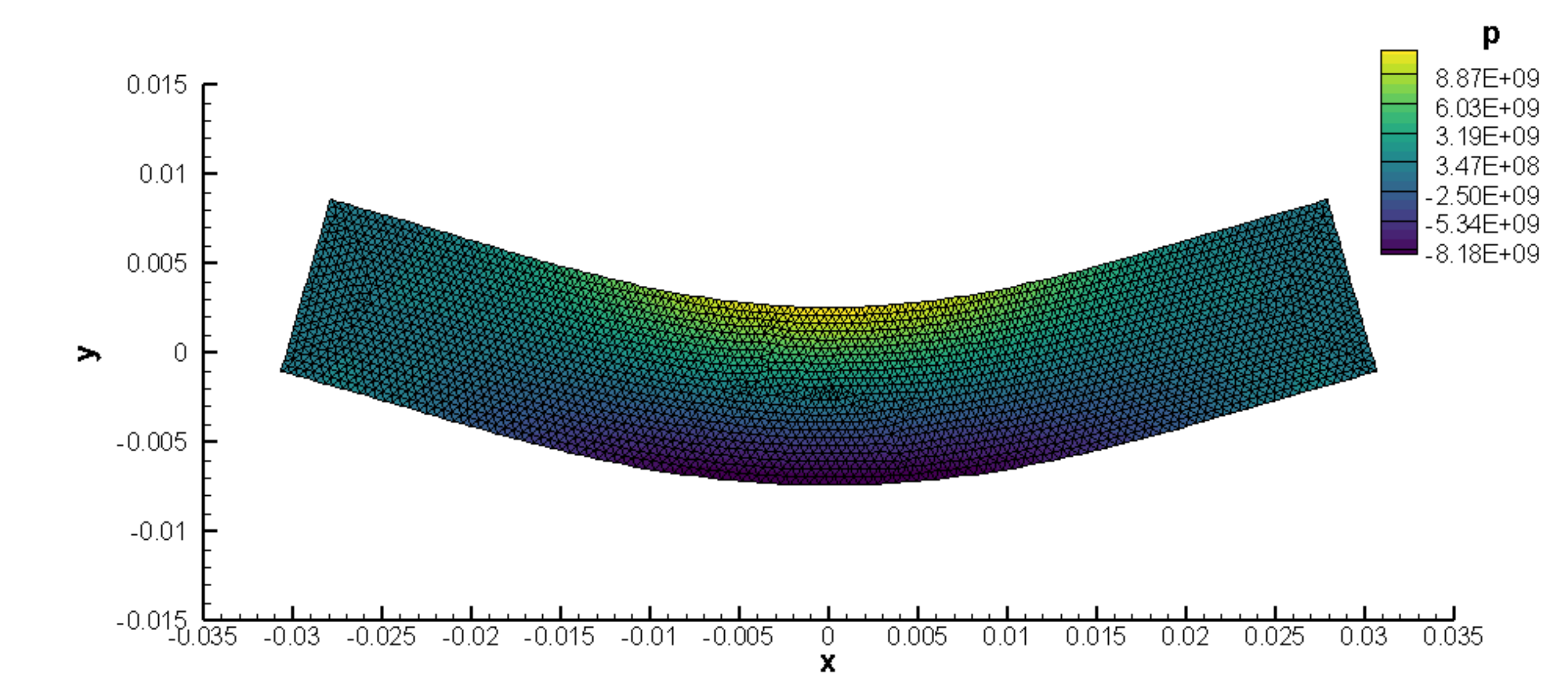}  &          
			\includegraphics[width=0.47\textwidth,draft=false]{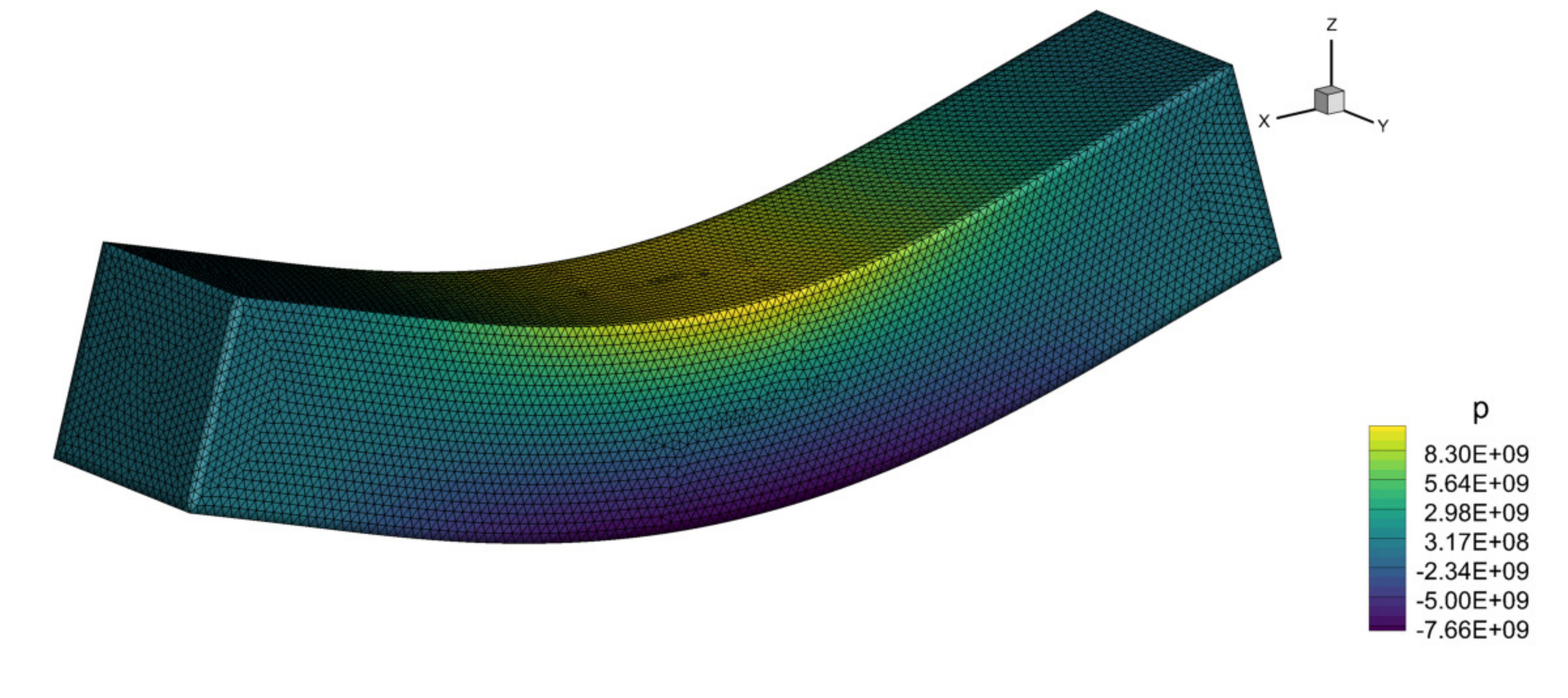} \\
			\includegraphics[width=0.47\textwidth,draft=false]{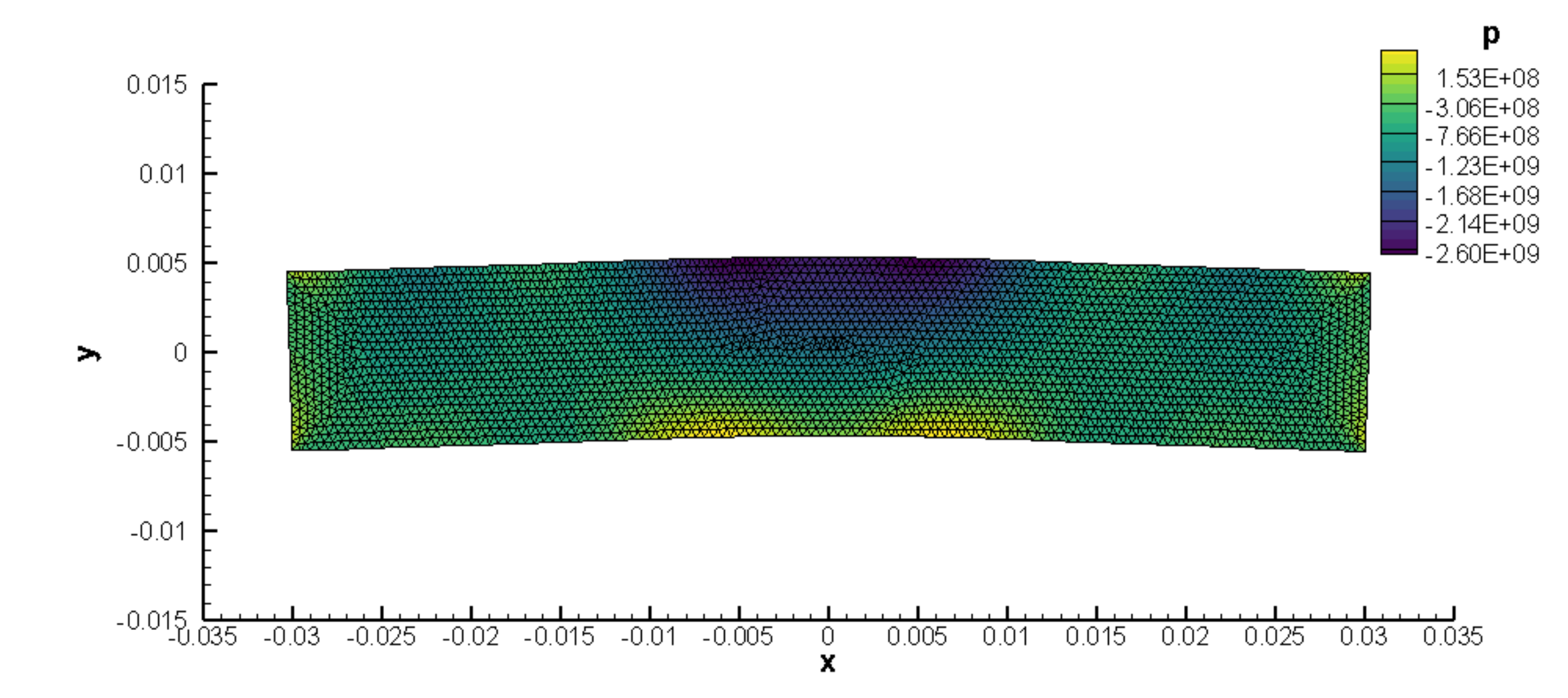}  &          
			\includegraphics[width=0.47\textwidth,draft=false]{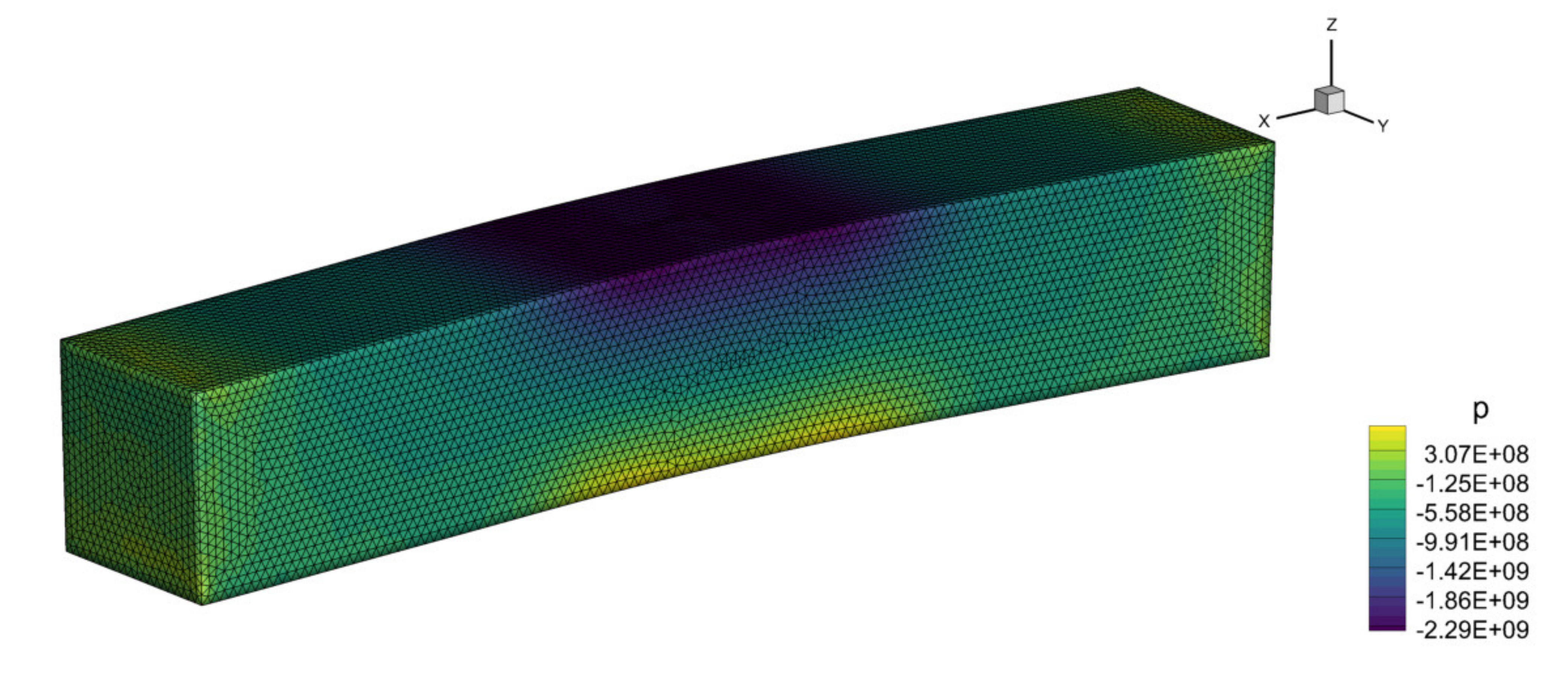} \\
		\end{tabular}
		\caption{Elastic vibration of a Beryllium plate. Numerical results at output times $t = 10^{-5}$ (top), $t = 2 \cdot 10^{-5}$ (middle) and $t = 3 \cdot 10^{-5}$ (bottom) in 2D (left) and 3D (right). }
		\label{fig.BePlate}
	\end{center}
\end{figure}

\begin{figure}[!htbp]
	\begin{center}
		\begin{tabular}{cc}
		\includegraphics[width=0.47\textwidth,draft=false]{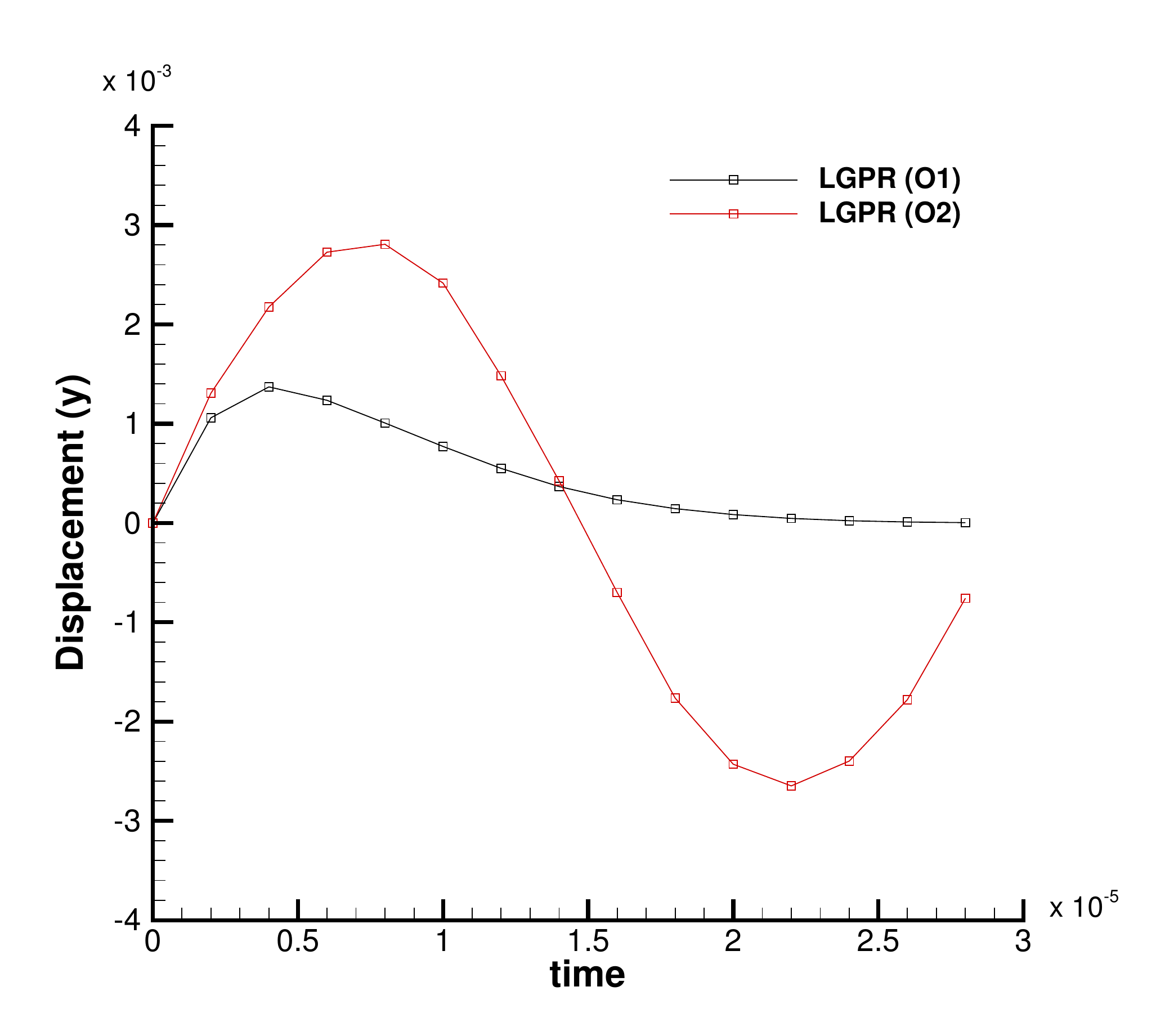}  &          
		\includegraphics[width=0.47\textwidth,draft=false]{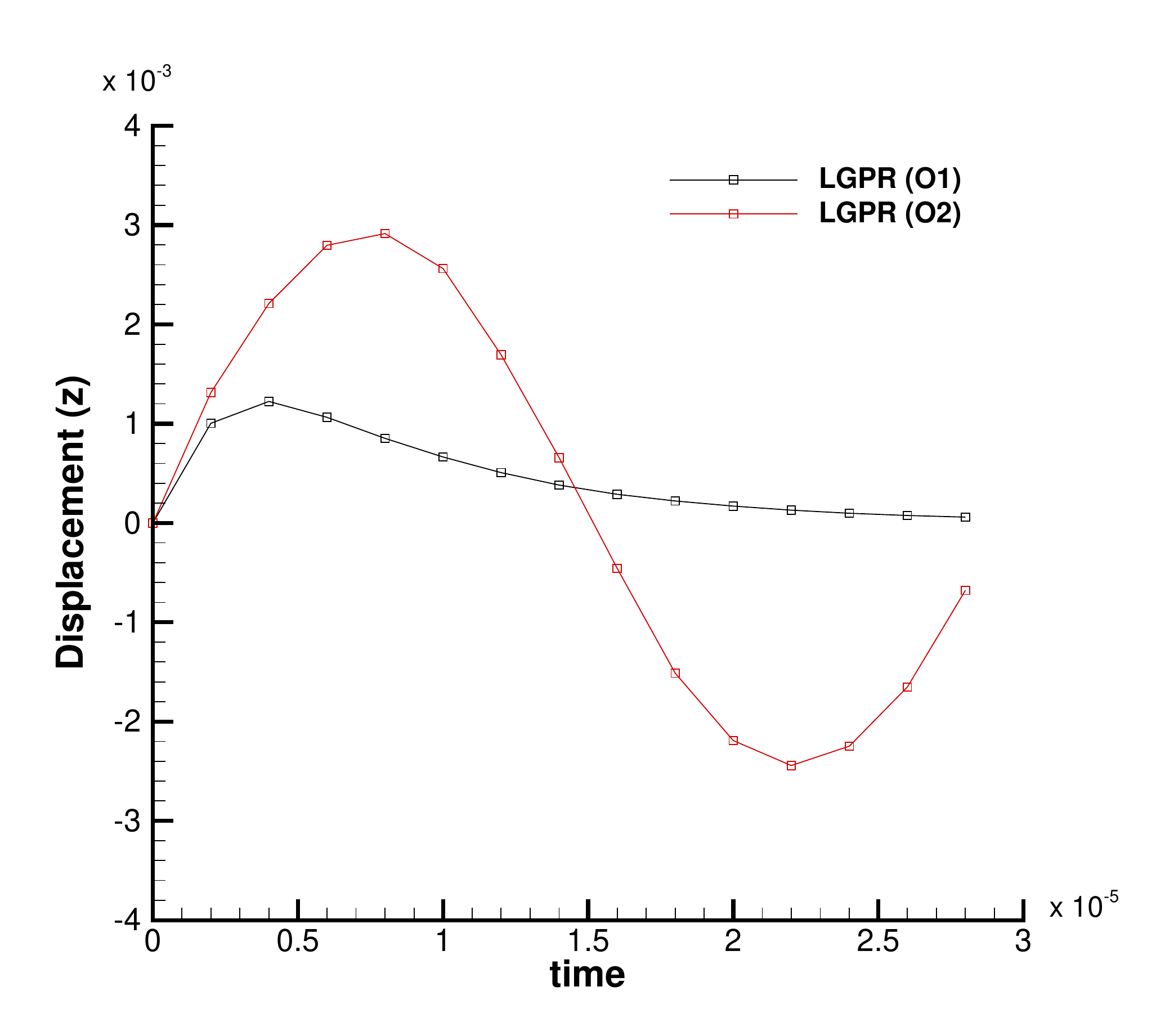} \\
    	\end{tabular}
		\caption{Elastic vibration of a Beryllium plate. Vertical displacement of the point initially located at $\x = (0,0,0)$ with first and second order LGPR scheme. }
			\label{fig.BePlate_disp}
	\end{center}
\end{figure}

This test case describes the elastic vibration of a beryllium plate or bar, see 
\cite{HyperHypo2019,CCL2020} for instance. Here we consider both the 2D and the 3D version. The 
computational domain is $\Omega^{2D}(0)=[-0.03;0.03]\times[-0.005;0.005]$ and 
$\Omega^{3D}(0)=[-0.03;0.03]\times[-0.005;0.005]\times[-0.005;0.005]$, thus the length of the bar 
is $L=0.06$. A characteristic mesh size of $h=1/200$ is used in 2D, while $h=1/100$ is adopted in 
3D. Free-traction boundary conditions are set everywhere and the final time of the simulation is 
chosen to be $t_f=3\cdot 10^{-5}$, so that approximately one oscillating period is completed. The 
hydrodynamics part of the energy is given by the Neo-Hookean EOS \eqref{eqn.NH}, 
with Young modulus 
$Y=3.1827 \cdot 10^{11}$ and Poisson ratio $\nu=0.0539$. The material is initially loaded via a 
perturbed initial velocity field $\vv^{2D}(0,\x)=(0,V_0(x),0)$ and $\vv^{3D}(0,\x)=(0,0,V_0(x))$ of 
the form
\bea
V_0(x) = A \omega \left[  a_1(\sinh(x')+\sin(x')) - a_2(\cosh(x')+\cos(x')) \right],
\eea
where $x'=\alpha(x+L/2)$, $\alpha=78.834$, $A=4.3369\times 10^{-5}$,
$\omega=2.3597\times 10^5$, $a_1=56.6368$ and $a_2=57.6455$. In this setting, the LGPR scheme \textit{exactly} collapses to the cell-centered finite volume scheme recently introduced in \cite{Boscheri2021}. Figure \ref{fig.BePlate} depicts the mesh configuration and the pressure distribution at three different output times during one flexural period, whereas in Figure \ref{fig.BePlate_disp} we show the time evolution of the vertical component of the velocity of the barycenter of the bar, i.e. the mesh point originally located at $\x_0=(0,0,0)$. The first and second order schemes are compared, demonstrating that second order accuracy in space and time is responsible of a remarkable reduction of numerical dissipation. The bar dissipates almost all initial kinetic energy for the first order scheme, which is not due to plastic deformations, since the material is a purely elastic solid, but only because of an excessive numerical dissipation in the nodal solver and the face-based fluxes. This demonstrates the advantages induced by a higher order scheme in space and time.
	
\subsection{Twisting column} \label{ssec.twcol}
\begin{figure}[!htbp]
	\begin{center}
		\begin{tabular}{cccc}
			\includegraphics[width=0.23\textwidth,draft=false]{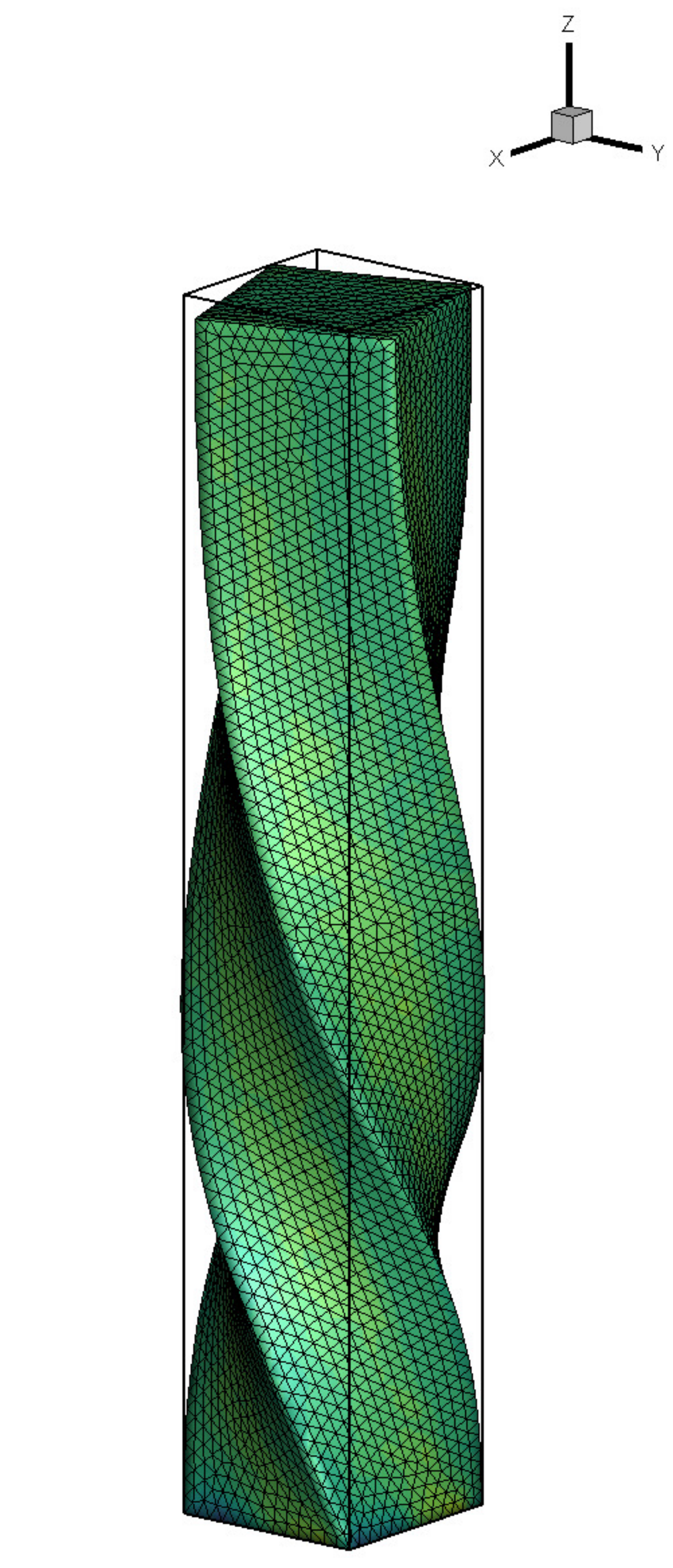} &          
			\includegraphics[width=0.23\textwidth,draft=false]{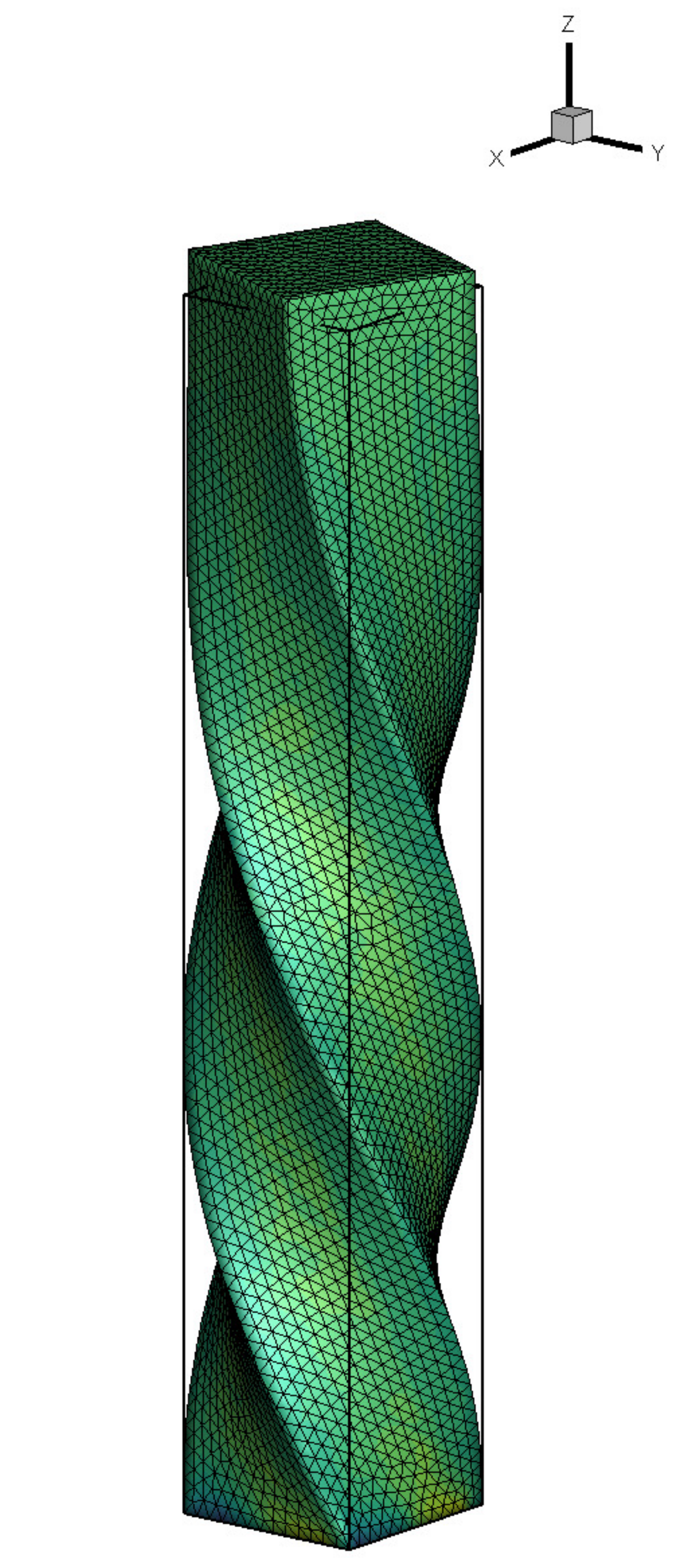} &
			\includegraphics[width=0.23\textwidth,draft=false]{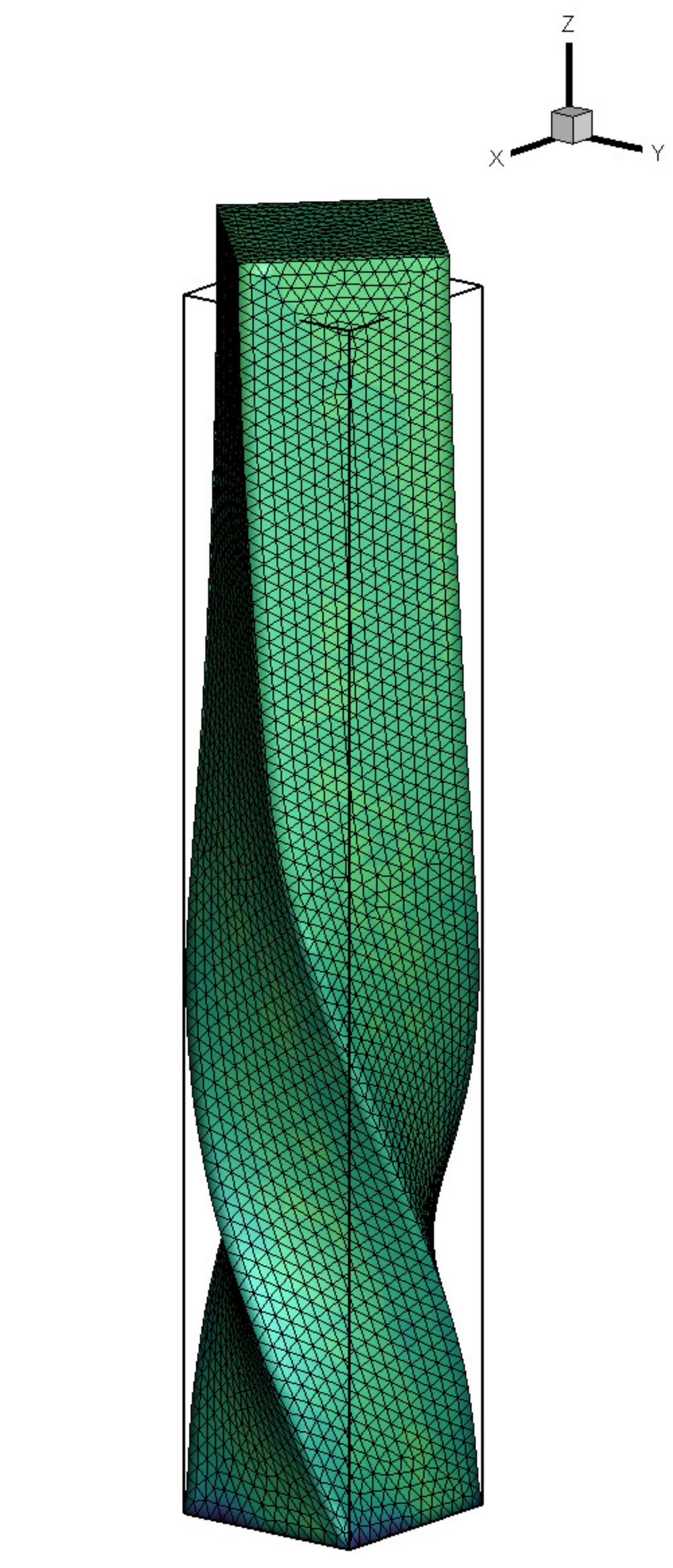} &
			\includegraphics[width=0.23\textwidth,draft=false]{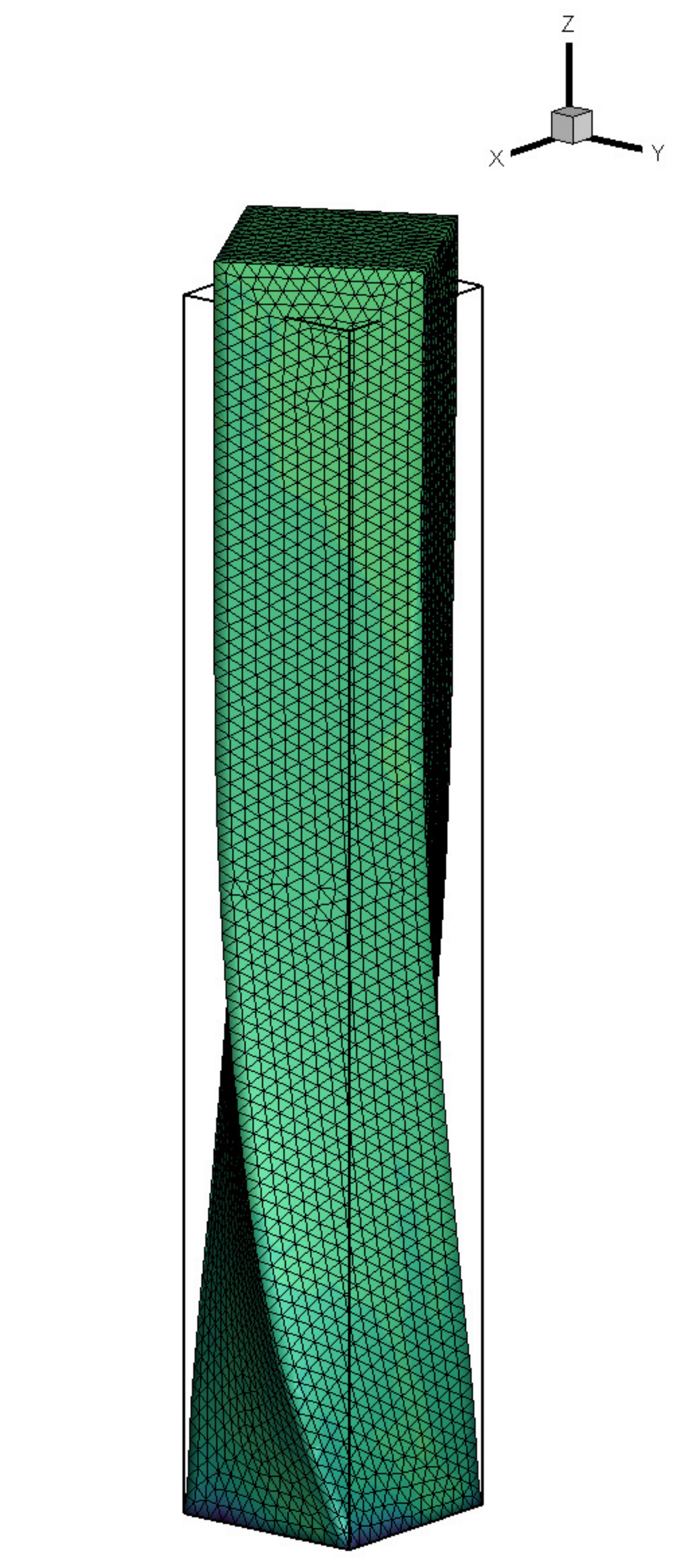} \\
			\includegraphics[width=0.23\textwidth,draft=false]{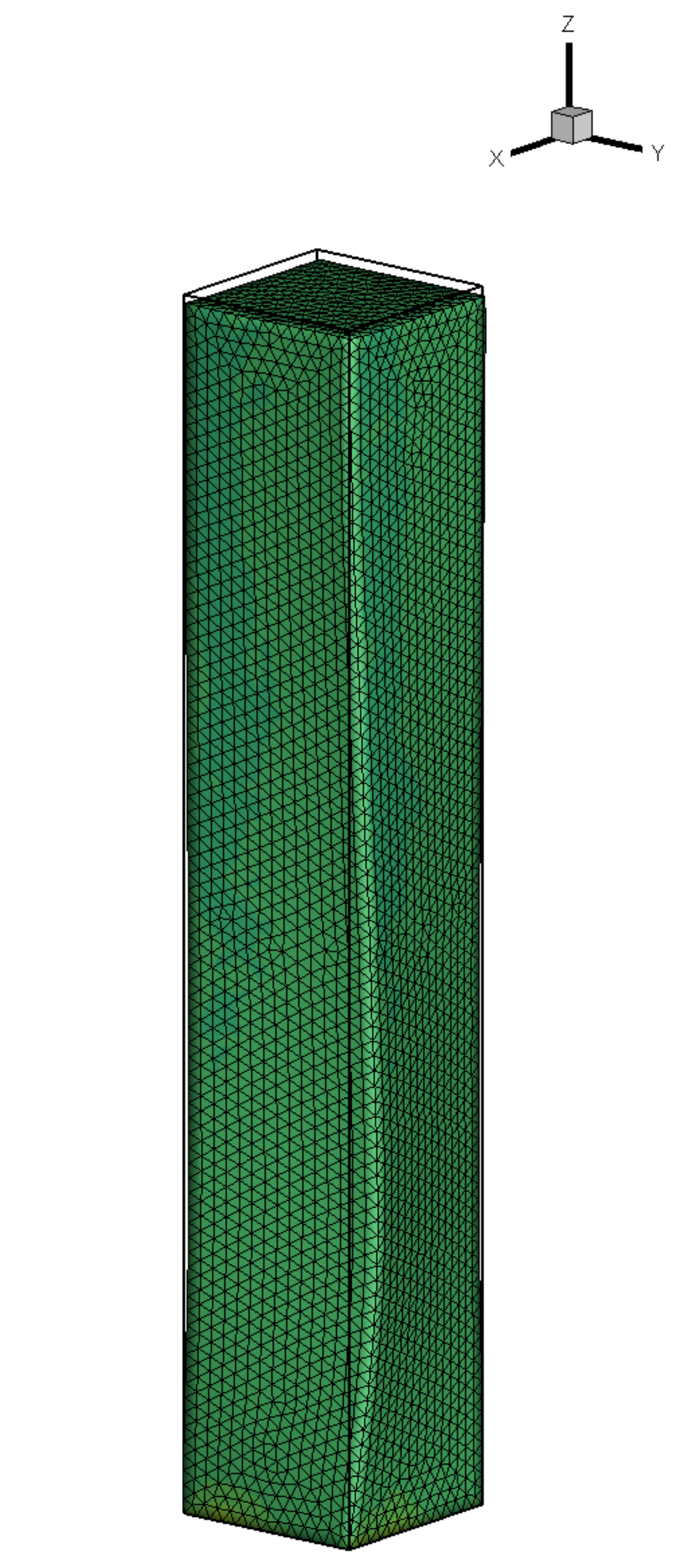} &          
			\includegraphics[width=0.23\textwidth,draft=false]{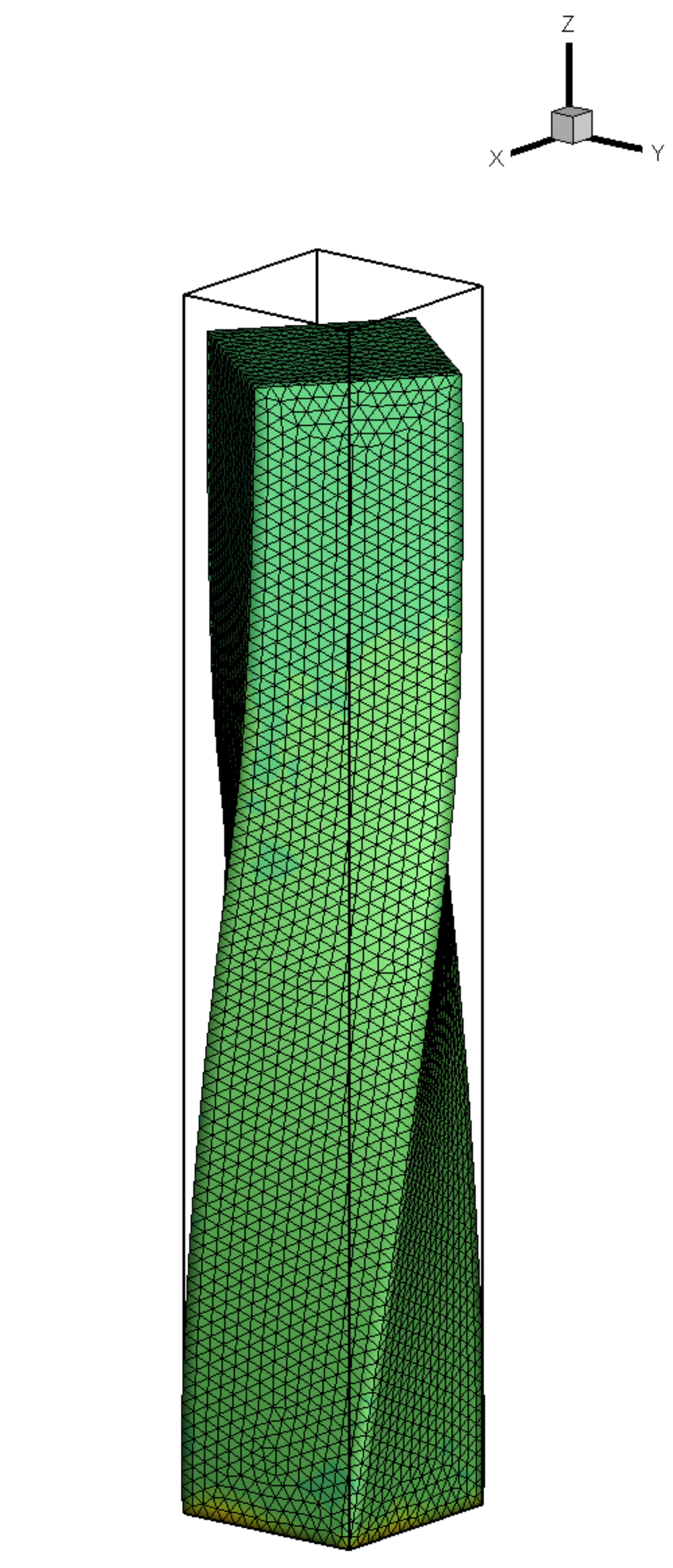} &
			\includegraphics[width=0.23\textwidth,draft=false]{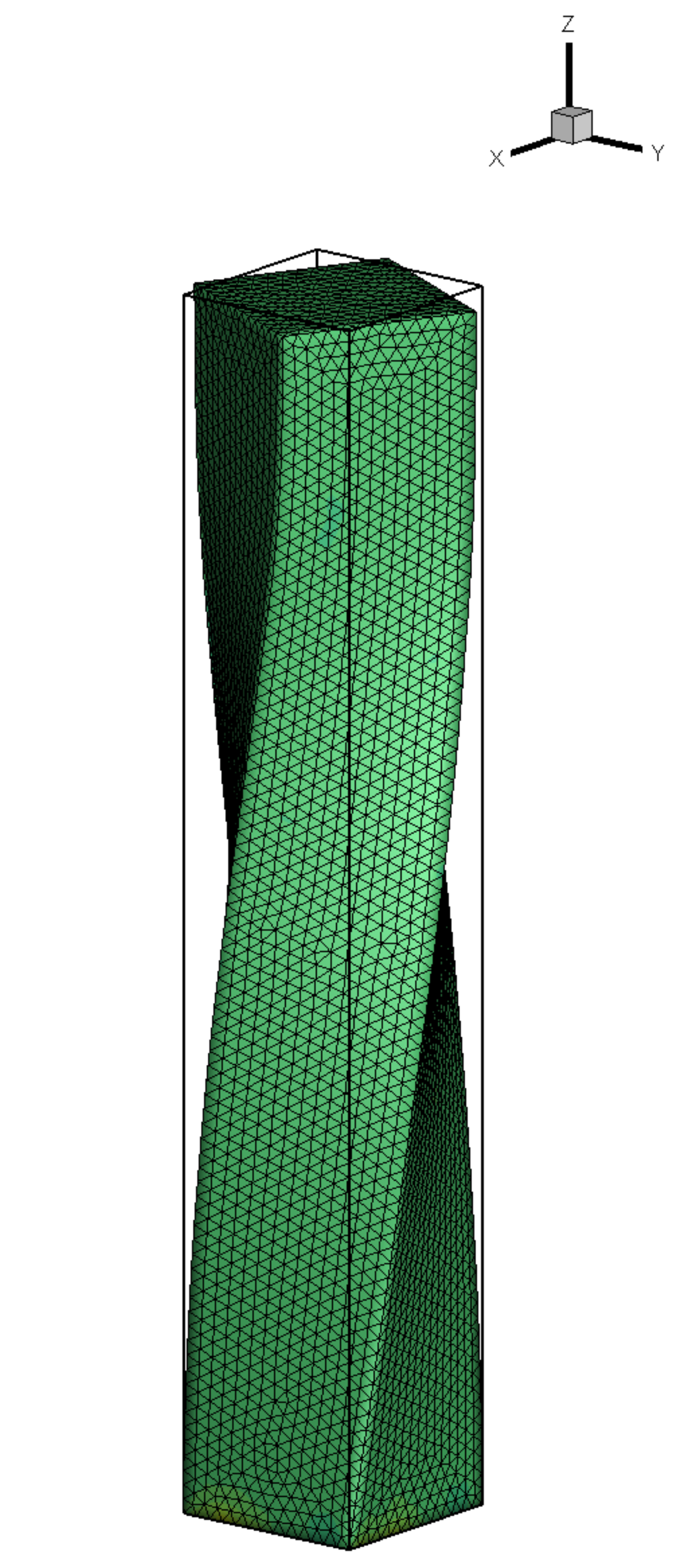} &
			\includegraphics[width=0.23\textwidth,draft=false]{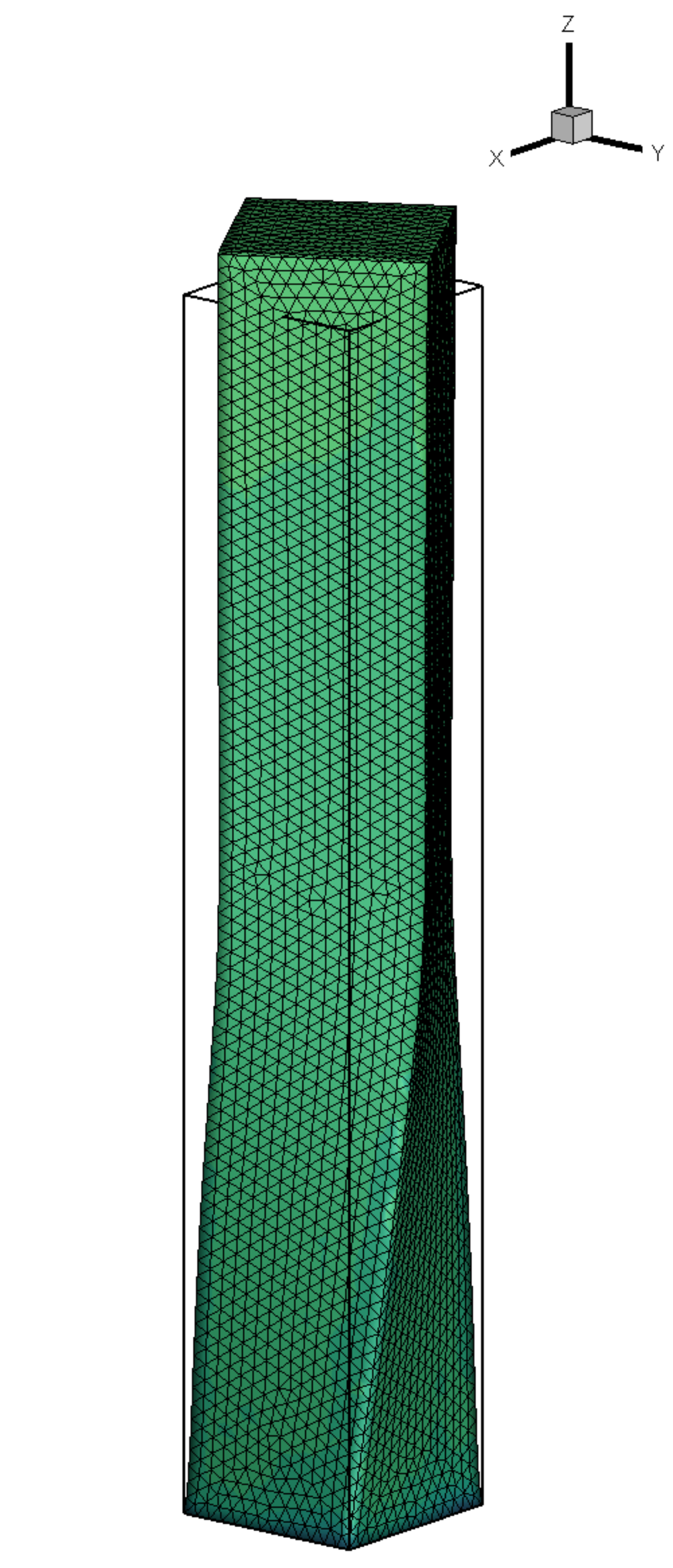} \\
			\multicolumn{4}{c}{\includegraphics[width=0.55\textwidth,draft=false]{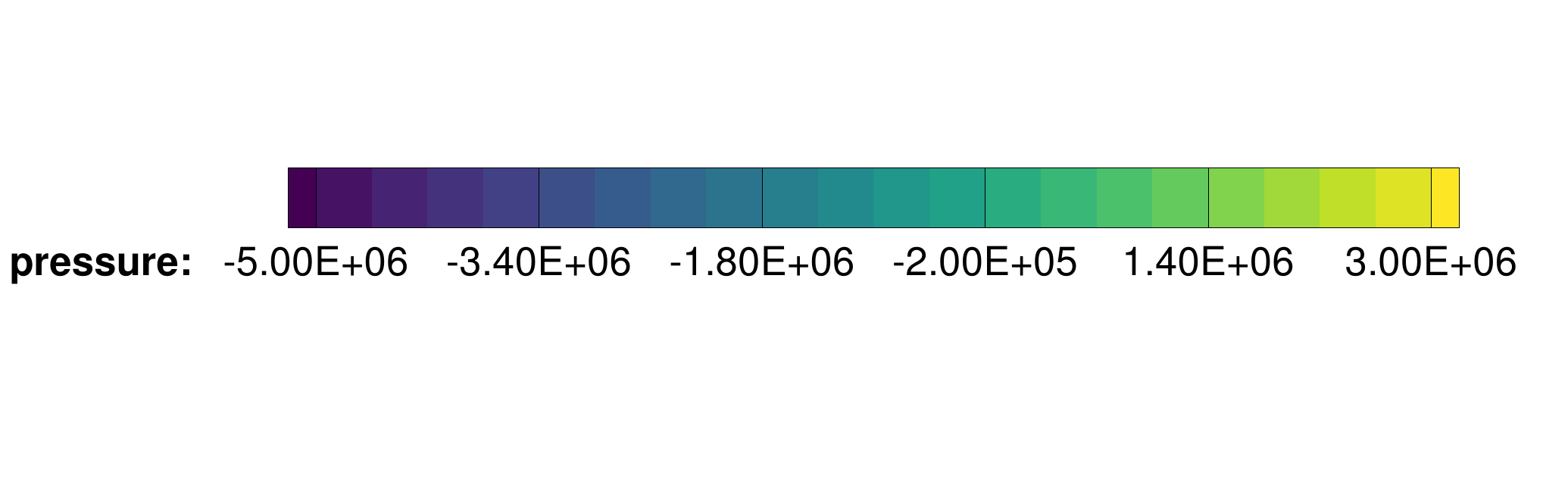}}
		\end{tabular}
		\caption{Twisting column with $\omega_0=100$. Column shape and pressure distribution at output times  $t = 0.00375$, $t = 0.075$, $t = 0.1125$, $t = 0.15$, $t = 0.1875$, $t = 0.225$, $t = 0.2625$ and $t = 0.3$ (from top left to bottom right). The shape is compared with respect to the initial configuration (hollow box). }
		\label{fig.TwistCol_u100}
	\end{center}
\end{figure}

\begin{figure}[!htbp]
	\begin{center}
		\begin{tabular}{cccc}
			\includegraphics[width=0.23\textwidth,draft=false]{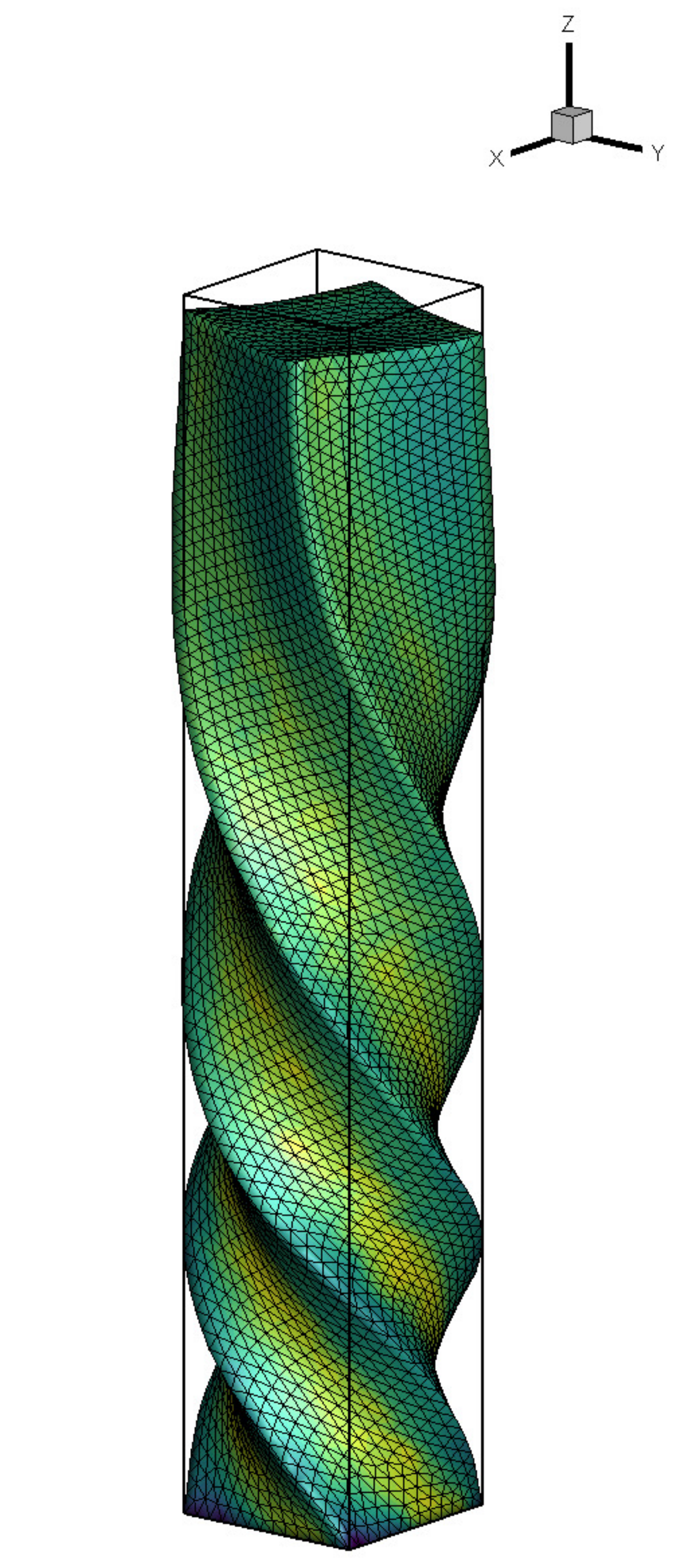} &          
			\includegraphics[width=0.23\textwidth,draft=false]{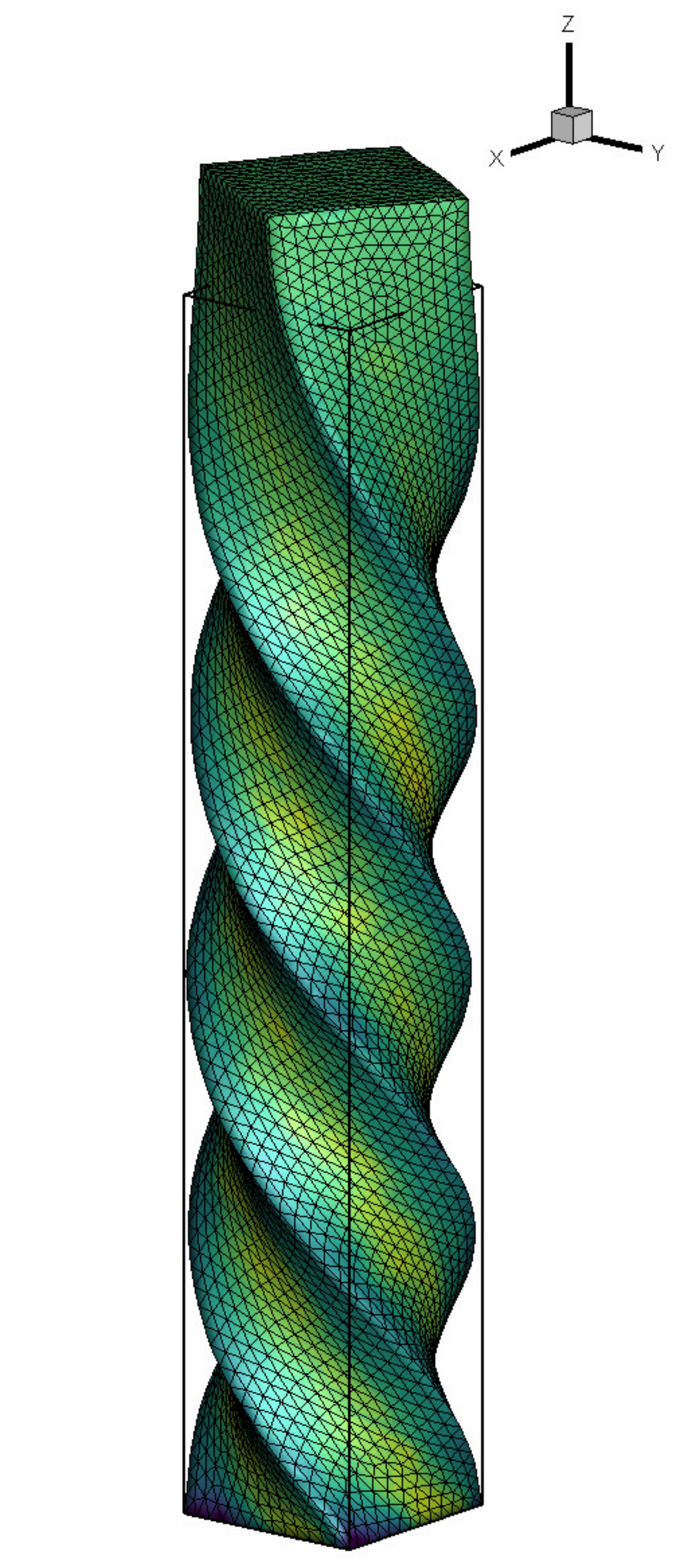} &
			\includegraphics[width=0.23\textwidth,draft=false]{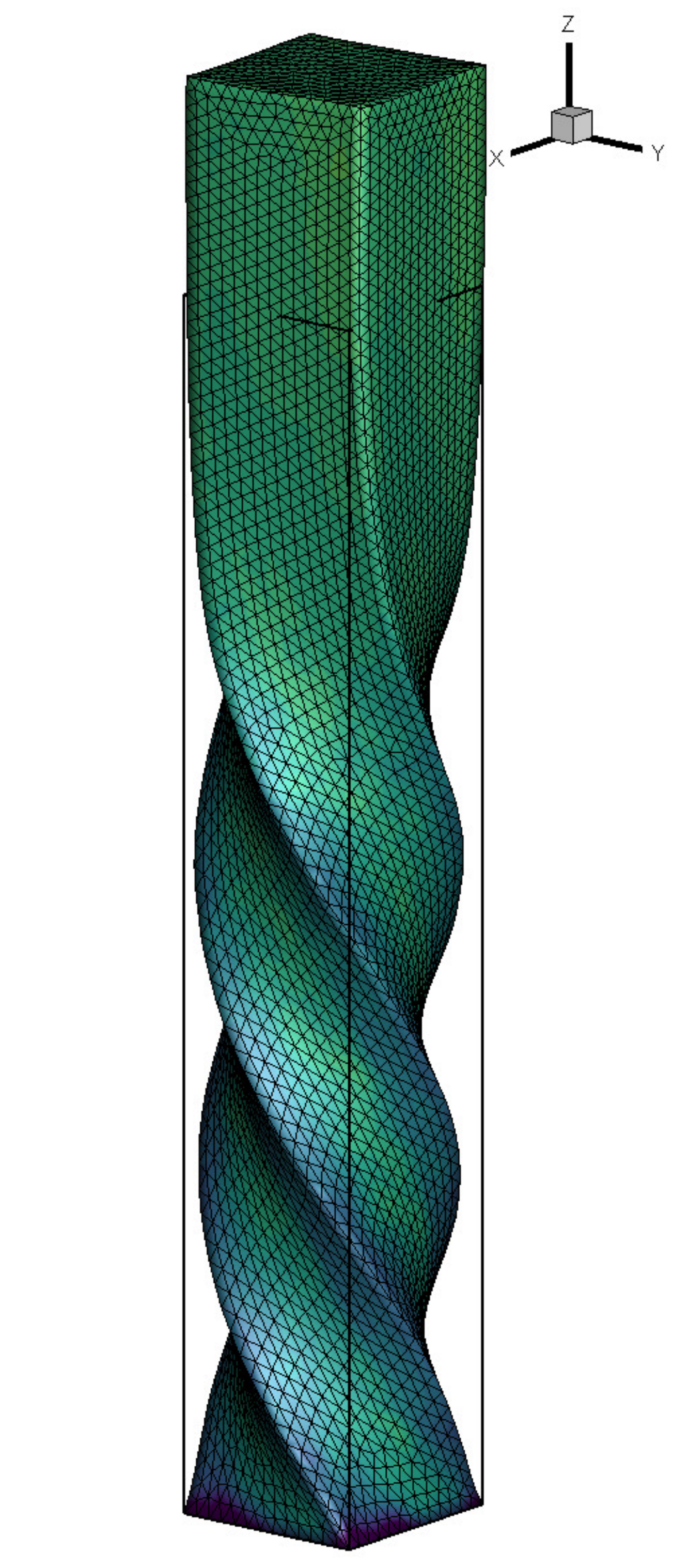} &
			\includegraphics[width=0.23\textwidth,draft=false]{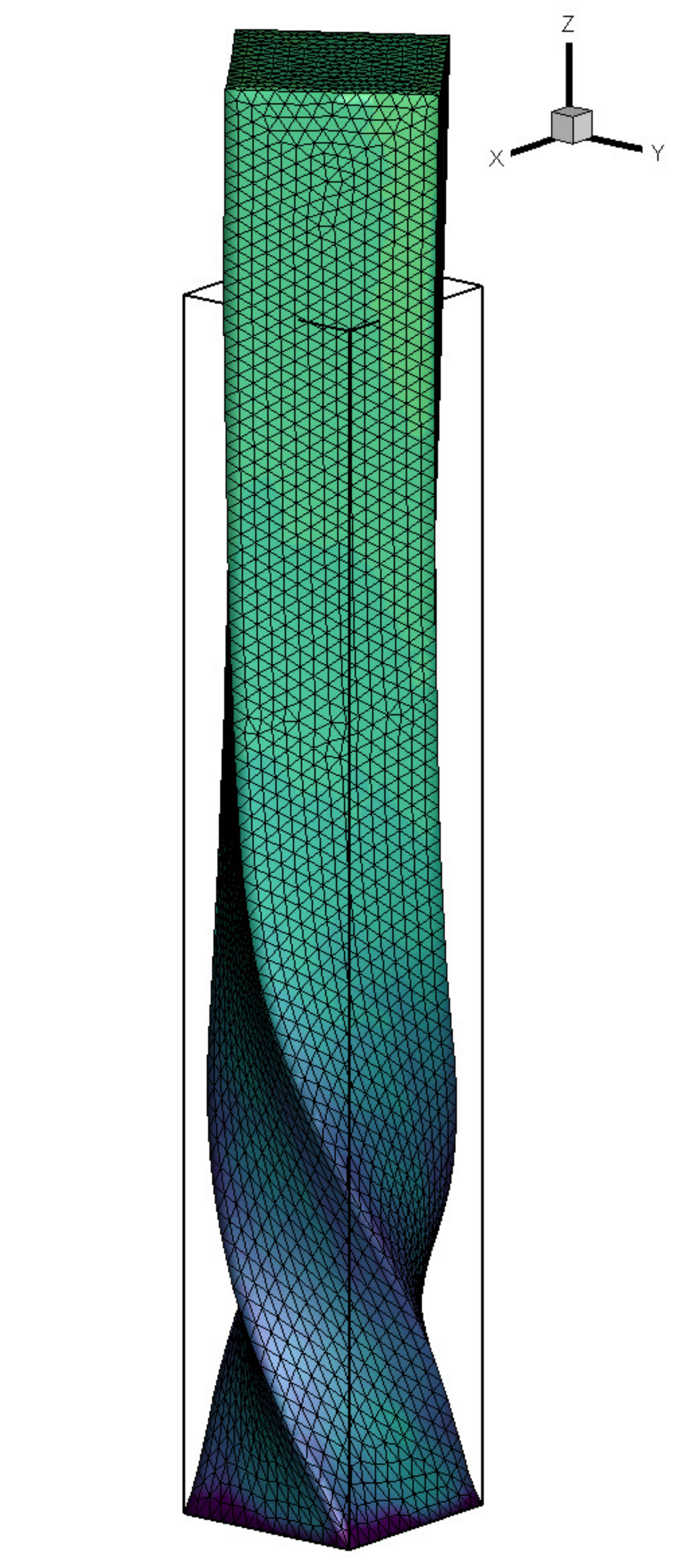} \\
			\includegraphics[width=0.23\textwidth,draft=false]{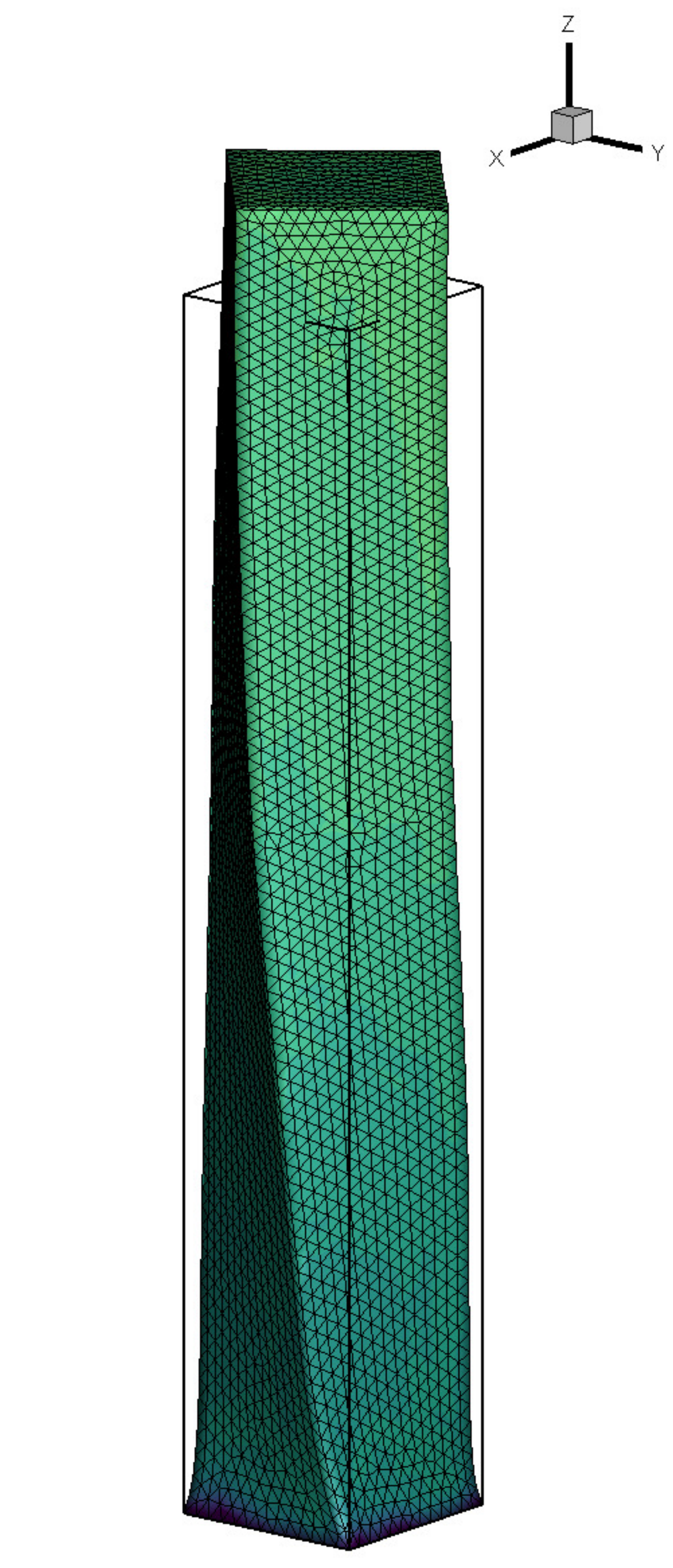} &          
			\includegraphics[width=0.23\textwidth,draft=false]{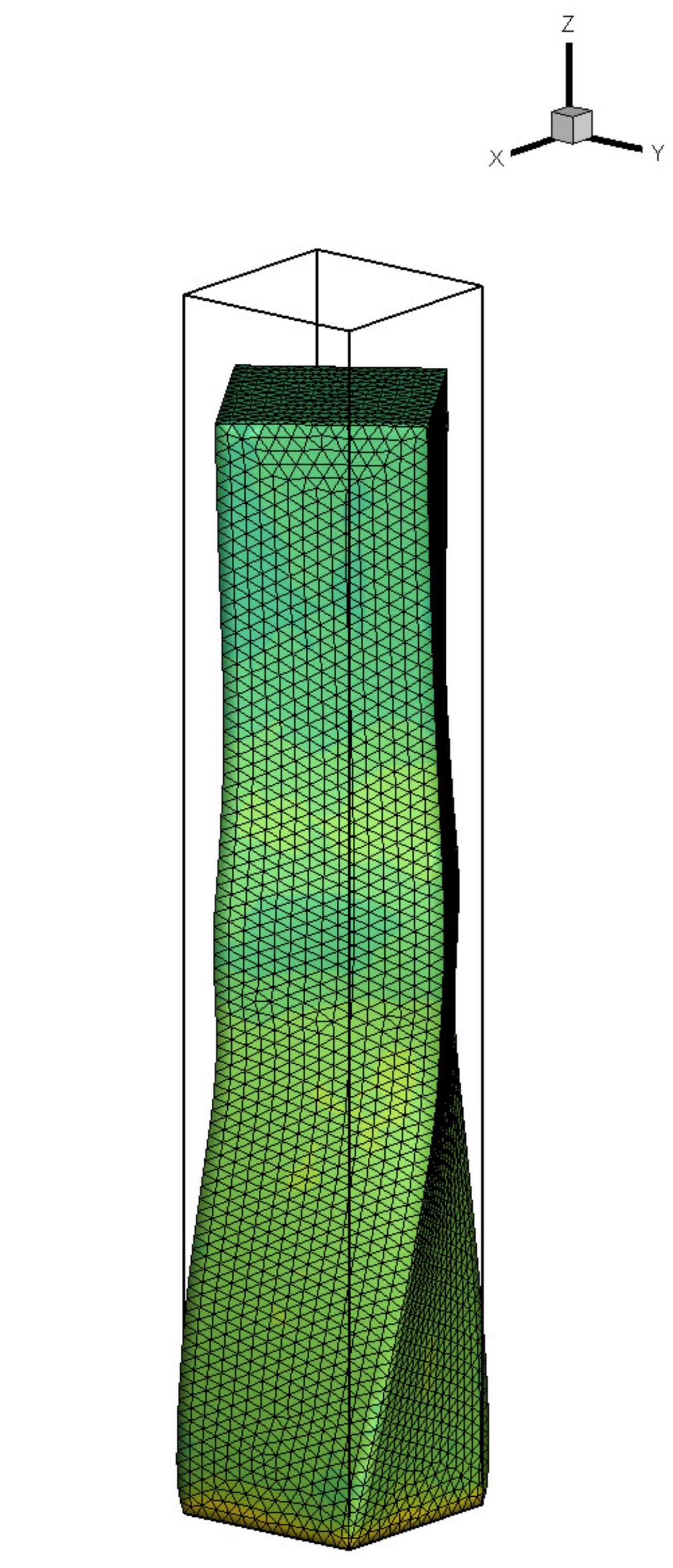} &
			\includegraphics[width=0.23\textwidth,draft=false]{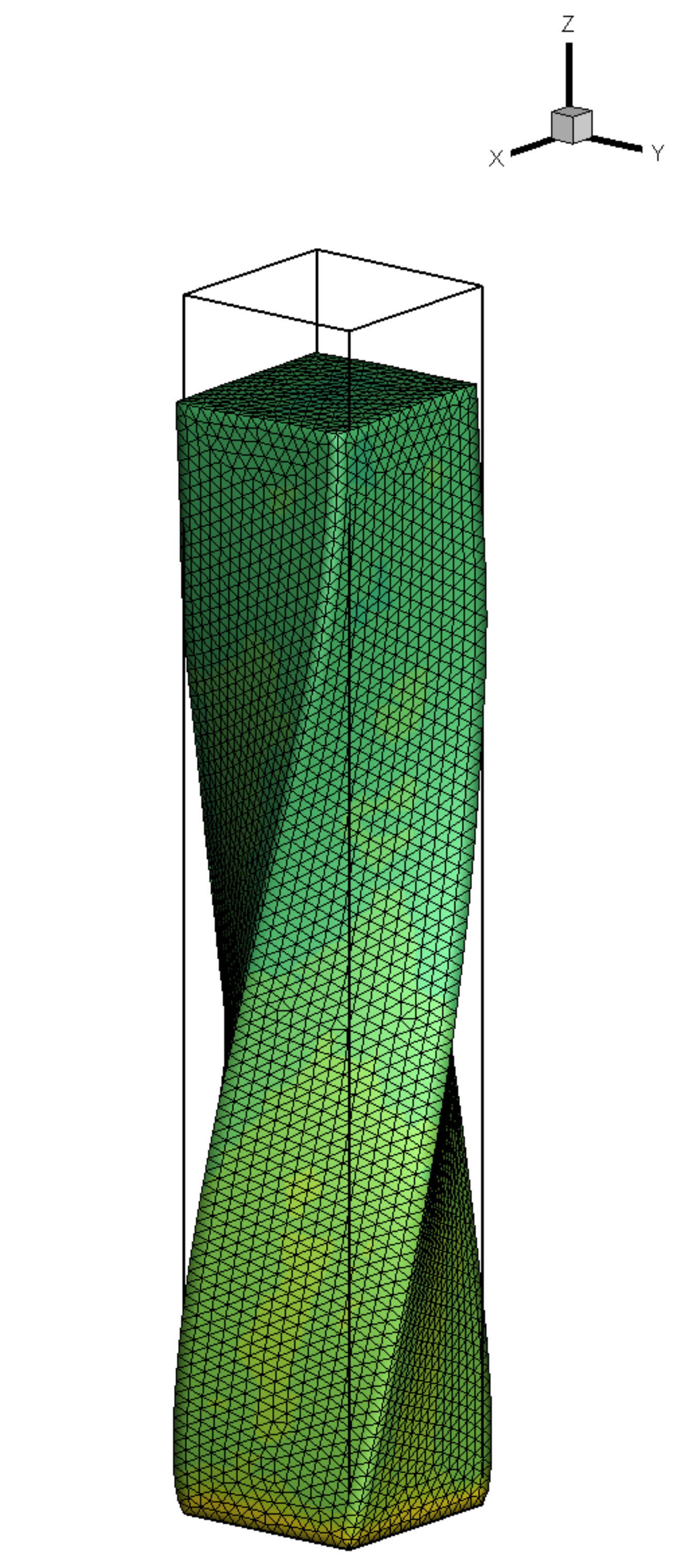} &
			\includegraphics[width=0.23\textwidth,draft=false]{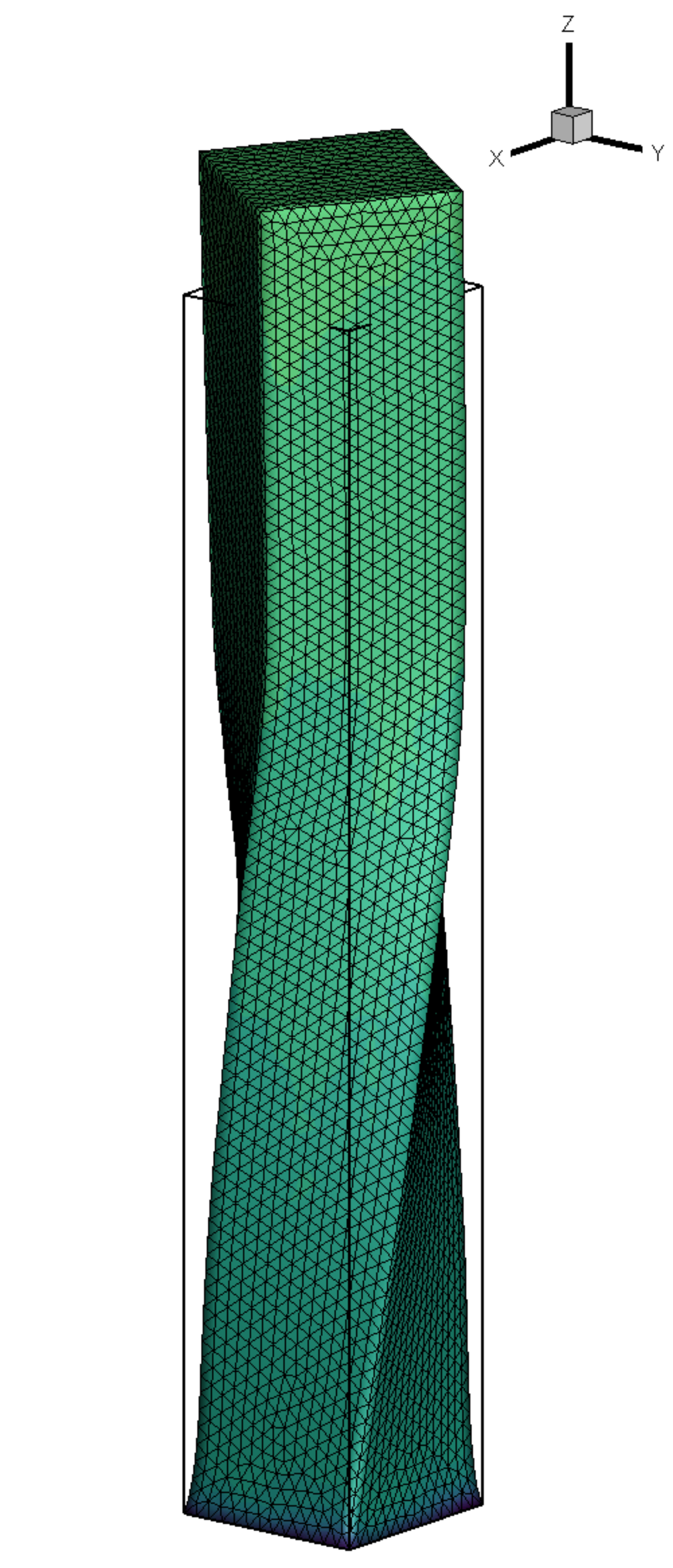} \\
			\multicolumn{4}{c}{\includegraphics[width=0.55\textwidth,draft=false]{TwistingColumn_legend_pressure}}
		\end{tabular}
		\caption{Twisting column with $\omega_0=200$. Column shape and pressure distribution at output times  $t = 0.00375$, $t = 0.075$, $t = 0.1125$, $t = 0.15$, $t = 0.1875$, $t = 0.225$, $t = 0.2625$ and $t = 0.3$ (from top left to bottom right). The shape is compared with respect to the initial configuration (hollow box). }
		\label{fig.TwistCol_u200}
	\end{center}
\end{figure}

Finally, a highly nonlinear scenario is considered by simulating a twisting column according to the test problem set up in \cite{Haider_2018} and the references therein. The computational domain is given by an initial unit squared cross section column of height $H = 6$, i.e. $\Omega(0)=[-0.5;0.5]\times[-0.5;0.5]\times[0;6]$. The $z=0$ face of the column is embedded into a wall, while the rest of the faces are assigned free-traction boundary conditions. An initial sinusoidal angular velocity field relative to the origin is given by
$\vv(0,\x) = \omega_0 \sin(\pi \frac{z}{2H}) ( y, -x, 0)^t$, while pressure is set to $p=0$. Two 
different magnitudes of the angular velocity are considered, namely $\omega_0=100$ and 
$\omega_0=200$. The main objective of this problem is to assess the capability of the proposed 
methodology to deal with the limit of incompressibility. A material with Neo-Hookean 
hydrodynamics EOS is used with Young modulus $Y = 1.7 \cdot 10^7$ and Poisson ratio $\nu = 0.45$. 
The simulation is run until $t_f=0.3$. Qualitatively one should observe a severe twist of the 
column which returns to its initial position. Driven by its own inertia, the bar keeps twisting 
until the final time. The mesh of the column is made of $N_E=119092$ tetrahedra with characteristic 
mesh size of $h=1/80$. Figure \ref{fig.TwistCol_u100} shows the pressure distribution at different 
output times for the case $\omega_0=100$, while Figure \ref{fig.TwistCol_u200} collects the results 
at the same output times for $\omega_0=200$. The initial column is represented as a hollow bar for 
comparison purposes and the expected behaviors are reproduced by the numerical simulation. Notice 
that there is no spurious oscillations nor nonphysical pressure distribution, thus the results 
obtained with the LGPR model are in agreement with the literature \cite{Boscheri2021,Gil2D_2014}. 
The pressure distribution is shown using 21 contour levels in the range $[-5;3] \cdot 10^6$.
	
\begin{figure}[!htbp]
	\begin{center}
		\begin{tabular}{cc}
			\includegraphics[width=0.47\textwidth,draft=false]{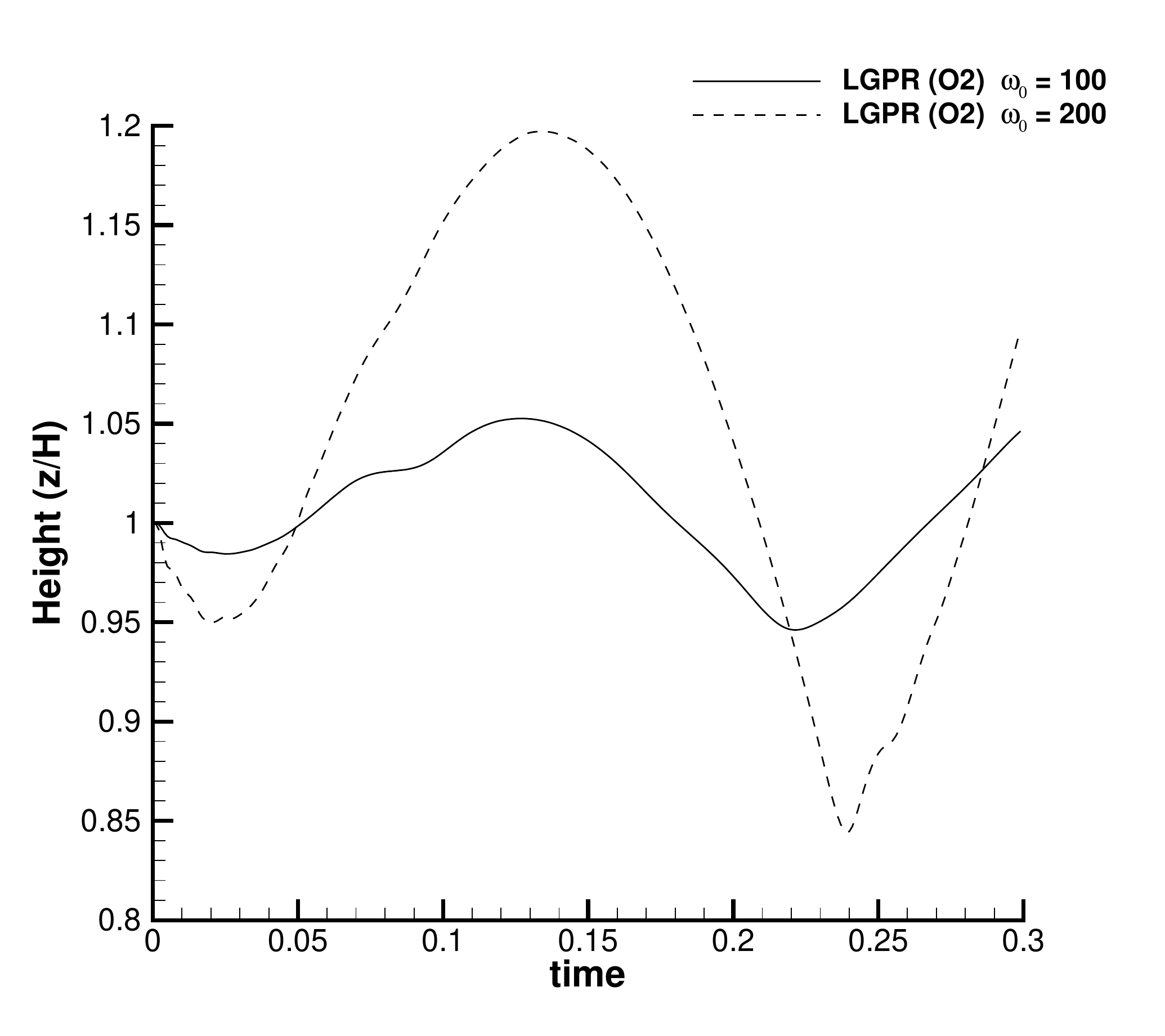}  &          
			\includegraphics[width=0.47\textwidth,draft=false]{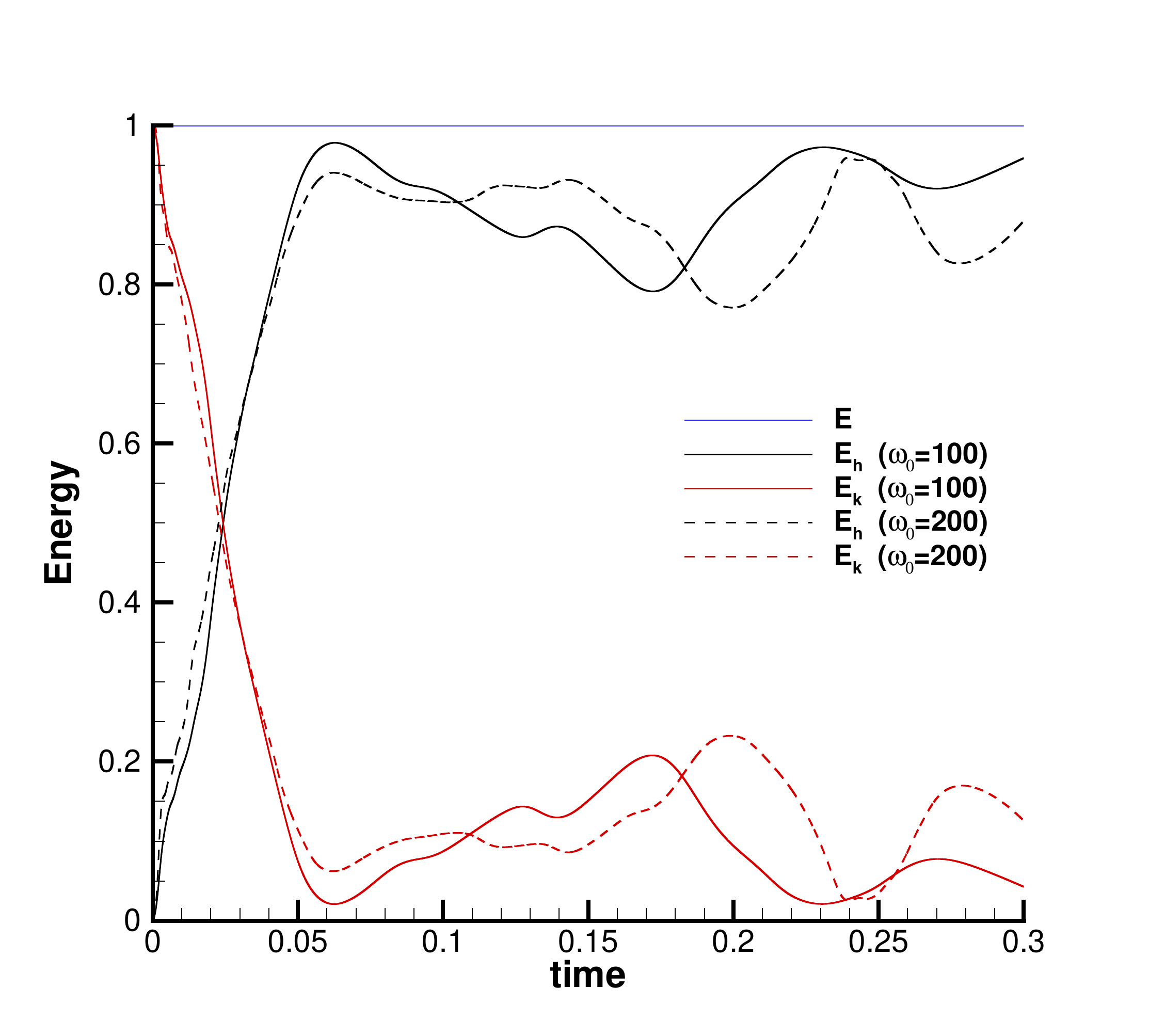} \\
		\end{tabular}
		\caption{Twisting column. Left: time evolution of non-dimensional height of the column measured at initial point $\x_T = (0,0,6)$. Right: analysis of internal ($E_h$), kinetic ($E_k$) and total ($E$) energy conservation normalized to unity for twisting velocity $\omega_0=100$ (solid lines) and $\omega_0=200$ (dashed lines). }
		\label{fig.TwistCol_z_en}
	\end{center}
\end{figure}

In Figure \ref{fig.TwistCol_z_en} we plot the time evolution of the dimensionless height of
the column measured at the point initially located at $\x_T = (0, 0, 6)$, showing that the case with $\omega_0=200$ exhibits a much stronger distortion and compression of the entire column. Furthermore, the analysis of energy conservation is also reported for both simulations, demonstrating that the total energy is fully conserved by the novel LGPR scheme \eqref{eqn.fv}. For comparison purposes, the total energy is normalized to unity for both simulations.
	
\section{Conclusions} \label{sec.concl}

We have presented a second-order in space and time updated Lagrangian IMEX scheme for the 
unified first-order hyperbolic formulation of continuum mechanics, also referred here as the 
Godunov-Peshkov-Romenski (GPR) model. The novelties of the paper can be summarized as follows. 
Firstly, for the first time, the GPR model was discretized with 
an updated Lagrangian scheme. Previously it was treated with either pure Eulerian 
\cite{DPRZ2016,SIGPR2021,Busto2020} or total Lagrangian  \cite{GodPesh2010} schemes, or 
even in the Arbitrary-Lagrangian-Eulerian (ALE) framework \cite{HyperHypo2019,LagrangeHPR}.
Secondly, since the GPR model is equipped with nonlinear relaxation source terms which may become stiff, 
the Lagrangian scheme is combined with a semi-analytical ad-hoc ODE solver, 
and second order extension is achieved via IMEX time discretization. 
This allows us to obtain the numerical solution consistent with the Navier-Stokes-Fourier limit for 
arbitrarily small relaxation times, thus satisfying the asymptotic preserving property. 
Moreover, the implicit treatment of the source terms is necessary to achieve better 
stability properties of the scheme in the stiff cases ($ \tau_1 \to 0 $ and $ \tau_2 
\to 0$). 
Thirdly, in this paper, we use a reduced version of the GPR model in which the evolution equation 
for the distortion field $ \Ae $ is substituted by the evolution equation of its symmetric part 
only (the metric tensor $ \Ge $). This simplification can be used if the rotational degrees of 
freedom of material particles encoded in the distortion field $ \Ae $ are negligible. 
The Lagrangian scheme is compatible with the GCL and the study of the asymptotic behavior 
of the first order fully discrete scheme is derived, demonstrating the asymptotic preserving property 
exhibited by the new scheme. Furthermore, because the GPR model provides a unified framework for the modeling of various material 
responses varying from ideal fluids to elastic solids, we have solved an extended series of 2D and 3D 
test cases covering examples from inviscid hydrodynamics, viscous heat conducting fluids, 
elastic and elasto-plastic solid mechanics. We emphasize that in all the presented test-cases for 
fluids and solids the \textit{same} set of governing equations \eqref{eqn.cl} and the \textit{same} scheme have been used. 
Whenever possible, we have compared the numerical solution obtained with the LGPR scheme 
against analytical solutions of inviscid Euler equations, Navier-Stokes-Fourier equations, and 
elasto-plasticity theory (e.g. see Section \ref{ssec.Shell}). A very good 
agreement between the solutions has been achieved, showing the capability of the scheme of dealing 
with different material properties. Finally, the new LGPR method collapses to already existing 
and well-established schemes when approaching the model limits, i.e. the hydrodynamics 
limit \cite{Maire2007} and the ideal hyperelasticity limit \cite{CCL2020,Boscheri2021}. 

Future work will consist in an extension of the present approach to general unstructured polygonal 
meshes as well as of generalization towards high-order (higher than 2) schemes which will require 
accounting for the mesh curvature via evolving the mesh characteristics (e.g. metric tensor) with 
high-order of accuracy. We also plan to discretize the original GPR model with the evolution 
equation for the distortion field $ \Ae $ and address the question of constructing structure 
preserving (e.g. curl-preserving) integration such as discussed in the pure Eulerian 
\cite{SIGPR2021} and total Lagrangian 
frameworks \cite{Gil2019}.
	
\section*{Acknowledgments}
WB and IP would like to thank the Italian Ministry of Instruction, University and Research (MIUR) to 
support this research with funds coming from PRIN Project 2017 No. 2017KKJP4X entitled "Innovative 
numerical methods for evolutionary partial differential equations and applications".	
SC acknowledges the financial support received by the Deutsche Forschungsgemeinschaft (DFG) under 
the project \textit{Droplet Interaction Technologies (DROPIT)}, grant no. GRK 2160/1. 

\bibliographystyle{plain}
\bibliography{biblio}

\begin{thebibliography}{100}

\bibitem{Aguirre2014}
Miquel Aguirre, Antonio~J Gil, Javier Bonet, and Aurelio Arranz.
\newblock {A vertex centred Finite Volume Jameson – Schmidt – Turkel ( JST
  ) algorithm for a mixed conservation formulation in solid dynamics}.
\newblock {\em Journal of Computational Physics}, 259:672--699, 2014.

\bibitem{AscRuuSpi}
U.~M. Ascher, S.~J. Ruuth, and R.~J. Spiteri.
\newblock {Implicit-explicit Runge-Kutta methods for time-dependent partial
  differential equations}.
\newblock {\em Appl. Numer. Math.}, 25:151--167, 1997.

\bibitem{StencilRec1990}
{T.J.} Barth and {P.O.} Frederickson.
\newblock Higher order solution of the euler equations on unstructured grids
  using quadratic reconstruction.
\newblock {\em 28th Aerospace Sciences Meeting}, pages AIAA paper no. 90--0013,
  January 1990.

\bibitem{BarthJespersen}
T.J. Barth and D.C. Jespersen.
\newblock The design and application of upwind schemes on unstructured meshes.
\newblock {\em AIAA Paper 89-0366}, pages 1--12, 1989.

\bibitem{BartonRom2010}
P~T Barton, D~Drikakis, and E~I Romenski.
\newblock {An Eulerian finite-volume scheme for large elastoplastic
  deformations in solids}.
\newblock {\em International Journal for Numerical Methods in Engineering},
  81:n/a--n/a, 2009.

\bibitem{Barton2019}
Philip~T. Barton.
\newblock {An interface-capturing Godunov method for the simulation of
  compressible solid-fluid problems}.
\newblock {\em Journal of Computational Physics}, 390:25--50, aug 2019.

\bibitem{Becker1923}
R.~Becker.
\newblock {Stosswelle und Detonation}.
\newblock {\em Physik}, 8:321, 1923.

\bibitem{Artzi}
M.~Ben-Artzi and J.~Falcovitz.
\newblock A second--order godunov--type scheme for compressible fluid dynamics.
\newblock {\em Journal of Computational Physics}, 55:1--32, 1984.

\bibitem{Benson1992}
D.J. Benson.
\newblock Computational methods in lagrangian and eulerian hydrocodes.
\newblock {\em Computer Methods in Applied Mechanics and Engineering},
  99:235--394, 1992.

\bibitem{MaireMM2}
M.~Berndt, J.~Breil, S.~Galera, M.~Kucharik, P.H. Maire, and M.~Shashkov.
\newblock {Two--step hybrid conservative remapping for multimaterial arbitrary
  Lagrangian–-Eulerian methods}.
\newblock {\em Journal of Computational Physics}, 230:6664--6687, 2011.

\bibitem{Bernstein60}
Barry Bernstein.
\newblock Hypo-elasticity and elasticity.
\newblock {\em Arch. Rational Mech. Anal.}, 6:89--104, 1960.
\newblock https://doi.org/10.1007/BF00276156.

\bibitem{Bonet2015}
J.~Bonet, A.~J. Gil, C.~Hean Lee, M.~Aguirre, and R.~Ortigosa.
\newblock A first order hyperbolic framework for large strain computational
  solid dynamics. part {I}: Total {L}agrangian isothermal elasticity.
\newblock {\em Computer Methods in Applied Mechanics and Engineering},
  283:689--732, January 2015.

\bibitem{Bonet2021}
J.~Bonet, C.~Hean Lee, A.~J. Gil, and A.~Ghavamian.
\newblock A first order hyperbolic framework for large strain computational
  solid dynamics. part {III}: Thermo-elasticity.
\newblock {\em Computer Methods in Applied Mechanics and Engineering},
  373:113505, January 2021.

\bibitem{BosRus}
S.~Boscarino and G.~Russo.
\newblock On a class of uniformly accurate {IMEX} {R}unge-{K}utta schemes and
  applications to hyperbolic systems with relaxation.
\newblock {\em SIAM J. Sci. Comput.}, 31:1926--1945, 2009.

\bibitem{BosARCME}
W.~Boscheri.
\newblock High order direct arbitrary-lagrangian–eulerian (ale) finite volume
  schemes for hyperbolic systems on unstructured meshes.
\newblock {\em Archives of Computational Methods in Engineering}, 24:751--801,
  2017.

\bibitem{BDLTV2020}
W.~Boscheri, G.~Dimarco, R.~Loub{\`{e}}re, M.~Tavelli, and M.H. Vignal.
\newblock {A second order all Mach number IMEX finite volume solver for the
  three dimensional Euler equations}.
\newblock {\em J. Comp. Phys.}, 415:109486, 2020.

\bibitem{BDT_cns}
W.~Boscheri, G.~Dimarco, and M.~Tavelli.
\newblock {An efficient second order all Mach finite volume solver for the
  compressible Navier-Stokes equations}.
\newblock {\em Computer Methods in Applied Mechanics and Engineering},
  374:113602, 2021.

\bibitem{Lagrange2D}
W.~Boscheri and M.~Dumbser.
\newblock {Arbitrary--{L}agrangian--{E}ulerian One--Step WENO Finite Volume
  Schemes on Unstructured Triangular Meshes}.
\newblock {\em Communications in Computational Physics}, 14:1174--1206, 2013.

\bibitem{Lagrange3D}
W.~Boscheri and M.~Dumbser.
\newblock {A direct Arbitrary-{L}agrangian-{E}ulerian ADER-WENO finite volume
  scheme on unstructured tetrahedral meshes for conservative and
  non-conservative hyperbolic systems in 3D}.
\newblock {\em J. Comput. Phys.}, 275:484--523, 2014.

\bibitem{SIGPR2021}
W.~Boscheri, M.~Dumbser, M.~Ioriatti, I.~Peshkov, and E.~Romenski.
\newblock {A structure-preserving staggered semi-implicit finite volume scheme
  for continuum mechanics}.
\newblock {\em Journal of Computational Physics}, 424:109866, jan 2021.

\bibitem{LagrangeHPR}
W.~Boscheri, M.~Dumbser, and R.~Loub\`ere.
\newblock {Cell centered direct Arbitrary-{L}agrangian-{E}ulerian ADER-WENO
  finite volume schemes for nonlinear hyperelasticity}.
\newblock {\em Computers and Fluids}, 134-135:111--129, 2016.

\bibitem{LAM2018}
W.~Boscheri, M.~Dumbser, R.~Loub\`ere, and P.-H. Maire.
\newblock A second-order cell-centered {L}agrangian {ADER-MOOD} finite volume
  scheme on multidimensional unstructured meshes for hydrodynamics.
\newblock {\em J. Comput. Phys.}, 358:103 -- 129, 2018.

\bibitem{Boscheri2021}
Walter Boscheri, Rapha{\"{e}}l Loub{\`{e}}re, and Pierre-Henri Maire.
\newblock {A cell-centered Lagrangian ADER-MOOD finite volume scheme on
  unstructured meshes for a class of hyper-elasticity models}.
\newblock apr 2021.

\bibitem{BOSCHERI2021110206}
Walter Boscheri and Lorenzo Pareschi.
\newblock High order pressure-based semi-implicit imex schemes for the 3d
  navier-stokes equations at all mach numbers.
\newblock {\em Journal of Computational Physics}, 434:110206, 2021.

\bibitem{Raviart.GRP.2}
A.~Bourgeade, P.~LeFloch, and P.A. Raviart.
\newblock An asymptotic expansion for the solution of the generalized riemann
  problem. {Part II}: application to the gas dynamics equations.
\newblock {\em Annales de l'institut Henri Poincar\'e (C) Analyse non
  lin\'eaire}, 6:437--480, 1989.

\bibitem{MaireMM1}
J.~Breil, S.~Galera, and P.H. Maire.
\newblock {Multi-material ALE computation in inertial confinement fusion code
  CHIC}.
\newblock {\em Computers and Fluids}, 46:161--167, 2011.

\bibitem{CCL2020}
J.~Breil, G.~Georges, and P.-H. Maire.
\newblock {3D cell-centered {L}agrangian second order scheme for the numerical
  modeling of hyperelasticity system}.
\newblock {\em Computer and Fluids}, 207:104523, 2020.

\bibitem{Busto2020}
Saray Busto, Simone Chiocchetti, Michael Dumbser, Elena Gaburro, and Ilya
  Peshkov.
\newblock {High Order ADER Schemes for Continuum Mechanics}.
\newblock {\em Frontiers in Physics}, 8(32), mar 2020.

\bibitem{Caramana1998}
E.J. Caramana, D.E. Burton, M.J. Shashkov, and P.P. Whalen.
\newblock The construction of compatible hydrodynamics algorithms utilizing
  conservation of total energy.
\newblock {\em Journal of Computational Physics}, 146:227--262, 1998.

\bibitem{SaltzmanOrg3D}
{E.J.} Caramana, {C.L.} Rousculp, and {D.E.} Burton.
\newblock {A compatible, energy and symmetry preserving {L}agrangian
  hydrodynamics algorithm in three-dimensional Cartesian geometry.}
\newblock {\em Journal of Computational Physics}, 157:89 -- 119, 2000.

\bibitem{Despres2009}
G.~Carr\'e, S.~Del Pino, B.~Despr\'es, and E.~Labourasse.
\newblock {A cell-centered {L}agrangian hydrodynamics scheme on general
  unstructured meshes in arbitrary dimension.}
\newblock {\em Journal of Computational Physics}, 228:5160--5183, 2009.

\bibitem{chengshu1}
J.~Cheng and C.W. Shu.
\newblock {A high order ENO conservative Lagrangian type scheme for the
  compressible Euler equations}.
\newblock {\em Journal of Computational Physics}, 227:1567--1596, 2007.

\bibitem{chengshu3}
J.~Cheng and C.W. Shu.
\newblock {A cell-centered Lagrangian scheme with the preservation of symmetry
  and conservation properties for compressible fluid flows in two-dimensional
  cylindrical geometry}.
\newblock {\em Journal of Computational Physics}, 229:7191--7206, 2010.

\bibitem{chengshu4}
J.~Cheng and C.W. Shu.
\newblock {Improvement on spherical symmetry in two-dimensional cylindrical
  coordinates for a class of control volume Lagrangian schemes}.
\newblock {\em Communications in Computational Physics}, 11:1144--1168, 2012.

\bibitem{cheng_jia_jiang_toro_yu_2017}
Jun-Bo Cheng, Yueling Jia, Song Jiang, Eleuterio~F. Toro, and Ming Yu.
\newblock A second-order cell-centered {L}agrangian method for two-dimensional
  elastic-plastic flows.
\newblock {\em Communications in Computational Physics}, 22(5):1224–1257,
  2017.

\bibitem{chiocchettimueller}
S.~Chiocchetti and C.~M{\"{u}}ller.
\newblock {A Solver for Stiff Finite-Rate Relaxation in Baer-Nunziato Two-Phase
  Flow Models}.
\newblock {\em Fluid Mechanics and its Applications}, 121:31--44, 2020.

\bibitem{Depres2012}
A.~Claisse, B.~Despr\'es, E.Labourasse, and F.~Ledoux.
\newblock {A new exceptional points method with application to cell-centered
  {L}agrangian schemes and curved meshes}.
\newblock {\em Journal of Computational Physics}, 231:4324--4354, 2012.

\bibitem{Despres_book_2017}
B.~Despr{\'e}s.
\newblock {\em Numerical Methods for {E}ulerian and {L}agrangian Conservation
  Laws}.
\newblock Frontiers in Mathematics. Springer International Publishing, 2017.

\bibitem{DepresMazeran2003}
B.~Despr\'es and C.~Mazeran.
\newblock {Symmetrization of Lagrangian gas dynamic in dimension two and
  multimdimensional solvers}.
\newblock {\em C.R. Mecanique}, 331:475--480, 2003.

\bibitem{Despres2005}
B.~Despr\'{e}s and C.~Mazeran.
\newblock {{L}agrangian gas dynamics in two dimensions and {L}agrangian
  systems}.
\newblock {\em Arch. Rational Mech. Anal.}, 178:327--372, 2005.

\bibitem{Dobrev1}
V.A. Dobrev, T.E. Ellis, Tz.V. Kolev, and R.N. Rieben.
\newblock {Curvilinear Finite elements for Lagrangian hydrodynamics}.
\newblock {\em International Journal for Numerical Methods in Fluids},
  65:1295--1310, 2011.

\bibitem{Dobrev2}
V.A. Dobrev, T.E. Ellis, Tz.V. Kolev, and R.N. Rieben.
\newblock {High Order Curvilinear Finite Elements for Lagrangian
  Hydrodynamics}.
\newblock {\em SIAM Journal on Scientific Computing}, 34:606--641, 2012.

\bibitem{Dobrev3}
V.A. Dobrev, T.E. Ellis, Tz.V. Kolev, and R.N. Rieben.
\newblock {High Order Curvilinear Finite Elements for axisymmetric Lagrangian
  Hydrodynamics}.
\newblock {\em Computers and Fluids}, 83:58--69, 2013.

\bibitem{SaltzmanOrg}
{J.K.} Dukovicz and B.~Meltz.
\newblock {Vorticity errors in multidimensional {L}agrangian codes.}
\newblock {\em Journal of Computational Physics}, 99:115 -- 134, 1992.

\bibitem{Dumbser2007693}
M.~Dumbser and M.~Kaeser.
\newblock Arbitrary high order non-oscillatory finite volume schemes on
  unstructured meshes for linear hyperbolic systems.
\newblock {\em Journal of Computational Physics}, 221:693 -- 723, 2007.

\bibitem{HYP2016}
Michael Dumbser, Ilya Peshkov, and Evgeniy Romenski.
\newblock {A unified hyperbolic formulation for viscous fluids and
  elastoplastic solids}.
\newblock In Christian Klingenberg and Michael Westdickenberg, editors, {\em
  Springer Proceedings in Mathematics and Statistics}, volume 237 of {\em
  Springer Proceedings in Mathematics and Statistics}, pages 451--463. Springer
  International Publishing, 2018.

\bibitem{DPRZ2016}
Michael Dumbser, Ilya Peshkov, Evgeniy Romenski, and Olindo Zanotti.
\newblock {High order ADER schemes for a unified first order hyperbolic
  formulation of continuum mechanics: Viscous heat-conducting fluids and
  elastic solids}.
\newblock {\em Journal of Computational Physics}, 314:824--862, jun 2016.

\bibitem{DPRZ2017}
Michael Dumbser, Ilya Peshkov, Evgeniy Romenski, and Olindo Zanotti.
\newblock {High order ADER schemes for a unified first order hyperbolic
  formulation of Newtonian continuum mechanics coupled with electro-dynamics}.
\newblock {\em Journal of Computational Physics}, 348:298--342, nov 2017.

\bibitem{Flanaghan_Belytchko_81}
D.~P. Flanagan and T.~Belytschko.
\newblock A uniform strain hexahedron and quadrilateral with orthogonal
  hourglass control.
\newblock {\em International Journal for Numerical Methods in Engineering},
  17(5):679--706, 1981.

\bibitem{Raviart.GRP.1}
P.~Le Floch and P.A. Raviart.
\newblock An asymptotic expansion for the solution of the generalized riemann
  problem. {Part I}: General theory.
\newblock {\em Annales de l'institut Henri Poincar\'e (C) Analyse non
  lin\'eaire}, 5:179--207, 1988.

\bibitem{Vilar2}
F.Vilar.
\newblock {Cell-centered discontinuous Galerkin discretization for
  two-dimensional Lagrangian hydrodynamics}.
\newblock {\em Computers and Fluids}, 64:64--73, 2012.

\bibitem{Vilar3}
F.Vilar, P.H. Maire, and R.~Abgrall.
\newblock {Cell-centered discontinuous Galerkin discretizations for
  two-dimensional scalar conservation laws on unstructured grids and for
  one-dimensional Lagrangian hydrodynamics}.
\newblock {\em Computers and Fluids}, 46(1):498--604, 2010.

\bibitem{Gavriluk08}
S.L. Gavrilyuk, N.~Favrie, and R.~Saurel.
\newblock Modelling wave dynamics of compressible elastic materials.
\newblock {\em Journal of Computational Physics}, 227:2941--2969, 2008.

\bibitem{3DHydroMaire}
G.~Georges, J.~Breil, and P.-H. Maire.
\newblock {A 3D GCL compatible cell-centered Lagrangian scheme for solving gas
  dynamics equations}.
\newblock {\em J. Comput. Phys.}, 305:921--941, 2016.

\bibitem{Gil2016}
A.~J. Gil, C.~Hean Lee, J.~Bonet, and R.~Ortigosa.
\newblock A first order hyperbolic framework for large strain computational
  solid dynamics. part {II}: Total {L}agrangian compressible, nearly
  incompressible and truly incompressible elasticity.
\newblock {\em Computer Methods in Applied Mechanics and Engineering},
  300:146--181, March 2016.

\bibitem{God1978}
S~K Godunov.
\newblock {\em {Elements of mechanics of continuous media}}.
\newblock Nauka, 1st russian edition, 1978.

\bibitem{GodPesh2010}
S.~K. Godunov and I.~M. Peshkov.
\newblock {Thermodynamically consistent nonlinear model of elastoplastic
  Maxwell medium}.
\newblock {\em Computational Mathematics and Mathematical Physics},
  50(8):1409--1426, aug 2010.

\bibitem{GodRom1972}
S.~K. Godunov and E.~I. Romenskii.
\newblock {Nonstationary equations of nonlinear elasticity theory in eulerian
  coordinates}.
\newblock {\em Journal of Applied Mechanics and Technical Physics},
  13(6):868--884, nov 1972.

\bibitem{GodRom2003}
S~K Godunov and E~I Romenskii.
\newblock {\em {Elements of continuum mechanics and conservation laws}}.
\newblock Kluwer Academic/Plenum Publishers, 2003.

\bibitem{GOUDREAU1982}
G.L. Goudreau and J.O. Hallquist.
\newblock Recent developments in large-scale finite element {L}agrangian
  hydrocode technology.
\newblock {\em Computer Methods in Applied Mechanics and Engineering},
  33(1):725--757, 1982.

\bibitem{Haider_2018}
J.~Haider, C.~Hean Lee, A.~J. Gil, A.~Huerta, and J.~Bonet.
\newblock An upwind cell centred total {L}agrangian finite volume algorithm for
  nearly incompressible explicit fast solid dynamic applications.
\newblock {\em Computer Methods in Applied Mechanics and Engineering}, 340:684
  -- 727, 2018.

\bibitem{Gil2019}
Osama~I. Hassan, Ataollah Ghavamian, Chun~Hean Lee, Antonio~J. Gil, Javier
  Bonet, and Ferdinando Auricchio.
\newblock {An upwind vertex centred finite volume algorithm for nearly and
  truly incompressible explicit fast solid dynamic applications: Total and
  Updated Lagrangian formulations}.
\newblock {\em Journal of Computational Physics: X}, 3:100025, jun 2019.

\bibitem{Howell2002}
B.P. Howell and G.J. Ball.
\newblock {A Free-Lagrange Augmented Godunov Method for the Simulation of
  Elastic–Plastic Solids}.
\newblock {\em Journal of Computational Physics}, 175(1):128--167, jan 2002.

\bibitem{HuShuTri}
C.~Hu and {C.W.} Shu.
\newblock Weighted essentially non-oscillatory schemes on triangular meshes.
\newblock {\em J. Comp. Phys.}, 150:97--127, 1999.

\bibitem{Jackson2019a}
Haran Jackson and Nikos Nikiforakis.
\newblock {A numerical scheme for non-Newtonian fluids and plastic solids under
  the GPR model}.
\newblock {\em Journal of Computational Physics}, 387:410--429, jun 2019.

\bibitem{Jackson2019}
Haran Jackson and Nikos Nikiforakis.
\newblock {A unified Eulerian framework for multimaterial continuum mechanics}.
\newblock {\em Journal of Computational Physics}, 401:109022, jan 2020.

\bibitem{KammLANL08}
J.~Kamm, J.~Brock, S.~Brandon, D.L. Cotrell, B.M. Johnson, P.~Knupp, T.G.
  Trucano, W.J. Rider, and V.G. Weirs.
\newblock Enhanced verification test suite for physics simulations codes.
\newblock {\em Technical Report LA-14379}, 2008.

\bibitem{SedovExact}
{J.R.} Kamm and {F.X.} Timmes.
\newblock On efficient generation of numerically robust sedov solutions.
\newblock {\em Technical Report LA-UR-07-2849, Los Alamos National Laboratory},
  2007.

\bibitem{kaeserjcp}
M.~K\"aser and A.~Iske.
\newblock {ADER} schemes on adaptive triangular meshes for scalar conservation
  laws.
\newblock {\em Journal of Computational Physics}, 205:486--508, 2005.

\bibitem{Kemm2020}
Friedemann Kemm, Elena Gaburro, Ferdinand Thein, and Michael Dumbser.
\newblock {A simple diffuse interface approach for compressible flows around
  moving solids of arbitrary shape based on a reduced Baer–Nunziato model}.
\newblock {\em Computers {\&} Fluids}, 204:104536, may 2020.

\bibitem{Kidder1976}
{R.E.} Kidder.
\newblock Laser-driven compression of hollow shells: power requirements and
  stability limitations.
\newblock {\em Nucl. Fus.}, 1:3 -- 14, 1976.

\bibitem{Kluth10}
G.~Kluth and B.~Despr{\'e}s.
\newblock Discretization of hyperelasticity on unstructured mesh with a
  cell-centered {L}agrangian scheme.
\newblock {\em Journal of Computational Physics}, 229(24):9092 -- 9118, 2010.

\bibitem{Gil2D_2014}
C.H. Lee, A.J. Gil, and J.~Bonet.
\newblock {Development of a stabilised Petrov-Galerkin formulation for
  conservation laws in {L}agrangian fast solid dynamics}.
\newblock {\em Comput. Methods Appl. Mech. Engrg.}, 268:40--64, 2014.

\bibitem{Morgan2019a}
E.J. Lieberman, X.~Liu, N.R. Morgan, D.J. Luscher, and D.E. Burton.
\newblock {A higher-order Lagrangian discontinuous Galerkin hydrodynamic method
  for solid dynamics}.
\newblock {\em Comp. Meth. Appl. Mech. Engrg.}, 353:467--490, 2019.

\bibitem{Morgan2019b}
E.J. Lieberman, N.R. Morgan, D.J. Luscher, and D.E. Burton.
\newblock {A higher-order Lagrangian discontinuous Galerkin hydrodynamic method
  for elastic-plastic flows}.
\newblock {\em Comput. Math. Appl.}, 353:318--334, 2019.

\bibitem{chengshu2}
W.~Liu, J.~Cheng, and C.W. Shu.
\newblock {High order conservative Lagrangian schemes with Lax–Wendroff type
  time discretization for the compressible Euler equations}.
\newblock {\em Journal of Computational Physics}, 228:8872--8891, 2009.

\bibitem{StagLag}
R.~Loub\`ere, P.H. Maire, and P.~V\'achal.
\newblock {A second--order compatible staggered Lagrangian hydrodynamics scheme
  using a cell--centered multidimensional approximate Riemann solver}.
\newblock {\em Procedia Computer Science}, 1:1931--1939, 2010.

\bibitem{LoubereSedov3D}
R.~Loub\`ere, {P.H.} Maire, and P.~V\'achal.
\newblock {3D staggered Lagrangian hydrodynamics scheme with cell-centered
  Riemann solver-based artificial viscosity.}
\newblock {\em International Journal for Numerical Methods in Fluids}, 72:22 --
  42, 2013.

\bibitem{phm109}
P.-H. Maire.
\newblock A high-order cell-centered {L}agrangian scheme for compressible fluid
  flows in two-dimensional cylindrical geometry.
\newblock {\em J. Comput. Phys.}, 228(18):6882--6915, 2009.

\bibitem{Maire2011a}
P.-H. Maire.
\newblock {A unified sub-cell force-based discretization for cell-centered
  {L}agrangian hydrodynamics on polygonal grids}.
\newblock {\em Int. J. Numer. Meth. Fluid}, 65:1281--1294, 2011.

\bibitem{Maire2011}
P.-H. Maire.
\newblock A high-order one-step sub-cell force-based discretization for
  cell-centered {L}agrangian hydrodynamics on polygonal grids.
\newblock {\em Computer and Fluids}, 46(1):341--347, 2011.

\bibitem{Maire_elasto}
P.-H. Maire, R.~Abgrall, J.~Breil, R.~Loub\`{e}re, and B.~Rebourcet.
\newblock {A nominally second-order cell-centered {L}agrangian scheme for
  simulating elastic-plastic flows on two-dimensional unstructured grids}.
\newblock {\em J. Comput. Phys.}, 235:626--665, 2013.

\bibitem{Maire2007}
P.-H. Maire, R.~Abgrall, J.~Breil, and J.~Ovadia.
\newblock {A cell-centered {L}agrangian scheme for two-dimensional compressible
  flow problems}.
\newblock {\em SIAM Journal on Scientific Computing}, 29:1781--1824, 2007.

\bibitem{phmbn09}
P.-H. Maire and B.~Nkonga.
\newblock {Multi-scale {G}odunov-type method for cell-centered discrete
  {L}agrangian hydrodynamics}.
\newblock {\em J. Comput. Phys.}, 228(3):799--821, 2009.

\bibitem{RomenskiMalyshev1987}
A~N Malyshev and E~I Romenskii.
\newblock {Hyperbolic equations for heat transfer. Global solvability of the
  Cauchy problem}.
\newblock {\em Siberian Mathematical Journal}, 27(5):734--740, 1987.

\bibitem{ContMechBook}
L.~Anand M.E.~Gurtin, E.~Fried.
\newblock {\em {The mechanics and thermodynamics of continua}}.
\newblock Cambridge University Press, 2009.

\bibitem{munz94}
C.D. Munz.
\newblock {On Godunov--type schemes for Lagrangian gas dynamics}.
\newblock {\em SIAM Journal on Numerical Analysis}, 31:17--42, 1994.

\bibitem{scovazzi1}
A.~L\'opez Ortega and G.~Scovazzi.
\newblock {A geometrically--conservative, synchronized, flux--corrected remap
  for arbitrary Lagrangian--Eulerian computations with nodal finite elements}.
\newblock {\em Journal of Computational Physics}, 230:6709--6741, 2011.

\bibitem{PR_IMEX}
L.~Pareschi and G.~Russo.
\newblock Implicit-explicit runge-kutta schemes and applications to hyperbolic
  systems with relaxation.
\newblock {\em J. Sci. Comput.}, 25:129--155, 2005.

\bibitem{HyperHypo2019}
Ilya Peshkov, Walter Boscheri, Rapha{\"{e}}l Loub{\`{e}}re, Evgeniy Romenski,
  and Michael Dumbser.
\newblock {Theoretical and numerical comparison of hyperelastic and hypoelastic
  formulations for Eulerian non-linear elastoplasticity}.
\newblock {\em Journal of Computational Physics}, 387:481--521, 2019.

\bibitem{nonNewtonian2021}
Ilya Peshkov, Michael Dumbser, Walter Boscheri, Evgeniy Romenski, Simone
  Chiocchetti, and Matteo Ioriatti.
\newblock {Simulation of non-Newtonian viscoplastic flows with a unified first
  order hyperbolic model and a structure-preserving semi-implicit scheme}.
\newblock {\em Computers {\&} Fluids}, 224:104963, jun 2021.

\bibitem{SHTC-GENERIC-CMAT}
Ilya Peshkov, Michal Pavelka, Evgeniy Romenski, and Miroslav Grmela.
\newblock {Continuum mechanics and thermodynamics in the Hamilton and the
  Godunov-type formulations}.
\newblock {\em Continuum Mechanics and Thermodynamics}, 30(6):1343--1378, nov
  2018.

\bibitem{HPR2016}
Ilya Peshkov and Evgeniy Romenski.
\newblock {A hyperbolic model for viscous Newtonian flows}.
\newblock {\em Continuum Mechanics and Thermodynamics}, 28(1-2):85--104, mar
  2016.

\bibitem{Torsion2019}
Ilya Peshkov, Evgeniy Romenski, and Michael Dumbser.
\newblock {Continuum mechanics with torsion}.
\newblock {\em Continuum Mechanics and Thermodynamics}, 31(5):1517--1541, 2019.

\bibitem{PTRSA2020}
Evgeniy Romenski, Ilya Peshkov, Michael Dumbser, Francesco Fambri, Fambri, and
  Francesco.
\newblock {A new continuum model for general relativistic viscous
  heat-conducting media}.
\newblock {\em Philosophical Transactions of the Royal Society A: Mathematical,
  Physical and Engineering Sciences}, 378(2170):20190175, may 2020.

\bibitem{Romenski1979}
E.~I. Romenskii.
\newblock {Dynamic three-dimensional equations of the Rakhmatulin
  elastic-plastic model}.
\newblock {\em Journal of Applied Mechanics and Technical Physics},
  20(2):229--244, 1979.

\bibitem{Rom1989}
E~I Romenskii.
\newblock {Hyperbolic equations of Maxwell's nonlinear model of elastoplastic
  heat-conducting media}.
\newblock {\em Siberian Mathematical Journal}, 30(4):606--625, 1990.

\bibitem{Rubin2019}
M~B Rubin.
\newblock {An Eulerian formulation of inelasticity: from metal plasticity to
  growth of biological tissues}.
\newblock {\em Philosophical Transactions of the Royal Society A: Mathematical,
  Physical and Engineering Sciences}, 377(2144):20180071, may 2019.

\bibitem{Rusanov:1961a}
V.~V. Rusanov.
\newblock {Calculation of Interaction of Non--Steady Shock Waves with
  Obstacles}.
\newblock {\em J. Comput. Math. Phys. USSR}, 1:267--279, 1961.

\bibitem{Hank2017}
J.~Massoni S.~Hank, N.~Favrie.
\newblock Modeling hyperelasticity in non-equilibrium multiphase flows.
\newblock {\em Journal of Computational Physics}, 330:65--91, 2017.

\bibitem{Sambasivan2013}
S.K. Sambasivan, M.J. Shashkov, and D.E. Burton.
\newblock {A finite volume cell-centered {L}agrangian hydrodynamics approach
  for solids in general unstructured grids}.
\newblock {\em Int. J. Numer. Meth. Fluid}, 72:770--810, 2013.

\bibitem{scovazzi2}
G.~Scovazzi.
\newblock {{L}agrangian shock hydrodynamics on tetrahedral meshes: A stable and
  accurate variational multiscale approach}.
\newblock {\em Journal of Computational Physics}, 231:8029--8069, 2012.

\bibitem{scovazzi3}
G.~Scovazzi, B.~Carnes, X.~Zeng, and S.~Rossi.
\newblock {A simple, stable, and accurate linear tetrahedral finite element for
  transient, nearly, and fully incompressible solid dynamics: a dynamic
  variational multiscale approach}.
\newblock {\em International Journal for Numerical Methods in Engineering},
  106:799--839, 2016.

\bibitem{chengshu5}
C.-W. Shu and J.~Cheng.
\newblock {A Third Order Conservative Lagrangian Type Scheme on Curvilinear
  Meshes for the Compressible Euler Equations}.
\newblock {\em Commun. Comput. Phys.}, 4:1008--1024, 2008.

\bibitem{godunov}
{S.K.}Godunov.
\newblock A finite difference method for the computation of discontinuous
  solutions of the equations of fluid dynamics.
\newblock {\em Mat. Sbornik}, 47:357--393, 1959.

\bibitem{Smith1999}
{R.W.} Smith.
\newblock {AUSM(ALE)}: a geometrically conservative arbitrary
  lagrangian--eulerian flux splitting scheme.
\newblock {\em Journal of Computational Physics}, 150:268–286, 1999.

\bibitem{tavellicrack}
M.~Tavelli, E.~Romenski, S.~Chiocchetti, A.-A. Gabriel, and M.~Dumbser.
\newblock {Space-time adaptive ADER discontinuous Galerkin schemes for
  nonlinear hyperelasticity with material failure}.
\newblock {\em Journal of Computational Physics}, page 109758, 2020.

\bibitem{TavelliELDI}
Maurizio Tavelli, Michael Dumbser, Dominic~Etienne Charrier, Leonhard
  Rannabauer, Tobias Weinzierl, and Michael Bader.
\newblock A simple diffuse interface approach on adaptive cartesian grids for
  the linear elastic wave equations with complex topography.
\newblock {\em Journal of Computational Physics}, 386:158--189, 2019.

\bibitem{Taylor}
G.I. Taylor.
\newblock The use of flat-ended projectiles for determining dynamic yield
  stress i. theoretical considerations.
\newblock {\em Proceedings of the Royal Society of London A: Mathematical,
  Physical and Engineering Sciences}, 194(1038):289--299, 1948.

\bibitem{toro3}
{V.A.} Titarev and {E.F.} Toro.
\newblock {ADER}: Arbitrary high order {Godunov} approach.
\newblock {\em Journal of Scientific Computing}, 17(1-4):609--618, December
  2002.

\bibitem{titarevtoro}
V.A. Titarev and E.F. Toro.
\newblock {ADER} schemes for three-dimensional nonlinear hyperbolic systems.
\newblock {\em Journal of Computational Physics}, 204:715--736, 2005.

\bibitem{Toro:2006a}
E.~F. Toro and V.~A. Titarev.
\newblock {Derivative Riemann solvers for systems of conservation laws and ADER
  methods}.
\newblock {\em Journal of Computational Physics}, 212(1):150--165, 2006.

\bibitem{toro.anomalies.2002}
E.F. Toro.
\newblock Anomalies of conservative methods: analysis, numerical evidence and
  possible cures.
\newblock {\em International Journal of Computational Fluid Dynamics},
  11:128--143, 2002.

\bibitem{ToroBook}
{E.F.} Toro.
\newblock {\em Riemann Solvers and Numerical Methods for Fluid Dynamics: a
  Practical Introduction.}
\newblock Springer, 2009.

\bibitem{Truesdell55}
C.~Truesdell.
\newblock Hypo-elasticity.
\newblock {\em Journal of Rational Mechanics and Analysis}, 4:83--1020, 1955.

\bibitem{Neumann1950}
J.~von Neumann and R.D. Richtmyer.
\newblock A method for the numerical calculations of hydrodynamical shocks.
\newblock 21:232--238, 1950.

\bibitem{Barton2020}
Tim Wallis, Philip~T. Barton, and Nikolaos Nikiforakis.
\newblock {A flux-enriched Godunov method for multi-material problems with
  interface slide and void opening}.
\newblock {\em Journal of Computational Physics}, 442:110499, oct 2021.

\bibitem{wil1}
M.L. Wilkins.
\newblock Calculation of elastic plastic flow.
\newblock In B.~Alder, S.~Fernbach, and M.~Rotenberg, editors, {\em Methods in
  Computational Physics}, volume~3, pages 211--263. Academic Press, New York,
  1964.

\bibitem{Morgan2019c}
T.~Wu, M.~Shashkov, N.R. Morgan, D.~Kuzmin, and H.~Luo.
\newblock {An updated Lagrangian discontinuous Galerkin hydrodynamic method for
  gas dynamics}.
\newblock {\em Comput. Math. Appl.}, 78:258--273, 2019.

\end{thebibliography}

\end{document}